\documentclass[11pt, a4paper, oneside]{book}

\usepackage[english]{babel}
\usepackage[T1]{fontenc}
\usepackage[utf8]{inputenc}
\usepackage[all]{xy}
\usepackage[]{geometry}
\usepackage{amsmath}
\usepackage{amssymb}

\usepackage{amsthm}
\usepackage{txfonts}
\usepackage{mathrsfs}
\usepackage{leftidx}
\usepackage{mathabx}
\usepackage{amscd}
\usepackage{pst-node}
\usepackage{tikz}
\usepackage{hyperref}
\usepackage{url}

\hypersetup{colorlinks,citecolor=black,filecolor=black,linkcolor=black,urlcolor=black}

\newcommand{\up}[1]{\textsuperscript{#1}}
\newcommand{\petitg}{g}
\newcommand{\thmunun}{\ref{t1.1}}
\newcommand{\propdeuxun}{\ref{p2.1}}
\newcommand{\propuncinq}{\ref{p1.5}}

\newcommand{\R}{\mathbb{R}}
\newcommand{\Q}{\mathbb{Q}}
\newcommand{\N}{\mathbb{N}}
\newcommand{\CC}{\mathbb{C}}
\newcommand{\M}{\mathscr{M}}
\newcommand{\D}{\mathscr{D}}
\newcommand{\C}{\mathscr{C}}
\newcommand{\Ss}{\mathscr{S}}
\newcommand{\HH}{\mathscr{H}}
\newcommand{\V}{\mathcal{V}}
\newcommand{\A}{\mathscr{A}}

\everymath{\displaystyle}

\renewcommand{\leq}{\leqslant}
\renewcommand{\geq}{\geqslant}

\newcommand{\addtotdm}[1]{\addcontentsline{toc}{section}{\\[-0.8em]\hspace*{-0.85cm}\textbf{#1}}}

\newtheorem{dfn}{Definition}
\newtheorem{lem}{Lemma}
\newtheorem{thme}{Theorem}
\newtheorem{dem}{Proof}
\newtheorem{prop}{Proposition}
\newtheorem{coro}{Corollary}

\theoremstyle{remark}
\newtheorem{rmq}{Remark}

\title{
  \bf\huge
  A Mathematicians' View\\[0.3em] of Geometrical Unification\\[0.3em] of General Relativity and Quantum Physics\\[1em]
  \normalsize(translated from french)
  }
\author{{\Large\centering{Michel Vaugon}}\\[0.5em]
  \sl Former Professor of Mathematics at Pierre et Marie Curie University\\
  \tt michel.vaugon@gmail.com
}
\date{September 2020}

\begin{document}

\maketitle

\section*{Abstract}

This document contains a description of physics entirely based on a geometric presentation: all of the theory is described giving only a pseudo-riemannian manifold $(\M, g)$ of dimension $n>5$ for which the tensor $g$ is, in studied
domains, almost everywhere of signature $(-, -, +,\dots, +)$. No object is added in this space-time, no general principle
is assumed. The properties we demand to some domains of $(\M, g)$ are only simple geometric constraints, essentially based
on the concept of `` curvature ''. These geometric properties allow to define, depending on considered cases, some objects
(frequently depicted by tensors) that are similar to the classical physics ones, they are however built here only from the  tensor $g$. The links between these objects, coming from their natural definitions, give, applying standard theorems
from the pseudo-riemannian geometry, all equations governing physical phenomenons usually described by classical theories, including general relativity and quantum physics. The purely geometric approach introduced here on quantum phenomena is profoundly different from the standard one. Neither Lagrangian or Hamiltonian is used. This document ends with a presentation of our approach of complex quantum phenomena usually studied by quantum field theory.

\newpage

\section*{Summary}

In this paper we offer a description of physics based entirely on a geometric presentation:
the whole theory is described from the simple data of a pseudo-Riemannian manifold $(\M,g)$ of dimension $n>5$
for
which the tensor $g$ is, in studied domains, almost everywhere of signature $(-,-,+,\dots,+)$. No
\begin{it} object \end{it} is added in this \begin{it} space-time\end{it}, no general principle is assumed. The
special properties
demanded on some domains of $(\M,g)$ are simple geometrical constraints essentially based on the
concept of curvature. These \begin{it} geometric properties \end{it} make it possible to define, according to the cases considered,
\begin{it} objects \end{it} (often represented by tensors)
which are similar to those of classical physics but which, here, are only constructed from the tensor $g$.
The links between these objects, which come from their natural definitions, make it possible to obtain, by
applying pseudo-Riemannian geometry standard theorems, all the equations that manage the description of
physical phenomena usually described by classical theories, including the theory of general relativity
and quantum physics. Neither \begin{it} Lagrangian \end{it} or \begin{it} Hamiltonian \end{it} is used. The domains of space-time $(\M,g)$
that are studied are locally diffeomorphic to $\Theta\times K$ where $\Theta$ is an open set of $\R^4$ and $K$ is a
compact manifold. The first negative sign of the signature of $g$ is associated to $\Theta$, the second negative sign
(corresponding to another notion of "time") is associated to $K$ and this one is an essential ingredient in the
description of electromagnetism.\\
Two different types of approximations make it possible to recover the usual equations of physics~:

 -A first type of approximation consists in neglecting some phenomena associated to the compact manifold $K$ and this
allows to recover the equations of non-quantum physics (general relativity including electromagnetism).

-A second type of approximation is to assume that the tensor $g$ is "bound" to the
Minkovski's metric on $\Theta$, but in this case we take into account precisely the characteristics of the compact manifold  $K$. In particular we
then recovers the qualitative and quantitative results usually obtained by classic quantum physics, without
using the axiomatic system of this theory. The purely geometrical view offered here on quantum phenomena is profoundly different from that of standard theories.

This text ends with a presentation of how we approach, with the theory presented here, the study of
complex quantum phenomena treated by quantum field theory.


\newpage
\section*{Foreword}
\addtotdm{Foreword}

What I am going to present here can be seen as an extension of what
was written in the manuscripts \cite{vaugon-1}, \cite{vaugon-2} and \cite{vaugon-3}
and in the paper \cite{steph-1}.
However, the presentation will be such that it is not necessary to consult
these for reading this paper.
\bigskip

\noindent This theory was developed with the friendly participation of:
\medskip

\begin{itemize}
\item[$-$] Stéphane Collion: University Doctor in
    Mathematics, Airline pilot at Air France.
    \smallskip
\item[$-$] Marie Dellinger: University Doctor in Mathematics,
    undergrad teacher at ENCPB.
    \smallskip
 \item[$-$] Zoé Faget: University Doctor in Mathematics,
    University Doctor in Computer Science, former Associate Professor at
    the University of Poitiers, undergrad teacher at Lycée Fenelon.
    \smallskip
  \item[$-$] Emmanuel Humbert: University Doctor of
    Mathematics, Professor and thesis Director at Tours-University.
    \smallskip
    \item[$-$] Benoît Vaugon: Mathematician, Physicist, Doctor in
    Computer Science \\ ENSTA-ParisTech.
    \smallskip
  \item[$-$] Claude Vaugon: Associate Professor of Mathematics at
    high school Jean de La Fontaine  (Château-Thierry).
  \end{itemize}
 \bigskip
 


\newpage
\section*{Introduction}
\addtotdm{Introduction}

This paper begins with a quick historical presentation of the considerations that led us
to the theory that will be presented here. They are linked to the (personal) view
of the physics of the last century that I summarize in three characteristic steps
which concern what is commonly called "classical physics" (non-quantum).

\paragraph{1 \up{era} stage} (before relativity theory) $ $
\smallskip

Space-time is modeled by $\R \times \R^3$ (or better by an affine space), $\R$ for time
considered absolute, $\R^3$ for space, provided with the Euclidean scalar product.
In this space "physical objects" exist, that are modeled either by "curves" 
(for particles for example), or by tensor fields: functions, vector fields, differential forms, etc.
(for "fluids", electric and magnetic fields, mass density or electric charge functions, for example). These physical objects are considered to have no influence on
the notions of time and distance given in an absolute way in $\R \times \R^3$.
These objects are governed by \textbf{laws} and some respect \textbf{principles}.
These laws are written associated to particular space-time observers
often referred to as "Galilean observers". We can mention: Newton's laws
for gravitation, Maxwell's laws for electromagnetism. The principles
admitted are: the invariance of the laws during the changes of Galilean observers,
homogeneity and isotropy of space. This model is essentially faulted
by the experimental observation of the absolute constancy of the speed of light
and the fact that Maxwell's equations are not invariant during changes
of Galilean observers. Which leads to the second step.

\paragraph{2 \up{th} stage} (theory of special relativity)
\smallskip

Space-time is modeled by $\R^4$ (or better, an affine space) equiped
with a Lorentz quadratic form: $q(t, x, y, z) = -c^2 t^2 + x^2 + y^2 + z^2$.
The notions of space and time are then strongly linked (time is no longer absolute).
This model makes it perfectly consistent that the speed of light is an absolute constant.
Physical objects are modeled as presented in the first step
and they still have no influence on the notions of time and distance given
by the  Lorentz quadratic form. Maxwell's laws now respect the new
principle of invariance by Lorentz transformations (which replace the Galileo's transformations). However, Newton's laws, used to describe gravitational phenomena,
become totally unsuited to this new space-time (with the Lorentz quadratic form). The basic idea that solves the problem is in the third step.

\paragraph{3 \up{th} stage} (theory of general relativity)
\smallskip

Space-time is now modeled by a 4-dimensional manifold $\M$ provided
with a Lorentz tensor $g$ whose signature, at each point of $\M$, is
$(-, +, +, +)$. In other words, the tangent space at each point of $\M$ is provided with
 a Lorentz quadratic form. Physical objects are still represented
by tensor fields defined on the manifold $\M$. To each object is associated its "energy-momentum" tensor which is a field of quadratic forms.

The fundamental law of physics is given by Einstein's equation which states
that the energy-momentum tensor, which characterizes objects found in
a domain of manifold, is equal to Einstein's curvature (for a good choice of "units").
So there is now
a close link between the notion of space-time and the physical objects themselves.
Electromagnetism is introduced very naturally in this space-time, and the "object"
caracteristing it is a $2$-differential form $F$ to which one associates its
energy-momentum tensor. The $2$-differential form $F$ is supposed to satisfy
the "Maxwell laws" whose differential operators are now derived themselves
from the Lorentz tensor $g$.

This representation of physics works very well with respect to
gravitation and electromagnetism (neglecting quantum effects). It
allows in particular the description of unexpected phenomena, represented by
"singularities" of space-time, such as big-bang, black holes, etc.
Of course, the models presented in the first and second stages
now appear as approximations of specific areas 
of general relativity. \textbf{We notice that in this theory, no "principle"
remains, homogeneity, isotropy, and more generally invariance under action
of some isometries groups, are only approximations
which allow to calculate but are obviously not general principles}.

\bigskip

In the 3 steps I have just described, the "laws of physics" are
given by differential equations which link the objects which
 intervene in space-time. Another point of view is to use
the \textbf{Lagrangian principle}: instead of giving laws axiomatically 
in the form of differential equations, we give, for a domain of
space-time containing physical objects, an \textbf{action} that will characterize
the behavior of these objects. The action is a map with real values,
the "variables" are the considered objects. The "axiom" then consists
to say that some of these objects form a stationary point of the action (minimize it,
for example). Mathematically, this express into the fact that these objects satisfy
differential equations (those that could have been given as axioms).
In summary, instead of directly giving the differential equations, we prefer
to give oneself an "action" from which one will deduce the equations.
The interesting aspect here is that often the expression of "action"
is more "aesthetic", even more "intuitive" than the differential equations
themselves (although in many cases, equations were found before actions). However, the notion of "Lagrangian"
has other advantages: it can sometimes considerably simplify  
calculations that are presented in some physical problems, but the importance given now
to the notion of Lagrangian (and its diversion in Hamiltonian) is essentially due
to the fact that this notion is necessary in the axiomatic system of physical theories
that deal with phenomena that are not described by "classical physics"
and that I will call the "quantum phenomena". In theories that try 
to describe quantum phenomena, I will classify: classical quantum mechanics,
quantum field theory, string theory, etc. These theories
have evolued (say for a century) parallel to classical physics theories that I presented in steps $2$ and $3$.
The procedures used to describe the (quantum) phenomena observed are
very different from those used in classical physics. If we start, too, from an "action" linked to the physical objects studied, the axiomatic system
used is very far from that of classical physics: we do not look
"stationary points" of the action, but we describe processes that allow,
from the action (modified), to obtain density of probabilities for the quantities
characteristics of these objects. In fact, these theories have been the subject of considerable work
these last decades and certainly can not be summed up in a few lines.
A more precise description would be of little interest to us because it is with an
\textbf{another point of view} that we will describe the quantum phenomena.
This point of view can be considered as an extension of the theory
of general relativity, however, many common points with the 
quantum fields theory will naturally appear as the informed reader will notice. \\
In this paper, the notion of Lagrangian-Hamiltonian
will be \textbf{completely let aside}. We replace it with the conceptually very different notion of "geometric type".

What I'm going to write now can be seen as the next step and we can, for the time being, disregard
from what I said about quantum phenomena.

\paragraph{4 \up{th} stage} $ $
\smallskip

Space-time is modeled by a manifold $\M$ of dimension $n > 5$,
provided with a pseudo-Riemannian tensor $g$ defined almost everywhere on $\M$
(a comment on this choice is given in the annex \ref{a3.8}).
\textbf{No physical object is "added" in this space-time}.
The "physical objects" studied, which correspond to the usual notions,
are only characteristics of the geometry of the pseudo-Riemannian manifold $(\M, g)$.
\textbf{No law, no principle is postulated}. The equations linking physical objects
(which are only defined from the geometry of $(\M, g)$) are just results
given by standard mathematical theorems on pseudo-Riemannian manifolds (often
consequences of Bianchi's identities in the case of "non-quantum" physics).
The pseudo-Riemannian manifold $(\M, g)$ is therefore supposed to be "totally anarchic".
Making physics in this "totally anarchic space-time" comes down to note
that some domains have particular geometric characteristics. These particular characteristics
define "objects" (which only make sense in this type of domain)
and the equations linking these objects are then mathematical consequences of their definitions.

Let's give a simple first example to clarify what we just said.
In this example, the reader may consider that $(\M, g)$ is the usual space-time
of the general relativity for which $dim \M = 4$ and $g$ is of signature $(-, +, +, +)$,
although it will not be the case later.
Suppose that at each point of a domain $\D$ of $\M$, the Einstein curvature (considered
as an endomorphism of the tangent space) admits a negative eigenvalue whose eigenspace is of dimension $1$ and timelike. Suppose also that it vanishes on the
subspace $g$-orthogonal to this eigenspace. We can then define, without ambiguities,
on this domain:
\begin{itemize}
 \item A function $\mu: \D \rightarrow \R$ such that, for each point $x$ of $\D$, $\mu(x)$ is
 the absolute value of the eigenvalue of the Einstein curvature.
 \item A vector field $X$, such that, for each point $x$ of $\D$, $X_x$ is the unit vector (ie. which satisfies $g(X, X) =
\text{-}1$)
 in the time-orientation of the eigenspace of dimension 1.
\end{itemize}
Such a pair $(\D, g)$ will be called a "fluid without pressure type domain" because, if one
reuses the language of the usual physics, the function $\mu$ will, by definition, be the
\textbf{energy density function} of the fluid. The flow of the vector field $X$ will be the
\textbf{fluid flow}. The choice of the nullity of the Einstein curvature on the subspace
$g$-orthogonal to the eigenspace will express the nullity of the pressure.

The reader can then verify (after some calculations) that by applying
Bianchi's second identity on Einstein curvature, one recover standard equations
of general relativity on fluids without pressure, which appear 
here as a simple mathematical consequence of the given definitions. (Of course,
this example will be repeated later in a more general context).

Here is a second example that describes an essential ingredient in the description of some
quantum phenomena. Here, the dimension of $\M$ is necessarily $> 5$ and the signature
of $g$ particular: we will say that a domain $\D$ of $(\M, g)$ is of type "oscillating metric
in a neutral potential" if the pseudo-Riemannian metric $g$ is conform to a metric
$g_0$ (ie. of the form $g = f g_0$ where $f: \D \rightarrow \R^+$) and is of
constant scalar curvature equal to that of $g_0$. The pseudo-Riemannian metric $g_0$
is a "reference" metric, the one with which the measurements are made,
 chosen so that, induced in the "apparent space" (of dimension $4$),
it coincides with the Minkovsky's metric (although we can generalize this).
This domain (slightly modified by "singularities") will represent what we call
in the usual language "particles in a vacuum". Contrary
to the first example, it is difficult at first sight to see the existing link
between such a domain $(\D, g)$ and the usual notion of particles,
and actually the theory that we are going to present will move away from
standard quantum theories on particles.

All this will be detailed in chapter 2 of this paper.
 We will show in particular that with this new point of view, we recover
 qualitative and quantitative description of classical quantum physics experiments by obtaining Klein-Gordon equations. Klein-Gordon equations 
 will give, in approximation, the standard Schrödinger equations that describe
the behavior of particles in a vacuum or in a potential.

The physics that I am going to present will be reduced to the search
for space-time domains that will be of (geometric) type "humanly
interesting".

As I said before, a \textbf{type} on a domain is a 
\textbf{geometric} constraint assumed on this domain (I just gave two examples).
It will be "humanly interesting" if it is "sufficiently deterministic".
This last notion can be mathematically defined in the following way
(it deserves more precise definitions than the one I am going to
give, but this is not the purpose of this paper (see Annex \ref{a3.9})):
a domain $(\Omega, g)$ is "sufficiently deterministic" if the knowledge
of some geometric quantities on a subdomain $\Omega'$ of $\Omega$
leads to the knowledge of these same quantities on~ $\Omega$.

In pseudo-Riemannian geometry, the results that make it possible to obtain
properties of this kind are called \textbf{rigidity theorems}.
In practice, this amounts to showing that the geometric constraints that define
the type give equations on the geometric quantities concerned that have
unique solutions for specified "initial conditions" (which are
the data of these quantities on $\Omega'$). We must not forget that
$(\Omega, g)$ contains the usual notion of "time" and that such results
of rigidity mean (in the common language) that identical initial conditions 
always give the same evolution in time, which justifies the "humanly interesting" aspect
of these domains.

It is important to note that the domains (with geometric constraints) that I am going to define are 
 generaly only \textbf{approximations}. If we go back to the first example,
we can consider that there are actually no domains that are exactly
"fluid without pressure" but that, under some circumstances, these
domains are a good approximation (used for example at very large
scale to describe the expansion of the universe in a space-time domain
containing a singularity of type "big-bang"). The second example,
which defines a domain of type "oscillating metric" can also,
of course, be considered as an approximation (the domain $\Omega$ of this
example can possibly be seen as a subdomain of the first example,
but in this case, the pseudo-Riemannian tensor chosen corresponds to a
"local approximation" totally different from the "global approximation"
of the first example).

In fact, the domains that I'm going to define naturally separate
in two classes that come from the fact that the experiments separate into
two distinct categories (we can also read annex \ref{a3.10+}):
\begin{itemize}
 \item Those whose \textbf{measurements} of the quantities concerned do not modify
 in no way (or have a negligible influence on) their unfolding.
 These experiments are commonly described by classical physics.
 \item Those whose \textbf{measurements} fundamentally change the course
 experiment. In this case, the description of the experiment is more subtel
 because it must include the measurement process itself. I will use the terminology
 "quantum phenomenon" when we look at this type of experiment in
  chapter \ref{part:deux}.
\end{itemize}

The domains that will correspond to the first category and describe the
"classical physics" will often be defined from the Einstein curvature: $(Ricc - \frac{1}{2}S g)$ where $Ricc$ is the Ricci curvature,
$S$ the scalar curvature and $g$ the pseudo-Riemannian tensor. This is
 natural because the second Bianchi identity simply says that the divergence of
Einstein curvature is nil. When this property of zero divergence
can be "transmitted" to a canonically defined vector field from
Einstein curvature (for a domain choice), the Stokes theorem
will give a \textbf{conservation law} for the magnitude associated to the flux of this vectors field, and this is an important feature that will allow to consider
that this domain is sufficiently deterministic. It's because of Einstein's curvature's divergence's nullity that the Einstein's curvature will appear
more naturally than the Ricci curvature in the definitions
presented in chapter \ref{part:un}.

For the domains that will correspond to the second category, 
it will be important to know precisely the pseudo-Riemannian tensor $g$
himself (and not just Einstein's curvature). First we will limit ourselves to the domains whose determination of the tensor $g$ comes 
solving down \textbf{linear} differential equations (chapter \ref{part:deux}),
which will be sufficient to make the connection with Klein-Gordon or
Schrödinger equations of standard quantum physics.

For these same domains, which correspond to the second category and describe
quantum phenomena, will intervene the subset of space-time $\M$
on which the pseudo-Riemannian tensor $g$ is not defined. Parts of this subset
(which will be sub-manifolds of $\M$ of dimension $< n$ and therefore of zero measure) will be
called \textbf{singularities} of the tensor $g$. (singularities of $g$ should not be confused
 with the "singularities" which correspond to the notions of "big-bang",
"black holes", etc., which are not considered subsets of $\M$).

No law will manage these $g$-singularities (space time is totally "anarchic"),
they will give "indeterminism" in quantum phenomena
(in contrast, the space-time in the neighborhood of these singularities will have a precise description). They will introduce properties of "localization"
for the domains corresponding to the notion of "particles".

\bigskip

The manifold $\M$ representing space-time will be of "big dimension"
($n$ probably $\geqslant 10$). The domains "sufficiently deterministic"
which will allow to recover usual physics results,
will be locally diffeomorphic to $\Theta \times S^1 \times W$ where $\Theta$ is an open set of
$\R^4$, $S^1$ is the circle and $W$ is a compact manifold. Moreover, we will see 
that the spin phenomena can be precisely described if we break down
$W$ into $S^3 \times V$ where $S^3$ is the standard sphere of dimension $3$.
The pseudo-Riemannian tensor $g$ will have almost everywhere a signature $(-, +, +, +, -, +, \cdots, +)$
which will give a "double" notion of time linked to the two $(-)$ signs 
(placed arbitrarily in first and fifth position). In many cases,
(and when $g$ is transported locally on $\Theta \times S^1 \times W$) this double time 
will be parameterized by $(t, u) \in \R \times S^1$, "$t$" can be assimilated to the usual time
 and "$u$" will be a completely new notion of time. We will see that all
the notions associated to electromagnetism will come from this new notion of time
parameterized by $u \in S^1$ (considered here as the fifth dimension) and this
whether for the domains corresponding to classical physics
or those describing quantum phenomena. In "quantum phenomena",
the concept of \textbf{mass} will be introduced as a \textbf{frequency} associated
to the notion of usual time $t \in \R$, and the concept of electric charge
as a \textbf{frequency} associated to the new notion of time $u \in S^1$
(which will force the electric charge to be an integer multiple of an 
elementary electric charge). It is this choice which will make it possible to recover, in approximation,
the results of standard quantum physics obtained from the Schrödinger equations that describe the behavior of "particles in a potential".
The fact that the signature of $g$ has exactly two signs $(-)$ will prevail in
this study (chapter 2, section \ref{s2.12}). Dimensions beyond 5
(which, locally, concern $W$) will be essential for the quantum phenomena 
(but they will also be very important in Chapter 1 (section \ref{s1.4}),
especially for "potential" type domains).\\
In fact, the need for theses choices of "dimension" and "signature" appears little by little
following the progress of the research starting from an entirely geometrical presentation
of physics: in the "classical physics" category presented in chapter 1,
the dimension $4$ of the space-time would have been sufficient for the domains representing the fluids
without electromagnetism (nothing new). For electromagnetism
(non-quantum) it has proved necessary, with our view on physics, that
the dimension of the space-time is at least 5. This point has been detailed in \cite{steph-1},
but in this article, the fifth dimension had been chosen "spacelike",
whereas it is now "timelike". If we had stayed there,
this choice of signature would technically not have had a great importance
(only a few signs change for the concept of electric charge in the obtained equations). It is in the description of potential areas and especially
in the study of quantum phenomena that the signature of $g$ has been imposed 
and that the dimension of space-time had to be increased.

With our point of view on physics, the study of notions that, 
in the usual language, correspond to "particles and their interactions",
will essentially be based on the \textbf{spectral theory associated to the compact manifold}
($S^1 \times W, g| _{S^1 \times W}$). This one is difficult and is determined
by the precise form of $(W, g|	_{W})$. We will stop in this paper
to the decomposition of $W$ in the form $W = S^3 \times V$ where $S^3$ is the standard sphere of dimension $3$, which will allow to describe the phenomena associated to
the notion of "spin" and, in particular, to address the phenomena of "quantum entanglement" and the "anomalous magnetic moment of the electron".

Whatever the domains studied in chapter 1 or 2, the dimension of space-time will be
considered equal to the same $n$ (although in some cases dimensions will be "neglected" but
this with a precise mathematical definition of course). \textbf{So we will have in particular
a perfectly unitary presentation of "classical physics" and
"quantum phenomena"}. \\
A summary of the theory presented in this paper is given in the annex \ref{a3.10+} and can help throughout the reading.

\newpage
\section*{Notations}
\addtotdm{Notations}

Notations are mainly those used in \cite{gourg} which sometimes
differ from those in \cite{vaugon-1}, \cite{vaugon-2} and \cite{vaugon-3}.

\noindent Let $(\M, g)$ be a pseudo-Riemannian manifold. We denote by:
\begin{center}
\begin{tabular} {lll}
 $Ricc_g$ &: & the Ricci curvature of $g$ \\
 $S_g$ &: & the scalar curvature of $g$ \\
 $Ein_g:= Ricc_g - \frac{1}{2}S_g g$ &: & the Einstein curvature of $g$ \\
 $G:= 2 Ein_g$ &: & twice the Einstein curvature (this will simplify \\
 & & formulas in many situations) \\
 $D$ &: & the covariant derivative associated with $g$, and when $T$ is \\
 & & a tensor field on $\M$ with coordinates $T_{(k)}^{(j)}$, \\
 & & the coordinates of $D T$ are denoted $\nabla_i T_{(k)}^{(j)}$.
\end{tabular}
\end{center}

\begin{itemize}
\item When $X$ is a vector field on $\M$, one denotes $D_X$ the covariant derivative associated with
this vector field.
\item When $f: \M \rightarrow \R$ is a function, $\nabla_g f$ or sometimes $\overrightarrow{grad f}$
refers to the \textbf{gradient} of $f$ ~~~ ($(\nabla_g f)^j = g^{i j}\partial_i f$).
\item $\nabla_g\cdotp X$ refers to the \textbf{divergence} of the vector field $X$ ~~~ ($\nabla_g\cdotp X
= \nabla_i
X^i$).
\item More generally, $\nabla_g\cdotp T$ refers to the \textbf{divergence} of the tensor field $T$ (for the
first index). Example: $(\nabla_g\cdotp T)^j = \nabla_i T^{i j}$ if $T \in \otimes^{(**)}\M$.
\item When $T$ is a totally covariant tensor (of coordinates $T_{i j k ...}$) and $g$ is a non degenerate symmetric bilinear form on a vector space $E$, its contravariant version (by $g$) is denoted $T^\sharp$.
\item When $T$ is totally contravariant ($T^{i j k ...}$), its covariant version (by $g$) is denoted $T^\flat$.
\item In the particular case where $B$ is a bilinear form (of coordinates $B_{i j}$) one denotes $\leftidx{^e}{B}$
the endomorphism of $E$ associated with $B$ (by $g$) (coordinates $B^i_{\hspace*{0.7mm}j} = g^{i k}B_{k j}$).
\item \textbf{Einstein's summation convention is used throughout this paper}
\end{itemize}

\newpage
\section*{Mathematical preliminary}

In all that follows, the universe is represented by a differential manifold $\M$ of dimension $n > 4$ provided with
a pseudo-Riemannian tensor $g$ (defined except on a zero measure subset). The tensor $g$ will have, in the domain
studied here, the signature $(-, +, +, +, -, +,\dots, +)$ almost everywhere.

\subsection*{Observation atlases}

The manifolds considered here are locally diffeomorphic to a product $\Theta \times K$ where
$\Theta$ is an open set of $\R^p$ and $K$ a compact manifold. It is useful to introduce a language
which takes into account this specificity. We will start by establishing a list of definitions that will specify the
terminology used throughout this paper. Part \textbf{A} of this preliminary will only deal with the
differential structure, Part \textbf{B} will concern the pseudo-Riemannian structure of the manifold.
\vspace{1mm}

\noindent \textbf{A.} Consider a $n$-dimensional manifold $\M$ of class $C^k$ where $k$ is large enough
to pose no problem regarding the objects defined next.

\begin{dfn}
 Let $K$ a compact manifold, $\Theta$ an open set of $\R^p$ and $\V$ an open set of $\M$.
 
- A $C^k$-diffeomorphism $\varphi: \V \rightarrow \Theta \times K$ will be called an \textbf{observation diffeomorphism}.

-A couple $(\V,\varphi)$, ~ where ~ $\varphi: \V \rightarrow \Theta \times K$ is an observation diffeomorphism, will be
called \textbf{a chart on $\M$ with values in $\Theta \times K$}.
\end{dfn}

\begin{dfn}
Let $\D$ an open set of $\M$ and $K$ a compact manifold. \\
A \textbf{$K$-observation atlas} on $\D$ is a family of charts $(\V_i,\varphi_i)_{i\in I}$ on $\D$ with values
in $\Theta \times K$ wich satisfies the following two properties:
\begin{enumerate}
 \item $\bigcup_{i\in I} \V_i = \D$
 \item \label{cond:2} for any $(i, j) \in I^2$, for any $x \in \V_i \cap \V_j~~ \varphi_i^{\text{-}1}(\{x_i^1\} \times K)
=
\varphi_j^{\text{-}1}(\{x_j^1\} \times K)$ where $x_i^1$ refers to the component of $\varphi_i(x)$ on $\Theta_i$.
\end{enumerate}
\end{dfn}

The property \ref{cond:2} is important, it will define unambiguously, in each point $x$ of $\D$,
a submanifold diffeomorphic to $K$ (associated to the data of an observation atlas).

\begin{prop}
\label{p0,1++}
 Let $\D$ a $n$-dimensional manifold and $K$ a compact $m$-dimensional manifold.
 There is equivalence between the following two assertions:
 \begin{enumerate}
  \item \label{hyp:1} There is a $K$-observation atlas on $\D$.
  \item There is a manifold $B$ of dimension $(n-m)$ and a submersion $\pi:\D \rightarrow B$ such that the triplet
  $(\D, B, \pi)$ is a fiber bundel with fiber $K$ and basis $B$.
 \end{enumerate}
(Proof is left to the reader).
\end{prop}

\begin{rmq}
 The compactness of $K$ is not necessary to obtain this proposition but it guarantees that the manifold $B$,
 constructed from the assumption \ref{hyp:1}, has a separate topology. If $K$ is not compact, just
assume further in the assertion \ref{cond:2} that the $\varphi_i^{\text{-}1}(\{x_i^1\} \times K)$ are subsets
closed in $\D$ (they are in $\V_i$) so that $B$ has a separate topology.
\end{rmq}

This proposition shows that what follows
 could be presented in the language of "bundles", but I think that this choice would not be very
natural since the "basis" of the bundles would not intervene (the basis would be of interest only
in very special cases).

The presentation chosen, in terms of "observation atlas" makes it possible to generalize the definitions
in the following manner. It is assumed that the manifold $K = K_1 \times  K_2$ where $K_1$ and $K_2$ are compact
(The reader will adapt what will follow in case $K = K_1 \times K_2 \times ... \times K_l$).

\begin{dfn}
  A \textbf{$K_1$-observation atlas} on $\D$ is a family of charts $(\V_i, \varphi_i)_{i\in I}$ on $\D$ with
values in $ \Theta_i \times K_1 \times K_2$ which satisfies the 
following two properties:
\begin{enumerate}
 \item $\bigcup_{i \in I} \V_i = \D$
 \item for any $(i, j) \in I^2$, for any $x \in \V_i \cap \V_j~~ \varphi_i^1(\{x_i^1\} \times K_1 \times \{x_i^3\}) =
\varphi_j^1(\{x_j^1\} \times K_1 \times \{x_j^3\})$ where $x_i^3$ is the component of $\varphi_i(x)$ on $K_2$.
\end{enumerate}
\end{dfn}

\begin{dfn}
 When $K_1$ is an oriented manifold, we will say that the $K_1$-observation atlas \\
 \textbf{keeps the orientation of
$K_1$} if: for any $(i, j) \in I^2$, for any $x \in \V_i \cap \V_j$, the orientation of the submanifold
$\varphi_i^{\text{-}1}(\{x_i^1\} \times K_1 \times \{x_i^3\})$ transported from that of $K_1$ by $\varphi_i$, is the
same as that transported by $\varphi_j$.
\end{dfn}

Of course, we write the same definitions for $K_2$.

\begin{dfn}
 A \textbf{$K_1$-$K_2$-observation atlas} on $\D$ is both a $K_1$-observation atlas and\\
 a $K_2$-observation atlas on $\D$.
\end{dfn}

\noindent \textbf{B.} We now consider a pseudo-Riemannian manifold $(\M, g)$ of dimension $n$ and class
$C^k$. Throughout this paper, the domains $\D$ of $\M$ will be locally diffeomorphic to $\Theta \times K$
where $\Theta$ is an open set of $\R^4$ and $K = S^1 \times W$, $S^1$ refers to the standard oriented circle and $W$ is a compact manifold. We therefore take the notation of the previous part \textbf{A} with $K_1 = S^1$ and $K_2 = W$.

Let $\D$ a domain of $\M$ and $\A$ a $S^1$-$W$-observation atlas on $\D$ which retains the orientation of $S^1$.
At each point $x$ of $\D$, the $1$-dimensional submanifold denoted $S_x^1$ (diffeomorphic to $S^1$) is then defined
by stating: \[S_x^1 = \varphi_i^{\text{-}1}(\{x_i^1\} \times S^1 \times \{x_i^3\})\] where $(\V_i, \varphi_i)$ is an observation 
diffeomorphism at $x$ of the atlas $\A$ and noting that $S_x^1$ does not depend on the choice of this observation diffeomorphism. The submanifold $S_x^1$ is, moreover, unambiguously oriented.

We define in the same way the submanifold of dimension $n-5$ denoted $W_x$ (diffeomorphic to $W$): \[W_x =
\varphi_i^{\text{-}1}(\{x_i^1\} \times \{x_i^2\} \times W)\]
Note that $x' \in S_x^1 \Leftrightarrow S_{x'}^1 = S_x^1$ and $x' \in W_x \Leftrightarrow W_{x'} =
W_x$.

\bigskip

The tangent space $T_x(\D)$ decomposes uniquely in the form: \[H_x \obot (T_x(S_x^1) \oplus T_x(W_x))\]
where $H_x$ is the $4$-dimensional subspace $g$-orthogonal to $T_x(S_x^1) \oplus T_x(W_x)$ \\ ($T_x(S_x^1)$ is not
assumed $g$-orthogonal to $T_x(W_x)$).

In the following $H_x$ will be called \textbf{the apparent vector space} at $x$ (associated to the observation atlas).

The field of spaces $H_x$ has no reason to be integrable, in other words, there is no reason that, through a point $x_0$, there passes a submanifold of dimension $4$ such that the tangent space at any point $x$ of this submanifold is
$H_x$.

\begin{dfn}
 A \textbf{$g$-observation atlas on $\D$} is a $S^1$-$W$-observation atlas on $\D$ that retains orientation
of $S^1$ and which furthermore satisfies the following properties associated to the pseudo-Riemannian tensor $g$:
\begin{enumerate}
 \item for any $x \in \D$ :
 \begin{itemize}
  \item the signature of $g|_{H_x}$ is $(-, +, +, +)$
  \item the signature of $g|_{T_x(S_x^1)}$ is $(-)$
  \item the signature of $g|_{T_x(W_x)}$ is $(+, \dots, +)$
 \end{itemize}
 \item \label{cond-local:2} for any $(i, j) \in I^2$, for any $x \in \V_i \cap \V_j~~~~
\varphi_i^*(\frac{\partial}{\partial t})_{\varphi_i(x)}$ ~ and ~ $\varphi_j^*(\frac{\partial}{\partial t})_{\varphi_j(x)}$
are timelike
and in the same time orientation (ie. $g(\varphi_i^*(\frac{\partial}{\partial t})_{\varphi_i(x)},
\varphi_j^*(\frac{\partial}{\partial t})_{\varphi_j(x)})<0$).
Here, $(\frac{\partial}{\partial t})_{\varphi_i(x)}$ refers to the tangent vector at $\varphi_i(x)$ associated
to the standard coordinate system $(t, x, y, z)$ of $\Theta_i \subset \R^4$ ~~ $(\varphi_i: \V_i \rightarrow \Theta_i
\times S^1 \times W)$.
\end{enumerate}
\end{dfn}

The property \ref{cond-local:2} defines a "classical" time orientation of each apparent space
$H_x$.

The following definitions use the classical process that introduces the concept of "complete" atlas (or
"saturated").

\begin{dfn}
 Two $g$-observation atlas on $\D$ are \textbf{equivalents} if their union is still a $g$-observation atlas.
\end{dfn}

\begin{dfn}
 Let $\A$ a $g$-observation atlas on $\D$, the \textbf{completed} (or the \textbf{saturated}) of $\A$ is the $g$-observation atlas formed from the union of all $g$-atlas which are equivalents to $\A$.
\end{dfn}

\begin{dfn}
 A $g$-observation atlas on $\D$ is \textbf{complete} (or \textbf{saturated}) if it is equal to its completed.
\end{dfn}

The domains that we will define later will be triplets $(\D, g, \A)$ where $\D$ is an open set of $\M$,
$g$ a pseudo-Riemannian tensor on $\D$ and $\A$ a \textbf{complete} $g$-observation atlas on $\D$. Types
of these domains will be specified by the data of \textbf{geometric constraints} imposed on $(\D, g)$.
The choice of the complete $g$-observation atlas will define the set of "observers" that one allows oneself for the
measurements of the quantities defined by the  domain (an "observer" is mathematically defined by an
observation diffeomorphism). As the atlas will be complete, it will leave a large choice of observer changes.
For example, if $(\V, \varphi)$ is a chart of this complete atlas $\A$ and if $\sigma_1: \Theta \rightarrow
\Theta' \subset \R^4$, $\sigma_2: S^1
\rightarrow S^1$, $\sigma_3: W \rightarrow W$ are three diffeomorphisms, then $(\V, \sigma_1 \times \sigma_2
\times \sigma_3 \circ \varphi)$ is a chart of this atlas (provided that
the "defined orientations" are retained).

\chapter{Non-Quantum Physics} \label{part:un}

\section{"Fluid" and "potential" type domains \label{s1.1}}

As we have just stated in the Mathematical preliminary, a "fluid" or a "potential" type domain is first a triplet $(\D, g, \A)$ where $\D$ is a domain of $\M$, $g$ a 
pseudo-Riemannian tensor and $\A$ a complete $g$-observation atlas on $\D$. We recall that at each point $x$ of $\D$
the tangent space $T_x(\D)$ decomposes in the form $H_x \obot (T_x(S_x^1) \oplus T_x(W_x))$.
The space $H_x$ (of dimension $4$) is the apparent space at the point $x$. The submanifold $S_x^1$ is diffeomorphic with
circle, oriented and timelike. Manifold $W_x$ (dimension $n-5$) is compact and space-like.

On $\D$, we define the vector field $Y$ taking for each $x$ of $\D$ the unique vector tangent to $S_x^1$,
in the orientation and such that $g_x(Y_x, Y_x) = -1$. \textbf{This vector field will be a fundamental object of
electromagnetism in $\D$}.


\section [Type "fluid" from the curvature of Einstein] {The domains of type "fluid" defined from the
Einstein Curvature \label{s1.2}}

The natural geometric constraint imposed on the domain so that it is of fluid type will be a
constraint given on the tensor $G:= 2 Ein_g$ which will allow to define \textbf{canonically} a vector field
$X_0$ on $\D$ timelike and such that $\forall x \in \D, ~~{X_0}_x \in H_x$. This vector field will be called
the \textbf{apparent vector field} of the fluid and the associated flow the \textbf{apparent flow} of the fluid.

This condition will also define canonically two functions $\mu$ and $\rho: \D \rightarrow \R$ which
will respectively represent the \textbf{energy density} and the \textbf{electric charge density}
fluid. The second Bianchi identity ($\nabla_g\cdotp  G = 0$) will then make it easy to obtain the conservation laws 
as well as the usual evolution equations, for a large class of fluids. This geometric constraint
 is the first stated in the definition \ref{def:4}. The second condition can be interpreted as
the fact that the quantum effects on electromagnetism are neglected (but this will not be apparent
before reading this paper completely). The following lemma shows that this second condition is only
the approximation which consists of a priori averaging the pseudo-Riemannian metric $g$ on the circles $S_x^1$
previously defined and thereby neglect the variations of $g$ on these circles. Of course, this second
condition will be dropped in chapter \ref{part:deux} because these are precisely the variations of $g$ on
"small compact manifolds" that describe quantum phenomena.

\begin{lem} \label{l1}
 We consider $Y$ the vector field tangent to the circles $S_x^1$ and normalized by $g(Y, Y)=-1$ defined
previously.
 We denote $\sigma$ the 1-parameter group of diffeomorphisms associated with $Y$. We define the "averaged" pseudo-Riemannian metric $\overline{g}$ by setting: \[\forall x \in \M,~~ \overline{g}_x = \frac{1}{l_x}
\int_{t_0}^{t_0 + l_x}\left(\sigma^*(t) g\right)_x dt\]
where $l_x$ is the length of the circle $S_x^1$ associated to $g$ ($\overline{g}_x$ does not depend on the choice of $t_0$ because
$\sigma_x(.)$ is periodic whith period $l_x$).

This way, $\overline{g}(Y, Y) = -1$ and, $\forall s \in \R,~~ \sigma^*(s) \overline{g} = \overline{g}$. In other
words, $Y$ is a Killing field for $\overline{g}$.
\end{lem}

\begin{dem}
 Prove that $\overline{g}(Y, Y) = \text{-}1$ is immediate, the rest is summarized by the following equalities \nolinebreak:
 \begin{displaymath}
 \begin{tabular} {rcl}
  $\displaystyle(\sigma^*(s)\overline{g})_x$
  & = & $\displaystyle\frac{1}{l_x}\int_{t_0}^{t_0+l_x}\sigma^*(s)(\sigma^*(t)g)_x~ dt$ \\
  & = & $\displaystyle\frac{1}{l_x}\int_{t_0}^{t_0+l_x}(\sigma^*(t+s)g)_x~ dt$ \\
  & = & $\displaystyle\frac{1}{l_x}\int_{t_0+s}^{t_0+s+l_x}(\sigma^*(t) g)_x~ dt$ \\[1.3em]
  & = & $\displaystyle\overline{g}_x$
 \end{tabular}
 \end{displaymath}
\end{dem}

\begin{dfn}
 \label{def:4}
 \textbf{A fluid type domain} is a triplet $(\D, g, \A)$ for which the following two conditions are
satisfied:
\begin{enumerate}
 \item \label{ei11} for any $x$ in $\D$ the endomorphism $~~\leftidx{^e}{G}|_{H_x}$ ~ admits an eigenspace $E_{-\mu}$ of dimension 1, 
timelike and whith eigenvalue $-\mu < 0$ ~~ (so that $\mu > 0$). ~~ (Here, $\leftidx{^e}{G}|_{H_x}$ is the endomorphism on
$H_x$ defined by: $\forall X \in H_x$ ~~ $\leftidx{^e}{G}|_{H_x}(X) = pr_{H_x} \leftidx{^e}{G}(X)$).
\item \label{ei12} The vector field $Y$ is a Killing field, in other words, the local diffeomorphisms generated
by the field $Y$ are $g$-isometries (see lemma \ref{l1}).
\end{enumerate}
\end{dfn}

\begin{rmq}
 \label{r2}
 When the condition \ref{ei11} is satisfied for a metric $g$, it remains satisfied for any metric in a
"neighborhood" of $g$ and its unique objective is, as already stated, to allow the canonical definition of the vector
field $X_0$ and the two functions $\mu$ and $\rho$. This is fundamentally different from the axiomatic Lagrangian principle  because the latter is characterized by the precise axiomatic data of "equalities" which define the
Lagrangian. It is recalled that the condition $\ref{ei12}$ has the unique purpose of neglecting the quantum effects associated to electromagnetism.
\end{rmq}


\subsection{The physical objects canonically defined in a fluid domain \label{ss1.1}}

\begin{enumerate}
 \item \textbf{The vector field $Y$} (already introduced) defined for any $x \in \D$ as the sole tangent vector at $x$
to $S^1_x$, in orientation, and such that $g(Y_x, Y_x) = -1$.
  \begin{itemize}
   \item \textbf{The associated 1-form} $Y^\flat$ (where $Y^\flat_i:= g_{ij} Y^j$).
   \item \textbf{The differential 2-form} characterizing electromagnetism, classically denoted $F$, is defined here
by $F= d (Y^\flat)$.
  \end{itemize}
 \item \textbf{The vector field $X_0$} (already presented) defined $\forall x \in \D$ as the sole vector ${X_0}_x$ of
the eigenspace $E_{-\mu} \subset H_x$, in orientation, and such that $g({X_0}_x, {X_0}_x) = -1$ (note that,
given the signature of $g|_{H_x}$, this timelike eigenspace is unique (see annex \ref{a3.1})). The
vector field $X_0$ will be called \textbf{the apparent vector field} of the fluid and the associated flow, \textbf{the 
apparent flow} of the fluid.
 \item \textbf{The function} $\mu: \D \rightarrow \R^+$ defined $\forall x \in \D$ by $\mu(x) = \mu_x$ where $-\mu_x$
is the eigenvalue associated with the eigenspace $E_{-\mu_x}$. This function is the \textbf{energy density function} of the 
fluid.
 \item \textbf{The function $\rho: \D \rightarrow \R$} defined $\forall x \in \D$ by $\rho(x) = G_x({X_0}_x,Y_x)$. This
function is \textbf{the electric charge density function} of the fluid.
 \item \textbf{The vector field $X:= X_0 + \frac{\rho}{\mu} Y$}, timelike, is the \textbf{vector field
of the fluid} ($X_0$ was only the apparent vector field of the fluid) and the associated flow, \textbf{the fluid flow}. (Of course,
$X = X_0$ if the electric charge density $\rho$ is zero).
 \item For all $x \in \D$, the \textbf{temporal space $\mathcal{T}_x$}, of dimension $2$, is the vector subspace
 of $T_x(\D)$ generated by ${X_0}_x$ and $Y_x$. The distribution (the field of temporal spaces) $\mathcal{T}$ is integrable since the
condition \ref{ei12} on the domains of type fluid has the consequence that $[X_0\ Y] = 0$ (cf.
annex \ref{a3.1}). The \textbf{temporal tube} $\tau_x$ is therefore defined as the integral submanifold of the distribution $\mathcal{T}$ passing through
$x$. The temporal tube $\tau_x$ is a submanifold of dimension $2$ \textbf{totally 
timelike} in the sense that all its tangent vectors are timelike. These tubes can be seen as the
generalization of fluid flow lines. They are oriented by the orientation given for $X_0$ and $Y$.
 \item In view of the above definitions, $\forall x \in \D$, the tensor $G_x$ induced in $\mathcal{T}_x$
is written in the form: $G|_{\mathcal{T}_x} = \mu(x) X^\flat \otimes X^\flat + \sigma(x) Y^\flat \otimes Y^\flat$ where $\sigma:\D
\rightarrow \R$ is a regular function. The field of bilinear forms $P = G - G|_{\mathcal{T}_x}$ will be called the
\textbf{fluid pressure} (it satisfies in particular: $P(X_0, X_0) = P(X_0, Y) = P(Y, Y) = P(X, X) = 0$).
The tensor field $G$ is written in the form:
 \begin{displaymath}
 \begin{tabular} {rcl}
  $\displaystyle G$
  & = & $\displaystyle \mu X^\flat \otimes X^\flat + \sigma Y^\flat \times Y^\flat + P$ \\
  & = & $\displaystyle \mu X^\flat_0 \otimes X^\flat_0 + \rho(X^\flat_0
  \otimes Y^\flat + Y^\flat \otimes X^\flat_0) + (\sigma + \frac{\rho^2}{\mu}) Y^\flat \otimes Y^\flat + P$
 \end{tabular}
 \end{displaymath}
\textbf{The apparent pressure} ${P_A}_x$ is the pressure $P$ induced in the apparent space $H_x$, that is, by
definition: $\forall (Z, Z') \in H^2_x$ ~~~~ $ {P_A}_x(Z, Z') = P_x(Z, Z')$, \\
$\forall Z \in T_x(\D) ~~\forall Z' \in T_x(S^1_x) \oplus T_x(W_x) ~~~~{P_A}_x(Z, Z') = 0$ and
${P_A}_x(Z, X_0) = 0$. \\
\textbf{The hidden pressure} is defined by ${P_H}_x:= P_x - {P_A}_x$.
The tensor field $G$ is also written in the form: $G = \mu X^\flat \otimes X^\flat + \sigma Y^\flat \otimes
Y^\flat + P_A + P_H$.
\end{enumerate}

\textbf{These objects are defined using only the geometry of the chosen type, ie. the
only pseudo-Riemannian tensor $g$ on $\D$. They are not independent of each other. Dependency links are
given by the sole mathematical properties of the pseudo-Riemannian manifolds (in particular Bianchi's second identity). As we   will verify, these dependencies are none other than the "physics laws" on the
fluids, in particular those of the standard general relativity. Here, no law, no principle is
added. The equations we are going to write are simple conclusions of the definitions we have just given.}


\subsection{General Equations for Fluids}

\begin{thme}
 \label{t1.1}
 In a fluid type domain, the following equalities are satisfied:
 \begin{enumerate}
  \item (generic of energy conservation)
  \[\nabla_g\cdotp  (\mu X) = \nabla_g\cdotp  (\mu X_0) = g(X_0, \nabla_g\cdotp  P)\]
  in coordinates:
  \[\nabla_i(\mu X^i) = \nabla_i(\mu X^i_0) = X_{0 j}\nabla_i P^{i j}\]
  Moreover:
  \[\mu^2X(\frac{\rho}{\mu}) = \mu g(Y, \nabla_g\cdotp  P) - \rho g(X_0, \nabla_g\cdotp  P)\]
  and:
  \[g(Y, \nabla_g\cdotp  P) = \nabla_g\cdotp  \leftidx{^e}{P}(Y)\]
  \item (generic of the conservation of the electric charge)
  \[\nabla_g\cdotp  (\rho X) = \nabla_g\cdotp  (\rho X_0) = \nabla_g\cdotp  (\leftidx{^e}{P}(Y)) = g(Y,
\nabla_g\cdotp  P)\]
  in coordinates:
  \[\nabla_i(\rho X^i) = \nabla_i(\rho X^i_0) = Y_j\nabla_i P^{i j}\]
  (Note: $\leftidx{^e}{P}(Y) = \leftidx{^e}{P_H}(Y)$).
 \vspace{2cm}
 
  \item (motion equations)
  \begin{enumerate}
   \item \[\mu D_x X = - \nabla_g\cdotp  P - g(X_0, \nabla_g\cdotp  P) X\]
   in coordinates:
   \[\mu X^i\nabla_iX^j =
-\nabla_iP^{i j} - {X_0}_j \nabla_i P^{i j}X^j\].
   \item \[\mu D_{X_0} X_0 = \rho \leftidx{^e}{F}(X_0) - pr_{\mathcal{T}^\bot}(\nabla_g\cdotp  P)\]
  \end{enumerate}
 \item (generic of Maxwell's second equation, the first is obvious since $F:= d(Y^\flat)$)\\
  \[\nabla_g\cdotp  F = \rho X_0 + \frac{1}{2}(F_{i j}F^{i j})Y - \leftidx{^e}{P}(Y)\]
  in coordinates:
  \[\nabla_iF^{i j} = \rho X^j_0 + \frac{1}{2}(F_{k l}F^{k l})Y^j - P^{i j} Y_i\]
 
 \end{enumerate}
\end{thme}

The proof of this theorem is quickly obtained using the pseudo-Riemannian manifolds standard properties (especially the second Bianchi identity) and is given in Annex \ref{a3.1}.

The equations given by the theorem \ref{t1.1} generalize the equations obtained in standard general relativity on
fluids. They are not "sufficiently deterministic" (too many unknowns relatively to the number of equations) and
are useful in this form only for the study of global behaviors on fluids. To be used 
more precisely, it is necessary to impose additional conditions on the geometry of the fluid type 
given. This must, of course, be considered as different approximation choices depending on the considered case
 (see Annex \ref{a3.9} for "approximations").

\subsection{Particular fluids}

The specific geometrical conditions that we are going to require now are intended to recover the exact form of the
known equations on fluids in standard general relativity written in dimension $4$. To us, the equations are written in dimension $n$, hense it's in the projection on the apparent spaces $H_x$ that the
comparison will be made. In fact, our goal is rather to show that the fluids considered in standard relativity
 cases are special cases of those just defined and will be discussed below. "Fluids" and "potentials" specific to the dimension $n > 4$ are particulary interesting. (For more details on the link between the mathematical model and the reality one can read annex \ref{a3.10+}).

\begin{dfn} \label{def:5}
~ \\ \vspace*{-1em}
\begin{enumerate}
 \item
 A \textbf{"perfect fluid"} type domain is a fluid type domain for which $\leftidx{^e}{P}(Y)=0$
~~ (actually, $\leftidx{^e}{P}(Y)=\leftidx{^e}{P}_H(Y)$ ~ where ~ $P_H$ ~ is the hidden pressure (defs. \ref{ss1.1.1}))
\item
A domain of type \textbf{"perfect isentropic fluid"} is a "perfect fluid" type domain for which,
at each point $x$ of $\D$, the pressure tensor $P_x$ is proportional to the tensor $(g_x-g_x|_{\mathcal{T}_x})$
~ (no space direction is preferred), ie:

$P=\tilde p(g+X_0^\flat\otimes X_0^\flat+Y^\flat\otimes Y^\flat)$

where $\tilde p:\D\rightarrow \R$ will be called \textbf{the pressure function}.
\item
A domain of type \textbf{"really-perfect fluid"} (or type "dust"), possibly charged
electrically, is a "perfect fluid" type domain for which ~ $\nabla_g\cdotp  P=0$.

When considering a domain of type \textbf{"truly-perfect fluid without electromagnetism"} it is assumed
more than $\rho=0$ and $F=0$.
\end{enumerate}
\end{dfn}

\begin{rmq} \label{r3}
 In the above definitions, the hypotheses of nullity of the tensors $\leftidx{^e}{P}(Y)$ and
$\nabla_g\cdotp  P$ can be replaced by the fact that these tensors are negligible compared to those in the following equations (specifying of course the notion of "negligibility").
\end{rmq}
The following proposition translates the previous theorem to the particular cases of fluids that we have just
 define.

\begin{prop} \label{p1.1}
~ \\ \vspace*{-1em}
\begin{enumerate}
 \item In a "perfect fluid" type domain the following equalities are satisfied:
\begin{enumerate}
  \item (generic of energy conservation)
  \[\nabla_g\cdotp  (\mu X) = \nabla_g\cdotp  (\mu X_0) = g(X_0, \nabla_g\cdotp  P)\]
  in coordinates:
  \[\nabla_i(\mu X^i) = \nabla_i(\mu X^i_0) = X_{0 j}\nabla_i P^{i j}\]
  Moreover:
  \[\mu^2X(\frac{\rho}{\mu}) =\mu^2X_0(\frac{\rho}{\mu}) = - \rho g(X_0,\nabla_g\cdotp  P)\]
  \item (electric charge conservation)
  \[\nabla_g\cdotp  (\rho X) = \nabla_g\cdotp  (\rho X_0) = 0\]
    in coordinates:
  \[\nabla_i(\rho X^i) = \nabla_i(\rho X^i_0) =0\]
  \item (motion equation)
  
  \[\mu D_x X = - \nabla_g\cdotp  P - g(X_0, \nabla_g\cdotp  P) X\]
   in coordinates: \\
   \[\mu X^i\nabla_iX^j =-\nabla_iP^{i j} - {X_0}_j \nabla_i P^{i j}X^j\].
   
which is also written in the form:

\[\mu D_{X_0} X_0 =\rho^2F(X_0)- \nabla_g\cdotp  P - g(X_0, \nabla_g\cdotp  P) X_0\]
\item (generic of Maxwell's second equation)
\[\nabla_g\cdotp  F = \rho X_0 + \frac{1}{2}(F_{i j}F^{i j})Y\]
  in coordinates:
  \[\nabla_iF^{i j} = \rho X^j_0 + \frac{1}{2}(F_{k l}F^{k l})Y^j\]
  
  (in particular: $\forall x\in\D ~~~pr_{H_x}(\nabla_g\cdotp  F)=\rho X_{0_x}$)
  
   (These properties are immediate consequences of theorem \ref{t1.1}).
   \end{enumerate}
   \item
   In a "perfect isentropic fluid" type domain, the following equalities are satisfied:
   \begin{enumerate}
    \item (generic of energy conservation)
    \[\nabla_g\cdotp  (\mu X_0)+\tilde p\nabla_g\cdotp  X_0=0\]
    Moreover:
    \[\mu^2X(\frac{\rho}{\mu})=\mu^2X_0(\frac{\rho}{\mu})=\rho\tilde p\nabla_g\cdotp  X_0\]
    \item (electric charge conservation)
    \[\nabla_g\cdotp  (\rho X)=\nabla_g\cdotp  (\rho X_0)=0\]
    \item (motion equation)
    \[(\mu+\tilde p)D_{X_0}X_0=\rho\leftidx{^e}{F}(X_0)-\nabla_g\tilde p-X_0(\tilde p)X_0\]
    
    Here we recover exactly the equations obtained in standard general relativity for isentropic charged fluids
 (we can see for example \cite{hawking}).

(The demonstration of these properties is quickly obtained after having verified that
$\nabla_g\cdotp  P=\nabla_g\tilde p+X_0(\tilde p)X_0+\tilde p(\nabla_g\cdotp  X_0)X_0+\tilde pD_{X_0}X_0$ ~~ then
~~ $g(X_0,\nabla_g\cdotp  P)=-\tilde p\nabla_g\cdotp  X_0$)
   \end{enumerate}
\item
In a "really-perfect fluid" type domain the following equalities are satisfied:
\begin{enumerate}
 \item (energy conservation)
 \[\nabla_g\cdotp  (\mu X)=\nabla_g\cdotp (\mu X_0)=0\]
 
 Moreover:
 \[X(\frac{\rho}{\mu})=X_0(\frac{\rho}{\mu})=0\]
 \item (electric charge conservation)
 \[\nabla_g\cdotp (\rho X)=\nabla_g\cdotp (\rho X_0)=0\]
 \item (motion equation)
 \[D_XX=0\]
 which is also written:
 \[\mu D_{X_0}X_0=\rho\leftidx{^e}{F}(X_0)\]
 In this case \textbf{$X$ is a geodesic field, whether the fluid is electrically charged or not}. (Of course, if the
electric charge density $\rho$ is zero $X=X_0$ is also a geodesic field).
\end{enumerate}
\end{enumerate}
\end{prop}
\begin{rmq} \label{r4}
 In the motion equation just given, $(D_{X_0}X_0)_x$ and $\leftidx{^e}{F}(X_0)_x$ are
g-orthogonal to $Y_x$ but not necessarily to $W_x$. If we want to guarantee the fact that $(D_{X_0}X_0)_x$ and
$\leftidx{^e}{F}(X_0)_x$ belong to the apparent space $H_x$ one can add the following condition (*):

(*) ~ Submanifolds $W_x$ are parallel along the geodesic circles $S^1_x$. Precisely: $\forall x\in\D
~~\forall Z\in T_x(W_x) ~~\forall x'\in S^1_x$, ~~ the parallel transport of $Z$ to $x'$ along the geodesic circle
$S^1_x$ is tangent to $W_x'$.

We then quickly verify that under this condition: $\forall x\in\D ~~\leftidx{^e}{F}(X_0)_x\in H_x$.

Condition (*) will be satisfied in the examples that we present next.
\end{rmq}


\section [Type "potential" from Einstein curvature] {Potential type domains defined from
 the Einstein curvature  \label{s1.3}}
\begin{dfn} \label{def:6}
 A \textbf{potential type domain} is a triplet $(\D,g,\A)$ that satisfies the following properties
\begin{enumerate}
\item $\forall x\in \D ~~~G_x|_{H_x}=0$ ~~ and ~~ $pr_{H_x}\leftidx{^e}{G}(Y)=0$
\item The field $Y$ is a Killing field. (See lemma \ref{l1} and its presentation).
\end{enumerate}

\end{dfn}

These potential type domains therefore appear as fluid type domains for which the
energy density, electrical charge density and apparent pressure $P_A$ are zero. We will notice
on the other hand there are no more canonically defined vector fields as were $X_0$ and $X$ for
fluids. Only those that define electromagnetism, ie. the field $Y$ (and then the 2-form $F$), remain as canonical objects.
 The tensor $G$ is then assimilated to the hidden pressure $P_H$ which thus verifies $\nabla\cdotp P_H=0$.

The following theorem \ref{t1.2} is immediately obtained using the proof of theorem
\ref{t1.1}'s part 4 on Maxwell equations.

\begin{thme} \label{t1.2}
In a potential type domain given by the definition \ref{def:6} the following equality is satisfied:

(Maxwell's second equation, the first being obvious since $F:=dY^\flat$)

\[\nabla_g\cdotp F=\frac{1}{2}	(F_{ij}F^{ij})Y-\leftidx{^e}{G}(Y)\]
In particular (since $pr_{H_x}\leftidx{^e}{G}(Y)=0$): \[pr_{H_x}\nabla_g\cdotp F=O\]
\end{thme}

The potential type domains are in practice very important because the knowledge of their geodesics gives,
in approximation, the curves of elementary objects (electrically charged or not) placed in these
potential when we consider that the incidence of these objects is negligible on the geometry of the potential type
. Indeed, if we introduce an elementary object in a potential type domain, this one
can then be considered as a domain of type "really-perfect fluid" in the subset where the density
energy is not zero. It is also very "localized in space". If we consider it outside of this
very localized domain, the geometry of the potential type domain is not modified, whereas the field $X$ of the fluid
(defined when $\mu\neq0$) is a geodesic field according to proposition \ref{p1.1} (2.3.c). The curves of the flow 
determined by $X$, which give the trajectory of the elementary object, are (in approximation)
geodesics of the potential type. In addition the ratio $\frac{\rho}{\mu}$ is equivalent to the ratio of the 
electric charge by the mass of the elementary object considered since it is "restricted in space". This principle
is classical in standard general relativity for elementary objects \textbf{not electrically charged} and is
used, for example, to determine the trajectory of planets in a Schwarzschild domain, the deviation of the
light, etc. What is remarkable is that this principle also applies now for elementary objects
\textbf{electrically charged}, but necessarily in a space of dimension $n\geqslant5$. In this case, the 
apparent trajectory is determined by the apparent field $X_0$ (non-geodesic), itself deduced from the geodesic field
 $X$. Of course, $X=X_0$ when the electric charge is zero.

Specific examples of calculations will be presented at the end of next section.
\begin{rmq} \label{r5}
 When a domain $(\D,g)$ is isometric to a domain of the form $(\D'\times V,g'\times g_V)$ where $V$ is a compact manifold
of dimension $k$ and $g'\times g_V$ is a product metric, we can define on the couple $(\D',g')$ (for
which $dim\D'=n-k$) notions of "fluid type" or "potential type" in an identical manner to what has been
done on $(\D,g)$. The "objects" defined from the Ricci curvature give very close notions on $(\D,g)$ or $(\D',g')$ since $g'\times g_V$ is a product metric, but it is not the same for those
defined by the Einstein tensor (or $G$) because the scalar curvature of $g_V$ is important
and $S_g$ is different from $S_{g'}$ when $S_{g_V}$ is not zero ($S_g=S_{g'}+S_{g_V}$). But the important notions are defined
from Einstein's tensor. For example, $(\D',g')$, of dimension $(n-k)$, can be of fluid type while $(\D,g)$ is not, or $(\D',g')$ is of fluid type
without pressure and $(\D,g)$ is of fluid type with pressure (or vice versa), etc. In addition, the objects defined
on fluids are different wether they are considered on $(\D,g)$ or on $(\D',g')$. It may therefore be smart, in
 special cases for which $(\D,g)$ is isometric to $(\D'\times V,g'\times g_V)$, to use the fluids
or the potentials defined on $(\D',g')$ rather than on $(\D,g)$. It will be necessary, of course, in this case to specify that they
are fluids or potentials defined on $(\D',g')$ whose characteristics are possibly different from
those given by $(\D,g)$. Of course, the theorems obtained previously are done on $(\D',g')$ of dimension
$(n-k)$ as well as on $(\D,g)$ of dimension $n$.
\end{rmq}


\section [Examples of domains given by the tensor $\protect\petitg$] {Examples of potential and fluid type domains given by
the\\
pseudo-Riemannian tensor $g$ himself \label{s1.4}}

Equalities given in theorem \ref{t1.1} are written using the pseudo-Riemannian tensor $g$.
Two equalities that have the same writing can describe fluids whose behavior are very different if the
pseudo-Riemannian tensors underlying are different. In fact, theorem \ref{t1.1}, written without precisions on the
tensor $g$, can only be applied to global studies on fluids. In order to describe
precise behavior, it is important to know how to define "fluid" or "potential" type
from the pseudo-Riemannian tensor $g$. This process is used in standard general relativity and the
 considered domains are often presented under the name of "exact solutions of the Einstein equation". We
can mention: the solutions of Schwarzschild-Kruskal, Reissner-Nordström, Kerr, Lemaître, etc. In all these
cases, $g$ is given explicitly.

All examples of "fluid" or "potential" type domains given in dimension 4 or dimension 5,
introduced in \cite{vaugon-1} or \cite{vaugon-2} can be immediately expressed into
what we present
here in dimension $n$. It suffices to define the domains of dimension $n$ which are locally isometric
to $(\Omega\times K,g\times g_K)$ where $(\Omega,g)$ is the domain considered in dimension 4 or 5 and $(K,g_K)$ is a
compact Riemannian manifold. The compact manifold $K$ is of the form $K=S^1\times W$ and the metric $g_K$ has for
signature $(-,+,+,\dots,+)$ if the dimension of $\Omega$ is 4.

Obviously these constructions have no mathematical interest, they only consist in bringing everything back to
dimension $n$ without modifying the local geometric properties.

Examples that we will present now are quite new and only make sense when $n>5$ and the
signature of the metric $g$ is of the form $(-,+,+,+,-,+,\dots,+)$. They will give in particular, with a
great simplicity, a very good approximation of the areas of classical physics which use 
electromagnetic and Newtonian potentials. These are the areas that will be reused later in the study of
quantum phenomena.
\bigskip

Before presenting the examples of potential type domains, we first specify some notations and
 definitions. These will be used throughout this paper.


\subsection{Some Notations and Definitions \label{ss1.1.1}}
\begin{enumerate}
 \item The circle $S^1(\delta)$ with radius $\delta$ is defined by setting $S^1(\delta)=\R/2\pi\delta\varmathbb{Z}$,
 
 which, using the surjection $\Pi:\R\rightarrow 2\pi\delta\varmathbb{Z}$, gives canonically:
 
 \textbf{an origin $P$} on $S^1(\delta)$ ~~ ($P:=\Pi(0)$),
 
\textbf{an orientation} (that of $\R$ induced by $\Pi$),
 
 \textbf{a coordinate} $u\in ]0 ~~2\pi\delta[$ for $\Pi(u)\in S^1(\delta)-\{P\}$,
 
 \textbf{a metric} $g_{S^1(\delta)}$ (that of $\R$ quotiented by $\Pi$).
 \item Let $\Theta$ be an open set of $\R^4$ and $\C=\Theta\times S^1(\delta)\times W$ where $W$ is a compact manifold (which
will then often be decomposed into the form $S^3(\rho)\times V$ and $(S^3(\rho),g_{S^3(\rho)})$ will be the Riemannian standard sphere of dimension 3)

$\C$ will be called \textbf{a type cell}.
\item \textbf{A standard coordinate system on $\C$} will be denoted: $(t,x^1,x^2,x^3,u,w)$ \\ where ~~ $(t,x^1,x^2,x^3)\in
 \Theta\subset \R^4$, ~~ $u\in S^1(\delta)$ ~ and ~ $w=(w^1,\dots,w^k)$ ) are the coordinates associated to the choice of 
a chart on $W$. The couple $(t,u)$, coordinates of the "double" time, will sometimes be denoted $(x^0,x^4)$.

The elements of $S^1(\delta)$ will be denoted u (of standard coordinate $u$). The
elements of $W$ will be denoted w (of standard coordinates $(w^1,\dots,w^k)$).
\item \textbf{A reference metric} on the cell $\C$ (which will then become a neutral potential metric)
is a pseudo-Riemannian metric $g_0$ written as:
\[g_0=g_\Theta\times(-g_{S^1(\delta)})\times g_W \]

where $g_\Theta$ is the Minkovski metric on $\Theta \subset\R^4$, ~~ $g_{S^1(\delta)}$ the standard Riemannian metric
 on $S^1(\delta)$ ~ and ~ $g_W$ a Riemannian metric on $W$.

The metric $g_0$ is written in a standard coordinate system:
\[g_0=-dt^2 +\sum_{k=1}^3(dx^k)^2-du^2+\sum g_W{_{ij}}dw^idw^j\]
It corresponds to the choices of geometric units since the coefficients of $dt^2$ and $du^2$ are
($-1$) and those of $(dx^k)^2$ are $(+1)$.
\end{enumerate}
\begin{rmq} \label{r6}
 As "manifold" the type cell could have been defined more simply by taking $S^1(1)$ instead of
$S^1(\delta)$ since by diffeomorphisms, radius do not matter. But then a "standard metric"
$g_0$ would have been defined by $g_0=-dt^2+\sum_{k=1}^3(dx^k)^2-\delta^2 du^2+\sum g_W{_{ij}}dw^idw^j$ to
obtain an identical result. In this case, the units of time on $\R$ and on $S^1(1)$ would have been
different. The choice mentioned in this paragraph seemed preferable.
\end{rmq}
\subsubsection{The notion of fundamental physical constant}
Many domains $(\D,g)$ that prove to be important have a metric $g$ expressed from the reference metric $g_0$. The fundamental constants are therefore chosen as characteristics of this reference metric $g_0$: \\ - The choice of the Minkovski metric for $g|_\Theta$ determines the speed of light (c = 1 in geometric unit). \\ -The radius $\delta$ of the circle $S^1(\delta)$ determines the elementary electric charge (see chapter 2 section \ref{s2.5}). \\ -Some dimensional characteristics of the compact manifold $(W,g|_W)$ determine the fine structure constant $\alpha\simeq1/137$ ) (see chapter 2 \ref{++2,3}). \\ -etc.

\subsection{Potential Metrics \label{ss1.2}}
Consider a chart $(\mathscr V,\zeta)$ of the observation atlas for which the type cell $\C$ is
 $\Theta\times S^1\times W$.\\
\textbf{The metrics representing potentials will be defined, for simplicity, on $\C=\Theta\times
S^1\times W$. In other words they are the metrics $g$ transported by $\zeta$ of that defined on $\mathscr V\subset
 \M$.}
\subsubsection{A- ~ Metrics representing neutral potentials}
\begin{dfn} \label{def:7}
 A pseudo-Riemannian metric $g_0$ defined on a type cell $\C=\Theta\times S^1\times W$ \textbf{represents
a neutral potential} if it is written in the form of a product metric:
\[g_0=g_\Theta\times(-g_{S^1(\delta)})\times g_W\]
where: $g_\Theta$ is the usual Minkovski's metric on $\Theta\subset\R^4$.

$g_{S^1(\delta)}$ is the standard metric of the circle $S^1$ of radius $\delta$.

$g_W$ is a Riemannian metric on the compact manifold $W$ such that the scalar curvature $S_{g_W}$ is
constant, ($S_{g_0}$ is then equal to this same constant).
\end{dfn}
The different neutral potentials are, for the moment, linked to the different possible choices of a 
compact Riemannian manifold $(W,g_W)$ with constant scalar curvature and to the choise of the radius $\delta$ of the circle $S^1$.

The signature of $g_0$ is in every point: $(-,+,+,+,-,\dots,+)$. The choice of Minkovski's metric $g_\Theta$
on $\Theta$ is due to the fact that one seeks, for the moment, only to recover the standard results of quantum theories
 (presented in Chapter 2).

As we will see, the metric $g_0$ will be considered, when describing specific experiments, as the metric
of the observer doing the measurements. A different choice of $g_0$ would be possible, for example to characterize a
"distortion" of the space-time in which the observer makes the measurements, but this would obviously 
strongly complicate calculations. It is essentially $g_\Theta$, ie. Minkovski's metric on $\Theta$, which will be used as the metric of an observer making the measurements.
\begin{rmq} \label{++r7}
 Metric $g_W$ may be different depending on the space-time domains that will be considered (but always with constant scalar curvature). This will be important in chapter 2 on quantum phenomena since, as we shall see, the notion of mass will depend on the scalar curvature of $(W,g_W)$.
\end{rmq}
\begin{rmq} \label{r7}
 The choice of $g_\Theta$ as Minkovski's metric allows us to have the following property:
 if $\Lambda:\R^4\supset\Theta\rightarrow\Theta'\subset\R^4$ is a standard Poincaré transformation, then:
$\sigma^* g_0=~g_0$ ~ when $\sigma:=\Lambda\times I_{S^1\times W}$ ~~ where ~~ $ I_{S^1\times W}$ is the 
"identity" function on $S^1\times W$. In other words, the notion of neutral potential is invariant by observers change
 that correspond to the map changes associated with Poincaré transformations.
In the rest of this paper, some important notions will not have this invariance property by Poincaré transformations, they will be notions related to particular coordinate systems representing observers. It should be noted that invariance by Poincaré transformations on $\Theta$ properties is certainly interesting when it is satisfied, but does not have a great conceptual importance within the framework of our theory. In this regard, it is recalled that in the context of standard general relativity invariance by Poincaré transformations properties appear only in very particular cases as a result of some approximations and can not therefore be considered as fundamental principles (whereas in quantum field theory this principle is considered fundamental).
\end{rmq}

\subsubsection{B- ~ Metrics representing active potentials \label{ss1.3}}
 If one gives a neutral potential $g_0$ on a type cell $\C$, any other metric $g$ on $\C$ is obviously written
 in the form $g=g_0+h$ ~ where ~ $h$ is a field of symmetrical bilinear forms. The field $h$ can be
considered as an endomorphisms field on $\C$, associated to $g_0$, denoted $\leftidx{^e}{h}$, (in coordinates,
$\leftidx{^e}{h}_j^i=g_o^{ik}h_{kj}$).

\textbf{The fundamental hypothesis verified by the metrics representing  potentials will be the NILPOTENCE of the endomorphism field} $\mathbf {\leftidx{^e}{h}}$  (see def. \ref{def:8}).
This hypothesis will allow to recover in particular, by simple calculations, many results obtained by standard classical and quantum physics which describe with a good precision the experimental results. Technically, the nilpotence of $\leftidx{^e}{h}$ makes it possible to carry out exact calculations (without approximations), in particular of the "inverse" of the metric tensor $g$ and the determinant of this denoted $|g|$ (cf. \ref{p++1,1} and expression (\ref{F1.2+})). That said, this nilpotence hypothesis can be considered as an approximation (which consists in neglecting the terms of the expression (\ref{F1.2+}) after some rank) and \textbf{it will be abandoned in section \ref{++2,3}} when it comes to studying the "anomalous magnetic moment of the electron", consequence of an experimental result with remarkable precision. \\

It is clear that $\leftidx{^e}{h}$ is a field of symmetrical endomorphisms associated to $g_0$ (see annex
\ref{a3.2}).
If the signature of $g_0$ was of the form $(+,+,\dots,+)$, the endomorphisms $\leftidx{^e}{h}$ would be
diagonalizable at any point $x$ of $\C$ and could therefore be nilpotent only by being identically null. This
is no longer the case if the signature of $g$ contains signs ($-$) and signs ($+$).

The following study 
show that \textbf{the sign ($-$) in the first position} in the signature of $g_0$ (which corresponds to the
standard Minkovski metric) allows the existence of \textbf{non nulls nilpotents} endomorphism fields $\leftidx{^e}{h}$  whose corresponding metrics $g=g_0+h$ will be in particular those which will give the \textbf{Newtonian} potential.

\textbf{The sign ($-$) in the fifth position} (which corresponds to the circle $S^1(\delta)$) allows
the existence of other \textbf{non-null nilpotent} endomorphism fields whose corresponding metrics will be
in especially those that will give the \textbf{electromagnetic potentials}.

\begin{rmq} The nilpotence property of $\leftidx{^e}{h}$ in the decomposition of $g$ in the form $g_0+h$ will allow to define interesting domains of types other than "potential" like, for example, the one presented in section \ref{s1.8}.
\end{rmq}

We will see that two objects are naturally attached to Newtonian potentials (and even to slightly more general potentials):

 a function $v:\C\rightarrow\R$ called the "potential function", and a lightlike vector field $X_1$
(cf proposition \ref{p1.2}).

Similarly, two objects are naturally attached to electromagnetic potential:

 a vector field $\Upsilon$ called the "electromagnetic potential vector field" (and we will recover here the
potential already defined from the vector field $Y$), and a lightlike vector field $X_2$ (cf.
proposition \ref{p1.3}).

In both cases, these two objects will completely characterize $h$.

The potentials will be qualified \textbf{active} if $h\neq0$. One can notice that for a metric $g_0$
representing a neutral potential, the components on $\Theta$ of the images of the geodesics are straight lines. It non't be the case for active potentials that we will present after having specified all these notions and given
some general properties.

\begin{dfn} \label{def:8}
 An endomorphism field $\leftidx{^e}{h}$ is \textbf{nilpotent of index $p\in \N$} if, for any
$x\in\C$, ~~ for any $q\geqslant p$, ~~the endomorphism $\leftidx{^e}{h}_x$ of $T_x(\C)$ satisfies
($\leftidx{^e}{h}_x)^q=0$ ~ and if exists ~ $x\in\C$ such that $(\leftidx{^e}{h}_x)^{p-1}\neq0$
\end{dfn}

\begin{dfn} \label{def:9}
A pseudo-Riemannian metric $g$ defined on a type cell $\C=\Theta\times S^1(\delta)\times W$
\textbf{represents an active potential} if it is written in the form $g=g_0+h$ ~ where ~ $g_0$ is a metric
representing a neutral potential and the endomorphism field $\leftidx{^e}{h}$ (associated to $g_0$) is nilpotent
of index $p\geqslant2$.
\end{dfn}

The properties set forth in the following proposition are direct consequences of the symmetry of $h$ and the nilpotence of $\leftidx{^e}{h}$.
They are proved in the annex \ref{a3.2}.

\begin{prop} \label{p++1,1}
When $g=g_0+h$ and $\leftidx{^e}{h}$ is nilpotent we have the following properties:
\begin{enumerate}
 \item $\forall x\in\C$, ~~ $\leftidx{^e}{h}$ is a symmetrical endomorphism for $g_0$, that is, \\ $\forall X$
and $Y\in T_x(\C)$, ~~ $g_0(\leftidx{^e}{h}_x(X),Y)=g_0(X,\leftidx{^e}{h}_x(Y))$.
\item $\forall q\in\N^*$, ~~ $\forall x\in\C$, ~~ $trace(\leftidx{^e}{h}_x)^q=0$.
\item for any basis of $T_x(\C)$, ~~ $det(g_{ij}(x))=det(g_0{_{ij}}(x))$.
\end{enumerate}
\end{prop}

A very important result of property 3 is the fact that the volume element $\eta$ associated
to $g$ ~~ ($\eta=\sqrt{det(g_{ij})}dx^0\wedge\dots\wedge dx^{n-1}$) ~ is the same as that associated to $g_0$. This result
will be very simplifying in the study of quantum phenomena that involves active potentials.

The following lemma, whose very simple proof is presented in the annex \ref{a3.2}, is a consequence of the signature choise of $g_0$ and the nilpotence of $\leftidx{^e}{h}$. The importance of the
signature of $g_0$ being exactly of the form $(-,+,+,+,-,+,\dots,+)$ appears here. This will be even more fundamental in the study of
quantum phenomena.

\begin{lem} \label{l2}
 Let $Y$ be the field of timelike unit vectors characterizing the electromagnetism already presented,
here transported on the cell $\C$.

Let $X_0$ be an unit timelike vector field $g_0$-orthogonal to $Y$ (there are an infinity of these vector fields taking into account
of the signature of $g_0$).

So, for any $x\in\C$ ~ the endomorphism $\leftidx{^e}{h}_x$ is null on the space $g_0$-orthogonal to the space
generated at the point $x$ by the $2p$ vector fields:

$Y,\leftidx{^e}{h}(Y),\dots,\leftidx{^e}{h}^{p-1}(Y),X_0,\leftidx{^e}{h}(X_0),\dots,\leftidx{^e}{h}^{p-1}(X_0)$.

The endomorphism field $\leftidx{^e}{h}$ is therefore entirely determined by its values on these $2p$ vector fields 
 (which, in general, are not independent).
\end{lem}
The following remark will justify the definition that follows.
\begin{rmq} \label{r8}
The $1$-differential form that characterizes electromagnetism has been defined in the general framework by $Y^\flat$
\textbf{where "$\flat$" is associated to the metric $g$}. In the case of an active potential where $g=g_0+h$, the vector field  associated with $Y^\flat$ \textbf{by $g_0$} (by $g$ it would obviously $Y$ itself) is none other than
$Y+\leftidx{^e}{h}(Y)$, ~ this since ~ $g_{ij}Y^j=g_0{_{ij}}Y^j+h_{ij}Y^j$ ~ and then
~ $g_0^{ki}g_{ij}Y^j=Y^k+\leftidx{^e}{h}^k_jY^j$.

In addition, $F:=d(Y^\flat)=d(g_{ij}Y^j)=d(h_{ij}Y^j)$ ~~ since ~ $d(g_0{_{ij}}Y^j)=0$.

The electromagnetic potential is therefore completely characterized by the vector field $\leftidx{^e}{h}(Y)$
\textbf{that we will later denote by $\Upsilon$}. We then give the following definition:
\end{rmq}
\begin{dfn} \label{def:10}
 An active potential is \textbf{without electromagnetism} if the vector field $\Upsilon:=\leftidx{^e}{h}(Y)=0$.
\end{dfn}
We will now focus on the particular cases of potential type domains that will allow us to
recover, among other things, all the standard results that describe the behavior of elementary objects,
electrically charged or not, in what is commonly called a "Newtonian potential" or an "electromagnetic potential".
It should be noted that these are exactly the same areas that will be reused in the description of the
quantum phenomena in chapter 2. \textbf{The nilpotence index corresponding to the two cases presented will be
$p=2$ and $p=3$}. Of course, since the dimension of $\M$ is $n$, the maximum nilpotency index of $\leftidx{^e}{h}$ is
$n-1$. The nilpotency index is therefore limited by the dimension of the compact manifold $W$. The choice of small
index of nilpotence can be interpreted as the fact that we neglect some effects linked, for
example, to a compact manifold $V_2$ in a decomposition of $W$ of the form $W=V_1\times V_2$.
\bigskip

\textbf{a- ~ Active potentials of index 2 without electromagnetism} (in particular Newtonian potentials).
\bigskip

The metric $g$ defined on the cell $\C=\Theta\times S^1(\delta)\times W$ is of form $g=g_0+h$ ~ where ~ $g_0$
is
the metric of a neutral potential and the endomorphism field $\leftidx{^e}{h}$ (associated to $g_0$) is \textbf{nilpotent
of index 2}.

As it is assumed that this potential is without electromagnetism (cf. def \ref{def:10}), one poses
$\Upsilon:=\leftidx{^e}{h}(Y)=0$.

According to lemma \ref{l2}, ~ $h$ is then entirely determined by its values on $Y$, $X_0$,
$\leftidx{^e}{h}(X_0)$ ~ where  $X_0$ is a timelike vector field, in time orientation,
$g_0$-orthogonal to $Y$, normed by $g_0(X_0,X_0)=-1$. ~ This field can be considered as an observation field.

Since the potential is active, ~ $h$ is assumed to be nonzero. On the other hand, it is easy to verify that:
$\leftidx{^e}{h}_x=0\Longleftrightarrow \leftidx{^e}{h}_x(X_0)=0$. ~ Indeed, $\leftidx{^e}{h}_x(Y)=0$ ~ and ~ $\leftidx{^e}{h}_x^2(X_0)=0$, if $\leftidx{^e}{h}_x(X_0)=0$ the lemma \ref{l2} shows that
$\leftidx{^e}{h}_x=0$. It is deduced in particular that the function $g_0(\leftidx{^e}{h}(X_0),X_0)$ is not
identically zero.

We then give the following definition:
\begin{dfn} \label{def:11}
 The non-zero function $v=-\frac{1}{2}g_0(\leftidx{^e}{h}(X_0),X_0)$ ~ is called \textbf{the potential function seen by
$X_0$} of the active potential domain.
(The coefficient $-\frac{1}{2}$ is put to recover the standard notion of Newtonian potential).
\end{dfn}

\begin{prop} \label{p1.2}
 If $g=g_0+h$ ~ is the metric representing an active potential of index 2 without electromagnetism, ~ then it exists, ~ on
the subset of $\C$ where the function $v$ is non-zero, \textbf{an unique vector field $X_1$} which satisfies:
\[h=-2vX_1^\flat\otimes X_1^\flat ~~~and~~~g_0(X_1,X_0)=1\]
We also have the following properties:

$g_0(X_1,X_1)=0$ ~~ ($X_1$ is a lightfield for $g_0$) ~~ and ~~ $g_0(X_1,Y)=0$.

(Recall that $X_1^\flat$ ~ is the $1$-form associated with $X_1$ \textbf{by $g_0$}).

So we have:
\[g=g_0-2vX_1^\flat\otimes X_1^\flat\]
\end{prop}

 \textbf{Proof}: We start by proving that:

\[g_0(\leftidx{^e}{h}(X_0),X_0)\leftidx{^e}{h}=\leftidx{^e}{h}(X_0)\otimes(\leftidx{^e}{h}(X_0))^\flat\]

This is easily obtained by noting that since $\leftidx{^e}{h}^2=0$, ~ the right hand side of this equality is
nilpotent of index $\leq2$, ~ then showing that the equality is true applied on $Y$, $X_0$ and
$\leftidx{^e}{h}(X_0)$  and therefore true everywhere according to lemma \ref{l2}.

Then we define, on the subset where $v$ is non-zero:
\[X_1:=-\frac{1}{2v}\leftidx{^e}{h}(X_0)\]
Hence:
\[g_0(X_1,Y)=-\frac{1}{2v}g_0(X_0,\leftidx{^e}{h}(Y))=0\]
\[g_0(X_1,X_0)=-\frac{1}{2v}g_0(X_0,\leftidx{^e}{h}(X_0))=1\]
And:
\[-2v\leftidx{^e}{h}=(2v)^2X_1\otimes X_1^\flat\]
Therefore:
\[\leftidx{^e}{h}=-2vX_1\otimes X_1^\flat\]
We get the uniqueness of $X_1$ using equalities ~ $\leftidx{^e}{h}=-2vX_1\otimes X_1^\flat$ ~ and
~ $g_0(X_1,X_0)=1$ which allow to write:

$\leftidx{^e}{h}(X_0)=-2vX_1$ ~~ hence ~~ $X_1=-\frac{1}{2v}\leftidx{^e}{h}(X_0)$.

Equality $g_0(X_1,X_1)=0$ ~ is immediate.
\begin{rmq} \label{r9}
 The potential function $v$ and the field $X_1$ depend on the choice of the observation field $X_0$.
 
 As we will see when calculating the Ricci curvature of this potential type (cf. proposition \ref{p1.4}), the
assumptions given in definition \ref{def:6} of a potential type domain defined from the Einstein curvature, can be satisfied if $\Delta_{g_0}v=0$ (see also remark \ref{r5} for the influence of
scalar curvature). We therefore give the following definition:
\end{rmq}
\begin{dfn} \label{def:12}
 A \textbf{Newtonian potential type domain} is an active potential type domain of index 2 without electromagnetism such
as the potential function $v$ satisfies $\Delta_{g_0}v=0$.
\end{dfn}
\textbf{b- ~ Electromagnetic potentials.}

  The metric $g$ defined on the cell $\C=\Theta\times S^1(\delta)\times W$ is of the form $g=g_0+h$ ~ where ~ $g_0$
is
the metric of a neutral potential and the endomorphism field $\leftidx{^e}{h}$ (associated to $g_0$) is \textbf{nilpotent
of index 2 or 3}.

The fact that this potential is "electromagnetic" is essentially characterized by the property:
$\leftidx{^e}{h}(Y)\neq0$ (see def. \ref{def:10}). However, it will be assumed that $g_0(\leftidx{^e}{h}(Y),Y):=h(Y,Y)=0$ ~
so that the vector field associated with $g_0$ at $Y^\flat$ (where here $\flat$ is associated to $g$) is the $g_0$-orthogonal sum of $Y$ and
$\leftidx{^e}{h}(Y)$ (see remark \ref{r8}). (We leave it to the reader to verify that, again, this hypothesis
can be interpreted as the fact that one neglects the quantum effects on electromagnetism).

To clarify that this potential is only electromagnetic and does not have a Newtonian component, it will be assumed that exist 
a timelike vector field $X_0$ such that $\leftidx{^e}{h}^2(X_0)=0$  and $g_0(\leftidx{^e}{h}(X_0),X_0)=0$
~ ~ ~ (lower assumptions than $\leftidx{^e}{h}(X_0)=0$, the first equality is
always true in the case of nilpotence index 2).

The reader will compare this with the definition of the active potential without electromagnetism and notice that
the roles of $Y$  and  $X_0$ are transposed (they are both timelike but for the different signs "$-$" of
the signature of $g_0$). We are further authorized here that the nilpotence is of index 3. An essential difference appears to be due to the fact that the field $Y$ is perfectly determined (taking into account the definition of an observation atlas) whereas $X_0$ is not and its choice is considered as that of an observation field. In the case of an electromagnetic potential that we define here, the objects that will characterize
$h$ do not depend on the choice of $X_0$ (thanks to the hypothesis $\leftidx{^e}{h}^2(X_0)=0$).

As already stated, we denote $\Upsilon:=\leftidx{^e}{h}(Y)$. \\ The analog of the proposition \ref{p1.2} is written in 
the following way:
\begin{prop} \label{p1.3}
 Let $g=g_0+h$ ~ be the metric representing an electromagnetic potential. Then it exists, on the subset of $\C$ where
$\Upsilon$ is not zero, \textbf{a unique vector field $X_2$} such that:
\[h=\Upsilon^\flat\otimes X_2^\flat+X_2^\flat\otimes \Upsilon^\flat\]
(where it is recalled that "$\flat$" in $X_2^\flat$ ~ and ~ $\Upsilon^\flat$ is associated to $g_0$).

This vector field has the following properties: \\
$g_0(X_2,X_2)=0$ ~~ ($X_2$ is a lightfield for $g_0$), ~ $g_0(X_2,Y)=1$ and $g_0(X_2,\Upsilon)=0$.

So we have:
\[g=g_0+\Upsilon^\flat\otimes X_2^\flat+X_2^\flat\otimes \Upsilon^\flat\]
\end{prop}

\textbf{Proof} The uniqueness of $X_2^\flat$ in decomposition $\Upsilon^\flat\otimes X_2^\flat+X_2^\flat\otimes
\Upsilon^\flat$, when $\Upsilon$ is not zero, is a simple tensor property due to the fact that $h$ is
symmetrical of rank 2.

The existence of such a decomposition is shown by considering two cases:
\begin{enumerate}
 \item It is assumed that $g_0(\Upsilon_x,\Upsilon_x):=g_0(\leftidx{^e}{h}_x(Y_x),\leftidx{^e}{h}_x(Y_x))\neq0$.
 
 (Here the nilpotency index of $\leftidx{^e}{h}_x$ is 3).
 
 We start by proving that, at point $x$ (omitted in the following lines):
 
\begin{eqnarray}
\leftidx{^e}{h}=\frac{1}{g_0(\Upsilon,\Upsilon)}(\Upsilon\otimes (\leftidx{^e}{h}^2(Y))^\flat+(\leftidx{^e}{h}
^2(Y))\otimes\Upsilon^\flat) \label{F1.1}
\end{eqnarray}
 
According to lemma \ref{l2}, $\leftidx{^e}{h}$ vanishes on the $g_0$-orthogonal space to: $Y$, $\leftidx{^e}{h}(Y)$,
$\leftidx{^e}{h}^2(Y)$, $X_0$, $\leftidx{^e}{h}(X_0)$. ~ One can also easily prove that the right hand side of
\ref{F1.1} also vanishes on this space.

To show \ref{F1.1}, all that remains is to prove that equality occurs when each hand side is applied
successively to $Y$, $\leftidx{^e}{h}(Y)$, $\leftidx{^e}{h}^2(Y)$, $X_0$, $\leftidx{^e}{h}(X_0)$, ~ which is easy
  using the fact that $\leftidx{^e}{h}$ is $g_0$-symmetric, ~ $\leftidx{^e}{h}^3=0$, ~ $h(Y,Y)=0$,
~ $\leftidx{^e}{h}^2(X_0)=0$ ~ and after showing that:
\[\leftidx{^e}{h}(X_0)=(h(Y,X_0)/g_0(\Upsilon,\Upsilon))\leftidx{^e}{h}^2(Y)\]
This last point is obtained the same method way considering the linear forms associated with each 
hand side of this equality and using the fact that $h(X_0,X_0):=g_0(\leftidx{^e}{h}(X_0,X_0)=0$.

The desired decomposition is then obtained by \ref{F1.1} by seting:
\[X_2:=(g_0(\Upsilon,\Upsilon))^{-1}\leftidx{^e}{h}^2(Y)\]
It is then immediate to prove that ~ $g_0(X_2,X_2)=0$, ~~ $g_0(X_2,Y)=1$ ~~ and ~~ $g_0(X_2,\Upsilon)=0$.
\item It is assumed that $g_0(\Upsilon_x,\Upsilon_x)=0$ ~ and ~ $\leftidx{^e}{h}_x(Y)\neq0$, ~ in other words
$\leftidx{^e}{h}_x(Y)$ is a lightfield (this is the case if the nilpotence index is 2).

Then, at the considered points $x$ (omitted in the following lines):
\[h(X_0,Y):=g_0(\leftidx{^e}{h}(X_0),Y)\neq0\]
(in particular $h(X_0)\neq0$)

Indeed, we have:
\[0=g_0(\leftidx{^e}{h}(Y),Y)=g_0(\leftidx{^e}{h}(Y),\leftidx{^e}{h}(Y))=g_0(\leftidx{^e}{h}(Y),\leftidx{^e}{h}^2
(Y))=g_0(\leftidx{^e}{h}(Y),\leftidx{^e}{h}(X_0))\]
If $g_0(\leftidx{^e}{h}(X_0),Y)=g_0(X_0,\leftidx{^e}{h}(Y)$ was zero, the linear form $g_0$-associated with
$\leftidx{^e}{h}(Y)$ would be zero according to lemma \ref{l2}, which is contrary to the hypothesis.

Using the same procedure as in the first case we then prove that:
\[\leftidx{^e}{h}=(h(X_0,Y))^{-1}(\leftidx{^e}{h}(Y)\otimes(\leftidx{^e}{h}(X_0))^\flat+\leftidx{^e}{h}
(X_0)\otimes(\leftidx{^e}{h}(Y))^\flat)\]

The desired decomposition is obtained by setting:
\[X_2=(h(X_0,Y))^{-1}\leftidx{^e}{h}(X_0)\]
It is then immediate to prove that: $g_0(X_2,X_2)=0$, ~~ $g_0(X_2,Y)=1$ ~~ and ~~ $g_0(X_2,\Upsilon)=0$.

\end{enumerate}
\begin{rmq} \label{r10}
 An active potential, of metric $g$, has been defined (for simplicity) on a type cell $\C=\Theta\times
S^1(\delta)\times W$, ~ but when $(\mathscr V,\zeta)$ is the chart of the observation atlas such that
$\zeta(\mathscr V)=\C$, ~ $g$ is the image $\zeta_*(g_\M)$ of the Riemannian metric defined on
$\M$.

Following the remark \ref{r7}, when $(\mathscr V',\zeta')$ is another chart of the observation atlas such that
$\zeta'(\mathscr V')=\C'\times S^1(\delta)\times W$, ~ when $\sigma:= \zeta\circ\zeta'^{-1}:\C'\rightarrow\C$ is
an isometry of the form $\Lambda\times I_{S^1\times W}$ ~ and ~ $\Lambda:\Theta'\rightarrow\Theta$ is a standard Poincaré transformation, the metric $g'=(\zeta'^{-1})^*(g)=\sigma^*g$, defined on $\Theta'\times S^1(\delta)\times W$,
is written as $g_0+\sigma^*h$.

In the context of an active potential without electromagnetism, the potential function $v'$ seen by $\sigma_*^{-1}X_0$,
"read" in $(\mathscr V',\zeta')$, is therefore equal to $v\circ\sigma$ ~ and ~ $X'_1=\sigma_*^{-1}X_1$.

In the context of an electromagnetic potential $\Upsilon'=\sigma_*^{-1}\Upsilon$ ~ and ~ $X'_2=\sigma_*^{-1}X_2$.
\end{rmq}


\section{Geodesics of potential type domains \label{s1.5}}
Study of the geodesics will be done in the usual way starting by calculating the Christoffel symbols of
 corresponding metrics in a standard coordinate system of the considered type cell.
\bigskip

Calculation of Christoffel's symbols requires the inverse of the matrix $(g)$ in the system of
chosen coordinates. The hypothesis of $p$-nilpotence of $\leftidx{^e}{h}$ allows to quickly obtain the matrix
$(g)^{-1}$ whose terms are conventionally denoted $g^{ij}$. Indeed, since $g=g_0+h$,
~~ $(g_0)^{-1}(g)=I+(\leftidx{^e}{h})$ ~ and, since $\leftidx{^e}{h}^p=0$, the inverse of the matrix $I+(\leftidx{^e}{h})$
is:
\[I-(\leftidx{^e}{h})+(\leftidx{^e}{h})^2+\dots+(-1)^{p-1}(\leftidx{^e}{h})^{p-1}\]
  Then, $(g)^{-1}(g_0)=I-(\leftidx{^e}{h})+\dots+(-1)^{p-1}(\leftidx{^e}{h})^{p-1}$, ~~ hense:
  \begin{eqnarray}
g^{ij}=g_0^{ij}-h^{ij}+h^i_kh^{kj}+\dots+(-1)^{p-1}h^i_{k_1}h^{k_1}_{k_2}\dots h^{k_{p-1}j}.\label{F1.2+}
\end{eqnarray}
\textbf{Warning}: here $(h^{ij})$ \textbf{is not} the inverse of the matrix ($h_{ij}$) ~ (which usualy is not
invertible) but is defined by $h^{ij}=g_0^{ik}g_0^{jl}h_{kl}$. At the opposite $(g^{ij})$ is indeed
the inverse of matrix $(g_{ij})$ ~~ (and is not $g_0^{ik}g_0^{jl}g_{kl}$). \\

We have if $p=2$:
\begin{eqnarray}
 g^{ij}=g_0^{ij}-h^{ij}. \label{F1.2}
\end{eqnarray}
If $p=3$:
\begin{eqnarray}
 g^{ij}=g_0^{ij}-h^{ij}+h^i_kh^{kj}. \label{F1.3}
\end{eqnarray}


\subsection{Geodesics of a potential type domain without electromagnetism \label{ss1.4}}
The type cell is of the form $\C=\Theta\times W$ ~ where here the circle $S^1(\delta)$ is considered to be a factor of the
compact manifold $W$, which does not intervene when there is no electromagnetism.

The open set $\Theta\subset\R^4$ will be of the form $\Theta=I\times\mathscr U$ ~ where ~ $I$ is an interval of $\R$ ~ and
~ $\mathscr U$ an open set of $\R^3$. This decomposition is justified by the fact that potential type domains without
electromagnetism, defined in proposition \ref{p1.2}, are not Poincaré-invariants when one chooses
$X_0=\frac{\partial}{\partial t}$ linked to the standard coordinate system.

According to proposition \ref{p1.2}, the pseudo-Riemannian metric $g$ is written:
\[g=g_0-2vX_1^\flat\otimes X_1^\flat \text{~~~where ~~~}g_0=g_\Theta\times g_W\]
The vector field $X_1$ (associated with $X_0$) satisfies:
\[g_0(X_1,X_1)=0, ~~~g_0(X_1,Y)=0, ~~~g_0(X_1,X_0)=1\]
In order to recover the precise results on the geodesics that describe the classical trajectory of an
elementary object in a potential without electromagnetism, we will assume the following hypothesis ($H_N$).
\bigskip

\textbf{Hypothesis $H_N$:}
\begin{enumerate}
 \item The potential function $v$ is defined on $\mathscr U$.
 \item The vector field $X_1$ is defined on $I\times W$ and is a Killing field.
\end{enumerate}
(Function $v$ and field $X_1$ can  naturally be considered as defined on $\C$: $v$ only depends on
variables of $\mathscr U$ and $X_1$ is tangent to $I\times W$)
\bigskip
  
It would be interesting to study the modifications made to the geodesics that we will describe when we do'nt assume
 the hypothesis $H_N$, especially in the case where the potential function is of the form $v=-\frac{m}{r}$ (cf.
section \ref{s1.7}) and estimate the disturbances on the usual conics.


\subsubsection{Christoffel symbols}
In a standard coordinate system of the cell $\C$, we denote:

  $\Gamma^k_{ij}$ ~ (resp $\tilde \Gamma^k_{ij})$ the Christoffel symbols of $g$ ~ (resp $g_0$).

  $T^k_{ij}$ the coordinates of the \textbf{tensor} ($\Gamma^j_{ij}-\tilde\Gamma^k_{ij})$.

We have, when $g=g_0+h$:
\begin{eqnarray}
 T^k_{ij}=\frac{1}{2}g^{kl}(\nabla_ih_{jl}+\nabla_jh_{il}-\nabla_lh_{ij}) \label{F1.4}
\end{eqnarray}

\textbf{where $\nabla_i$ are associated to $g_0$}.

(This result is easily shown by taking a normal coordinate system for $g_0$).

Here $h=-2vX_1\otimes X_1$.

\textbf{In the following, $X_1$ will be denoted $X$ to simplify writing}.

We have: $g^{kl}=g_0^{kl}+2vX^kX^l$

And: $\nabla_j(h_{il})=-2\nabla_j(vX_iX_l)=-2((\nabla_jv)X_iX_l+v(X_i\nabla_jX_l+X_l\nabla_iX_j))$

So, developing \ref{F1.4} and using the fact that $X_1$ is a Killing field (ie. $\nabla_iX_j+\nabla_jX_i=0$), we get:
\[T^k_{ij}=-(g_0^{kl}
+2vX^kX^l)((\nabla_jv)X_iX_l+(\nabla_iv)X_jX_l-(\nabla_lv)X_iX_j+2v(X_i\nabla_jX_l+X_j\nabla_iX_l))\]

Since $X(v)=0$ and $X^lX_l=0$:
\begin{eqnarray}
 T^k_{ij}=-X^k((\nabla_jv)X_i+(\nabla_iv)X_j)+(\nabla^kv)X_iX_j-2v(X_i\nabla_jX^k+X_j\nabla_iX^k) \label{F1.5}
\end{eqnarray}

And we recall that: $\Gamma^k_{ij}=\tilde\Gamma^k_{ij}+T^k_{ij}$.

We deduce, since $X^kX_k=0$ (therefore $X_k\nabla_iX^k=0)$ and $X(v)=0$:

$X_k\Gamma^k_{ij}=X_k\tilde\Gamma^k_{ij}$

But, $\nabla_iX_j=\partial_iX_j-X_k\tilde\Gamma^k_{ij}$, it follows, since $\nabla_iX_j+\nabla_jX_i=0$:

$2X_k\tilde\Gamma^k_{ij}=\partial_iX_j+\partial_jX_i$

Hence:
\begin{eqnarray}
 X_k\Gamma^k_{ij}=\frac{1}{2}(\partial_iX_j+\partial_jX_i) \label{F1.6}
\end{eqnarray}


\subsubsection{Study of the important geodesics}
We consider a geodesic $x:\R\supset I\rightarrow\C$.

$\forall s\in I$ ~~~ $x(s)=(x^0(s),x^1(s),\dots,x^{n-1}(s))$.

For $k$ ~ from ~ $0$ to $n-1$ ~~ and ~~ $\forall s\in I$:
\begin{eqnarray}
 {x^k}''_{(s)}+\Gamma^k_{{ij}_{(x(s))}}{x^i}'_{(s)}{x^j}'_{(s)}=0 \label{F1.7}
\end{eqnarray}
\hspace{0.5cm} Using \ref{F1.6} we can deduce:

$X_{k_{(x(s))}}{x^k}''+(\partial_iX_j)_{(x(s))}{x^i}'_{(s)}{x^j}'_{(s)}=0$

In other words:

$\frac{d}{ds}(X_{k_{(x(s))}}{x^k}'_{(s)})=0$

Hence:

$X_k{x^k}'=K$ ~~ where ~~ $K$ is a constant

(The parameterization index $s$ is no longer specified here and in the following)
\bigskip

We are now interested in the components of geodesics in the apparent space $\mathscr U$.

\textbf{For $k$ ~ from ~ $1$ to $3$}, ~~ $X^k=0$ ~~ and ~~ $\tilde\Gamma^k_{ij}=0$, ~~ then, according to \ref{F1.5}:

$\Gamma^k_{ij}=(\nabla^kv)X_iX_j$ ~~~ (where here $\nabla^kv=\partial_kv$).

The geodesic equation then gives:

\textbf{For $k$ from ~ $1$ to $3$}, ~ ${x^k}''=-(\nabla^kv)X_iX_j{x^i}'{x^j}'=-K^2\nabla^kv$.

An affine parameter change in the geodesic allows us to choose
$K=1$ (we will not be interested in particulars geodesics
 for which
$K=0$).

We finally get:
\bigskip

$({x^1},{x^2},{x^3})''_s=-(\frac{\partial v}{\partial x^1},\frac{\partial v}{\partial x^2},\frac{\partial
v}{\partial x^3})_{x(s)}=-(\nabla_{(x^1,x^2,x^3)}v)_{x(s)}$
\bigskip

\textbf{Which is nothing else than the Poisson equation} in classical physics when $v$ is the Newtonian potential
 and $(x^1(s),x^2(s),x^3(s))$ describes the trajectory of a material point in such a potential, \textbf{but
 considering that $s$ is the time parameter which here does not correspond to $x^0=t$}. The parameter $s$ can be
interpreted as the proper time associated with the image of the corresponding geodesic.

Note that if we assume the velocity, associated with a geodesic, very small compared to the
speed of light, ie. if $\forall k\neq0$ ~~~ ${x^k}'(s)=\circ(1)$ ~ ($x^0$ is the time variable), then ${x^0}'(s)=1+\circ(1)$. Indeed, as $X_k{x^k}'(s)=1$, ~~ ${x^0}'(s)=1-\sum_{k\neq0}X_k{x^k}'(s)$ ~
since $X_0=1$ and more $\sum_{k\neq0}X_k^2=1$.
This means that, in this case, the parameter $s$ is very close to the time $x^0$ given by the coordinates system, which corresponds to the usual non-relativistic approximation.
\begin{rmq} \label{r11}
 It can be verified that, in the framework of the last assumptions \\
$g(\overrightarrow{\text{v}}(s),\overrightarrow{\text{v}}(s))=g_0(\overrightarrow{\text{v}}(s),\overrightarrow{\text{v}}
(s))+h(\overrightarrow{\text{v}}(s)
, \overrightarrow{\text{v}}(s))$ (where $\overrightarrow{\text{v}}$ is the tangent vector to the geodesic), which is
necessarily
a constant $C_0$, will be very close to $-1$ if we assume that the potential $v=\circ(1)$. We could have chosen the
parametrization of the geodesic so that $C_0=-1$ ~ (classical normalization), the constant $K$ would then be
 different from $1$ while still being very close.
\end{rmq}


\subsection{Geodesics of an electromagnetic potential type domain \label{ss1.5}}
The type cell is of the form $\C=\Theta\times S^1(\delta)\times W$ (here the circle $S^1(\delta)$ is important).

We consider the vector field $\Upsilon=\leftidx{^e}{h}(Y)$ which defines the electromagnetic potential and the vector field
 $X_2$ given by proposition \ref{p1.3}. According to this one, the pseudo-Riemannian metric is written:

$g=g_0+\Upsilon^\flat\otimes X_2^\flat+X_2^\flat\otimes\Upsilon^\flat$ ~~ where ~~ $g_0=g_\Theta\times(-g_{S^1})\times g_W$.

The vector field $X_2$ satifies:

$g_0(X_2,X_2)=0$, ~~~ $g_0(X_2,Y)=1$, ~~~ $g_0(X_2,\Upsilon)=0$.

In order to recover the precise results on the geodesics that describe the classical trajectory of an
elementary \textbf{electrically charged} object in an electromagnetic field, we will assume the hypothesis
$H_E$ similar to hypothesis $H_N$ previously studied.
\bigskip

\textbf{Hypothesis $H_E$:}
\begin{enumerate}
 \item a- ~ The vector field $\Upsilon$ is defined on $\Theta$
 \item b- ~ The vector field $X_2$ is defined on $S^1(\delta)\times W$ ~ and is a Killing field.
\end{enumerate}
(The fields $\Upsilon$ and $X_2$ can be naturally considered as defined on $\C$: $\Upsilon$ is tangent to
$\Theta$ and depends only on the variables of $\Theta$, ~ $X_2$ is tangent to $S^1(\delta)\times W$ and depends only on the
variables of $S^1(\delta)\times W$).


\subsubsection{Christoffel symbols}
We start from expression \ref{F1.4} given in the previous study.

Here $h=\Upsilon^\flat\otimes X_2^\flat+X_2^\flat\otimes\Upsilon^\flat$.

\textbf{In the following $X_2$ will be denoted $X$ to simplify writing.}

We have: $g^{kl}=g_0^{kl}-h^{kl}+h^{km}h_m^l$, ~~ in other words:

$g^{kl}=g_0^{kl}-(\Upsilon^kX^l+\Upsilon^lX^k)+(\Upsilon^kX^m+\Upsilon^mX^k)(\Upsilon_mX^l+\Upsilon^lX_m)$

Hence, since $X^mX_m=0$ ~~ and ~~ $X^m\Upsilon_m=0$:

$g^{kl}=g_0^{kl}-(\Upsilon^kX^l+\Upsilon^lX^k)+(\Upsilon^m\Upsilon_m)X^kX^l$

Then:

$T^k_{ij}=\frac{1}{2}(g_0^{kl}-(\Upsilon^kX^l+\Upsilon^lX^k)+(\Upsilon^m\Upsilon_m)X^kX^l)(\nabla_jh_{il}+\nabla_ih_{jl}
-\nabla_lh{ij})$

By developing and using the fact that $\nabla_iX_j+\nabla_jX_i=0$, we find:

$(*):=(\nabla_jh_{il}+\nabla_ih_{jl}-\nabla_lh_{ij})=X_iF_{jl}+X_jF_{il}
+2(\Upsilon_i\nabla_jX_l+\Upsilon_j\nabla_iX_l)+X_l(\nabla_j\Upsilon_i+\nabla_i\Upsilon_j)$

where $F_{ij}=\nabla_i\Upsilon_j-\nabla_j\Upsilon_i=\partial_i\Upsilon_j-\partial_j\Upsilon_i$ ~ was
components of the 2-differential form

$F:=dY^\flat=d\Upsilon^\flat$.

Note that, according to the assumption ($H_E$), $X^lF_{jl}=0$.

Then, since $X^lX_l=0$ ~ (then $X^l\nabla_iX_l=0$):

$X^l(*)=0$

Hence, since $\Upsilon^lX_l=0$:
\begin{eqnarray}
 T^k_{ij}=\frac{1}{2}(X_iF_j^{~~k}+X_jF_i^{~~k}
)+\Upsilon_i\nabla_jX^k+\Upsilon_j\nabla_iX^k+X^k(\nabla_j\Upsilon_i+\nabla_i\Upsilon_j-\Upsilon^l(*)) \label{F1.8}
\end{eqnarray}

It follows:

$X_kT^k_{ij}=0$

Hence

$X_k\Gamma^k_{ij}=X_k\tilde\Gamma^k_{ij}$

Then, as for the potential without electromagnetism, one deduces:

$2X_k\tilde\Gamma^k_{ij}=\partial_iX_j+\partial_jX_i$

Then:
\begin{eqnarray}
 X_k\Gamma^k_{ij}=\frac{1}{2}(\partial_iX_j+\partial_jX_i) \label{F1.9}
\end{eqnarray}


\subsubsection{Study of the important geodesics}
We consider a geodesic $x:\R\supset I\rightarrow\C$.

$\forall s\in I$ ~~~ $x(s)=(x^0(s),x^1(s),\dots,x^{n-1}(s))$.

For $k$ ~ from ~ $0$ to $n-1$ ~~ and ~~ $\forall s\in I$:
\begin{eqnarray}
 {x^k}''_{(s)}+\Gamma^k_{{ij}_{(x(s))}}{x^i}'_{(s)}{x^j}'_{(s)}=0 \label{F1.10}
\end{eqnarray}
\hspace{0.5cm} With \ref{F1.9} we can deduce:

$X_{k_{(x(s))}}{x^k}''+(\partial_iX_j)_{(x(s))}{x^i}'_{(s)}{x^j}'_{(s)}=0$

In other words:

$\frac{d}{ds}(X_{k_{(x(s))}}{x^k}'_{(s)})=0$

Hence:

$X_k{x^k}'=K$ ~~ where ~~ $K$ is a constant

(The parameterization index $s$ is no longer specified here and in the following)
\bigskip

We are now interested in the first four components of geodesics that correspond to the "classical" space-time
 $\Theta$.

According to \ref{F1.8} and since $\Gamma^k_{ij}=\tilde\Gamma^k_{ij}+T^k_{ij}$,

\textbf{For $k$ from ~ 0 ~ to ~ 3}:

$\Gamma^k_{ij}=\frac{1}{2}(X_iF_j^{~~k}+X_jF_i^{~~k})$

Then equation \ref{F1.10} of the geodesics is written:

${x^k}''+\frac{1}{2}(X_iF_j^{~~k}+X_jF_i^{~~k}){x^i}'{x^j}'=0$

Which means,

\textbf{For $k$ from ~ 0 ~ to ~ 3}:
\begin{eqnarray}
 {x^k}''+KF_i^{~~k}{x^i}'=0 \label{F1.11}
\end{eqnarray}

The geodesic is conventionally parameterized so that $g(x'(s),x'(s))=-1$ ~ and \\ ~ $x'(s)|_{T_{x(s)}(\Theta)}$ is
in the time orientation given by $x^0$.

We denote $\overrightarrow{\text{v}}_{(s)}=({x^0}'(s), {x^1}'(s), {x^2}'(s), {x^3}'(s))$ the vector defined by the four
first components of the vector $x'(s)$ tangent to the geodesic.

Equation \ref{F1.11} is then written since $F_i^{~~k}=-F^k_{~~i}$:

$\overrightarrow{\text{v}}'_{(s)}=K\leftidx{^e}{F}_{x(s)}(\overrightarrow{\text{v}})_{(s)}$

where $\leftidx{^e}{F}$ is the endomorphism field associated with $F$ associated to $g_0$.

\textbf{We recover the classical expression (of special relativity) which gives the equation of motion of a mass particle $m$ and charge
$q$ in an electromagnetic field $F$ when defining $K=\frac{q}{m}$ and when $s$ is
actually the proper time of the particle (see, for example, \cite{gourg} (17-61) p.554). The result that
the trajectory of an elementary object is described by geodesics is known in general relativity but only
in the context of gravitation. We have just shown here that for us this result is still valid in the context
of electromagnetism}.

What we have just written is only a verification, in this example, of the principle of calculating the trajectory of an elementary object
electrically charged in a potential discussed in section \ref{s1.3}.

Here $K=\frac{q}{m}=X_k{x^k}'$ is a feature of geodesic over compact manifold $S^1(\delta)\times W$
since $X_k=0$ for $k$ from 0 to 3.

\begin{rmq} \label{r12}
 As the parameterization of the geodesic is chosen so that $g(x'(s),x'(s))=-1$ and since
$g_{ij}=g_{0_{ij}}+\Upsilon_iX_j+X_i\Upsilon_j$ we have:

$\forall s\in I$
~~~ ${{x^0}'}^2_{(s)}+{{x^4}'}^2_{(s)}=1+\sum_{k\neq0et4}{{x^k}'}^2_{(s)}+2K\Upsilon_{i(x(s))}{x^i}'_{(s)}$

If we assume that for $i$ from 0 to 3: $K\Upsilon_i=\circ(1)$ ~ and that, ~~ for $k\neq0$: ${x^k}'(s)=\circ(1)$ ~ (which
means in particular that the speed determined by the geodesic is very small compared to the speed of the
light) then ${x^0}'(s)=1+\circ(1)$ (considering the choice of parameterization in the time orientation associated to $x^0$)
and the parameter $s$ is very close to the time $x^0$ given by the coordinate system, which corresponds to
the usual non-relativistic approximation.
\end{rmq}


\section [Ricci curvature and scalar curvature] {Ricci curvature, scalar curvature and different properties
 of potential type domains
\label{s1.6}}


\subsection{Active Potentials without electromagnetism with index 2 \label{ss1.6}}
We recall that, according to proposition \ref{p1.2}, ~~ $g=g_0-2vX_1^\flat\otimes X_1^\flat$ ~~ where  $v$ is the potential function and $X_1$ satisfies:

$g_0(X_1,X_1)=0$, ~~~ $g_0(X_1,Y)=0$, ~~~ $g_0(X_1,X_0)=1$.

It is further assumed here, as in the study of geodesics (see \ref{ss1.4}), that the hypothesis $H_N$ is satisfied, and,
 to not use terms related to $D_{g_0}X_1$ in the results that follow, we
will assume that
$D_{g_0}X_1=0$, which is a stronger assumption than the fact that $X_1$ is a Killing field. These assumptions
can be seen as the fact that some quantum effects on the field $X_1$ are neglected.

We then obtain the following result.

\begin{prop} \label{p1.4}
~ \\ \vspace*{-1em}
\begin{enumerate}
  \item $R_{{icc}_g}^\sharp=R_{{icc}_{g_0}}^\sharp-(\Delta_{g_0}v)X_1\otimes X_1$ \hspace{1cm} (Here
~ $\Delta_{g_0}:=-\nabla^k\nabla_k$)

where $R_{{icc}_{g_0}}$ is entirely determined by $R_{{icc}_{{g_0}_W}}$ since $R_{{icc}_{{g_0}_\Theta}}=0$.
\item $S_g=S_{g_0}+2vR_{{icc}_{g_0}}(X_1,X_1)$

where $S_{g_0}=S_{{g_0}_W}$.
\item $X_1$ is also a lightfield for the metric $g$ and $D_gX_1=0$.
\item $D_gY=0$, in particular $Y$ is a Killing field and a geodesic field for the metric $g$ (this is
obviously true for the metric $g_0$).
 \end{enumerate}
\end{prop}
The proof of this proposition is detailed in annex \ref{a3.3}.
\begin{rmq} \label{r13}
 If $R_{{icc}_{g_0}}(X_1,X_1)=0$, the scalar curvature $S_g$ remains equal to that of the neutral potential $S_{g_0}$. This
is the case, for example, when $W$ is a product Riemannian manifold $V_1\times V_2$ ~ where ~ ($V_2,g_{V_2}$) is an
Einstein manifold  and when $X_1$ is tangent to $V_2$ ~ (one can have dim $V_1=0$). Indeed, we have
$R_{{icc}_{g_0}}(X_1,X_1)=C^{te}g_0(X_1,X_1)=0$.

In fact, in standard experiments the potential function $v$ is $\ll 1$ (in geometric units).
Given the normalization of $X_1$ by $g_0(X_1,X_0)=1$, we can write in this case:
$|2vR_{{icc}_{g_0}}(X_1,X_1)| \ll |S_{g_0}|$ ~ and $S_g$ remains very close to $S_{g_0}$ which is constant by definition.
\end{rmq}
\subsection{Electromagnetic potentials \label{ss1.7}}
According to proposition \ref{p1.3}: $g=g_0+\Upsilon^\flat\otimes X_2^\flat+X_2^\flat\otimes\Upsilon^\flat$
 where the vector field $\Upsilon$ is the electromagnetic potential and $X_2$ satisfies:

$g_0(X_2,X_2)=0$, ~~~ $g_0(X_2,Y)=1$, ~~~ $g_0(X_2,\Upsilon)=0$.

We suppose here, as for the study of geodesics (cf. \ref{ss1.5}), that hypothesis $H_E$ is satisfied, even more, in order to anticipate on section \ref{s2.13} where one takes into account the effects of "spin", we assume the following hypothesis $H'_E$. \\
\textbf{Hypothesis $H'_E$}
\begin{enumerate}
 \item The vector field $\Upsilon$ is defined on $\Theta\times S^3(\rho)$.
 \item The vector field $X_2$ is defined on $S^1(\delta)\times V$.
\end{enumerate}
This when the type cell is of the form $\C=\Theta\times S^1(\delta)\times S^3(\rho)\times V$ \\ (see section
\ref{s2.13}).

Whether under the assumption $H_E$ or $H'_E$, we will assume, as for proposition \ref{p1.4} with the field
$X_1$, that $D_{g_0}X_2=0$.
\begin{prop} \label{p1.5}
 We denote $F$ the $2$-differential form of electromagnetism:
$F=d(Y^{\flat_g})=d(\Upsilon^\flat)$.
 \begin{enumerate}
  \item $R_{{icc}_g}^\sharp=R_{{icc}_{g_0}}^\sharp+\frac{1}{2}(f(X_2\otimes X_2)-(\nabla_{g_0}\cdotp  F)\otimes
X_2-X_2\otimes (\nabla_{g_0}\cdotp  F))$

where the function $f:=\frac{1}{2}g_0^{ik}g_0^{jl}F_{kl}F_{ij}$ ~~ and ~~ $\nabla_{g_0}\cdotp  F$ ~ is the vector field 
"divergence of $F$" associated to $g_0$ ($\nabla_{g_0}\cdotp  F$ has for components $\nabla_i(g_0^{ik}F_k^j)$).
\item $S_g=S_{g_0}+(\Upsilon^k\Upsilon_k)R_{icc_{g_0}}(X_2,X_2)$.
\item $X_2$ is also a lightfield for $g$ ~ and ~ $D_g{X_2}=0$.
\item $Y$ is a Killing field and a geodesic field for $g$ (this is obviously true for $g_0$).
\end{enumerate}
\end{prop}

The proof of this proposition is detailed in annex \ref{a3.4}.

\begin{rmq} \label{r14}
 If $R_{icc_{g_0}}(X_2,X_2)=0$, the scalar curvature
$S_g$ remains equal to that of the neutral potential $g_0$ (cf. remark \ref{r13}). Again, in the standard experiments, $|(\Upsilon^k\Upsilon_k)R_{icc_{g_0}}(X_2,X_2) |\ll |S_{g_0}|$ ~ and ~ $S_g$ therefore remains very close to
$S_{g_0}$.
\end{rmq}


\section{Remarks on the Newtonian potential domain \label{s1.7}}
The cell considered here is of the form $\C=\Theta\times W$
with ~~ $\Theta=I \times \mathscr U \subset \R^4$.

According to proposition \ref{p1.2} the pseudo-Riemannian metric $g$ satisfies $g=g_0-2vX_1^\flat\otimes X_1^\flat$
~ and one chooses here $X_0=\frac{\partial}{\partial t}$ related to the standard coordinate system of the cell $\C$. We
further assumes that the assumption used for proposition \ref{p1.4} is still valid.

This area deserves a more precise study than the one we present, in particular
when the function $|v|$ is not $\ll 1$ (ie. when the Newtonian approximation is no longer valid). We will only make a few basic remarks.

Consider for example the case for which, in the standard coordinate system, \\ $v=-\frac{m}{r}$ ~ where ~ $m$ ~ is a
positive constant and ~ $r=(\sum_{k=1}^3{x^k}^2)^{1/2}$. \\ The domain has "spherical symmetry in
space" for the usual coordinates $(x_1,x_2,x_3)$ of $\mathscr U$ ~ (but keeping in mind that the others
"dimensions" are fundamental).

We have here $\Delta v=0$ ~~ (def. \ref{def:12}) ~~ and we have seen that (cf. paragraph \ref{ss1.4}), for the metric
$g=g_0+\frac{2m}{r}X_1^\flat\otimes X_1^\flat$ ~~ geodesics ~ $x_{(s)}=(x^0_{(s)},\dots,x^{n-1}_{(s)})$ ~ satisfy
(at least the ones we are interested in): ~~ $(x^1_{(s)},x^2_{(s)},x^3_{(s)})''=(-\nabla v)_{x(s)}=mr^{-2}_{x(s)}$.

We deduce, as in standard Newtonian physics, that the components on $\mathscr U$ of these geodesics haves their
images that are conics for which "$O$" is a focus, and that Kepler's laws remain valid, \textbf{but considering
the parameter $s$ and not the time $x^0=t$ of the coordinate system}, $s$ can be very different
of $t$ if the function $|v|$ is not $\ll 1$.

Since $g(\frac{\partial}{\partial t},\frac{\partial}{\partial t})=-1+\frac{2m}{r}$, the vector field
$\frac{\partial}{\partial t}$ ($=\frac{\partial}{\partial x^0}$) is timelike when $r>2m$, ~~ lightlike when
$r=2m$ ~~ and spacelike when $0<r<2m$. The critical distance is here $2m$ which corresponds exactly to the 
Schwarzschild radius and this suggests to compare this domain of "Newtonian potential type" that has just been studied, to a
 domain of "Schwarzschild type" whose cell will be defined as the following couple $(\C,g_S)$:\\
 $\C=\Theta\times W$
~~ ~~ where ~~ $\Theta=\R\times{\R^3}^*$ ~~ and ~~ $g_S$ is the product metric: $g_S:=g_{S_\Theta}\times
g_W$ ~ for which $g_{S_\Theta}$ is \textbf{the classical Schwarzschild metric}  which is written for
coordinates $(t,r,\varphi,\phi)$ ~ on ~ $\R\times]2m\text{~~}+\infty[\times S^2\sim\R\times {\R^3}^*$:
\[g_{S_\Theta}(t,r,\varphi,\phi)=(-1+\frac{2m}{r})dt^2+(1-\frac{2m}{r})^{-1}dr^2+r^2(d\varphi^2+\sin^2\varphi
d\phi^2)\]
Of course, there is mathematically no interest in making the product of $g_{S_\Theta}$ by $g_W$, it is
simply to bring back the classical Schwarzschild domain of dimension $4$ to the domains that we studied here of
dimension $n$. (We can also consider more generally the extended Schwarzschild domain for $0<r<2m$).

We now compare some properties of type (1): "Newtonian potential"

with those of type (2): "Schwarzschild" that we have just defined.
\begin{enumerate}
 \item The two Ricci curvatures: $R_{icc_g}$ ~~ and ~~ $R_{icc_{g_S}}$ are identical, both equal to
$R_{icc_{g_0}}$ ~ determined by $R_{{icc}_{{g_0}_W}}$ since $R_{{icc}_{{g_0}_\Theta}}=0$ (cf.
proposition \ref{p1.4}).
\item For $r \gg 2m$, according to proposition \ref{p1.2}, the components of the geodesics on $\Theta$ of type (1) that
we have described and those timelike of type (2), give back, with a very good precision, the trajectories of an
elementary object around a "mass" $m$ with spherical symmetry (in the space $\mathscr U$ calculated with the standard 
Newtonian physics.
\end{enumerate}
Note also that the coefficient $(-1+\frac{2m}{r})$ of $dt^2$ of the tensor $g$ of the type (1) is the same as
that of the tensor $g_S$ of the type (2), however, for the type (1) the potential ~ $\frac{2m}{r}$ ~ disturb the
compact manifolds without touching $\mathscr U$, whereas for type (2) the potential disturb $\mathscr U$ without
touching $W$.

These two types therefore appear as two special "extreme" cases and one can probably describe a family
of Newtonian potential types for which the potential disturbs both $\mathscr U$ and $W$ and such
that the properties 1. and 2. are preserved. These domains may have different properties when $r$ is not $\gg 2m$.


\section [A "really perfect static fluid"] {An example of a "realy perfect static fluid" type domain \label{s1.8}}
This static domain will be without electromagnetism and the apparent pressure will be zero (but not the hidden pressure).

The interest of presenting such an example is essentially the fact that such a domain can not exist in the context
of standard general relativity. The possibility of its existence here is due to the hidden pressure (which,
obviously, does not make sense when $dim \M=4$)

For the simplicity of the following calculations which will follow, we will present a very specific case.
\bigskip

The fluid type domain considered is a triplet $(\D,\tilde g,\A)$ (see definition \ref{def:4}) where $(\D,\tilde g)$
is isometric to $(\C,g)$ defined as follows:

$\C=Z\times W$ ~~~ where ~ $Z=\Theta\times S^1(\epsilon)$ ~~ and ~~ $\Theta=I \times \mathscr U\subset R^4$

$W=S^1(\delta)\times V$ ~~~ where $V$ is a compact manifold, but here this decomposition of $W$ will not be used
because the domain is assumed without electromagnetism and the usual circle $S^1(\delta)$ does not intervene.

(Circle $S^1(\epsilon)$ in the $Z$ decomposition of $\C$ was introduced only to simplify the following calculations).

One note $(t,x^1,x^2,x^3,v,w)$ the standard coordinate system of the cell $\C$~~ (the coordinate $v$ being 
that of $S^1(\epsilon)$).

The metric $g$ is defined by $g=g_Z\times g_W$, ~~ the signature of $g_W$ is $(-,+,\dots,+)$.

The metric $g_Z$ is defined by $g_Z=g'_0+\beta\otimes X_1^\flat+X_1^\flat\otimes\beta$

where ~~ $g'_0=g_\Theta\times g_{S^1(\epsilon)}$ ~~ ($g_\Theta$ being the Minkovski metric and ~ $g_{S^1(\epsilon)}$ ~ the
standard Riemannian metric on $S^1(\epsilon)$).

$\beta$ is a $1$-differential form defined on $\mathscr U$: $\beta=adx^1+bdx^2+cdx^3$ ~ where ~ $a,b,c:\mathscr
U\rightarrow\R$ are three regular functions.

$X_1$ is the vector field defined on $Z$ such that ~ $X_1^\flat=dt+dv$ ~~ (where $\flat$ is associated to $g'_0$).

($a,b,c$ ~ and ~ $X_1$ can, of course, be considered defined on the cell $\C$)
\bigskip

In the standard coordinate system, the components of $g$ do not depend on $t$, it is in this sense that
\textbf{the domain is qualified as static}.

The tensor $g$ is to be compared with that of a domain of "electromagnetic potential type", unlike
 that $X_1^\flat=dt+dv$ replaces $X_2^\flat=du+dv$, (this one would have given exactly an example of electromagnetic potential, here $u$ is the coordinate of $S^1(\delta)$). Since the considered functions are independent of
$t$, one can use the results obtained for an electromagnetic potential type domain  where
~ $(\Upsilon_0,\Upsilon_1,\Upsilon_2,\Upsilon_3)$ ~ becomes ~ $(0,a,b,c)$. ~ The roles of $t$ and $u$ are
transposed by the choice of $X_1$ instead of $X_2$. In fact, a calculation software gives quickly the
next result:

We denote ~ $(A,B,C):=rot(a,b,c)$ ~ which means ~ $A=(\frac{\partial b}{\partial x^3}-\frac{\partial c}{\partial x^2})$,
~~ $B=(\frac{\partial c}{\partial x^1}-\frac{\partial a}{\partial x^3})$, ~~ $C=(\frac{\partial a}{\partial
x^2}-\frac{\partial b}{\partial x^1})$.

We denote ~ $(\mathcal A,\mathcal B,\mathcal C):=rot(A,B,C)=rotrot(a,b,c)$.

We obtain for the matrix $(R_{{icc}_{g_Z}})$ of the Ricci curvature for the metric $g_Z$ associated to the system of
standard coordinates of $Z$ ~ (order matrix 5):
\bigskip

\[ (R_{{icc}_{g_Z}}) = \frac{1}{2}\begin{pmatrix}
                                         
  (A^2+B^2+C^2)&\mathcal A&\mathcal B&\mathcal C&(A^2+B^2+C^2) \\                                    
   \mathcal A&0&0&0&\mathcal A\\
   \mathcal B&0&0&0&\mathcal B\\
   \mathcal C&0&0&0&\mathcal C\\
   (A^2+B^2+C^2)&\mathcal A&\mathcal B&\mathcal C&(A^2+B^2+C^2)       
  \end{pmatrix}\]
\bigskip

And for the scalar curvature: $S_{g_Z}=0$
\bigskip

Since $g$ is a product metric $g:=g_Z\times g_W$, the Ricci curvature for $g$ is now determined by the Ricci curvature for $g_W$. Furthermore ~ $S_g=S_{g_Z}+S_{g_W}=S_{g_W}$.

We now assume that the scalar curvature $S_{g_W}$ is zero (this assumption is simplifying for the fluid type
 presented here but we can do without it (see remark \ref{r5}). Tensor $G$ then satisfies
 $G=2R_{{icc}_g}$.

To obtain a "really perfect" fluid one adds the following conditions:

$(\mathcal A,\mathcal B,\mathcal C)=rotrot(a,b,c)=0$ ~~ and ~~ $\mu:=A^2+B^2+C^2\neq0$.

Then ~ $G_{g_Z}=\mu X_1^\flat\otimes X_1^\flat$ ~~ where ~~ $\mu=|rot(a,b,c)|^2=A^2+B^2+C^2>0$.

In the standard coordinate system, the $g_Z$ matrix is:
\[ (g_Z) =\begin{pmatrix}
 -1&a&b&c&0\\
 a&1&0&0&a\\
 b&0&1&0&b\\
 c&0&0&1&c\\
 0&a&b&c&1
\end{pmatrix}\]

For any $x\in \C$, the apparent space $H_x$ is of dimension 4. It is easy to verify that the four vector fields
: $X_0:=\frac{\partial}{\partial t}$, ~~ $X_1:=a\frac{\partial}{\partial t}+\frac{\partial}{\partial
x^1}-a\frac{\partial}{\partial v}$ ~~ $X_2:=b\frac{\partial}{\partial t}+\frac{\partial}{\partial
x^2}-b\frac{\partial}{\partial v}$, ~~ $X_3:=c\frac{\partial}{\partial t}+\frac{\partial}{\partial
x^3}-c\frac{\partial}{\partial v}$ ~~ form in each point $x$ a $g_Z$-orthonormed basis of $H_x$.

As $G(X_0,X_0)=\mu$ ~~ and ~~ $G(X_i,X_j)=0$ ~ when $(i,j)\neq(0,0)$, ~ the matrix of $G_{H_x}$ \textbf{in this
basis} is:
\[ (G_{H_x}) =\begin{pmatrix}
 \mu&0&0&0\\
 0&0&0&0\\
 0&0&0&0\\
 0&0&0&0
\end{pmatrix}\]

According to the definitions given (see \ref{ss1.1.1}), the field ~ $X_0=\frac{\partial}{\partial t}$ ~ is the apparent field
of the fluid (here $X_0=X$ is also the fluid field) and $\mu$ the energy density function.

The apparent pressure $P_A$ is zero, the hidden pressure $P_H$ is the pressure $P$ itself and $P=G-\mu
X_0^\flat\otimes X_0^\flat$.

We then verify that $\nabla_g\cdotp  P=0$:

Since $\nabla_g\cdotp  G=0$,~ just verify that ~
$\nabla_g\cdotp  (\mu X_0^\flat\otimes X_0^\flat)=0$, ~ ie. $\nabla_i(\mu T^{ij})=0$ ~~ where ~~ $X_0\otimes
X_0:=T^{ij}\partial_i\partial_j$ ~ with ~ $T^{00}=1$ ~~ and ~~ $T^{ij}=0$ ~ whether ~ $(i,j)\neq(0,0)$.

We have:

$\nabla_i(\mu T^{ij})=\partial_i(\mu T^{ij})
+\mu(T^{lj}\Gamma^i_{il}+T^{il}\Gamma^j_{il})=\mu(T^{0j}\Gamma^i_{i0}+\Gamma^j_{00})$ ~ since ~ $\partial_i(\mu
T^{ij})=\partial_0(\mu T^{00})=\partial_0\mu=0$ ~ ($\mu$ does not depend on $t$).

$\Gamma^j_{00}=-g^{jl}(\partial_l g_{00})=0$ ~~ because ~~ $g_{00}==-1$

$\Gamma^i_{i0}=\frac{1}{2}g^{il}(\partial_ig_{l0}+\partial_0g_{li}-\partial_lg_{i0})=\frac{1}{2}g^{il}(\partial_ig_{l0}
-\partial_lg_{i0})=0$ ~~ because ~~ $g^{il}$ is symmetric in $(i,l)$ ~ and ~ $(\partial_ig_{l0}
-\partial_lg_{i0})$ ~~ antisymmetric.

The domain $(\C,g,\A)$ is thus a domain of type "really perfect fluid" (verifying other properties is easy), its apparent pressure is zero and the vector field of the fluid
$X_0=X=\frac{\partial}{\partial t}$ ~ is a geodesic field (for $g$ and $g_0$).
\vspace{2cm}

The potential type domains that have been presented, for which the pseudo-Riemannian tensor $g$ is given
explicitly, will be particularly important in the study of quantum phenomena that will be addressed
now, but it will not be the geodesics that will intervene, the process will be very different.

\chapter{Quantum phenomena}
\label{part:deux}

\vspace{-11 mm}

\section{Introduction}
\vspace{-3mm}

Pursuing what was done in the first chapter, we will now define domains
of space-time $(\M,g)$ which will represent what, classically, we call "particles in a
potential".
As the point of view that we are going to have on these notions is fundamentally different from that of standard quantum physics
, the term "particle" will no longer be appropriate and we will sometimes refer to it only to link
what we are going to present and the classical theories.

The essential difference between the areas that we will define and those presented in the first chapter will be due
to the fact that the \textbf{variations} of the pseudo-Riemannian tensor $g$ (defined on a cell
$\mathscr {C} = \Theta \times S^1 \times W$), associated to compact manifold
 $S^1 \times W $, will become the essential elements of the description of quantum phenomena.
In particular we will no longer consider that the field of vectors $Y$ tangent to the circles $S^1_x$ is a Killing field. If we had taken into account precisely in the first chapter the possible variations of $g$ on $S^1
\times W$, many additional terms would have appeared in the equations and would have rendered them
unmanageable, in the language of standard physics this would have meant that the equations would have described both classical and quantum phenomena. The approximations that we will make by defining the
geometric types that correctly describe
the quantum phenomena will make it possible to obtain usable equations, but will make their
"nature" different from those used in the first chapter.

For the reader to have a clear idea of the structure of this second chapter, we start with an
exhaustive presentation of what will follow:

The domains $(\D, g)$ that we will define will have, as pseudo-Riemannian metrics,  metrics $g$
constructed from a reference metric $g_0$ which will be that of a neutral potential that has already been defined
in the previous chapter. This metric $g_0$, defined on a cell
($\mathscr {C} = \Theta \times S^1(\delta) \times W$),
may be considered associated with an observer doing the
measures. The fact that, by definition, $g_0|_\Theta$ is Minkovski's metric specifies the chosen approximation.
We could use for $g_0|_\Theta$ a metric that takes into account a possible distortion of space-time
$\Theta$ associated to the observer, but this would make things more complicated and we will just seek to recover the standard results of quantum physics considering that $g_0|_\Theta$ is the Minkovski metric.

\textbf{The fundamental geometric type that will describe all quantum phenomena will simply be a domain
$(\D,g)$ of "constant scalar curvature and conform to a potential"}.

In other words, the pseudo-Riemannian tensor $g$ will be of the form $fg_{\mathcal P}$ where $g_{\mathcal P}$ is a metric
of potential as defined in the first chapter (from $g_0$), $f$ will be a non negative real function, and the
scalar curvature of $g$ will be constant, equal to that of $g_{\mathcal P}$ (itself equal to that of $g_0$).

The function $f$ will then become the important object of these domains and it will be easer to consider the function
$a:=f^{4/n-2}$ (where $n=dim\M$) because the definition of the chosen geometric type gives a differential equation "for
$a$" simpler than the equivalent equation written with $f$. It turns out that the non-linear part of
the differential equation "for $a$" will express, in language of classical physics (which will no longer be appropriate
here) \textbf{the interaction of particles with each other} in the considered field. As most of what we will
study will be devoted to phenomena that neglect the interactions of particles between them, it is
a \textbf{linear approximation} of the general equation we will be using (at least up to the section
\ref{s2.16}). The areas for which we will use
 this linear approximation will be called \textbf{domains of oscillating metric} because the 
fact that the metric $fg_{\mathcal P}$ remains at constant scalar curvature will, in most cases, lead to
strong oscillations for the function $f$.

The characteristics of the function $f$ (and therefore of the function $a$) imposed by the selected geometric type will
give back the standard notions
 of mass, electric charge, spin, etc.

\textbf{The linear equations obtained will therefore be necessary conclusions of the definition of the fundamental geometric type
selected (and the linear approximation). No law, no principle will be added. These linear equations
 will be of "Klein-Gordon" type (different according to the considered potentials). 
Schrödinger equations used in standard quantum physics (non-relativistic) will appear as approximations (which we will specify) of Klein-Gordon equations obtained. This will show later that one
recovers in particular the qualitative and quantitative results of the generic experiments of classical quantum physics
 (diffraction, Young's slits, potential deviation, Stern-Gerlach experiment, quantum entanglement, etc.)}.

In fact, the domains of oscillating metric type will not be enough to completely describe the notion of
"particles in a potential" and therefore the previously mentioned experiments. Therelack at this stage a notion of
localization (which, for example, will specify the notion of impact on a screen). It will be necessary to define
more
precisely domains that will be associated to the domains "particles in a potential" of standard physics. For that it will be enough "to add" in a domain $(\D, g)$ "with constant scalar curvature conform to a potential"
a subset $\Ss$ of $\D$ of zero measure (see section \ref{s2.11}). The connected components of $\Ss$ will be
called \textbf{the elementary singularities of $\D$}. In fact this subset can be considered as the one
where the tensor $g$ is not defined,
but this point of view also requires to clarify the behavior of $g$ in the neighborhood of $\Ss$. This last point
will not be useful for the study of the quantum phenomena presented in this paper but will have to be developed for the study of
more complex phenomena. The reader can therefore be content to consider the elementary singularities as
simple subset
 in $\D$. No law will be given on the behavior of these singularities and it will be assumed
only that these are "equiprobably distributed" relative to the metric $g$ (but not to $g_0$) induced in the spacelike submanifolds of maximum dimension $(n-2)$.

It should not be considered that one elementary singularity  represents a "particle" in the usual sense.
The good point of view is rather the following:

A singularity \textbf{in an oscillating metric type domain} whose characteristics of the function
$f ~(=a^{(n-2)/4})$
correspond to those of an electron could be associated to the standard notion of electron. If the characteristics of the
function $f$ correspond to those of a photon, the singularity \textbf{in this domain} will be associated to the standard notion
 of photon. When the singularity is in a potential type domain (for which $f=C^{te}$) it
will be undetectable, etc.

This point of view is fundamentally different from any standard quantum (or non-quantum) theory:
characteristics (mass, electric charge, spin, etc.) will be those \textbf{of the oscillating metric} in the neighborhood
in which the singularity is. These considerations express into the fact that the dynamics of the
particles in a potential will be entirely managed by the oscillating metric. The study of this will be the 
deterministic part of the description of the considered physical phenomena, since the differential equations deduced
will only have one solution for well-defined boundary conditions. The random part will be linked 
to the set $\Ss$ of the elementary singularities whose probability of presence in a domain will give, by
definition,
that of the position of the particles.

To clarify what we have just said we use these same notions in the context of a simple experiment expressed
in the classical language of physics:

Particles thrown with a known velocity vector enter a box in which vacuum has been created and where
may exist a potential. Particles leave impacts on a screen on the opposite side of
the inlet port of the particles. This experiment therefore occurs in a space-time domain $\D$ where, apart from the
domain inside the box (written $\D_B$) and during the experiment, the metric will be assumed to be that of a
 neutral potential $g_0$ (in dimension $n$ with two minus signs in the signature). The domain ($\D-\D_B,g_0$)
may be considered as the domain of an observer who makes measurements associated to the experiment (measurement of
position of the impacts on the screen for example). These measurements will necessarily be made relative to the metric
$g_0$ of the observer. For example, the measurement of the position of the impacts on the screen will be made by "transparency"
on the outside of the screen where the metric is actually $g_0$, or, after the experiment, by dismounting the
box and there again the measurement of the trace of the impacts will be made with the metric $g_0$ (in fact it is 
$g_0|_\Theta$ which is used, ie. Minkovski's metric). It is impossible for an observer to make 
any measurement inside the box during the experiment as the introduction of a physical object into
it would necessarily completely disrupt the experiment. With the point of view that is presented here, this 
simply means that the presence of another object in the box (oscillating metric with singularities, for example), 
completely changes the equations and the study then becomes that of another domain of space-time.

The domain $(\D_B,g)$ representing the inside of the box during the time of the experiment will therefore be of type
"particles in a potential". The conformal factor $a:= f^{(n-2)/4}$ of the metric conform to $g_{\mathcal P}$
corresponding
to the oscillating metric of this type of "particles" will satisfy a linear differential equation (relative to
 the metric $g_0$). Accuracy of the solution will derived from boundary conditions on the domain
$\D_B$. One of the conditions is, of course, the one that specifies that the particles enter in the box with a 
given velocity vector. For this it is necessary to introduce a third domain $\D'\subset{\D-\D_B}$ where is rigorously defined
 a type "particles in a neutral potential moving at a velocity $v$ in one direction". This particular type will be important because
it will serve as an initial condition for many experiments. Other boundary conditions
will reflect the fact that there are no reflections of the particles on the inside edges of the box, etc. Problems
associated to boundary conditions will be the same as those in classical quantum physics
with the state function. The oscillating metric will therefore be sufficiently determined. On the other hand, as we have already said,
the singularities will only be "equiprobably distributed" relative to the metric $g$.

To visualize the experiment, the reader can imagine that the inside of the box is "elastic" (as well as
the screen that we will consider having a thickness). During the preparation of the experiment the metric is everywhere
$g_0$. During the experiment, during the impact of the particles on the screen, the effective metric is the oscillating metric
 $g$ which is "$g_0$ disturbed by the conformal metric change", ie. $fg_0$. The screen is
 deformed
relatively to $g_0$. The position of singularities (and thus impacts on the screen) is equiprobable on the distorded screen
 (the effective metric is $g$). When measuring, after the experiment, the metric is $g_0$ (the screen has its
initial form), the position of the impacts is no longer equiprobable on the screen and it is clear that the probability law
on the position of impacts (for $g_0$) is completely determined by the distortion of space-time at
the inside of the box during the impacts, ie. by the oscillating metric $g$, and we foresee the close link between
the volume element for $g$ and the "presence probability density" of the singularities relative to the metric $g_0$.

In fact, all that we have just described happens in a domain of dimension $n>5$ that we can consider to be
$\Theta\times W$ where $\Theta$ is a open set of $\R^4$ and $W$ a compact manifold, equiped with a metric having a
signature $(-, +, +, +, -, +, \cdots, +)$. \\ $\Theta$ can be seen as the standard space-time and $W$ as the extra
small dimensions  (including some sort of time). The notion of "box", in such a space,
needs to be defined, but this does not cause any problems since it is enough to define it classically in $\Theta$, and
then make the Cartesian product with $W$. The box is considered here only as a boundary domain. All
tensors and differential operators that intervene in the description of this experiment depend, among
other things, of the "additional compact manifold variables". This is fundamental and it is what distinguishes
as we have already said, the description of the quantum phenomena from those described in the first chapter.
\begin{rmq} \label{++r2,1}
 More generally, a domain $(\D,g,\Ss)$ of type "particles in a potential" can be considered of the form $((\D_0\cup_{k=1}^m \D_k),g,\Ss)$ where the metric $g|_{\D_0}=g|_\mathcal{P}$ is the metric of the potential alone and the $g_k:=g|_{\D_k}$ are oscillating metrics to be specified ($\D_0$ and $\D_k$ are assumed to be connected and $\Ss$ defines singularities). Said more simply, such a domain consists of $m$ "bubbles" of oscillating metrics (the $(\D_k,g_k)$) in a potential type domain. Be careful, a bubble does not be view like a particle and a bubble can contain many singularities. When the oscillating metrics of each bubble $(\D_k,g_k)$ are of the same type (which is often the case in the context of standard experiments), the study of the physical phenomenon (for example the probabilistic distribution of the impacts on a screen) is performed on only one domain $(\D_k,g_k)$ because the (probabilistic) result is the same for each bubble.
\end{rmq}

Classification of oscillating metrics (which thus becomes the equivalent of the classification of particles into
classical physics) will be reduced to the choice of particular terms in the spectral decomposition of the function $a$
associated to the
Dalembertian operator of the compact manifold ($W,g_W$). Notions usually called "interaction" (between particles
or (and) with potentials) classically managed by Lagrangian-Hamiltonian choices and then by a procedure of
"quantification", will be managed here only by the very simple axiomatic system that has already been specified (constant scalar curvature domain conform to a potential) and by the choice of the terms of the spectral decomposition. No law, no
principle will be added. The study of the spectral decomposition of a function defined on ($W,g_W$) is complex and is
obviously completely associated to the precision of the data of ($W,g_W$). We will present here essentially the study of
particle interactions with a potential (in other words, the results given by the linear approximation),
this one including the notion of "spin" and thus naturally the description of Stern -Gerlach type experiments. The particular phenomena that involve the notion of spin that are "quantum entanglement" and the "anomalous magnetic moment of the electron" are treated in sections \ref{+2.2} and \ref{++2,3}. The description of more complex quantum phenomena (discussed
currently by quantum field theory) will be briefly presented in section \ref{s2.17}.

Let's turn now to the precise statements of the definitions, the results obtained and their proofs.


\section{The fundamental geometric type \label{s2.1}}
The geometric type that we will define now will allow the complete study of what is classically called
"particles in a potential". The sole application of its definition will be enough to recover the standard results
 of quantum physics on the subject. We will not go much further for this study in this paper, but this
geometric domain will also allow to address the phenomena currently treated by the quantum field theory (which we will write "Q.F.T" later) and this with a remarkable simplicity of the axiomatic system. We specify
for the reader familiar with Q.F.T, that the following definitions are more closely linked to the notion of
"fields" of the Q.F.T. The "quanta" of these fields will be associated to the "singularities" defined in the
section \ref{s2.11}.

The potentials that we will consider are those defined in (\ref{ss1.2}). They will be of constant scalar curvature. The pseudo-Riemannian metric $g$ of a domain $\D$ that will characterize the fundamental geometric type
 will be defined by two properties:
\begin{enumerate}
 \item metric $g$ will be conform to a potential metric
$g_\mathcal P$, that is, it will be of the form $fg_{\mathcal P}$ where $f$ is a non negative real function.

\item The scalar curvature $S_g$ will be kept equal to $S_{g_0}$ (so constant).
\end{enumerate}

It is recalled that a conformal transformation of a pseudo-Riemannian metric $g$ is a transformation that does not
privileged
direction (whether of space or time), it consists in multiplying, for each point $x$ of $\M$, the quadratic form
 $g_x$ of the tangent space $T_x(\M)$ by a real number $f(x) \geq 0$.

It turns out to be simpler, when we make a conformal metric change $fg_{\mathcal P}$, to write $f$ in the form
${\left|a\right|}^{4/n-2}$ (where $n = \text{dim} \M$) because the differential equation "for $a$" deduced from the previous properties, has a simpler form than that written "for $f$".

The precise definition of the geometric domain which, for us, will describe all the quantum phenomena is therefore
the next one.

\begin{dfn}
\label{d2.1}
 \textbf{A domain with constant scalar curvature, conform to a potential} is a domain $\D$ of $\M$ that
satisfies the following two properties:
\begin{enumerate}
 \item For any $x\in\D$ there is a chart
($\mathscr{V},\varphi$) around $x$ of the $g$-observation atlas, and a regular function
$a:\varphi(\mathscr{V})=\mathscr {C}=\Theta \times S^1 \times W\rightarrow \R$ such that
 $\varphi:(\mathscr{V},g_{\M})\rightarrow (\mathscr{C},{\left|a\right|}^{4/n-2}g_{\mathcal P})$ is an
isometry.

Here $g_{\mathcal P}$ is a metric representing a potential.
 \item The scalar curvature $S_{g_{\M}}$ is equal to $S_{g_{\mathcal P}}$ (and so is a constant equal to
$S_{g_0}$ for the considered potentials).
 \end{enumerate}
 \end{dfn}
 Condition $2.$ of this definition is actually a condition of \textbf{normalization} of the function $a$ in the
conformal metric change. Multiplying $a$ by a positive constant $\lambda$ causes the scalar curvature
$S_{g'}$,  where $g'=|a|^{4/n-2}g_\mathcal P$,  to be multiplied by $\lambda^{-4/n-2}$. We could have
choosen another form of normalization that leaves $S_{g'}$ constant but different from $S_{g_{\mathcal P}}$. The
condition $2.$ seemed to me the simplest choice. The importance of the choice of normalization will appear in the section
\ref{s2.12}.

\begin{rmq} \label{r1'}
This remark is to be compared with the remark \ref{r5} of the first chapter. We can define a wider class
of domains "with constant scalar curvature conform to a potential" when the cell $\C$ of the definition
\ref{d2.1} is of the form $\mathscr {C} = \Theta \times S^1 \times V_1\times V_2$ and the potentials are defined by
a product metric $g_{\mathcal P} = g'_{\mathcal P}\times g_{V_2}$ where $g'_{\mathcal P}$ is a potential metric
on $\mathscr{C'} = \Theta \times S^1 \times V_1$, this considering conformal transformations of
$g'_{\mathcal P}$ and not $g_{\mathcal P}$. Precisely, condition $1.$ of definition \ref{d2.1} is replaced
by: For any $x\in\D$ there is a chart ($\mathscr{V},\varphi$) around $x$ of the $g$-observation atlas, there exists
a regular function $a:\mathscr {C'}=\Theta \times S^1 \times V_1\rightarrow \R$ such that
$\varphi:(\mathscr{V},g_{\M})\rightarrow (\mathscr{C},({\left|a\right|}^{4/n-2}g'_{\mathcal P})\times g_{V_2})$ is
an isometry. Note that with this definition, some space directions are preferred. Of course, all the results we will get from definition \ref{d2.1} can be applied to this new
definition considering that the manifold is of dimension $(n-dim V_2)$. The fact that
${\left|a\right|}^{4/n-2}g'_{\mathcal P}\times g_{V_2}$ is a product metric is then simple to manage. However,
the scalar curvature $S_{g_{\mathcal P}}$ is different from $S_{g'_{\mathcal P}}$ if $S_{V_2}$ is $\not=0$, (since
$S_{g_{\mathcal P}} 	= S_{g'_{\mathcal P}}+S_{V_2}$), and this may be important because we will see that
the notion of "mass" will be influenced by the scalar curvature (def. \ref{d2.8}). Extending the definition
\ref{d2.1}
as we have just presented it, will therefore allow more possibilities associated to the notion of mass.
\end{rmq}


\section [The fundamental equation ...] {The fundamental equation of a domain "with constant scalar curvature, conform to a potential "\label{s2.2}}
We know that the  transformation law of the scalar curvature, when a conformal metric change is made
$g'={\left|a\right|}^{4/n-2}g$, is given by Yamabe's equation:

\begin{eqnarray}
\frac{4(n-1)}{n-2}\Box_ga+S_ga=S_{g'}{\left|a\right|}^{4/n-2}a\label{F0}
\end{eqnarray}

where $S_g$ (resp $S_{g'}$) denotes the scalar curvature relative to $g$ (resp $g'$) and
$\Box_g:=-{\nabla_g}^i{\nabla_g}_i$
is the (geometric) Dalembertian operator associated to $g$.

Here the function $a$ can change sign. On the set where $a$ vanishes the curvature relative to $g'$ is not
defined. On this set, equation \ref{F0} is reduced to $0 = 0$.

Given definition \ref{d2.1}, the fundamental equation satisfied by the function $a$, for a domain of "constant scalar     curvature conform to a potential", is therefore the following:

\begin{eqnarray}
\dfrac{4(n-1)}{n-2}\Box_{g_{\mathcal P}}a+S_{g_{\mathcal P}}a=S_{g_\mathcal P}{\left|a\right|}^{4/n-2}a\label{F0'}
\end{eqnarray}
where $S_{g_{\mathcal P}}$ is actually $S_{g_0}$ which is constant for the considered potentials.

When $S_{g_{\mathcal P}}\not=0$, this equation is nonlinear, which makes its use very delicate,
however we will see that in many cases corresponding to classical experiments, the function $a$ will be
$\ll 1$
(relative to the normalization given in the definition \ref{d2.1} - $2.$). Equation \ref{F0'} can then be
approximated by the equation where the right hand side is zero. This is specified in the next section. (If the function is 
 close to 1 equation \ref{F0'} can also be approximed by a linear equation, this situation will be presented in \ref{4-ss-3}).

\section [The linear approximation] {The linear approximation. \\
  Oscillating metric in potential domains \label{s2.3}}

We are now interested in the domains of "constant scalar curvature, conform to a potential" for
which solutions of the fundamental equation \ref{F0'} are very close to the solutions of the associated
\textbf{linear} equation:
$$\dfrac{4(n-1)}{n-2}\Box_{g_\mathcal P}a+S_{g_0}a=0$$
this for specified boundary conditions (fundamental for the validity of this approximation). We shall
see below (section \ref{s2.12}) that this approximation corresponds to experiments for which, in the
classical language, we neglect the interaction of particles between them (but not with the potential of course). The fact
that the function $a$ remains $\ll 1$ will be the translation of a "low particle density" in the considered experiment. The non-linearity of equation \ref{F0'} ensures that, when $a \ll 1$, the term
$S_{g_\mathcal P}{\left|a\right|}^{4/n-2}$ can be
absorbed by the term $S_{g_\mathcal P}$ of the left hand side of the equation by modifying it very little. Of course,
this simple fact is not mathematically sufficient for the solutions of the linear equation to be close to
those of the equation \ref{F0'}, things have to be specified. (Annex \ref{ea3.4} presents a very simple example of approximation of solutions of a nonlinear equation by those of an associated linear equation, in
the spirit of what has just been presented, which can help understanding the process).

We give the following definition.

\begin{dfn}
\label{d2.2}
 An oscillating metric domain \textbf{in a potential} is a domain "with constant scalar curvature,
conform to a potential" for which the condition $2.$ of definition \ref{d2.1} is replaced by the fact that the
function $a$ satisfies:

\begin{eqnarray}
\Box_{g_\mathcal P}a+Sa=0\label{F1}
\end{eqnarray}
where we set $S=\frac{4(n-1)}{n-2}S_{g_0}$
\end{dfn}
\textbf{Equation \ref{F1} thus becomes the fundamental equation of an oscillating metric domain in a
potential}.

The "oscillating metric" terminology is justified by the fact that in most cases, as we shall see
later, definition \ref{d2.2} shows that the function $a:\mathscr {C}=\Theta
\times S^1 \times W\rightarrow \R$ must have strong oscillations in $t\in\R$ and possibly in $u\in S^1$.
The frequency of the oscillations in $t\in\R$ will be associated to the notion of "mass" and the oscillations in $u\in
S^1$ to the notion of "electric charge" (defs. \ref{d2.6} and \ref{d2.8}).

The main objective of this second chapter is to recover the  
classical quantum physics results that describe
standard experiments (diffraction, Young's slits, influence of potentials, Stern-Gerlach experiment,
quantum entanglement, etc.).
For these results, the interactions of the particles with each other are neglected (we can consider them "one by
one" except for the phenomena of quantum entanglement) and it is therefore equation \ref{F1} that we are going to
use.

\begin{rmq}
 The fundamental equation \ref{F1} can be considered as the linear approximation of a more general equation than
\ref{F0'}. It is not necessary in definition \ref{d2.1} to assume that the conformal metric change
 preserves the scalar curvature of ${g_\mathcal P}$ (nor even leaves it constant). It is sufficient, for the validity of
\ref{F1}, to bound $S_{g'}$, ~ where $g'=|a|^{4/n-2}g_\mathcal P$, this to normalize the function $a$. Requiring the constancy of the scalar curvature will be of interest only to make the parallel with some studies
presented in Q.F.T.
\end{rmq}


\section{Classification of oscillating metrics \label{s2.4}}
As we have already said, a domain with constant scalar curvature, conform to a potential is characterized by two
"objects":

\begin{enumerate}
\item The potential given by: $g_{\mathcal P} = g_0 + h$ ~~~ (see \ref{ss1.2}).
\item The function $a:\mathscr{C} = \Theta \times S^1(\delta) \times W\rightarrow \R$.
\end{enumerate}

For any $x\in\Theta$, the function $a_x(.):\mathscr {C}= S^1(\delta) \times W\rightarrow \R$, given by
$a_x(u,v):=a(x,u,v)$, is defined on the \textbf{compact} pseudo-Riemannian manifold $S^1(\delta) \times W$. It
therefore admits a spectral decomposition associated to the Dalembertian operator $\Box_{g_0{_{S^1(\delta) \times W}}}$. In fact, given the signature $(-,+,+,\dots,+)$ of $g_0|_{S^1(\delta) \times W}$, it will be better to consider independently
the spectral decompositions of $a_{(x,u)}(.): W\rightarrow \R$ and $a_{(x,w)}(.):S^1(\delta) \rightarrow \R$
associated to the \textbf{Riemannian Laplaciens operator} of $(W,g|_W)$ and $(S^1(\delta),g_0|_{S^1(\delta)})$. This will be specified
in the next paragraph.

The principle of decomposition that we will use is the same as that which consists in decomposing "sounds" in
"pure sounds", which mathematically are expressed by the Fourier series decomposition of
periodic functions defined on $\R$. Periodic functions defined on $\R$ identify with functions
defined on a circle $S^1$ and the later is a compact manifold. Compactness is fundamental for validity
of a discrete spectral theory. Obtaining discrete quantities in the following theory, which
justifies the term "quantum", will come from the compactness of the manifold $S^1(\delta) \times W$ in the cell.
Of course, the spectral decomposition that we will use is much more complex than for the case of a circle. The classification of the oscillating metrics that will be deduced from the spectral decomposition
is similar to the classification of particles in standard physical theories.

We begin by recalling some results of spectral theory which concern the Riemannian compact manifolds.
These results will be applied to the manifolds $(W,g_0|_W)$ and $(S^1(\delta),g_0|_{S^1(\delta)})$.


\subsection{Some results of spectral theory on Riemannian compact manifolds} \label{ssn1}

\textbf{Spectral Theorem}: We consider a \textbf{compact} Riemannian manifold $(V,g)$ and the associated (geometric) Laplacian operator: $\Delta:=-{\nabla}^i{\nabla}_i$.

\begin{enumerate}
 \item The Laplacian eigenvalues form an increasing sequence of non negative real numbers converging to $+\infty$:

$0 = \lambda_0 < \lambda_1 < \lambda_2 <\dots < \lambda_n < \dots$
\item For each eigenvalue $\lambda_i$, the corresponding eigenspace $E_{\lambda_i}$ is of finite dimension and,
for any $i\neq j$, $E_{\lambda_i}$ and $E_{\lambda_j}$ are orthogonal for the standard scalar product of
$L^2(V,g)$.

($\dim E_0=1,~E_0$ being the set of constant functions on $V$)
\item \textbf{The algebraic sum of the eigenspaces $E_{\lambda_i}$ is dense in $C^\infty(V)$ with the uniform topology. In particular, the Hilbert space $L^2(V,g)$ has an hilbertian basis of eigenfunctions}.
\end{enumerate}

The following simple proposition is about the Dalembertian operator on a pseudo-Riemannian manifold.

\begin{prop} \label{p2.1'}
We consider $k$ compact Riemannian manifolds $(V_i,g_i)$.

On the product manifold $V=V_1\times...\times V_k$ one defines the \textbf{pseudo}-riemannian tensor: \\
$g = (-g_1) \times (-g_2) \times \dots (-g_p) \times (g_{p+1}) \dots \times (g_k)$.

The associated Dalembertian operator is defined by $\Box=-{\nabla^i}{\nabla_i}$. ~~~ (If $p=0,~~\Box=\Delta$).

We consider the  eigenspaces $E_{\lambda_1},..,E_{\lambda_k}$ respectively associated with the Laplacians operator of 
$(V_1,g_1), \dots ,(V_k,g_k)$.

Then the set $F$ consisting of "finite sums of products" $f_1f_2 \dots f_k$, where the $f_i~\in~E_{\lambda_i}$,
is a
subspace of the eigenspace $E_\lambda$ associated to the Dalembertian operator $\Box$ of the 
pseudo-Riemannian manifold $(V,g)$ associated with the eigenvalue
$\lambda = -\lambda_1\dots-\lambda_p+\lambda_{p+1}+\dots+\lambda_k$.

(Functions $f_i$ are here considered to be defined from $V$ to $\R$, the abuse of notation used is to identify
$f_i:{V_i}\rightarrow \R$ and $f_i\circ p_i:V\rightarrow \R$ (where $p_i:V\rightarrow V_i$ is the canonical projection), this
abuse will be frequent in the following).
\end{prop}

In other words, \textbf{$F$ canonically identifies with the tensor product $E_{\lambda_1}\otimes\dots \otimes
E_{\lambda_k}$} that
 will be considered as a vector subspace of the eigenspace $E_\lambda$.

\begin{rmq}
It can be proven that if $\lambda$ is an eigenvalue of the Dalembertian operator of $(V,g)$ and if there exists a \textbf{unique}
$k$-uple
$(\lambda_1,\dots,\lambda_k)$ of eigenvalues for $(V_i,g_i)$ such that
$\lambda=-\lambda_1\dots-\lambda_p+\lambda_{p+1}+\dots+\lambda_k$, then the eigenspace $E_\lambda$ is exactly
$E_{\lambda_1}\otimes\dots \otimes E_{\lambda_k}$. On the other hand, care should be taken as when the 
pseudo-Riemannian manifold $(V,g)$ is not Riemannian, there may be specific spaces associated to the Dalemberian operator that
are of infinite dimension.
\end{rmq}
The manifold $V$ that we consider here is $S^1(\delta)\times W$ and the pseudo-Riemannian metric $g_0|_V$ is the
product metric $(g_0|_{S^1(\delta)}\times g|_W)$ with signature $(-,+,..,+)$. The
function
$a:\mathscr {C}=\Theta \times {S^1(\delta)} \times W\rightarrow \R$, which specifies the 
oscillating metric, thus satisfies the following properties according to the spectral theorem:
\begin{enumerate}
 \item For any $(x,u)\in\Theta\times{S^1(\delta)}$, the function $a_{(x,u)}(.):W\rightarrow\R$ admits the
spectral decomposition $a_{(x,u)}(.)=\sum_{i=1}^{\infty}\varphi_i{_{(x,u)}}\alpha_i(.)$ where $(\alpha_i)$ is an 
orthonormal hilbertian basis of eigenfunctions of $(W,g)$.

\item For any $i\in \N^*$, for any $x\in\Theta$, the function
$\varphi_{i,x}(.):S^1(\delta)\rightarrow\R$ admits the spectral decomposition (Fourier decomposition):

$\varphi_{i,x}(u)=\sum_{j=0}^\infty \zeta_{1,i,j}{_{(x)}}\cos(2\pi ju/\delta)+\zeta_{2,i,j}{_{(x)}}\sin(2\pi
ju/\delta)$.
\end{enumerate}
\textbf{Classifying oscillating metrics will consist in taking only some terms of the spectral decomposition of
the function $a$}. This will allow, among other things, to define characteristic constants for each chosen elementary
oscillating metric. These constants will give the notions of "mass", "electric charge",
"spin", etc., associated with the considered elementary oscillating metric, which will then correspond, when we
have clarified the notion of "singularities" (section \ref{s2.11}), to the standard notions linked to "particles".

The choice of the definition that we will present now is essentially guided by the fact that we want
to recover the standard classifications on particles. It takes into account that, given the signature of $g$,
$S^1(\delta)$ does not play the same role as $W$.


\subsection{Elementary oscillating metrics}
\textbf{Elementary} oscillating Metrics can be seen as "pure sounds" in the decomposition of a
" periodic sound ".

We consider an eigenspace of the Dalembertian operator $\Box_{g_0{_{S^1(\delta) \times W}}}$ of the form
\\ $E_{\lambda,\mu}:=E_{S^1(\delta)}(\lambda)\otimes E_W(\mu)$.

\begin{dfn} \label{d2.3}
 \textbf{An elementary oscillating metric domain in a potential associated with $E_{\lambda,\mu}$} is
an oscillating metric domain in a potential (definition \ref{d2.2}) for which the function
$a:\C\rightarrow\R$ satisfies: for any $x\in\Theta$, \\ $$a_{(x)}(.)\in E_{\lambda,\mu}$$.
\end{dfn}

If $\lambda>0$, the eigenspace $E_{S^1(\delta)}(\lambda)$ is of dimension 2.

It is of dimension 1 when $\lambda=0$.

We shall often choose the natural basis:
\begin{enumerate}
 \item

 $(\alpha_1,\alpha_2)$ ~~ if $\lambda>0$, ~~~~ where ~ $\alpha_1(u):=\cos(\sqrt{\lambda}~u)$ and
$\alpha_2(u):=\sin(\sqrt{\lambda}~u)$.
 \item $(\alpha_1=1)$ when $\lambda=0$.
\end{enumerate}
(We will then identify $E_{S^1(\delta)}(\lambda)$ to 
$\varmathbb{C}$ if $\lambda \neq 0$ ~~ and ~~ $E_{S^1(\delta)}(0)$ to $\R$).

\medskip
If we choose a basis $(\beta_1,\dots,\beta_k)$ of $E_W(\mu)$ ~~ ($L^2_W$-orthonormal, for example), the
function $a$ then satisfies: for any $x\in\Theta$,
\begin{eqnarray}
  a_x(.)=\sum_{i=1}^k(\varphi_{1,i}{_{(x)}}~\alpha_1+\varphi_{2,i}{_{(x)}}~\alpha_2)~\beta_i\label{f4'}
\end{eqnarray}
When the basis of the eigenspaces are specified, the function $a$ is entirely determined by the
$2k$ functions ~ $\varphi_{1,i}$ and ~ $\varphi_{2,i}$ ~ ($k$ real functions if $\lambda=0$).

We will see in section \ref{s2.5'} that, using natural approximations, we can limit the
number of functions $\varphi_{1,i}$ and ~ $\varphi_{2,i}$ that determine the function $a$ in decomposition \ref{f4'}
and thus refine the classification of elementary oscillating metrics. Before that, we present the first
important features of elementary oscillating metrics.


\section [Associated constants] {Constants associated with elementary oscillating metric \label{s2.5}}
We introduce here in particular the notion of \textbf{electric charge} and \textbf{mass} of an elementary oscillating metric. Of course, the terminology is chosen so that these notions then correspond exactly to those of 
standard physics.

According to definition \ref{d2.3} the function $a:\C\rightarrow\R$ which characterizes the elementary oscillating metric
associated with $E_{\lambda,\mu}$ satisfies: for any $x\in\Theta$ 
$$a_x(.)\in E_{\lambda,\mu}:=E_{S^1(\delta)}(\lambda)\otimes E_W(\mu)$$
\textbf{The important constants we are going to define are none other than those built from the eigenvalues 
$\lambda$ and $\mu$ of $E_{\lambda,\mu}$}. They will therefore be invariant by change of charts in the observation atlas
 built from the Poincaré transformations (change of charts "leave fixed", by definition, the
compact manifold $S^1(\delta)\times W$ so eigenspace $E_{\lambda,\mu}$).


\subsubsection{Absolute electric charge}
A natural basis of the eigenspace $E_{S^1(\delta)}(\lambda)$, which has been used in decomposition \ref{F4'} of
$a_x(.)$, is given by:

$(\alpha_1,\alpha_2)$ if $\lambda>0$ ~~ (where $\alpha_1(u)=\cos({\sqrt\lambda}u)$ and $\alpha_2(u)=\sin({\sqrt\lambda}u)$)

and by ($\alpha_1=1$) when $\lambda=0$.

\begin{dfn} \label{d2.6}
 \textbf{The electric charge frequency} of the elementary oscillating metric associated with $E_{\lambda,\mu}$, is the
 non negative constant $Q^+:=\sqrt{\lambda}$.
\end{dfn}
To recover the standard notion of "electric charge" expressed in the S.I units, we give the following definition:

\begin{dfn} \label{d2.7}
 \textbf{The absolute electric charge} of the elementary oscillating metric associated with $E_{\lambda,\mu}$ is the
constant $q^+:=\hbar Q^+$ ~~ where $\hbar$ is the Planck constant.
\end{dfn}
Definition of \textbf{relative} electric charge (positive or non positive) will be presented later for  more
restrictive oscillating metrics (def. \ref{d2.14}).

The Planck constant appears here only as a factor that reduces a frequency (expressed in geometrical units) to
a S.I unit of electric charge.

Regarding this matter, remember that, in geometric units, a mass is a "length", an electric charge is
a "length", Planck's constant is "the square of a length", a frequency is the "inverse of a
length", etc.


\subsubsection{Mass}
The use of the term "mass" (as well as that of "electric charge") will become clearer 
when translating fundamental equation \ref{F1} into Klein-Gordon equations and then
Schrödinger equations (thms. \ref{2.1} and \ref{2.2}). The definition given here is specific to elementary oscillating metrics and, contrary to the notion of absolute electric charge, several variants will be considered for extensions
 to more general frames but will be briefly discussed in this paper (see section \ref{s2.17}).

\begin{dfn} \label{d2.8}
 \textbf{The mass frequency} of the elementary oscillating metric associated with $E_{\lambda,\mu}$ is the non negative constant $M$ which satisfies:
$$M^2=S+\mu-\lambda$$
where ~~ $S:=\frac{n-2}{4(n-1)}S_{g_0}$, ~~ $S_{g_0}$ is the constant scalar curvature of the metric $g_0$ (and of the
metrics $g_\mathcal P$ associated with potentials).
\end{dfn}

The mass frequency is obviously defined here only when $S+\mu-\lambda$ is non negative (the cases for
which $S+\mu-\lambda$ is negative will be studied when the notion of "life time" intervenes, but this subject
will only be briefly discussed in this paper (see section \ref{s2.17})).

To recover the notion of "mass" expressed in the S.I units, we give the following definition.

\begin{dfn} \label{d2.10}
 \textbf{The mass} of the elementary oscillating metric associated with $E_{\lambda,\mu}$ is the constant $m:=\hbar
c^{-1}M$ ~~ where $c$ is the speed of light.
\end{dfn}

Note that the scalar curvature $S_{g_0}$ appears in the definition of the mass added to the value
$\mu-\lambda$: this is to recover, as we have already said, the standard notion of mass of classical physics.
But we will see (section \ref{s2.15}) that this fact can be associated to the notion of "Higgs field", (the scalar curvature, when it is positive "gives" mass to oscillating metrics for which $\mu-\lambda$ is negative).

\begin{rmq} \label{r2.3}
 By its definition, the electric charge frequency is an integer multiple of $\delta^{-1}$ ~ where $\delta$ is the
radius of the circle $S^1(\delta)$. The absolute electric charge $q^+$ is therefore an integer multiple of an elementary charge.
On the other hand, the mass of such an oscillating metric has a set of possible values much more complex. (It is a discrete set for a fixed manifold ($S^1(\delta)\times W,g_0|_{S^1(\delta)\times W)}$)).
\end{rmq}

When the compact manifold $W$ is a product $V_1\times\dots\times V_k$, the
eigenvalues of the eigenspaces $E_{V_k}(\mu_k)$ are characteristics of the oscillating metric that can
be important. This will be the case for a decomposition $W=S^3(\rho)\times V$ for which the eigenspaces
$E_{S^3(\rho)}(\gamma)$ will give the classification in terms of "spin" of the elementary oscillating metrics. This
will be detailed in section \ref{s2.13}.


\section [Refinement of Classification] {Refinement of oscillating metric classification.\\
Elementary oscillating metrics of order 1 and 2 in a potential \label{s2.5'}}

The relevance of the physical results obtained through the study of elementary oscillating metrics is, of course, associated to the choise of the compact Riemannian manifold $(W,g_W)$. This choice is determined little by little, guided
by the desire to describe as finely as possible the experimental results. To be able to write 
sufficiently precise theorems it is necessary to know how to neglect the "unknown part" of the manifold $(W,g_W)$.
\vspace{5mm}

To recover results of standard quantum physics that do not take into account the notion of spin, we
neglect (with precise definition (def. \ref{d2.4})) the "quantum effects" associated to the manifold
 $(W,g_W)$ and only the "quantum effects" associated to $(S^1(\delta),g_{0_{S^1}})$ will be considered. \\
Elementary oscillating metrics corresponding to this situation will be said of \textbf{order 1} (def.
\ref{d2.9} -1.) And we will see that it will be enough then to use only \textbf{two} real functions (or equivalently a
complex function) to completely determine the equation that satisfies the function $a$.
\vspace{3mm}

To recover the results that take into account the notion of spin, we will consider that the manifold
$(W,g_W)$ is a product:

$(W,g_W)=(S^3(\rho)\times V,g_{0_{S^3}}\times g_V)$ \\ where $(S^3(\rho),g_{0_{S^3}})$
is the standard Riemannian sphere of dimension 3 and radius $\rho$. A precise description of this situation
can be made by neglecting the quantum effects associated to $(V,g_V)$ and considering only those
associated to $(S^1(\delta)\times S^3(\rho), g_{0_{S^1\times S^3}})$. \\
Elementary oscillating metrics corresponding to this situation will be said \textbf{of order 2} (def.
\ref{d2.9} -2.). The number of real functions that will be sufficient to determine the equation satisfied by the function
$a$ will be limited here by the dimension of the eigenspaces of the Laplacian operator on the sphere $(S^3(\rho),g_{0_{S^3}})$.
\vspace{3mm}

An increasingly fine description of quantum effects can be pursued by continuing to specify the
properties of the manifold $(W,g_W)$. This precision will not necessarily be obtained by decomposing the manifold $V$ under
product form, but if this is the case, we will define elementary oscillating metrics of order $k>2$.
\vspace{5mm}

Definition \ref{d2.4} specifies what has just been presented, i.e. the way by which
we "neglects" some "quantum effects" associated to the manifold $(W,g_W)$. It concerns the metrics representing the
potentials and it applies in the general context for which the compact manifold $(W,g_W)$ decomposes in the form
$(V_1,g_{V_1})\times (V_2,g_{V_2})$. Oscillating metrics of order 1 just mentioned (specified in 
definition \ref{d2.9} -1.) are those for which $V_1$ is reduced to a singleton (uninteresting) and $V_2=W$. \\
Oscillating metric of order 2 (specified in definition \ref{d2.9} -2.) are those for which $V_1=S^3(\rho)$ ~~ (and
$V_2$ is written $V$).
\vspace{5mm}

We consider a  cell $\C=\Theta\times S^1(\delta)\times W$ with $W=V_1\times V_2$ where $V_1$ and
$V_2$ are two compact Riemannian manifolds ($V_1$ may be of 0-dimension). \\
Let:

- $g_0$ be a neutral potential metric on $\C$: $g_0=g_\Theta\times (-g_{S^1(\delta)})\times g_{V_1} \times
g_{V_2}$.

- $g_{\mathcal P}$ be a metric representing an active potential: $g_{\mathcal P}=g_0+h$ ~~ (subsection \ref{ss1.3}).

- $E_{V_2}(\mu_2)$ be an eigenspace of the Laplacian operator on the manifold $(V_2,g_{V_2})$ with eigenvalue $\mu_2$.

\begin{dfn} \label{d2.4}
 A potential (of metric $g_{\mathcal P}$) \textbf{is neutral on $E_{V_2}(\mu_2)$ (or, more simply neutral on $V_2$)} if:  for any 
functions $\varphi:\Theta\times S^1(\delta)\times V_1\rightarrow\R$ ~~ and ~~ $\beta\in E_{V_2}(\mu_2)$:
\begin{eqnarray}
\Box_{g_\mathcal P}(\varphi\beta)=(\Box_{g_\mathcal P}\varphi+\mu_2\varphi)\beta \label{F4'}
\end{eqnarray}
($\varphi$ and $\beta$ can be considered defined on $\C$)
\end{dfn}
Obviously, equality \ref{F4'} is always satisfied for the metric $g_0$, this being a "product" metric.
It is also the case for any potential metric $g_\mathcal{P}$ if $\mu_ 2=0$. This definition reflects the fact that
the metric $g_\mathcal{P}$ has no more influence on the eigenspace $E_{V_2}(\mu_2)$ than the metric $g_0$
(or at least the difference in influence is negligible).
\vspace{5mm}

If one considers an elementary oscillating metric domain in a potential associated with $E_{\lambda,\mu}$ (definition \ref{d2.3}), for which ~~ $W=V_1\times V_2$ ~~ and $E_{\lambda,\mu}=E_{S^1(\delta)}(\lambda)\otimes
 E_{V_1}(\mu_1)\otimes E_{V_2}(\mu_2)$ ~~~ (with $k_1:=\text{dim} V_1$ and $k_2:=\text{dim} V_2$), the function $a$ is 
written in the form:
 \begin{eqnarray}
  a=\sum_{i=1}^{k_2}\varphi_i\beta_i\label{F5}
 \end{eqnarray}
where the $\beta_i$ form a basis of $E_{V_2}(\mu_2)$ and the functions $\varphi_i:\Theta\times S^1(\delta)\times
V_1\rightarrow\R$
satisfy: \\ for any $ x\in\Theta$, ~~~ $\varphi_{i_x}(.)\in E_{S^1(\delta)}(\lambda)\otimes E_{V_1}(\mu_1)$
~~~~ ($E_{V_1}(\mu_1)=\R$ if dim $V_1=0$).

Fundamental equation \ref{F1} is written:
\begin{eqnarray}
 {\Box_{g_\mathcal{P}}(\sum_{i=1}^{k_2}\varphi_i\beta_i)+S\sum_{i=1}^{k_2}\varphi_i\beta_i}=0\label{F6}
\end{eqnarray}

\textbf{When the metric $g_{\mathcal{P}}$ is neutral on $E_{V_2}(\mu_2)$}, equation \ref{F6} is equivalent,
according to
\ref{F4'}, to $k$ \textbf{identical} equations:
\begin{center}
 $\Box_{g_\mathcal{P}}(\varphi_i)+(\mu_2+S)\varphi_i=0$
\end{center}
to describe the function $a$, it suffices, when $g_\mathcal {P}$ is neutral on $E_{V_2}(\mu_2)$,
to study a single equation of the form:
\begin{eqnarray}
\Box_{g_\mathcal{P}}(\varphi)+(\mu_2+S)\varphi=0\label{F7}
\end{eqnarray}
where here, the function $\varphi$ is defined on $\Theta\times S^1(\delta)\times V_1$.
This gives important simplification when ~ $k_1:= \text{dim} V_1=1$, which is the case when $g_\mathcal{P}=g_0$
since $g_0$ is neutral on $W$ according to definition \ref{d2.4}, but this will also be the case for active potentials
 of order 1, as will be specified in the following paragraph.
\vspace{5mm}

The potentials that we will consider are essentially the following active potentials:\\
1- without electromagnetism \\
2- electromagnetic (subsection \ref{ss1.3}).

As already announced, the experimental results concerning the electromagnetic potentials will be
correctly described (spin effect included) when the cell decomposes in the following form:
\begin{eqnarray}
\mathscr {C} = \Theta \times S^1(\delta)\times S^3(\rho)\times V\label{F8}
\end{eqnarray}
where $S^3(\rho)$ is the standard 3-dimensional Riemannian sphere of radius $\rho$.

\textbf{When we neglect the effects of spin}, we will assume that the electromagnetic potential
$g_\mathcal{P}$ is neutral on $E_{{S^3(\rho)}\times{V}}(\mu)$ (in which case decomposition $W=S^3(\rho)\times{V}$ is
unnecessary). This for a domain of type "elementary oscillating metric associated with
$E_{\lambda,\mu}=E_{S^1(\delta)}\times E_{S^3\times{V}}(\mu)$"~ (def. \ref{d2.3}).

\textbf{When considering the effects of spin}, we assume that the potential $g_\mathcal{P}$ is neutral
only on $E_V(\mu_2)$.
(In the first case we can consider that we have "averaged" a priori the potential $g_\mathcal{P}$ of the second
case on $S^3(\rho)$).

Potentials "without electromagnetism" will often be considered neutral on $E_{S^3\times{V}}(\mu)$.
\bigskip

Active potential "without electromagnetism" and "electromagnetic" have been described precisely by
 propositions \ref{p1.2} and \ref{p1.3}. The following proposition, whose (very simple) proof is given in
Annex \ref{a3.5}, specifies conditions that make these potentials "neutral" on $E_{V_2}(\mu)$ when the
cell is of the form $\mathscr {C} = \Theta \times S^1(\delta)\times V_1\times V_2$.
\begin{prop} \label{p2.1}
~ \\ \vspace*{-1em}
\begin{enumerate}
  \item
  If, for an active potential without electromagnetism (proposition \ref{p1.2}), the vector field $X_1$ vanishes
on $E_{V_2}(\mu_2)$ ~ (ie.  if for any $ \beta\in E_{V_2}(\mu_2),~X_1(\beta)=0$), then this potential is neutral
on $E_{V_2}(\mu_2)$.
\item
If, for an electromagnetic potential (proposition \ref{p1.3}), the two vector fields $\Upsilon$ and $X_2$
vanishes on $E_{V_2}(\mu_2)$, then this potential is neutral on $E_{V_2}(\mu_2)$.
 \end{enumerate}
\end{prop}

The preceding considerations lead to refining the classification of elementary oscillating metrics. We
gives in the following paragraph the precise definition of elementary oscillating metrics of order 1 and 2.


\subsubsection{Elementary oscillating metrics of order 1 and 2 in a potential}
The considered cell is $\mathscr {C} = \Theta \times S^1(\delta)\times W$ ~ and ~ $W$ is decomposed into the form
$W=S^3(\rho)\times V$ for order 2.

A neutral potential metric is of the form: \\
$g_0= g_\Theta\times (-g_{S^1(\delta)})\times g_W$ ~ and ~ $g_W=g_{S^3(\rho)}\times {g_V}$ for order 2.
\newpage

\begin{dfn} \label{d2.9}
~ \\ \vspace*{-1em}
\begin{enumerate}
 \item
 \textbf{An elementary oscillating metric domain of order 1 in a potential} is an 
oscillating metric domain in a potential (def. \ref{d2.2}) for which the function $a:\C\rightarrow \R$ is of the
form:

$a=\varphi\beta$ ~~ where ~~ $\beta\in E_W(\mu)$ ~~ and
$\varphi:\Theta\times S^1(\delta)\rightarrow \R$ satisfies: \\
for any $ x\in\Theta$, ~~~ $\varphi_x(.)\in E_{S^1(\delta)}(\lambda)$. \\
Furthermore, the potential is neutral on $E_W(\mu)$.
\item
\textbf{An elementary oscillating metric domain of order 2 in a potential} is an 
oscillating metric domain in a potential (def. \ref{d2.2}) for which the function $a:\C\rightarrow \R$ is of the
form:

$a=\phi\beta$ ~~ where ~~ $\beta\in E_V(\mu)$ ~~ and $\phi:\Theta\times S^1(\delta)\times
S^3(\rho)\rightarrow \R$ satisfies: \\
for any $ x\in\Theta$, ~~~ $\phi_x(.)\in E_{S^1(\delta)}(\lambda)\otimes E_{S^3(\rho)}(\gamma)$. \\
Furthermore, the potential is neutral on $E_V(\nu)$.
\end{enumerate}
\end{dfn}

Either for order 1 or 2, this definition corresponds to definition \ref{d2.3} when we choose for the
function $a$ \textbf{only one} term in the decomposition \ref{F5}. (It is recalled that this is justified by the fact
that the constraint on the potentials gives identical equations of the form \ref{F7} for each of the terms of the
sum \ref{F5}).

\begin{rmq} \label{r2.1}
We could have defined oscillating metrics "of order 0" in the particular case where $\lambda=0$ but we have chosen to consider them as being of order 1. \\
Since $E_{S^1(\delta)}(0)$ naturally identifies with $\R$ and
$E_{S^1(\delta)}(\lambda)$ to $\CC$ when $\lambda>0$ (as will be specified in section \ref{s2.8}), the function
$\varphi$ of the above definition identifies to a \textbf{real} function defined on $\Theta$ when
$\lambda=0$ and a \textbf{complex} function if $\lambda>0$. Similarly, since $E_{S^1(\delta)}(\lambda)\otimes
E_{S^3(\rho)}(\gamma)$ identifies with the complexified $E_{S^3(\rho)}^{\varmathbb C}(\gamma)$, the function $\varphi$ of an
elementary oscillating metric of order 2 identifies with a function with values in
$E_{S^3(\rho)}^{\varmathbb C}(\gamma)$ defined on $\Theta$.

 As already mentioned, the process can be continued to define oscillating metrics
of order $k>2$ corresponding to finer decompositions of $W$ than that of the form $S^3(\rho)\times V$.
 
\end{rmq}
\begin{rmq} \label{r2.2}
 The reader familiar with Q.F.T can begin to make the following parallel with the theory presented here:
\begin{enumerate}
\item Elementary oscillating metrics of order 0 are similar to real scalar fields of the Q.F.T, those of order 1 to 
complex scalar fields, those of order 2 to spinorial fields, etc.
\item "Fields" quanta of  Q.F.T will be similar to singularities presented in section \ref{s2.11}.
\item The Fock space of a "particle system" is replaced by the space of regular functions:\\
 $\{a:\mathscr {C} = \Theta \times S^1(\delta) \times W\rightarrow \R\}$.
\item No Lagrangian, Hamiltonian and quantization process will be used here. Everything will be "managed" by
the sole fundamental equation \ref{F0'} which is a simple consequence of the geometric type chosen.
\item Precise description of experiments that will be deduced from the theory presented here will be profoundly different
to that given in classical quantum physics and in Q.F.T. It will be essentially linked to the "deformation" of
space-time $\M$ whose interpretation can only be specified after section \ref{s2.12}.
\end{enumerate}
\end{rmq}


\section [Important Examples] {Important examples of elementary oscillating metrics in neutral potential
 \label{s2.6}}
The particular examples of oscillating metrics that we will describe in this section are those which, in the
language of classical physics, represent the "flow of particles" moving at a constant velocity
$\overrightarrow{v}$ (therefore necessarily in a neutral potential and in the context of the linear approximation
). They will allow, for example, to write
precisely "boundary conditions" for the description of experiments whose principle is based on "sending 
particles" thrown at a velocity $\overrightarrow{v}$ in a physical system.These examples will also allow to
naturally introduce the concept of \textbf{relative} electric charge and this will then be extended to a class of much more general oscillating metrics (section \ref{s2.7}). We start by considering the case of elementary oscillating metrics
 of order 1 (def. \ref{d2.9}), extension to the more general case of elementary oscillating metrics
associated with $E_{\lambda,\mu}$ (def. \ref{d2.3}) does not bring any difficulties and will be presented at the end of this section.


\subsection{Homogeneous elementary oscillating metrics of order 1 in neutral potential}\label{4-ss-1}
According to definition \ref{d2.9}, the function $a$ which characterizes an elementary oscillating metric of order 1 in a
potential is of the form $a=\varphi\beta$ ~ where $\beta\in E_W(\mu)$ and $\varphi:\Theta\times S^1(\delta)\rightarrow
\R$ satisfies~: for any $ x\in\Theta ~~~\varphi_x(.)\in E_{S^1(\delta)}$.

It is assumed that this oscillating metric is \textbf{"homogeneous-stationary"}, which is expressed by saying that the
function $a$ does not depend on space variables $(x^1,x^2,x^3)$ of $\Theta$.

This property is obviously not invariant by the Poincaré transformations on $\Theta$, the word
"stationary" is relative to the choice of a chart in the observation atlas.

The function $\varphi$ is then of the form:
$$\varphi(t,u)=\varphi_1(t)\cos(Q^+u)+\varphi_2(t)\sin(Q^+u)$$
where $t=x^0$, $u=x^4$ and $Q^+:=\sqrt \lambda$ is the electric charge frequency (def. \ref{d2.6}).

The fundamental equation \ref{F1} is written here:
$$\Box_{g_0}a+Sa=0$$
$$\text{where} ~~ S=\frac{n-2}{4(n-1)}S_{g_0}~~ \text{and} ~~~\Box_{g_0}=\partial^2/\partial
t^2-\sum_{k=1}^{3}{\partial^2/{(\partial x^k)^2}}+\partial^2 /{\partial  u}^2+\Delta_{{g_0}_W}$$

This gives, using the fact that $\Delta_{{g_0}_W}(\beta)=\mu\beta$:
$$\frac {\partial^2{\varphi_1}}{\partial t^2}\cos(Q^+u)+\frac {\partial^2{\varphi_2}}{\partial
t^2}\sin(Q^+u)+(S+\mu-(Q^+)^2)\varphi=0$$
Which is equivalent to:

$$\frac {\partial^2{\varphi_1}}{\partial t^2}=-M^2\varphi_1 ~~~\text{and}~~~\frac {\partial^2{\varphi_2}}{\partial
t^2}=-M^2\varphi_2$$
since $M^2=S+\mu-(Q^+)^2$

In other words, oscillating metrics that satisfy the previous properties have their functions $a$ (seen in the
chart of the observation atlas) which are of the form: \\
$a=\beta((A_1\cos(Mt)+A_2\sin(Mt))\cos(Q^+u)+(B_1\cos(Mt)+B_2\sin(Mt))\sin(Q^+u))$ \\
where $A_1,A_2,B_1,B_2$ are constants.

This is written, after transformations:
$$a=\beta(C_1\cos(Mt+Q^+u)+C_2\sin(Mt+Q^+u))+\beta(C_3\cos(Mt-Q^+u)+C_4\sin(Mt-Q^+u))$$
where $C_1=1/2(A_1-B_2), C_2=1/2(A_2+B_1), C_3=1/2(A_1+B_2), C_4=1/2(A_2-B_1)$

Functions $a$ are written in the form ${a^+}+{a^-}$ with:

$$a^+=\beta(C_1\cos(Mt+Q^+u)+C_2\sin(Mt+Q^+u))$$
$$a^-=\beta(C_3\cos(Mt-Q^+u)+C_4\sin(Mt-Q^+u))$$
And one will notice that there is no ambiguity for the definition of $a^+$ and $a^-$ provided that $M$ is nonzero.

When one define $Q=Q^+$ for $a^+$ and $Q=-Q^+$ for $a^-$, functions $a^+$ and $a^-$
are written in the same form:
$$\beta( C\cos(Mt+Qu)+C'\sin(Mt+Qu))$$
In this particular example, constant $Q$ will be defined as \textbf{the relative electric charge} of the oscillating metric, positive for $a^+$, negative for $a^-$, and this is only well defined if $M>0$. (We can link $a^+$ and
$a^-$ to the notion of "particle" and "antiparticle")

We therefore give the following definition:
\begin{dfn} \label{d2.11} 
\textbf{An homogeneous-stationary elementary oscillating metric domain relative to $(\V,\zeta)$, of order 1}, is a
domain $\V$ of type "elementary oscillating metric of order 1 in a neutral potential" for which, in the
chart $(\V,\zeta)$ function $a$ is of the form:

$a=\varphi\beta$ where $\beta\in E_W(\mu)$ and $\varphi:\Theta\times S^1(\delta)\rightarrow \R$ satisfies:
$$\varphi(t,x^1,x^2,x^3,u)=C\cos(Mt+Qu)+C'\sin(Mt+Qu)$$
(and therefore does not depend on $x^1,x^2,x^3$). \\
(Here $Q$ is the relative electric charge and can be positive or negative).
\end{dfn}
To define an homogeneous elementary oscillating metric \textbf{that moves at a constant velocity
$\overrightarrow{v}$ relative to a chart ($\V$ ', $\zeta'$) of the observation atlas}, it suffices to
consider a Poincaré transformation: $\wedge:\Theta\rightarrow\Theta'$ which corresponds to an observer associated to
$\Theta'$ moving at a velocity $-\overrightarrow{v}$ relative to $\Theta$. We then define the
$g_0$-isometry \\
$\sigma:\C=\Theta\times S^1(\delta)\times W\rightarrow\mathscr{C'}=\Theta'\times S^1(\delta)\times W$ by setting
$\sigma= \wedge\times I_d$ ~~ where $I_d$ is the identity map on $S^1(\delta)\times W$.

The reader can check that the function $a\circ\sigma$, which corresponds to the function of the oscillating metric
 "seen" in the chart ($\V$ ', $\zeta'$), when $a$ is that "seen" in $(\V,\zeta)$, is of the
form $a\circ\sigma=\varphi'\beta$ ~~ where $\beta\in E_W(\mu)$ and \\
$\varphi'(t,(x^k),u)=C\cos(M't-\sum_{k=1}^3\lambda_kx^k+Qu)+C'\sin(M't-\sum_{k=1}^3\lambda_kx^k+Qu)$. \\
The velocity vector $\overrightarrow{v}$ is written $(1/M')(\lambda_1,\lambda_2,\lambda_3)\in \R^3$ ~~ and ~~
$\sum_{k=1}^3\lambda_k^2<M'^2$. \\
(Here the tangent space at any point of $\Omega$ is canonically identified with
$\R^3$ when $I\times\Omega= \Theta\subset\R^4$). \\
The constant $M'>0$ is called \textbf{the relativistic mass seen in the chart ($\V$ ', $\zeta'$)}. The fact that
 ~~ $\Box_{g_o}a+Sa=0$ ~~ shows that $M=(1-\lvert\overrightarrow{v}\rvert^2)^{1/2}M'$ ~~ where $M$ is the mass
(at rest) already defined.

We therefore give the following definition:
\begin{dfn} \label{d2.12}
 An homogeneous elementary oscillating metric (order 1) in a neutral potential \textbf{has a constant propagation velocity
 $\overrightarrow{v}$ relative to a chart $(\V,\zeta)$ of the observation atlas} if the function
$a=\varphi\beta$ is such that: \\
$\varphi(t,(x^k),u)=C\cos(M't-\sum_{k=1}^3\lambda_kx^k+Qu)+C'\sin(M't-\sum_{k=1}^3\lambda_kx^k+Qu)$. \\
The velocity vector is then $\overrightarrow{v}=(1/M')(\lambda_1,\lambda_2,\lambda_3)\in \R^3$.

\end{dfn}
When $\overrightarrow{v}=0$ we obviously recover the definition \ref{d2.11}.

\subsection{The more general case of elementary oscillating metrics associated with $E_{\lambda,\mu}$}  \label{ss+2}
If one chooses a basis $(\beta_1,\dots,\beta_k)$ in $E_W(\mu)$, the function $a$ representing such an oscillating metric  decomposes in the form \ref{F4'}:
$$a(x,u,w)=\sum_{i=1}^k(\varphi_{1,i}(x)\cos(Q^+u)+\varphi_{2,i}(x)\sin(Q^+u))~\beta_i(w)$$
Each term of this sum can be considered as the function $a_i$ of an elementary oscillating metric of
order 1.

\begin{dfn} \label{d2.13}
 An elementary oscillating metric type domain associated with $E_{\lambda,\mu}$ is \textbf{homogeneous
stationary for a chart $(\V,\zeta)$} if for any $i$ from 1 to k, the functions
$a_i$ correspond to those given
in definition \ref{d2.11}, all with the same value of $Q$, (when $Q=\sqrt{\lambda}=Q^+$ the relative electric charge
 is positive, when $Q=-\sqrt{\lambda}=-Q^+$ it is negative).
\end{dfn} We also generalize definition \ref{d2.12}.

These definitions obviously do not depend on the choice of the basis $(\beta_1,\dots,\beta_k)$.


\section{Relative electric charge \label{s2.7}}
Here we extend the notion of relative electric charge, which was introduced in the previous subsection, to the more general framework of elementary oscillating metrics associated with $E_{\lambda,\mu}$ (def. \ref{d2.3}) for which the notion of absolute electric charge has been specified in definitions \ref{d2.6} and \ref{d2.7}. Remember here that the cell is $\C=\Theta\times S^1(\delta)\times W$ and that the function $a:\C\rightarrow\R$ satisfy: for any $ x\in\Theta ~~~a_x(.)\in E_{\lambda,\mu}$.

\begin{dfn} \label{d2.14}
 An elementary oscillating metric domain associated with $E_{\lambda,\mu}$ in a potential \textbf{has a 
well-defined (relative) electric charge} if $Q^+=0$, \textbf{or} if the vector field defined on $\Theta$ by:
\begin{eqnarray}
\int_{S^1\times W} \frac{\partial a}{\partial u}~~(\overrightarrow{grad}_{g_0{_\Theta}}a)~~dv_{g_0{_{S^1\times
W}}}\label{F9}
\end{eqnarray}
is a \textbf{timelike} vector field which is furthermore:
\begin{enumerate}
\item Either everywhere in the time orientation of $\Theta$.
\item Either everywhere opposed to the time orientation of $\Theta$.
\end{enumerate}
In the first case, \textbf{the relative electric charge frequency $Q$ is defined by $Q^+$}, in the second case \textbf{by
$-Q^+$}.

(We define as previously \textbf{the electric charge $q$} by setting $q=\hbar Q$)
\end{dfn}
In the integral \ref{F9}, the function $\frac{\partial a}{\partial u}:\C\rightarrow\R$ is none other than the
function $Y(a)$ where $Y$ is the vector field that defines electromagnetism, and the field
$\overrightarrow{grad}_{g_0{_\Theta}}a$ is the vector field tangent to $\Theta$ defined by: \\
$\overrightarrow{grad}_{g_0{_\Theta}}a:=(\partial_0a)\partial_0-\sum_{k=1}^3(\partial_ka)\partial_k$ ~~
where ~~ $\partial_0:=\frac{\partial}{\partial x^0}=\frac{\partial}{\partial t}$ ~~ and ~~
$\partial_k:=\frac{\partial}{\partial x^k}$.
\vspace{5mm}

By its "intrinsic" definition, the relative electric charge is invariant by Poincaré transformations
on
$\Theta$ that keep the time orientation.

The reader can verify that the relative electric charge given by definition \ref{d2.14}
corresponds to that given in the particular case of homogeneous oscillating metrics (defs. \ref{d2.12} and
\ref{d2.13}) and that, when $M=0$, the vector field defined by \ref{F9} is spacelike or lightlike,
relative electric charge is therefore not well defined in this case.


\section [Associated canonical functions] {Canonical functions associated with elementary oscillating metrics
of order 1 and 2 \label{s2.8}}
It is recalled that the considered cell for these oscillating metrics is $\mathscr {C} = \Theta \times
S^1(\delta)\times W$ and that $W$ is decomposed in the form $W= S^3(\rho)\times V$ for order 2. \\
A neutral potential metric is of the form $g_0= g_\Theta\times (-g_{S^1(\delta)})\times g_W$ where
$g_W=g_{S^3(\rho)}\times {g_V}$ for order 2.
\vspace{5mm}

According to definition \ref{d2.9}, the function $a$ that characterizes the oscillating metric satisfies for order 1:

$$a=\varphi\beta$$
where ~~ $\beta\in E_{S^3\times V}(\mu)$ ~~ and ~~ $\varphi:\Theta\times S^1\rightarrow\R$ satisfies:
for any $ x\in\Theta, ~~~~~\varphi_x(.)\in E_{s^1(\delta)}(\lambda)$.

We have: ~~ for any $ (x,u)\in \Theta\times S^1(\delta)$,
\begin{eqnarray}\label{F10}
\varphi_x(u)=\varphi_1(x)\cos(Q^+u)+\varphi_2(x)\sin(Q^+u).
\end{eqnarray}
And for order 2:

$$a=\phi\beta$$
where ~~ $\beta\in E_V(\nu)$ ~ and ~ $\phi:\Theta\times S^1(\delta)\times
S^3(\rho)\rightarrow\R$ satisfies:

for any $ x\in\Theta, ~~\phi_x(.)\in E_{S^1(\delta)}(\lambda)\otimes E_{S^3(\rho)}(\nu)$.

We have: ~~ for any $ (x,u,s)\in \Theta\times S^1(\delta)\times S^3(\rho)$,
\begin{eqnarray}\label{F11}
\phi_x(u,s)=\phi_{1,x}(s)\cos(Q^+u)+\phi_{2,x}(s)\sin(Q^+u).
\end{eqnarray}
 \textbf{We will define canonical functions from the preceding functions $\varphi$ and $\phi$. These
functions will contain all the essential information of the function $a$ and will be easier to use. They will also be interesting because beeing defined on "apparent space" $\Theta\subset\R^4$ and because they
can be compared to functions involved in standard quantum theories}.

\begin{dfn} \label{d2.15}
~ \\ \vspace*{-1em}
 \begin{enumerate}
  \item We consider an " elementary oscillating metric domain \textbf{of order 1} in a potential ". \\
\textbf{The canonical function} associated with this domain is the \textbf {complex} function $a_c:\Theta\rightarrow\CC$ ~~
given, when $Q^+\neq 0$, ~ by ~ $a_c(x):=\varphi_1(x)+i\varphi_2(x)$ \\ where ~ $\varphi_1$ and $\varphi_2$ are defined
by \ref{F10}, \\ and given simply by $a_c=\varphi_1$ if $Q^+=0$ ~~ (it is, in the latter case, a real function).
\item We consider an "elementary oscillating metric domain \textbf{of order 2} in a potential ". \\
\textbf{The canonical function} associated with this domain is the function 
$a_c:\Theta\rightarrow E_{S^3(\rho)}^{\CC}(\gamma)$ \\ given,  when ~ $Q^+\neq0$, ~ by
$a_c(x):=\phi_{1,x}+i\phi_{2,x}\in E_{S^3(\rho)}(\gamma)+iE_{S^3(\rho)}(\gamma)$ \\ (where ~~ $E_{S^3(\rho)}^{\CC}(\gamma)$ is
the "complexified"
$\CC\otimes E_{S^3(\rho)}(\gamma)=E_{S^3(\rho)}(\gamma)+iE_{S^3(\rho)}(\gamma)$ ~ and ~
$\phi_{1,x}$ and $\phi_{2,x}$ are defined by \ref{F11}), \\ and given simply by $a_c(x):=\phi_{1,x}$ if $Q^+=0$.
\end{enumerate}
\end{dfn}
We can, of course, extend this definition to elementary oscillating metrics of order $k>2$ when the compact manifold
 $V$ is broken down into product.

In fact, canonical functions $a_c$ are directly associated to the existence of the following canonical isomorphisms
that we will use during the proofs:

We consider the eigenspaces $E_{S^1(\delta)}(\lambda)$. For eigenvalues $\lambda\neq0$ all these spaces
are of dimension 2 and possess a canonical basis made up of the two functions "$\cos(\sqrt \lambda u)$" and
"$\sin(\sqrt \lambda u)$" expressed in the standard coordinate system of $S^1(\delta)$ oriented
(see \ref{ss1.1}).
\\ When $\lambda\neq0$, the canonical isomorphism $\CC_\lambda:E_{S^1(\delta)}(\lambda)\rightarrow\CC$ is defined
by:
\begin{eqnarray}\label{F12}
 \CC_{\lambda}(A\cos(\sqrt \lambda (.))+B\sin(\sqrt \lambda(.)):=A+iB
\end{eqnarray}
In the case of order 1 oscillating metrics, the canonical function $a_c$ thus satisfies: \\ for any $ x\in\Theta,
~~a_c(x)=\CC_\lambda(\varphi_x(.))$.

The isomorphism $\CC_\lambda$ then naturally induces an isomorphism:
\begin{eqnarray}\label{F13}
\CC_{\lambda,\gamma}:E_{S^1(\delta)}(\lambda)\otimes E_{S^3(\delta)}(\gamma)\rightarrow E_{S^3(\rho)}^{\CC}(\gamma)
\end{eqnarray}
In the case of order 2 oscillating  metrics, the canonical function $a_c$ satisfies: \\ for any $ x\in\Theta,
~~a_c(x)=\CC_{\lambda,\gamma}(\phi_x(.))$. \\
We then quickly prove that for order 1 or 2:
\begin{eqnarray}\label{F14}
 a=Re(\beta e^{-iQ^+u}a_c)
\end{eqnarray}


\section [Klein-Gordon equations] {Klein-Gordon equations (for different potentials) obtained as
simple consequences of
the definition of "geometric type" \label{s2.9}}
The important result of this section is stated in the following theorem \ref{2.1}. It concerns the 
elementary oscillating metric domain of order 1 in a potential. In fact, those of order 2, which take into account the notion of "spin", are treated identically, but for the sake of clarity we will present this case only in the
section \ref{s2.13}. \\
\textbf{Equations obtained are only translations of the fundamental equation $\Box_{g_{\mathcal P}}a+Sa=0$}.

 We begin by recalling and clarifying what are the "potential" type domains that we will use and that we
have already studied (see \ref{ss1.2} and \ref{s1.5}):

We consider a chart $(U,\zeta)$ of the observation atlas such that $\zeta(U)$ is the cell
$\C=\Theta\times S^1(\delta)\times W$ whose coordinates will be written $(x,u,w)$ with $x=(t,x^1,x^2,x^3)\in
\Theta\subset\R^4$, the coordinates of "time" $x^0$ and $x^4$ are written here $t$ and $u$.
\begin{enumerate}
 \item \textbf{Neutral potential}
 \vspace{4mm}
 
 The metric is $g_0= g_\Theta\times (-g_{S^1(\delta)})\times {g_W}$ ~~ where ~~ $g_\Theta$ is the Minkovski's metric on
$\Theta \subset \R^4$, ~~ $g_{S^1(\delta)}$ is the standard metric of $S^1(\delta)$, ~~ $g_W$ is a Riemannian metric on the compact manifold $W$ with constant scalar curvature.
\item \textbf{Active potential without electromagnetism}
\vspace{4mm}

The metric is, as we have seen (cf \ref{ss1.2}), of the form $g=g_0+h$ where, for any point $P$, the endomorphism
$\leftidx{^e}{h_P}$ is nilpotent of index 2. We have (proposition \ref{p1.2}):

$h=-2vX_1^\flat\otimes X_1^\flat$, ~~ $g_0(X_1,\frac{\partial}{\partial t})=1$, ~~ $g_0(X_1,Y)=0$, ~~ $g_0(X_1,X_1)=0$.

Here $v$ is the potential function and $X_0$ was chosen equal to $\frac{\partial}{\partial t}$ which is linked to the chart of
the observation atlas.

It is assumed that the hypothesis $H_N$ set for the study of geodesics (see \ref{ss1.4}) holds (with
$D_{g_0}X_1=0$),
that $S_g=S_{g_0}$ (see remark \ref{r13}) and that the metric $g$ is neutral on $E_W(\mu)$ (def. \ref{d2.4})
when the oscillating metric is associated with $E_{\lambda,\mu}:=E_{S^1}(\lambda)\otimes E_W(\mu)$.

These last points are recapitulated in \textbf{the following hypothesis $H_{1,N}$}:
\begin{enumerate}
 \item $S_{\mathcal P}=S_{g_0}$.
 \item The potential function $v$ is defined on $\mathcal U$ when $\Theta=I\times\mathcal U\subset \R\times\R^3$.
 \item The vector field $X_1$ is defined on $I\times W$, ~~ $D_{g_0}X_1=0$ and $X_1$ vanishes on $E_W(\mu)$.
\end{enumerate}
\item \textbf{Electromagnetic potential} \\
The metric is, as we have seen (cf \ref{ss1.2}), of the form $g=g_0+h$ where, for any point $P$, the endomorphism
$\leftidx{^e}{h_P}$ is nilpotent of index 2 or 3. We have (proposition \ref{p1.3}):

$h=\Upsilon^\flat\otimes X_2^\flat+X_2^\flat\otimes\Upsilon^\flat$ ~~ $g_0(X_2,Y)=1$ ~~ $g_0(X_2,\Upsilon)=0$ ~~ $g_0(X_2,
X_2)=0$.

It is assumed that hypothesis $H_E$ set for the study of geodesics (cf \ref{ss1.4}) holds (with
$D_{g_0}X_2=0$), that $S_g=S_{g_0}$ (see remark \ref{r14}) and that the metric $g$ is neutral on $E_W(\mu)$ (def.
\ref{d2.4}).

These last points are recapitulated in \textbf{the following hypothesis $H_{1,E}$}:
\begin{enumerate}
 \item $S_{\mathcal P}=S_{g_0}$.
 \item The vector field $\Upsilon$ is defined on $\Theta$.
 \item The vector field $X_2$ is defined on $S^1(\delta)\times W$, ~~ $D_{g_0}(X_2)=0$ and $X_2$ vanished on
$E_W(\mu)$.
\end{enumerate}
\end{enumerate}

Hypothesis $H_{1,N}$ and $H_{1,E}$ can be seen as "approximations" that reflect the fact that
some quantum effects on the compact manifold $W$ are neglected. Not assuming these assumptions would bring
additional perturbative terms in the equations given by theorem \ref{2.1}. It can be interesting
to study these disturbances, however the problem is more delicate than it seems at first sight because, if one goes into the details of the proof of theorem \ref{2.1}, one realizes that having chosen
the linear approximation (section 2.3) imposes constraints on the potentials (realized in particular by the
hypotheses
$H_{1,N}$ and $H_{1,E}$). To completely remove the assumptions $H_{1,N}$ and $H_{1,E}$ would therefore require to use
 fundamental equation \ref{F0'} and not its linear approximation \ref{F1}, which makes the study much more
difficult.

\begin{thme} \label{2.1}
We consider an "elementary oscillating metric domain of order 1 in a potential" (def. \ref{d2.9}).
  
Then, in the 3 cases of considered potential, the
canonical function $a_c$ satisfies the following "Klein-Gordon" equations:
\begin{enumerate}
 \item \textbf{In neutral potential}.
 \begin{eqnarray}\label{F15}
  \Box_\Theta a_c+ M^2a_c=0      
 \end{eqnarray}
where ~~ $\Box_\Theta=\frac{\partial^2}{(\partial t)^2}-\sum_{k=1}^3\frac{\partial^2}{(\partial x^k)^2}$ ~~ and ~~ $M$ is the
mass frequency of the oscillating metric.
\item \textbf{In a potential without electromagnetism under the assumption $H_{1,N}$}.
\begin{eqnarray}\label{F16}
\Box_\Theta a_c+ M^2a_c-2v\frac{\partial^2a_c}{(\partial t)^2}=0
\end{eqnarray}
where ~~ $v$ is the potential function (def. \ref{def:11})
\item \textbf{In an electromagnetic potential under the assumption $H_{1,E}$}.
\begin{eqnarray}\label{F17}
 \sum_{j=0}^3\varepsilon_j(i\frac{\partial}{\partial x^j}+Q^+\Upsilon^j)^2a_c+M^2a_c=0
\end{eqnarray}
where ~~ $\varepsilon_j=g_{0jj}$ ~~ ie: ~~ $\varepsilon_0=-1$ ~~ and ~~ $\varepsilon_1=\varepsilon_2=\varepsilon_3=+1$.

(The square of a differential operator is understood, of course, as the composition by itself).
\end{enumerate}
\end{thme}
Proof of this theorem is detailed in annex (\ref{a3.6}).

\begin{rmq} \label{r2.4}
 Equations \ref{F15} and \ref{F17} are invariant by chart changes that come from Poincaré transformations. This is not the case for equation \ref{F16} because the "field
"$X_0$ was chosen equal to $\frac{\partial}{(\partial t)}$ linked to a chart of the observation atlas, and
this one is not invariant by Poincaré transformations.
\end{rmq}
\textbf{It is important to note that equations given by this theorem are necessary conclusions of the sole
definition of a domain "with constant scalar curvature conform to a potential" as part of the linear approximation. No law, no principle has been added to obtain these equations}.
\bigskip

\textbf{Remark on the Dirac equation}

 Klein-Gordon equations given by theorem \ref{2.1} are of the form $D(a_c)=0$ where $D$ is an \textbf{order 2} differential operator. In neutral potential, $D=\Box_\Theta+M^2$. Dirac equation is obtained by
"factoring" the differential operator $D$ into two \textbf{order 1} differential operators (see subsection \ref{ssAv}). Having
1st order differential operators is fundamental in standard quantum physics, essentially to make coherent
the probabilistic axiomatic system of the latter. This is not the case for theory that we present here
as will be specified in section \ref{s2.12}. We are therefore not interested, from a conceptual point of view, in
introducing
equations of order 1, and it is equations of theorem \ref{2.1} which will remain fundamental. \\ However, for a
 simply "computational" point of view, it can be useful, in the study of the solutions of the Klein-Gordon equations,
to decompose the differential operator of order 2 into "composition" of operators of order 1. This makes it possible to introduce the
notion of "exponential operator" and provides a powerful tool (often used in Q.F.T) that can be interesting
even for us.

That said, we do not have enough elements at the end of this section to give the physical interpretation of
results of theorem \ref{2.1} to explain and describe completely experiments
usually studied by standard quantum physics, it will be necessary to wait for the end of section \ref{s2.12}.

Equations given by theorem \ref{2.1} are, of course, valid when $M=0$, although in this case the 
relative electric charge is generally not defined (but it does not intervene in these equations). 
Oscillating metric of "zero mass" have very different physical characteristics from 
 positive ones. Section \ref{s2.14} will clarify things and also consider oscillating metrics of
zero mass of order 2.

When $M$ is positive, the next section will show that theorem \ref{2.1} gives an approximation to
 Schrödinger equations of standard quantum physics on the state function that allow the study of the
behavior of "particles in a potential". The link between "the state function" that we will define here and
that of the usual quantum physics will appear only later.


\section [The state function] {The state function of an "oscillating metric in a\\
potential domain". \\
 Schrödinger equations found as approximations of\\
 Klein-Gordon equations \label{s2.10}}
The state function will be defined from the canonical function $a_c$ (def. \ref{d2.15}) but it will only make sense in the context of elementary oscillating metrics for which the electric charge is well defined (def.
\ref{d2.14}), in particular the mass will be assumed to be positive. Equations which we will obtain will be valid in a more restricted framework than that supposed for theorem \ref{2.1}.

The terminology "state function" was chosen because equations satisfied by this one will be very close to the  
Schrödinger equations (in different potentials) concerning the "state function" of 
classical quantum physics. However, the interpretation of the "probabilistic" results associated to the state function (as well
that to the canonical function $a_c$) will be fundamentally different and will be specified in section \ref{s2.12}.
\begin{dfn} \label{d2.16}
 We consider an "elementary oscillating metric domain in a potential" whose mass $m$ is 
 positive and charge $q$ well defined (def. \ref{d2.14}).

\textbf{The state function} associated with the function $a$ that characterizes the oscillating metric (seen in a chart of
the observation atlas) is the function:

$\varPsi:\Theta\rightarrow\CC$ ~~ for order 1,

$\varPsi:\Theta\rightarrow E_{S^3(\rho)}^{\CC}(\gamma)$ ~~ for order 2,

which satisfies:

$\varPsi=e^{iMt}a_c$ ~~ if the electric charge is positive,

$\varPsi=e^{iMt}\overline{a_c}$ ~~ if the electric charge is negative.

where the conjugate $\overline{a_c}$ is defined in $\CC$ for order 1 and in $E_{S^3(\rho)}^{\CC}(\gamma)$ for
order 2.

(If the electric charge is zero, these two equalities are identical:
 $\varPsi=e^{iMt}\varphi$, since then $a_c=\overline{a_c}=\varphi$).
\end{dfn} 
We notice that, according to \ref{F14}:
\begin{eqnarray}\label{F1'}
 a=\beta(Re(\varPsi e^{-i(Mt+Qu)}))=\beta(\varPsi_1\cos(Mt+Qu)+\varPsi_2\sin(Mt+Qu))
\end{eqnarray}
where ~~ $\varPsi=\varPsi_1+i\varPsi_2$.
\begin{rmq} \label{r2.5}
 The idea that leads to the definition of the state function that we have just given is to factorize in the
function $a$ that characterizes the oscillating metric, the maximum of oscillations in $t$. Indeed, as we have seen,
the pseudo-Riemannian metric $g$ of an elementary oscillating metric satisfies a form of wave equation
and the frequency of this wave is precisely the mass frequency $M$. In writing $a_c=e^{-iMt}\varPsi$
we make the variations of $\varPsi$ in $t$ as small as possible by factoring $a_c$ by
$e^{iMt}$ which only involves the fundamental constant $M$. This only really makes sense when
the considered oscillating metric is close to an homogeneous oscillating metric (def. \ref{d2.13}) (when
it is stationary it is easily proven that $\varPsi$ no longer depends on $t$), which is the case in the
description of generic experiments of standard quantum physics.
\end{rmq}
\begin{rmq} \label{r2.6}
 The modulus of $\varPsi$ is equal to that of $a_c$. In standard quantum physics this means that the "density of 
probability of presence of a particle" is also given by $|a_c|$. This equality will be important for us during the "probabilistic interpretation" given in section \ref{s2.12}.
\end{rmq}
The following theorem is only a simple rewriting of theorem \ref{2.1} when the canonical function $a_c$ has been replaced
 by the function $\varPsi$ given in definition \ref{d2.16} (here we restrict ourselves to oscillating metrics of order 1,  those of order 2 will be studied in section \ref{s2.13}
. Verification of the equations obtained is left to the reader.

\begin{thme} \label{2.2}
 Under the same assumptions as those of theorem \ref{2.1} and when the electric charge is well defined, the
state function $\varPsi$ satisfies the following equations:
\begin{enumerate}
 \item \textbf{In neutral potential}.
 \begin{eqnarray}\label{F2}
2iM\frac{\partial\varPsi}{\partial t}=\Delta\varPsi+\frac{\partial^2\varPsi}{\partial t^2} 
 \end{eqnarray}
where $\Delta$ is the Laplacian operator: $\Delta=-(\sum_1^3\frac{\partial^2}{(\partial x^k)^2})$ and $M$ the
mass frequency.
\item \textbf{In an active potential without electromagnetism}.
 \begin{eqnarray}\label{F3}
 2iM\frac{\partial\varPsi}{\partial t}=\Delta\varPsi+2vM^2\varPsi+(1-2v)\frac{\partial^2\varPsi}{\partial t
^2}+4iMv\frac{\partial\varPsi}{\partial t}
 \end{eqnarray}
where ~~ $v$ ~~ is the potential function.
\item \textbf{In an electromagnetic potential}.
 \begin{eqnarray}\label{F4}
 2iM\frac{\partial\varPsi}{\partial t}=\sum_{j=O}^3\varepsilon_j(i{\frac{\partial}{\partial
x^j}}+Q\Upsilon^j)^2\varPsi-2MQ\Upsilon^0\varPsi
\end{eqnarray}
where ~~ $\varepsilon_j=g_{0jj}$ ~~ ie: ~~ $\varepsilon_0=-1$ ~~ and ~~ $\varepsilon_1=\varepsilon_2=\varepsilon_3=+1$, \\
$\Upsilon$ is the electromagnetic potential, $Q$ the relative electric charge frequency.
\end{enumerate}
\end{thme}
\begin{rmq} \label{r2.7}
 As already mentioned in remark \ref{r2.4} concerning the Klein-Gordon equations, equations
\ref{F2} and \ref{F4} are invariant by chart changes that come from Poincaré transformations and
this is not the case for equation \ref{F3}. Of course, none of these equations will be invariant when one
will have removed the "negligible" terms to recover the classical Schrödinger equations.
\end{rmq}

It is interesting to rewrite the equations obtained in theorem \ref{2.2} by returning to the more usual units 
in which the standard Schrödinger equations are written. The "time" that has been written $t$ in
the system of geometric units will be written $\tau$ in the S.I system, and the "time" $u$ will be written $u'$.
The mass $m$ (in kg) and the electric charge (in C) are defined, as we have seen, from the mass frequency $M$ and
 electric charge frequency $Q$ by: $m:=\hbar c^{-1}M$ and $q:=\hbar Q$.

If we introduce the $\tau$-mass frequency (in $s^{-1}$) and the $\tau$-electric charge frequency $Q'$ (in
$s^{-1}$) by setting $M':=cM$ and $Q':=cQ$, the previous equalities are written:
 \begin{eqnarray}\label{F5''}
  mc^2=\hbar M' ~~~{\text et} ~~~ qc=\hbar Q'
 \end{eqnarray}
Similarly, it is natural to introduce the Newtonian potential $v'$ (in $m^2s^{-2}$) by setting $v':=c^2v$ (where $v<1$ is
without unit) and the electric potential $\phi:=c\Upsilon_0=-c\Upsilon^0$ since it is "$\Upsilon_0~ dt$" which
intervenes in the
definition of the 1-form $(\leftidx{^e}{h}(Y))^\flat$.

Equations \ref{F2}, \ref{F3}, \ref{F4} obtained by theorem \ref{2.2} are then written, when one denotes
again $\varPsi$ the function depending of the variables $(\tau,x^1,x^2,x^3)$:
\begin{eqnarray}\label{F6'}
 i\frac{\partial\varPsi}{\partial
\tau}=\frac{\hbar}{2m}\Delta\varPsi+\frac{\hbar}{2{mc^2}}\frac{\partial^2\varPsi}{\partial \tau^2}
\end{eqnarray}
\begin{eqnarray}\label{F7'}
 i\hbar\frac{\partial\varPsi}{\partial
\tau}=\frac{{\hbar}^2}{2m}\Delta\varPsi+v'm\varPsi+\frac{1}{c^2}(\frac{\hbar^2}{2m}(1-2v'/c)\frac{\partial^2\varPsi}
{\partial
\tau^2}+2i\hbar v'\frac{\partial\varPsi}{\partial \tau})
\end{eqnarray}
\begin{eqnarray}\label{F8'}
 i\hbar\frac{\partial\varPsi}{\partial \tau}=\frac{{\hbar}^2}{2m}\sum_{k=1}^3(i{\frac{\partial}{\partial
x^k}}+\frac{q}{\hbar}\Upsilon_k)^2\varPsi+q\phi\varPsi-\frac{1}{2mc^2}(q\phi-i\hbar{\frac{\partial}{\partial
\tau}})^2\varPsi
\end{eqnarray}

The reader will notice that, when the last term of the right hand side has been removed, the three equations obtained are
identical to classical Schrödinger equations for the state function of a particle in a potential. The
essential difference is the fact that the differential equations of theorem \ref{2.2} are of order 2 in $t$.

The last term of equation \ref{F6'}, \ref{F7'}, \ref{F8'} written in the S.I system include the coefficient
$\frac{1}{c^2}$, which suggests that it can be easily "neglected" to recover the standard Schrödinger equation. This point must obviously be clarified and this is the purpose of the following subsection.


\subsection{The $\varepsilon$-approximations}
We specify here the mathematical conditions which make it possible to justify the approximations which bring back the
equations given by theorem \ref{2.2} to the corresponding classical Schrödinger equations. Conditions
are classified into two categories: the $\varepsilon$-approximations for the state function and the
$\varepsilon$-approximations for the potentials.
\bigskip

For the sake of clarity, we start by writing what we call the "classical Schrödinger's equations"
 and then we specify their differences with the equations \ref{F2}, \ref{F3}, \ref{F4} obtained in 
theorem \ref{2.2}.
\newpage

\textbf{The classical Schrödinger equations}

\begin{enumerate}
 \item In neutral potential.
\begin{eqnarray}\label{F9'}
 2iM\frac{\partial\varPsi}{\partial t}=\Delta\varPsi
\end{eqnarray}
\item In an active potential without electromagnetism.
\begin{eqnarray}\label{F10'}
2iM\frac{\partial\varPsi}{\partial t}=\Delta\varPsi+2vM^2\varPsi 
\end{eqnarray}
\item In an electromagnetic potential.
\begin{eqnarray}\label{F11'}
  2iM\frac{\partial\varPsi}{\partial t}=\sum_{k=1}^3(i{\frac{\partial}{\partial
x^k}}+Q\Upsilon^k)^2\varPsi-2MQ\Upsilon^0\varPsi
\end{eqnarray}
\end{enumerate}

To recover equations \ref{F2}, \ref{F3}, \ref{F4}, terms to be added are successively the following ones:
\begin{eqnarray}\label{F12'}
 1: \frac{\partial^2\varPsi}{\partial t^2}, ~~~~~2: (1-2v)\frac{\partial^2\varPsi}{\partial t
^2}+4iMv\frac{\partial\varPsi}{\partial t}, ~~~~~3: -(i{\frac{\partial}{\partial
t}}+Q\Upsilon^0)^2\varPsi
\end{eqnarray}


\subsubsection{The $\varepsilon$-approximations for the state function}
The intuitive idea is that these approximations will be valid when the velocities involved will be very small 
relative to the speed of light. In the context of the elementary oscillating metrics, the notion of velocity has been
specified only for homogeneous oscillating metrics in a neutral potential (def. \ref{d2.12}). So we start, on
this example, to compute the state function $\varPsi$ to determine the links that exist between the derivatives of $\varPsi$
, $\varPsi$ itself and the propagation velocity associated to a chart $(U,\zeta)$ in the observation atlas.

The function $a$ associated with an homogeneous oscillating metric in a neutral potential is, by definition, of the form
$a=\varphi\beta$ where: \\
$\varphi(t,(x^k),u)=C\cos(M't-\sum_{k=1}^3\lambda_kx^k+Qu)+C'\sin(M't-\sum_{k=1}^3\lambda_kx^k+Qu)$

and ~ $\beta\in E_W(\mu)$ ~ with $\sum_{k=1}^3 \lambda_k^2 <M'^2$.

The velocity vector $\overrightarrow{v}$ has components $\frac{1}{M'}(\lambda_1,\lambda_2,\lambda_3)$. $M'$ is the relativistic mass frequency and is associated to $M$ by: $M=\sqrt{1-|\overrightarrow{v}|^2}M'$.

A simple calculation shows that: \\
$\varphi(t,(x^k),u)=(C\cos\alpha+C'\sin\alpha)\cos(Qu)+(C'\cos\alpha-C\sin\alpha)\sin(Qu)$

where $\alpha=M't-\sum_{k=1}^{3}\lambda_kx^k$.

We deduce that the associated canonical function (def. \ref{d2.15}) is, if $Q>0$:

$a_c=(C\cos\alpha+C'\sin\alpha)+i(C'\cos\alpha-C\sin\alpha)=(C+iC')e^{-i(M't-\sum_{k=1}^{3}\lambda_kx^k)}$

and is equal to the conjugate if $Q<0$.

The state function $\varPsi$ (def. \ref{d2.16}) thus satisfies (wathever the sign of $Q$): \\
$\varPsi=(C+iC')e^{-i((M'-M)t-\sum_{k=1}^{3}\lambda_kx^k)}$

Then: \\ $|\varPsi|^2=C^2+C'^2, ~~|\frac{\partial\varPsi}{\partial t}|=(M'-M)|\varPsi|,
~~|\frac{\partial^2\varPsi}{\partial t^2}|=(M'-M)|\frac{\partial\varPsi}{\partial t}| \\=(M'-M)^2|\varPsi|,
~~|\triangledown_{x_1,x_2,x_3}\varPsi|^2=(\sum_{k=1}^{3}\lambda_k^2)|\varPsi|=|\Delta_{x_1,x_2,x_3}\varPsi|.$

As $|\overrightarrow{v}|^2=\frac{1}{M'^2}\sum_{k=1}^{3}\lambda_k^2$ ~~ and
~~ $M'-M=((1-|\overrightarrow{v}|^2)^{-1/2}-1)M=(1/2|\overrightarrow{v}|^2+\circ(|\overrightarrow{v}|^2))M$ ~~~
we can deduce:
\begin{eqnarray}\label{F13'}
 |\frac{\partial\varPsi}{\partial t}|=(1/2|\overrightarrow{v}|^2+\circ|\overrightarrow{v}|^2)M|\varPsi|,
~~~~|\frac{\partial^2\varPsi}{\partial
t^2}|=(1/4|\overrightarrow{v}|^4+\circ|\overrightarrow{v}|^4)M^2|\varPsi|
\end{eqnarray}
The idea is then to consider that elementary oscillating metrics (for which the notion of velocity does not have
been defined) which will be used in the experiments for which the approximations that we are going to specify are
valid (diffraction, Young's slits, Stern-Gerlach experiment, quantum entanglement, etc.) are only
"low disturbances" of
homogeneous oscillating metrics. Relations \ref{F13'} as well as terms given by \ref{F12'} thus inspire the
next definition.
\begin{dfn} \label{d2.17}
 We consider an "elementary oscillating metric domain of order 1 in a potential" for which the
associated state function, associated to a chart of the observation atlas, is written $\varPsi$.

We will say that \textbf{$\varPsi$ satisfies the $\varepsilon$-approximations} if:

There is a real number ~ $0<\varepsilon<1$ ~ and two functions $\mu$ and $\mu':\Theta\rightarrow\CC$ such that:
\begin{eqnarray}\label{F14'}
\frac{\partial\varPsi}{\partial t}=\varepsilon^2M\mu\varPsi,
~~~~\frac{\partial^2\varPsi}{\partial t^2}=\varepsilon^4M^2\mu'\varPsi
\end{eqnarray}
\vspace{-1cm}
  
where, on $\Theta$:
\begin{enumerate}
 \item $|\mu|<1$ ~~ and ~~ $|\mu'|<1$
 \item There are two constants $C_1$ and $C_2$ such that: ~~~ for any $k$ from $0$ to $3$,
 
 $|\frac{\partial\mu}{\partial x^k}|<C_1$ ~~ and ~~ $|\frac{\partial\mu'}{\partial x^k}|<C_2$
\end{enumerate}
In geometric units, $C_1$ and $C_2$ are inverse lengths, their choices are specified in the lines that
follow.
\end{dfn}

If we use the fact that $\frac{\partial^2\varPsi}{\partial t^2}=\varepsilon^4M^2\mu'\varPsi$, equation \ref{F2}
is written:

$$2iM\frac{\partial\varPsi}{\partial t}=\Delta\varPsi+\varepsilon^4M^2\mu'_1\varPsi+i\varepsilon^4M^2\mu'_2\varPsi$$

where we put $\mu'=\mu'_1+i\mu'_2$.

The term $\varepsilon^4\mu'_1$ therefore appears exactly like the Newtonian potential $2v$ of equation \ref{F10'} and the
term $\varepsilon^4M^2\mu'_2$ as one of the magnetic potential terms $Q\frac{\partial\Upsilon_k}{\partial x^k}$ of
 equation \ref{F11'}.

Equation \ref{F2} then corresponds to a classical Schrödinger equation with a Newtonian potential and a
magnetic potential.

If $\varepsilon \ll 1$ and $C_1$, $C_2$ are sufficiently small relative to the characteristic constants of
experiment, these potentials can be considered negligible, and the solution $\varPsi$ of
\ref{F2} very close to that given by \ref{F9'} (for specified boundary conditions).

This principle of comparison between the solutions can be applied without difficulty to equations \ref{F3} and \ref{F4}
using the fact that $\frac{\partial\varPsi}{\partial t}=\varepsilon^2M\mu\varPsi$ and the constraints that
we will now impose on potentials.
\begin{rmq} \label{r2.8}
 The solution of a Schrödinger equation "with active potentials" depends on the "derivatives" of these potentials,
this is why the second condition has been applied to the derivatives of $\mu$ and $\mu'$.

The values of $\varepsilon$ that we choose as well as the constants $C_1$ and $C_2$ to make terms \ref{F12}
negligible, depend on the considered experiment. (Potential may have a negligible effect on results
of an experiment if the time of this one is short, but not to have it for this same type of experiment on a "long" time).
\end{rmq}


\subsubsection{The $\varepsilon$-approximations for potentials}
As will be seen in section \ref{s2.12}, it will be necessary to recover the probabilistic interpretation of the state function
$\varPsi$, that volume elements associated to potential metrics, induced in "spacelike sub-manifolds", are close to those of the metric $g_0$, itself induced in these submanifolds.
This is reflected in the fact that, in the standard coordinate system, the determinant of the matrix $M_{g_{\mathcal P_r}}$ is close to
that of $M_{g_0{_r}}$ when $g_{\mathcal P_r}$ is the potential metric "induced in the spacelike submanifold". The metric $g_{\mathcal P}$ of the potential is generally of the form $g_{\mathcal P}=g_0+h$. In the
examples of potentials that have been given, $detM_{g_{\mathcal P_r}}$ is close to $detM_{g_0{_r}}$ if, in the system
of standard coordinates: $|h_{ij}| \ll 1$, which is written for the active potential without electromagnetism: $|v| \ll 1$ and
for the electromagnetic potential: $|\Upsilon_k| \ll 1$ for $k$ from $0$ to $3$.

On the other hand, if we look at the Schrödinger equations \ref{F10} and \ref{F11}, the "comparable" terms are:
\begin{eqnarray}\label{F15'}
 M\frac{\partial\varPsi}{\varPsi}, ~~vM^2, ~~Q^2\Upsilon^2_{1,2,3}, ~~MQ\Upsilon_0.
\end{eqnarray}
Similarly, according to expression 3. \ref{F12}, ~ $Q\frac{\partial\Upsilon_0}{\partial t}$ is comparable to
$Q^2\Upsilon_0^2$. \\ If we consider valid the $\varepsilon$-approximations for $\varPsi$ given by \ref{F14'}, it 
is natural, so that all terms of expression \ref{F15'} are of the same order in $\varepsilon$, to give
the following definition.
\begin{dfn} \label{d2.18}
 We consider an "elementary oscillating metric domain of order 1 in a potential" for which the
associated state function, associated to a chart of the observation atlas, is written $\varPsi$.

We will say that \textbf{the potentials satisfy the $\varepsilon$-approximations} if there is a real number $0<\varepsilon<1$
such that the following inequality sequences 1. and 2. are satisfied:
\begin{enumerate}
 \item $|v|<\varepsilon^2, ~~|Q\Upsilon_0|<M\varepsilon^2, ~~|Q\Upsilon_{1,2,3}|<M\varepsilon$.
 \item There is a constant $C$ such that, for $k$ from
$0$ to $3$: $|Q\frac{\partial^2\Upsilon_0}{\partial x^k \partial t}|<CM^2\varepsilon^4$, ~ whose choice is based on the
same principle as used for the constants $C_1$ and $C_2$ in definition \ref{d2.17}.
\end{enumerate}
\end{dfn}
The three inequalities in Part 1 of this definition are sufficient to estimate the difference between the elements
associated to the metrics $g_{\mathcal P}$ and $g_o$ induced in "spacelike submanifolds"
rigorously defined in the following section.

Inequalities given in definitions \ref{d2.17} and \ref{d2.18} allow to show, when $\varepsilon \ll 1$
and the constants $C_1$, $C_2$ are well chosen, that terms \ref{F12} appear in equations \ref{F2},
\ref{F3}, \ref{F4} exactly as standard potential terms in Schrödinger's equations and that these
are "negligible".

In the following, we will simply say that an elementary oscillating metric domain of order 1 in a
potential "\textbf{satisfies the $\varepsilon$-approximations} if inequalities given in definitions
\ref{d2.17} and \ref{d2.18} are satisfied with $\varepsilon \ll 1$.


\subsubsection
{The rigorous use of $\varepsilon$-approximations.}

We consider an experiment for which the space-time domain is an "elementary oscillating metric domain
of order 1 in a potential ".

Potentials are assumed to satisfy the $\varepsilon$-approximations (def. \ref{d2.18}) with $\varepsilon \ll 1$ (this
is a datum of the experiment) and that, for given boundary conditions, equation \ref{F2} or \ref{F3}
or \ref{F4} admits an unique solution $\varPsi$ (this equation is of order 2 in $t$). If $\varPsi$
satisfies the $\varepsilon$-approximations (def. \ref{d2.17}) with $\varepsilon \ll1$, then the
stated considerations
previously allow to say that $\varPsi$ is very close to the solution $\varPsi'$ of Schrödinger equation
\ref{F9} or \ref{F10} or \ref{F11} for boundary conditions compatible with the previous ones and
which ensure uniqueness of $\varPsi'$ (Schrödinger equation is of order 1 in $t$).

The method of using $\varepsilon$-approximations that we have just presented is in fact purely theoretical
because, in practice, exact determination of $\varPsi$ is, except for trivial cases, very difficult to obtain.
The main purpose of this section was actually to show that the description of experiments that with our
point of view on physics, uses equations of theorem \ref{2.2}, will give the same description as that obtained by classical quantum physics.


\section{Singularities \label{s2.11}}
It is necessary, in order to accurately describe the experiments that classically use the notion of
"particles", to introduce a notion of localization for domains of oscillating metric type.
This localization experimentally corresponds to the phenomena of "impacts" on a screen, "traces" in a
bubble chamber, drift, wire, etc.

We begin in the first part of this section, by defining the "simple" objects (singularities) that
suffice to describe qualitatively and quantitatively the results of classical experiments in quantum physics
 (particles in a potential, diffraction, Young's slits, Stern-Gerlach type experiments, quantum entanglement
, etc.). As
this text will not go further in the theory that is sufficient for this description, the reader can be satisfied with the
first part of this section and go directly to the next section. However, for the purpose of studying
phenomena much more delicate that are, for example, what is classically called "interactions 
between particles" we will, in the second part of this section, give a 
more precise point of view on the notion of singularity. Currently, the description of these phenomena is addressed by
Q.F.T. The parallel between the theory presented here and the Q.F.T was succinctly presented in remark
\ref{r2.2}.


\subsection{Singularities seen "simply"}
Consider a domain ($\D,g$) representing a geometric type (domain with constant scalar curvature
in a potential, for example) for which the signature of $g$ is everywhere $(-,+,+,+,-,+,\dots,+)$.

\begin{dfn} \label{d2.19}
 A triplet $(\D,g,\Ss)$, where $\Ss$ is a non-empty and zero-measure part of $\D$, will be called \textbf{a domain
with singularities}. The set $\Ss$ will be called \textbf{the singular subset of $\D$}.
\end{dfn}
\begin{dfn} \label{d2.20}
 Consider a domain with singularities $(\D,g,\Ss)$ and $\HH$ a spacelike submanifold of dimension $n-2$
(maximum dimension).

The connected parts of $\HH\cap \Ss$ will be called \textbf{the elementary singularities of $\HH$}.
\end{dfn}
It is tempting to say that the elementary singularities of $\HH$ are the "particles" (in the classical sense) in $\HH$. However, no law will be given on the behavior of these elementary singularities
 and it will be assumed only that the distribution of these in $\HH$ is random with respect to the
metric $g$; this notion will of course be specified in the following section.

Note also that the characteristic quantities usually associated with "particles" (mass, electric charge, spin, momentum, etc.) are for us, as we have already said, characteristics \textbf{of the metric
$g$ in the domain $\D$} (which will usually be "oscillating metric in a potential") and are absolutely 
not associated with singularities. This is profoundly different from the views of other standard physical theories.

It should be noted that elementary singularities are not necessarily points of $\HH$ but only
connected parts whose diameter of each is perfectly defined relative to $g_\HH$ which is, at each point of
$\HH$, a scalar product since $\HH$ is spacelike.

The notion of singularity that we have just presented was introduced by adding an "extra object" (the singular part
 $\Ss$) to the structure of the pseudo-Riemannian manifold $(\D,g)$. However, it is conceivable that this
singular part is simply a characteristic of the pseudo-Riemannian tensor in $\D$, written
now $g_{\Ss}$, and that this one is not defined on the part $\Ss$ of $\D$. The pseudo-Riemannian tensor
$g_{\Ss}$ would then be very close to the tensor $g$ given in definition \ref{d2.19}, except on a
neighborhood of the part $\Ss$ for which the knowledge of the asymptotic behavior of $g_{\Ss}$ (relative to a reference metric $g_0$) would make it possible to approach precisely the study of some interaction phenomena. This point of view is justified in the next paragraph. It shows that oscillating metric type domains can be very close to
fluid type domains presented in section \ref{s1.2} and that collapses of parts of these
domains in singularities can be justified by calculations similar to those used in general relativity
which describe the collapse of some fluid type domains into singularities (black holes for
example). The tensor $g_{\Ss}$ can then be considered as the result of the "evolution" of the tensor $g$
representing an oscillating metric after some local collapses (the word "evolution" requires
to be specified). Singular part $\Ss$ of $\D$ presented in definition \ref{d2.19} is then the part on
which the tensor $g_{\Ss}$ is not defined.


\subsection{Singularities seen as oscillating metric collapses}
Consider an "homogeneous oscillating metric domain in a neutral potential" having a propagation of
constant velocity (def. \ref{d2.12}). The function $a$ is of the form $a=\varphi\beta$ ~~ where ~
~ $\varphi=C\cos(M't-\sum_{k=1}^3\lambda_kx^k+C')$ ~~ and ~~ $\beta\in E_W(\mu)$ ~~ (it is assumed here for simplicity that
$Q^+=0$).

We are interested in the particular case for which $\beta=C^{te}$ which corresponds to a domain of type "oscillating metric"   associated with the notion of "Higgs field" as will be seen in section \ref{s2.15}.

The pseudo-Riemannian metric $g$ is therefore of the form $g=|\varphi|^{4/n-2}g_0$.

 According to the fundamental equation of an oscillating metric: $\Box_{g_0}\varphi+S\varphi=0$. \\ We thus have:
 $(\sum_{k=1}^3\lambda_k^2-M'^2)\varphi+S\varphi=0$ ~~~ hence: $S=M'^2-\sum_{k=1}^3\lambda_k^2$.
 
 The following proposition makes the link between domains of "oscillating metric type" and of "fluid type", it is proved by  starting from the standard expression of the Ricci curvature after a
conformal metric change $g=|a|^\frac{4}{n-2}g_0$: ~~~ on the open set where $a$ does not cancel,

$R_{ij}(g)=R_{ij}(g_0)+\frac{2n}{n-2}a^{-2}\nabla_ia\nabla_ja-2a^{-1}\nabla_i\nabla_ja-\frac{2}{n-2}a^{-2}
(\nabla^ia\nabla_ia+a\nabla^i\nabla_ia)g_{0_{ij}}$

where the covariant derivatives $\nabla_i$ are associated to $g_0$.

\begin{prop} \label{p2.2}
 For the "oscillating metric" type domain that has just been presented, Ricci curvature satisfies:
 \begin{enumerate}
  \item If $S:=\frac{n-2}{4(n-1)}S_{g_0}>0$:
  
  $R_{icc}^{\sharp}(g)=R_{icc}^{\sharp} (g_0)+S\alpha_1X_0\otimes X_0+S\alpha_2g_0$
  
  where ~~ $\alpha_1:=2(\frac{n}{n-2}\tan^2 z+1)$ ~~ and ~~ $\alpha_2:=\frac{2}{n-2}(\tan^2 z-1)$ \\
with ~ $z:=M't-\sum_{k=1}^3\lambda_kx^k+C'$.

$X_0=S^{-1/2}(M'\partial_t-\sum_{k=1}^3\lambda_k\partial_k)$ ~~ and ~~ $g_0(X_0,X_0)=0$
\item If $S_{g_0}=0$:

$R_{icc}^{\sharp}(g)=R_{icc}^{\sharp}(g_0)+S\alpha_1X\otimes X$

where $X=M'\partial_t-\sum_{k=1}^3\lambda_k\partial_k$ ~~ and ~~ $g_0(X,X)=g(X,X)=0$.

Moreover, $D_gX=0$, in particular $X$ is a geodesic field for $g$.
\end{enumerate}
\end{prop}
Note that, when $S>0$, the Ricci curvature is of the form:

 $R_{icc}^{\sharp}(g)=\alpha'_1 X'_0\otimes X'_0-S\alpha'_2Y\otimes Y+S\alpha'_2(Y\otimes Y+X'_0\otimes
X'_0+g)+R_{icc}^{\sharp} (g_0)$

where we set $X'_0=(\cos^{2/n-2}z)X_0$.

So, $g(X'_0,X'_0)=-1, ~~R_{icc}^{\sharp} (g_0)(X'_0,X'_0)=0$, ~~ (we also have $g(Y,Y)=-1$).

As $S_g=0$ since ~ $\Box_{g_0}\varphi+S\varphi=0$, ~~ the tensor $G$ satisfies $G=2R_{icc}(g)$.

\textbf{The domain is then of type "fluid"} according to definition \ref{def:4}. The energy density function is
$S\alpha'_1$, it is positive. The apparent vector field of the fluid is $X'_0$ (here non-geodesic). The apparent pressure
 is $S\alpha'_2(Y\otimes Y+X'_0\otimes X'_0+g)$. The hidden pressure is $R_{icc}(g_0)$.

When $S=0$, the domain is of type "lightlike fluid" since $g(X,X)=0$ and in this case $X$ is a geodesic field. This last domain is comparable to that of type "active potential without electromagnetism" presented in
\ref{ss1.3} for which $R_{icc}^{\sharp}(g)=(\Delta_{g_0}v)X_1\otimes X_1 +R_{icc}^{\sharp}(g_0)$ ~~ where $v$ is the
potential function, $X_1$ is a lightfield and $D_g X_1=0$. Note, however, that the vector field $X_1$
was tangent to $I\times W$ while the field $X$ is tangent to the "apparent" space $\Theta$.

Collapse phenomena of a massing fluid into a singularity (a black hole for example) are described in
classical general relativity. They can be seen mathematically as a consequence of Raychaudhuri's theorem
which says that, under some assumptions, the "expansion $\varTheta$" of a geodesic vector field $X$
characterizing the fluid, goes to $-\infty$ in a finite time. We will see that this same phenomenon is described
simply (even in large dimension) for any light-geodesic field $X$, when the Ricci curvature of the
considered domain is of the form $R^\sharp_{icc}(g)=\alpha X\otimes X+P$ with $R_{icc}(X,X)\geq0$, which is the case
for oscillating metrics just described or for an active potential without electromagnetism.

We consider the classical characteristics of a geodesic vector field $X$ that are: expansion
$\varTheta$, ~~ vorticity $\omega$, ~~ strain $\sigma$ whose definitions are not discussed here.

 Raychaudhuri theorem, in dimension $n$ and for a signature of $g$ of the form \\
$(-,+,+,+,-,\dots,+)$, says that:
if $\beta$ is a lightlike geodesic parameterized by $s$, such that $\dot{\beta}(s)=X_{\beta(s)}$, then
the following egality is satisfied:
\begin{eqnarray}\label{F34}
 (\varTheta_{\beta(s)})'_s=
-R_{icc}(X_{\beta(s)},X_{\beta(s)})
+2\omega_{ij}\omega^{ij}_{\beta(s)}-2\sigma_{ij}\sigma^{ij}_{\beta(s)}-\frac{1}{n-2} \varTheta^2_{\beta(s)}
\end{eqnarray}

It is easy to deduce the following proposition:
\begin{prop} \label{p2.3}
 It is assumed that the vorticity $\omega$ of the field $X$ is zero and that for a value $s_0$ of the parameter: $\varTheta_
{\beta(s_0)}<0$. \\ Then there exists $s_1$ satisfying $s_0<s_1\leq s_0+(n-2)|\varTheta(\beta(s_0)|^{-1}$ such that the
expansion $\varTheta_{\beta(s)}$ goes to $-\infty$ when $s$ goes to $s_1$ by lower values.
\end{prop}
\textbf{Proof}: As $R^\sharp_{icc}(g)=\alpha X\otimes X+R^\sharp_{icc}(g_0)$ and as $X$ is 
lightlike, ~~ $R_{icc}(g)(X_{\beta(s)},(X_{\beta(s)})=0$ ~~ since ~~ $R_{icc}(g_0)(X,X)=0$. ~~ So, as the vorticity
is zero, Raychaudhuri's theorem shows that:

$(\varTheta_{\beta(s)})'_s\leq \frac{-1}{n-2}\varTheta^2(\beta(s))$.

When $\varTheta_{\beta(s)}\neq0$ ~~ we have: $(\varTheta^{-1}(\beta(s))'_s\geqslant\frac{1}{n-2}$, ~~ hence, ~~ for ~~
$s\geqslant s_0$: \\ $\varTheta^{-1}_{\beta(s)}-\varTheta^{-1}_{\beta(s_0)}\geqslant\frac{1}{n-2}(s-s_0)$.

As $\varTheta^{-1}_{\beta(s)}$ is increasing for $s$ and is negative for $s=s_0$, ~~ $\varTheta^{-1}_{\beta(s)}$
goes to $0$ when $s$ goes to $s_1\leq s_0+(n-2)|\varTheta(\beta(s_0)|^{-1}$.

This proposition describes, under the condition of null vorticity, a possible
"collapse" into a singularity for the geodesics generated by the field $X$, this for a finite value of
parameter $s$, provided that the expansion $\varTheta$ takes negative values. In fact, equality \ref{F34}
shows that this phenomenon can take place under much weaker assumptions, and it is not necessary to assume
 nullity of vorticity or $R{icc}(X,X)$ since only inequality $(\varTheta_{\beta(s)})'_s\leq
\frac{-1}{n-2}\varTheta^2(\beta(s)$ is used. This suggests that the phenomenon of "collapse" of an
oscillating metric in a singularity takes place for much more general oscillating metrics than those where
$R^\sharp_{icc}(g)$ is of the form $\alpha X\otimes X+P$. Proof of proposition \ref{p2.3} shows that
it's the
fact that the expansion $\varTheta$ is negative in some domain that "causes" the phenomenon of collapse. It
is conceivable that only subset of the oscillating metric domain collapses into a singularity and that we get
then an "oscillating metric domain with singularities".

Specifying the behavior of the metric $g$ near a singularity will only matter for the study of more complex phenomena. For  physical phenomena described in this paper up to section
\ref{s2.17}, only the "localization" of the singularities (random relatively to $g$) will be used. This is the subject of the next section.


\section{The probabilistic part \label{s2.12}}
In the axiomatic system that we have chosen to represent space-time, no principle, no law governs the singularities. The fact that, in the generic experiments of standard
quantum physics (diffraction, Young's slits, etc.) the singularities do not appear equiprobably
distributed on "the screen"~~ (we consider that it is the singularities that leave traces on the screen)
\textbf{is not due to laws that would govern these singularities, but only to the distortion of the pseudo-Riemannian metric
 $g$, in the subset of space-time where they are, relative to the metric $g_0$ used during the
measurement of the results of the experiment}. As expected, we will obtain probabilistic results to describe results of these experiments and, to some approximations, these results will be identical to those
of standard quantum physics, without using its axiomatic system.
\bigskip

We consider a domain with singularities $(\D,g,\Ss)$ (def. \ref{d2.19}).

\textbf{The fact that no law governs the singularities of $\D$ is mathematically expressed in the following way}:

Any spacelike submanifold $\HH$ of $\D$ of dimension $n-2$ (ie. of maximum dimension) satisfies the two
following properties:
\begin{enumerate}
 \item If $\HH_2 \subset \HH_1$ are two subsets of $\HH$ with finite volumes associated to the Riemannian metric
 $g_{\HH}$. \\
If "$\varsigma$" is an elementary singularity in $\HH_1$ (def. \ref{d2.10}),
then the probability that "$\varsigma$" is in $\HH_2$ is:
\begin{eqnarray}\label{F35}
 p=\text{vol}_{g_{\HH}}(\HH_2)~/~\text{vol}_{g_{\HH}}(\HH_1)
\end{eqnarray}
In other words, the probability density of presence of the singularity "$\varsigma$" in $\HH_2$ is the uniform 
density given by the $(n-2)$-differential form on $\HH_1$ defined by\\
$\sigma=(\text{vol}_{g_{\HH}}(\HH_1))^{-1}dv_{g_\HH}$ where $dv_{g_\HH}$ classically refers to the $(n-2)$-Riemannian volume form on $(\HH,g_{\HH})$.
  
(When $(v^1,\dots,v^{n-2})$ is a coordinate system on $\HH$, the volume form is written
: $dv_{g_\HH}=\sqrt{det g_\HH}dv^1\dots dv^{n-2}$ where
$det g_\HH$ is the determinant of the $g_\HH$ matrix in the coordinate system).
\item If $\varsigma_1,\dots,\varsigma_N$ are $N$ elementary singularities in $\HH_1$, then the probability that $k$
of them are in $\HH_2$ is:
\begin{eqnarray}\label{F36}
 p(N,k)=\displaystyle\binom{N}{k}p^kq^{N-k}
\end{eqnarray}
where $p=\text{vol}_{g_{\HH}}(\HH_2)~/~\text{vol}_{g_{\HH}}(\HH_1)$ ~~ and ~~ $q=1-p$

In other words, the associated probability law is a binomial law. This highlights the fact that the presence of an
elementary singularity in $\HH$ has no influence on the presence of others.
\end{enumerate}
Note that the properties 1. and 2. are not linked to the choice of a chart of the observation atlas.
\bigskip

Properties 1. and 2. exactly correspond to the experiment consisting in randomly throwing
punctual objects on a surface $\HH_1$ and calculating the probabilities of arrival of these objects on a surface
$\HH_2\subset\HH_1$ (but for us, of course, $\HH_2$ and $\HH_1$ are of dimension $n-2$ and the metric is
$g_{\HH_1}$). \\
It is essentially the property 1 that will be used, property 2 will only be used when we introduce the notion of "density of singularities". \\
In fact, it will be sufficient to consider that properties 1 and 2 are valid only in the space-time domains that correspond to the concerned experiments.
\bigskip

We now consider a chart of the observation atlas for which the cell $(\C,g)$ is that of an
oscillating metric in a potential (def. \ref{d2.2}): $\C=\Theta\times S^1(\delta)\times W$ and
$g=|a|^\frac{4}{n-2}g_\mathcal P$. ($\Theta$ will be assumed of the form $I\times \mathscr U\subset \R\times\R^3$.)

We "fix" $(t,u)\in \R\times S^1(\delta)$ ~~ ("double" time.)

Let $\HH_{2_{(t,u)}}\subset \HH_{1_{(t,u)}}\subset \C$ be of the form:

$\HH_{2_{(t,u)}}=\{t\}\times \omega\times\{u\}\times W\subset\C$

$\HH_{1_{(t,u)}}=\{t\}\times \Omega\times\{u\}\times W\subset\C$

where $\omega\subset\Omega\subset \R^3$ are two domains in $\R^3$.

To visualize things, the reader can assume that $\Omega\subset\R^3$ is a domain representing a screen
(with a thickness) and $\omega$ a subdomain of this screen during a standard experiment which consists of studying the
distribution of the impacts of "particles" on the screen, when these are sent through a potential
$\mathcal P$.

$\HH_{2_{(t,u)}}$ and $\HH_{1_{(t,u)}}$ are spacelike submanifolds of dimension $n-2$ for the two metrics
$g$ and $g_0$.

According to \ref{F35}, when "$\varsigma$" is an elementary singularity in $\HH_{1_{(t,u)}}$, the probability
that it lies in $\HH_{2_{(t,u)}}$ is:
\begin{eqnarray}\label{F37}
 p_{(t,u)}=\int_{\HH_{2_{(t,u)}}}dv_{g_{\HH_{1_{(t,u)}}}}~/~\int_{\HH_{1_{(t,u)}}}dv_{g_{\HH_{1_{(t,u)}}}}
\end{eqnarray}
In a standard coordinate system of $\C$, this is written:
\bigskip

$p_{(t,u)}=\int_{\HH_{2_{(t,u)}}}\sqrt{det(g_{\HH_{1_{(t,u)}}})}dx^idv^j~/~\int_{\HH_{1_{(t,u)}}}\sqrt{det(g_{\HH_{1_{(t
,u)}}})}dx^idv^j$
\bigskip

 (here ~~ $dx^idv^j:=dx^1dx^2dx^3dv^1\dots dv^{n-5}$)
 \bigskip
 
 hence:
 \begin{eqnarray}\label{F38}
p_{(t,u)}=\int_{\HH_{2_{(t,u)}}}a^2_{(t,u)}\sqrt{det(g_{{\mathcal P}_{\HH_{1_{(t,u)}}}})}dx^idv^j~/
\int_{\HH_{1_{(t,u)}}}(idem)chart
 \end{eqnarray}
\textbf{WE WILL NOTICE HERE THE IMPORTANCE OF THE SIGNATURE OF \boldmath $g$ \unboldmath  WHICH COUNTS TWO SIGNS \boldmath "$-$" \unboldmath AND \boldmath $(n-2)$ \unboldmath SIGNS
\boldmath "$+$" \unboldmath because it is this fact that makes it possible to obtain the exponent "2" on the function \boldmath $a$ \unboldmath  since the dimension of \boldmath $\HH_1$ \unboldmath is then
equal to \boldmath $(n-2)$ \unboldmath and therefore: \boldmath $\sqrt{det(|a|^{\frac {4}{n-2}}g_{{\mathscr
P}_{\HH_{1_{(t,u)}}}})}=a^2\sqrt{det(g_{{\mathcal P}_{\HH_{1_{(t,u)}}}})}$ \unboldmath.\\
(the exponent "2" on the function \boldmath $a$ \unboldmath 
is particularly important in the rest)}.

It is now assumed that the metric $g_{\mathcal P}$ is "sufficiently" close to $g_0$, in other words that the
"potential functions" that occur in $g_{\mathcal P}$ are $\ll 1$. We can then write:
\begin{eqnarray}\label{F39}
 p_{(t,u)}\simeq\int_{\HH_{2_{(t,u)}}}a^2_{(t,u)}\sqrt{det(g_{0_{\HH_{1_{(t,u)}}}})}dx^idv^j~/~\int_{\HH_{1_{(t,u)}}}
( idem)
\end{eqnarray}
An estimation of the error in the "equality" ~ $\simeq$ ~ can be made using the definition of
$\varepsilon$-approximations on potentials (def. \ref{d2.18}), but we won't develop this here.
\bigskip

Let us now examine the case of an \textbf{elementary} oscillating metric of order 1 (def. \ref{d2.9}) for which,
by
definition:

$a=\varphi\beta$ ~~ where ~~ $\varphi:\Theta\times S^1(\delta))\rightarrow\R$ satisfies:
$\varphi=\varphi_1\cos(Q^+u)+\varphi_2\sin(Q^+u)$

Here $\varphi_1$ and $\varphi_2$ are two real functions defined on $\Theta$ ~ and ~ $\beta\in E_W(\mu)$.

We recall that the canonical function $a_c$ is defined by $a_c=\varphi_1+i\varphi_2$.

(The case of elementary oscillating metrics of the 2nd order is treated in a similar way but will be presented, for
the clarity of the exposition, only in the section concerning the "spin". The process can be generalized without
difficulties when the order $k>2$ considering the decomposition of $W$ into a product of compact manifolds).

As $g_0$ is a "product" metric, \ref{F39} is written:
\bigskip

$p_{(t,u)}=(\int_\omega\varphi^2_{(t,x,u)}dv_{g_{0|_\omega}})(\int_W\beta^2dv_{g_{0|_W}}) ~/
~(\int_\Omega\varphi^2_{(t,x,u)}dv_{g_{0|_\Omega}})(\int_W\beta^2dv_{g_{0|_W}})$
\bigskip

where the integral of $\beta^2$ on $W$ can be deleted.

We thus obtain:
\begin{multline}\label{F40}
p_{(t,u)}\simeq\int_\omega(\varphi^2_1{_{(t,x,u)}}\cos^2(Q^+u) +
\varphi^2_{2_{(t,x,u)}}\sin^2(Q^+u)+\\
2\varphi_1\varphi_2{_{(t,x,u)}}\cos(Q^+u)\sin(Q^+u))dx^1dx^2dx^3~ / ~\int_\Omega (idem)
\end{multline}

In the particular case where $Q^+=0$, the canonical function $a_c$ is the function $\varphi_1$.
The previous "equality" is then written:
\bigskip

$p_{(t,u)}=p_{(t)}\simeq\int_\omega\varphi^2_1{_{(t,x)}}dx^i ~/ ~\int_\Omega\varphi^2_1{_{(t,x)}}dx^i$
\bigskip

hence:
\bigskip

$p_{(t)}\simeq\int_\omega|a_c|^2_{(t,x)}dx^i~/ ~\int_\Omega|a_c|^2_{(t,x)}dx^i $
\bigskip

When $Q^+$ is positive, it is interesting to consider the "averages" over the "time"
$u\in S^1(\delta)$, that is equivalent to assume that the variations in $u$ of $p_{(t,u)}$ given in \ref{F40} are not
visible in practice.

It is therefore assumed that, for an observer associated to the standard coordinate system of the cell $\C$, the density of
probability of the presence of a singularity $\varsigma$ in $\HH_2{_{(t)}}:=\{t\}\times\omega\times S^1(\delta)\times
W$ is uniform with respect to time $u\in S^1(\delta)$, which is expressed as follows:

When, for $t\in\R$, "$\varsigma$" is an elementary singularity in $\HH_1{_{(t)}}:=\{t\}\times\Omega\times
S^1(\delta)\times W$, the probability that it lies in $\HH_2{_{(t)}}$ is:
\begin{eqnarray}\label{F41}
 p_{(t)}\simeq\int_{\HH_2{_{(t)}}}\sqrt{det(g_{\HH_{1_{(t,u)}}})}dudx^idv^j ~/
~\int_{\HH_1{_{(t)}}}\sqrt{det(g_{\HH_{1_{(t,u)}}})}dudx^idv^j
\end{eqnarray}
In the case of an elementary oscillating metric of order 1 and with the same approximations as those described
previously, \ref{F41} becomes (starting from \ref{F40}):
\vspace{0cm}

$p_{(t)}\simeq\int_{\omega\times S^1(\delta)}(*)dx^idu~ / ~\int_{\Omega\times
S^1(\delta)} (*)dx^idu$
\vspace{1mm}

where $(*)=\varphi^2_1{_{(t,x,u)}}\cos^2(Q^+u) +\varphi^2_{2_{(t,x,u)}}\sin^2(Q^+u)+
2\varphi_1\varphi_2{_{(t,x,u)}}\cos(Q^+u)\sin(Q^+u)$
\vspace{1mm}
hence, since: $\int_{S^1(\delta)}\cos^2(Q^+u)du=\int_{S^1(\delta)}\sin^2(Q^+u)du=\pi\delta$\\
and $\int_{S^1(\delta)}\cos(Q^+u)\sin(Q^+u)du=0$:
\vspace{2mm}

$p_{(t)}\simeq\int_{\omega}(\varphi_1^2+\varphi_2^2)_{(t,x)}dx^i ~/ ~\int_{\Omega}(\varphi_1^2+\varphi_2^2)_{(t,x)}dx^i$
\vspace{5mm}

One recovers the expression that one had obtained in the particular case for which $Q^+$ was zero:
\begin{eqnarray}\label{F42}
 p_{(t)}\simeq\int_{\omega}|a_c|^2_{(t,x)}dx^i ~/ ~\int_{\Omega}|a_c|^2_{(t,x)}dx^i
\end{eqnarray}

Of course, a singularity $\varsigma$ in ${\HH_1{_{(t)}}}$ (resp. ${\HH_2{_{(t)}}}$) is seen "in practice"
in $\Omega$ (resp. $\omega$) at time $t$ by the observer associated to the coordinate system.

Probabilities given by equalities \ref{F36} may also be extended to the case "average over $S^1(\delta)$", the reader may precisely formalise this result.
\vspace{2mm}

Consider the particular case where the elementary oscillating metric has a well-defined electric charge (def.
\ref{d2.14}). The state function $\varPsi$ is then itself well defined (def. \ref{d2.16}) and we have
$|\varPsi|=|a_c|$.
Probability $p_{(t)}$ given by \ref{F42} can therefore also be written:
\begin{eqnarray}\label{F43}
 p_{(t)}\simeq\int_{\omega}|\varPsi|^2_{(t,x)}dx^i ~/ ~\int_{\Omega}|\varPsi|^2_{(t,x)}dx^i
\end{eqnarray}

Here we recover one of the main axioms of classical quantum physics, although the
denominator $\int_{\Omega}|\varPsi|^2_{(t,x)}dx^i$ may depend on "$t$". This fact is not a conceptual problem
 in the physic theory presented here in contrast to standard quantum physics, let's remember that
"singularity" is not a notion equivalent to that of "particle" (see introduction to chapter 2).

If we assume the validity of the $\varepsilon$-approximations, we have seen that the state function $\varPsi$ satisfies (in
approximation) classical Schrödinger equations which are of the form $i\frac{\partial \varPsi}{\partial
t}=H(\varPsi)$ where $H$ is an Hermitian operator, and if $\varPsi$ satisfies good boundary conditions we classicaly deduce
that $\partial({\int_\Omega|\varPsi|^2}) ~/~\partial t=0$ and that then $\int_\Omega|\varPsi|^2$ does not depend
on "$t$". We can then, as in standard quantum physics, \textbf{normalize} the state function $\varPsi$
so that $\int_\Omega|\varPsi|^2=1$ and one recovers exactly with \ref{F43} the considered standard axiom. This process
is valid only when $\varPsi$ satisfies a Schrödinger equation (of order 1 in $t$) and not for a Klein-Gordon equation which is of order 2 in $t$, but, as we have already said, this poses no conceptual problem for the physics presented here.

The important axiom of quantum physics given by \ref{F43} "normalized" is for us only a consequence of the
fact that no law governs the singularities on $\M$ which is expressed mathematically by \ref{F35}, and,
in the case of oscillating metrics in a potential, give \ref{F38}. The probability given by \ref{F38} depends on the
potentials, which in practice seriously complicates things, but in the context of the
$\varepsilon$-approximations, dependence on these potential is negligible and we thus recover, for
 experiments studied by standard quantum physics, identical results.


\subsection{The notion of "density of singularities"}\label{4-1}
The concept introduced in this subsection is in fact not really essential for the description of standard 
quantum physics experiments, but it will become useful in the study of quantum phenomena currently described
by Q.F.T. This notion of "density of singularities" will allow, among other things, to address quantitatively the
problem of determining the conditions for which the linear fundamental equation \ref{F1} is a good
approximation for the fundamental nonlinear equation \ref{F0'}.

We consider, as before, a spacelike submanifold $\HH$ (in a domain $\D$) of dimension $n-2$.
When the number of elementary singularities in $\HH$ is "large enough" we express that no law
 governs the singularities in $\D$ (more precisely, we express the equiprobability of presence of singularities in $\HH$
relative to $g_\HH$) by assuming the following:

There is a constant $D$ (the density) such that, for any $N\in\N$ and any domain $\HH_N\subset\HH\subset\D$
which contains $N$ elementary singularities, we have:
\begin{eqnarray}\label{F44}
 N(1+\varepsilon(N))=Dv_{g_{\HH_N}}
\end{eqnarray}
where $v_{g_{\HH_N}}$ denotes the volume associated to the Riemannian metric $g_\HH$ in the domain $\HH_N$, ~~ and ~~ $lim
~~ \varepsilon (N)=0$ when $N$ goes to infinity.

We consider a sequence of domains $(\HH_N)_{N>N_0}$ such that $v_{g_{\HH_N}}$ goes to infinity
with $N$, and $\HH_1$ a subdomain of $\HH_N$.

According to \ref{F36}, when $\varsigma_1,\dots,\varsigma_N$ are the $N$ elementary singularities in $\HH_N$, the
probability that $k$ of them are in $\HH_1$ is:
\vspace{2mm}

$p(N,k)=\displaystyle\binom{N}{k}p^kq^{N-k}$
\vspace{2mm}

where $p=v_1/v_N$, ~~~ $q=1-p$ ~~~ when we have set $v_N:=v_{g_{\HH_N}}$.

Then, using relation \ref{F44}, we deduce that:
\vspace{0mm}

$p(N,k)=\frac{N!}{k!(N-k)!}(v_1/v_N)^k(1-v_1/v_N)^{N-k}$
converges to $\frac{(Dv_1)^k}{k!}e^{-Dv_1}$
\vspace{1mm}

when $N$ goes to infinity.

We can therefore consider that, when there is a large number of elementary singularities in $\HH$, the probability that there are $k$ in a subdomain $\HH_1$ of $\HH$ is:
\begin{eqnarray}\label{F44'}
 p(k)=\frac{(Dv_1)^k}{k!}e^{-Dv_1} 
\end{eqnarray}
(We can notice that we have $\sum_{k=0}^\infty 	p(k)=1$)
\vspace{2mm}

When $(\C,g)$ is a cell of type "oscillating metric in a potential" for which
$g=|a|^{\frac{4}{n-2}} g_{\mathcal P}$, and when admitting the $\varepsilon$-approximations such that $g_{\HH(t)}$
is "close" to $|a|^{\frac{4}{n-2}}g_0{_{\HH(t)}}$ as before, relation \ref{F44'} shows that at
time "$t$" for an observer associated to a standard coordinate system of $\C$, the probability that
 $k$ elementary singularities are in $\HH_1(t)$ ~ is ~ $ p(k)=\frac{(Dv_1(t))^k}{k!}e^{-Dv_1(t)}$.

(In fact, these singularities are "seen" in $\omega$ when $\HH_1(t)=\{t\}\times \omega\times S^1(\delta)\times
W$)

 where here $v_1(t)=\int_{\omega\times S^1(\delta)\times W}a^2dx^idudv^j$.
 \vspace{2mm}
 
 Everything happens as if the density of elementary singularities \textbf{seen by the observer associated to the system of
coordinates} (which makes the measurements with $g_0$) was given by the function $Da^2$.

If we assume that the function $a^2$ is small enough compared with the constant $1$ that linear equation
\ref{F1} satisfied by the function $a$ is a good approximation of equation \ref{F0'}, this is no longer
necessarily the case for the function $(\lambda a)^2$ where $\lambda$ is a constant "sufficiently large", althougt the function $\lambda a$ satisfies the linear equation \ref{F1}. From the previous results this says
simply that, when the density of singularities becomes large (it grows in $\lambda^2)$, the linear equation
ceases to be a good approximation and one can no longer neglect the non-linear part of equation \ref{F0'}. It
corresponds to the classical physics interpretation that when the particle density becomes large, we
can no longer neglect the interaction of particles with each other, but we must understand that, with our point of view of
 physics, we do not assume any kind of interaction between the singularities, it's only the oscillating metric,
through the function $a$, which takes into account this phenomenon. Of course, the density
of singularities that appear in an experiment depends on the initial (or boundary) conditions of that
experiment.
\begin{rmq} \label{++r2,4}
In the context of a "particles in a potential type domain" presented in the remark \ref{++r2,1} in the introduction of chapter 2, the notion of "density of singularities" that we have just presented applies to each "bubble" $(\D_k,g_k)$ for which the singularities are "identifiable" because these domains carry the characteristics of corresponding oscillating metrics (mass, electric charge, spin, etc.). On the other hand, the singularities of the domain $\D_0$ are not necessarily identifiable because the metric $g_\mathcal P$ is that of a potential. It follows that the notion of "density of identifiable singularities" in the whole domain $\D$ is also associated to the "bubble density" in $\D$.
 
\end{rmq}


\section{The spin \label{s2.13}}
 Description of the experimental results associated to "spin" phenomena, usually treated by 
standard quantum physics, turns out to be natural to us if we consider that, in the cell of the form
$\C=\Theta\times S^1(\delta)\times W$, the compact manifold $W$ is a product manifold of the form
$S^3(\rho)\times V$ where $S^3(\rho)$ is the standard sphere of dimension $3$, radius $\rho$, with the standard Riemannian metric $g_{S^3(\rho)}$ induced by the $\R^4$ Euclidean metric (to describe the reference metrics $g_0$).

Study of $(S^3(\rho),g_0|_{S^3(\rho)})$ and in particular that of eigenspaces of the Laplacian operator is very
important for us and is specified in Annex \ref{sa3.7}. We recall the essential points in the two
subsections that follow.


\subsection{Eigenvalues and eigenspaces of the Laplacian operator on the sphere $S^3(\rho)$ as well as those associated to the
Hopf fibration}
Eigenspaces $E_{S^3(\rho)}(\gamma_p)$ of the (geometric) Laplacian operator ($\Delta:=-{\triangledown_i}{\triangledown^i}$)
have for eigenvalues: $\gamma_p=\rho^{-2}p(p+2)$ ~~~~ $p\in \N$

We will classify them according to the increasing values of $\gamma_p$ by noting them:
\begin{eqnarray}\label{F45}
 E_0,E_1,\dots,E_p,\dots
\end{eqnarray}
In fact, $p/2$ will correspond to the "spin index" of standard quantum physics.

It is important to note that when \textbf{$p$ is even}, each eigenspace $E_p$ contains a linear subspace
 $E'_p$ that can be identified with the eigenspace $E_{S^2(\rho/2)}(\gamma_p)$ of the Laplacian operator of the standard sphere $S^2$
of radius $\rho/2$ (Annex \ref{sa3.7}).

The odd indices of $E_p$ correspond to the eigenspaces linked only to $S^3$ whereas the even indices are
associated to the eigenspaces of $S^2$. It is interesting to note that, \textbf{if $p$ is even}, eigenfunctions
$\beta\in E_p$ are invariant by the "antipodal" isometry $\sigma$ of $S^3(\rho)$ ~~ ($\beta\circ\sigma=\beta$),
and change sign ~~ ($\beta\circ\sigma=-\beta$) \textbf{if $p$ is odd}.

\subsection{Three vector fields that parallelize $S^3$ and the endomorphisms of their canonically associated  eigenspaces}
The sphere $S^3$ is a parallelizable manifold (only the spheres $S^1,S^3,S^{15}$ are).

We classically define $S^3(\rho)$ by setting:

$S^3(\rho):=\{(x_1,x_2,x_3,x_4	)\in \R^4 /~ \sum_{k=1}^4x_k^2=\rho^2 \}$

Let $L_1,L_2,L_3$ be the three vector fields defined on $\R^4$ by:
\begin{multline}\label{F46}
 L_1=-x_3\partial_1-x_4\partial_2+x_1\partial_3+x_2\partial_4\\
 L_2=-x_4\partial_1+x_3\partial_2-x_2\partial_3+x_1\partial_4\hspace{8.6cm}\\
 L_3=x_2\partial_1-x_1\partial_2-x_4\partial_3+x_3\partial_4\hspace{8cm}
\end{multline}
where $\partial_i$ is for $\frac{\partial}{\partial x_i}$.

 At each point $x=(x_1,x_2,x_3,x_4)$ these three vector fields are orthogonal to the radial vector
$x_1\partial_1+x_2\partial_2+x_3\partial_3+x_4\partial_4$ and are orthogonal to each other (these are, in a way, the
more "simple" that we can write in the coordinate system and have these properties).

We write $L_k|_{S^3}$ the three vector fields on $S^3(\rho)$ "restrictions" of $L_k$. At each
point
$x\in S^3(\rho)$ the three vectors $L_k|_{S^3}$ form an orthogonal basis of the tangent space $T_x(S^3(\rho))$.
($S^3$ is therefore parallelizable).

$L_1$, $L_2$, $L_3$ naturally define (as differential operators) three endomorphisms of
$C^\infty(\R^4)$ since for any $ f\in C^\infty(\R^4)$ ~~ $L_k(f)\in C^\infty(\R^4)$.

Similarly, $L_1|_{S^3}$, $L_2|_{S^3}$, $L_3|_{S^3}$ define, by restriction, three endomorphisms
of $C^\infty(S^3(\rho))$.

The following proposition is fundamental. His proof is given in Annex \ref{ss3.2}.
\begin{prop} \label{p2.4}
 The eigenspaces $E_{S^3(\rho)}(\gamma)$ written $E_p$ ~ ($p\in \N$) specified in \ref{F45} and spaces
$E'_q$ ~ ($q\in 2\N$) \textbf{are invariant} by each of the three endomorphisms defined by $L_k|_{S^3}$.
\end{prop}
This result gives the following definition:
\begin{dfn} \label{d2.21}
 For each $E_p$ (resp. $E'_q$), the three endomorphisms written $S_1$, $S_2$, $S_3$ (without reference to $p$ or
$q$) defined by $L_k|_{S^3}$ induced in $E_p$ (resp. $E'_q$) will be called \textbf{the canonical endomorphisms of
eigenspaces of $S^3(\rho)$}.
\end{dfn}
As we have already seen (prop \ref{p2.1'}), the eigenspaces on $S^1(\delta)\times
S^3(\rho)$ associated to the pseudo-Riemannian metric $(g_0|_{S^1})\times (g_0|_{S^3})$ can be identified with
$E_{S^1}(\lambda)\otimes E_{S^3}(\gamma_p)$. The corresponding eigenvalues are:\\ $(\gamma_p-\lambda)$ ~ where ~
$\gamma_p=\rho^{-2}p(p+2)$ ~ and ~ $\lambda={Q^+}^2$.

When $\lambda\neq 0$, the eigenspace $E_{S^1}(\lambda)$ can be identified with $\CC$ by isomorphism ~ $\CC_\lambda$~(see \ref{F12}) \\ and ~ $E_{S^1}(\lambda)\otimes E_{S^3}(\gamma_p)$ with $\CC\otimes E_p$ ~ (see \ref{F13}) ~ which is 
 the complexification of $E_p$ written $E^{\CC}_p$.

For each $E_p$ (resp. $E'_q$) the three endomorphisms $S_k$ (def. \ref{d2.21}) naturally extend over the
complexification $E^{\CC}_p$ (resp. $E'^{\CC}_q$) by setting:

for any $ (\varphi_1+i\varphi_2)\in E_p+iE_p=E^{\CC}_p$
~~~~~ $S^{\CC}_k(\varphi_1+i\varphi_2):=S_k(\varphi_1)+iS_k(\varphi_2)$

However, only in order to recover exactly the endomorphisms involved in standard quantum physics
 for the phenomena associated to spin, one gives the following definition where one introduces the coefficient "$\textstyle{\frac{i}{2}}$".
\begin{dfn} \label{d2.22}
For each $E_p$ (resp $E'_q$), the three endomorphisms of $E^{\CC}_p$ (resp. $E'^{\CC}_q$) defined by 
$\hat{S}_k:=\textstyle{\frac{i}{2}} S^{\CC}_k$ will be called \textbf{the canonical endomorphisms of the complexified eigenspaces of
$S^3(\rho)$}.
\end{dfn}
(Note: in the version of this paper, definitions of $\hat{S}_k$ and the gyromagnetic constant presented next, differ from previous versions  by a constant factor, so that they now correspond exactly to those of classical quantum physics , it follows that in some equations some coefficients differ from those of previous versions). \\
The reader will be able to verify that the three endomorphisms $\hat{S}_1$, $\hat{S}_2$, $\hat{S}_3$
have the following properties which are none other than those satisfied by the observables of kinetic momentum (but for a 
factor $\hbar$) in classical quantum physics: \\
$\hat{S}_2\hat{S}_3-\hat{S}_3\hat{S}_2=i\hat{S}_1$, ~~ $\hat{S}_1\hat{S}_3-\hat{S}_3\hat{S}_1=-i\hat{S}_2$,
~~ $\hat{S}_1\hat{S}_2-\hat{S}_2\hat{S}_1=i\hat{S}_3$


\subsection{Domain of type "Oscillating metric with spin in a potential"} \label{ss2.15+}
Here we will focus on domain of type "elementary oscillating metric \textbf{of order 2} in a
potential" (def. \ref{d2.9}).

Experimentally, the particular physical effects associated to the notion of spin appear in areas where
"electromagnetism" is present. We will therefore mainly study the case of oscillating metrics with
spin in an \textbf{electromagnetic} potential.

To clarify things we give the following definition:
\begin{dfn} \label{d2.23}
  \textbf{A domain of type "elementary oscillating metric with spin in an electromagnetic potential"} is
a domain of type "elementary oscillating metric of order 2" (def. \ref{d2.9}) for which the 
pseudo-Riemannian metric $g_\mathcal P$ of the electromagnetic potential has the specific form presented in the next paragraph.
\end{dfn}


\subsubsection{Specific form of the ectromagnetic potential $g_\mathcal P$}
Recall that the cell $\C$ is of the form $\C=\Theta\times S^1(\delta)\times S^3(\rho)\times V$, and
according to proposition \ref{p1.3}, ~~ $g_\mathcal P=g_0+sym(\Upsilon^\flat\otimes X^\flat)$.

In the study of elementary oscillating metrics of order 1, $\Upsilon$ was chosen as a vector field defined on $\Theta$ (that is, defined to $\C$ but tangent to $\Theta$ and depending only on the variables
of $\Theta$). This can be interpreted as the fact that we neglect the "quantum effects" associated to $W$ (but
not to $S^1(\delta)$).

We now consider that $\Upsilon$ is a vector field defined on $\Theta\times S^3(\rho)$ (we neglect the
quantum effects associated to $V$ but not to $S^1(\delta)\times S^3(\rho)$).

The vector field $\Upsilon$ decomposes naturally in the form $\Upsilon=A+C$ where $A$ is the component of
$\Upsilon$ tangent to $\Theta$ and $C$ the tangent component to $S^3(\rho)$.

The vector field $A$ will be considered to represent the classical electromagnetic potential and it will be assumed
that it depends only on the variables of $\Theta$ (this last point is obviously to be interpreted as an
"approximation"). It is written in the standard coordinate system: $A=\sum_{i=0}^3 A^i\partial_i$ ~ where ~
$\partial_i:= \frac{\partial}{\partial x^i}$.

We denote $C^k$ the components of the vector field $C$ tangent to $S^3(\rho)$ associated to the basis of the three vector fields $L_1|_{S^3}$,
$L_2|_{S^3}$, $L_3|_{S^3}$ which parallelize $S^3(\rho)$. We have:\\
 $C=\sum_{k=1}^3C^kL_k|_{S^3}$.

It should be noted that, in "geometric units", the components $C^k$ are "inverse lengths"
(since the $L_k$ have components of the form $x_k\partial_k$) while the components $A^j$ are "without
unit".\\

\textbf{The specific form of the electromagnetic potential will be mainly due to the particular choice of the components
$C^k$  which we will specify in the lines that follow}.
 It should be noted that we keep the fact that the nilpotence index of the endomorphism associated to $sym(\Upsilon^\flat\otimes X^\flat)$ is at most 3, the hypothesis of nilpotence will be deleted in section \ref{++2,3}, but in this section we are content to simply recover the standard results.\\
 When the metric g is of the form $g=|a|^{\frac{4}{n-2}} g_{\mathcal P}$, where:
 
 $a = \phi \beta$~~ whith $\phi: \Theta\times S^1(\delta)\times S^3(\rho)\rightarrow\R$~~ and $\beta\in E_V(\nu)$,
 
 $g_\mathcal P=g_0+sym(\Upsilon^\flat\otimes X_2^\flat)$~~(cf proposition \ref{p1.3})~~whith $\Upsilon=A+C$ previously specified,\\
 then the fundamental equation \ref{F1} is expressed (under the assumptions of theorem \ref{2.3}), by (cf dem. theorem \ref{2.3} part 3 Annex \ref{a3.7}):

\begin{eqnarray}\label{F49''}
\sum_{j=0}^3\varepsilon_j(i\frac{\partial}{\partial x^j}+Q^+\Upsilon^j)^2a_c+M^2a_c+2
Q^+\sum_{k=1}^3C^k\hat S_k(a_c)+{Q^+}^2\rho^2|C|^2a_c=0
\end{eqnarray}

where ~~ $\varepsilon_j=g_{0jj}$, ~~$|C|^2:=\sum_{k=1}^3{C^k}^2$, ~ $\rho$ is the radius of the sphere $S^3(\rho)$, ~ the $\hat S_k$ are the canonical isomorphisms of the complexified eigenspaces of $S^3(\rho)$
~ (def. \ref{d2.22}) ~~ and ~~ $\hat S_k(a_k):\Theta\rightarrow E^{\CC}_p$ is defined for any $ x\in \Theta$ ~ by ~ $\hat
S_k(a_c)(x):=\hat S_k(a_c(x))$.\\

\textbf{In the spin 1/2 case} for which $E^{\CC}_p=E^{\CC}_1$, \textbf{if we choose $C^k=B^k$} where the function $B^k$ are the components of the magnetic field (in the coordinate system):\\
$B^1=\frac{\partial A^3}{\partial x^2}-\frac{\partial A^2}{\partial x^3}$, ~~ $B^2=\frac{\partial A^1}{\partial
x^3}-\frac{\partial A^3}{\partial x^1}$, ~~ $B^3=\frac{\partial A^2}{\partial x^1}-\frac{\partial A^1}{\partial
x^2}$,\\
\textbf{then}, the 3 first terms of the left hand side of equation \ref{F49''} \textbf{can be factorized} and this allows to recover \textbf{standard Dirac equation} of the electron in a magnetic field. This is detailed in the subsection \ref{ssAv}.\\
\textbf{To generalize the study to spin other than 1/2 we then choose the vector field $C$ of the form:}\\
$$C:={\textstyle\frac{\varrho}{2}}\sum_{k=1}^3B^kL_k|_{S^3}$$
Where $\varrho$ is a constant and the coefficient $\textstyle\frac{1}{2}$ is choosen so that $\varrho$ corresponds to the standard gyromagnetic constant (the Landé factor) of classical quantum physics ($\varrho=2$ for the electron). The three functions $B^k$ are the three components of the magnetic field ($\varrho$ is a constant "without unit" and the $B^k$ have "units" as inverse lengths).\\

 The vector field $\Upsilon$ defined in this way is associated to the choice of the coordinate system, it is not invariant by  Poincaré transformations on $\Theta$ since only the magnetic field intervenes in the expression of $C$.\\
 \textbf{We will therefore consider in the following of this section that we place ourselves in a coordinate system for which the electric field defined by the potential $A$ is zero}.\\
 This condition of nullity of the electric field is rather restrictive, it did not appear in the section \ref{s2.9}, but it should be noted that  experimental results associated to the notion of spin (experiment of Stern-Gerlach, measurement of the magnetic moment of electron, etc.) are studied theoretically with this restriction (see also remark \ref{r7} of section \ref{ss1.2}).

\begin{dfn} \label{d2.24}
 The real number $\varrho$ will be called \textbf{the gyromagnetic constant} of the domain "oscillating metric with spin in an electromagnetic potential ".
\end{dfn}
(see remark \ref{r2.10} on this topic).
\begin{rmq} \label{r2.9}
 The image of the orthogonal basis field ($L_1|_{S^3}$, $L_2|_{S^3}$, $L_3|_{S^3}$) by an isometry $\sigma$ of
the Euclidean space $\R^4$ (induced on $S^3(\rho)$) is still a field of orthogonal basis that can be written
($L^\sigma_1|_{S^3}$, $L^\sigma_2|_{S^3}$, $L^\sigma_3|_{S^3}$). \textbf{$\Upsilon$ may actually be chosen more
usually of the form $\Upsilon=A+C_\sigma$} ~~ where ~~ $C_\sigma:={\textstyle{\frac{\varrho}{2}}}\sum_{k=1}^3B^kL^\sigma_k|_{S^3}$. This is not of great importance since in the final results on the simple spin measures that we will obtain, $S^3(\rho)$ will no longer intervene precisely, however this fact will gain importance when studying phenomena of "quantum entanglement" in section \ref{+2.2}.
\end{rmq}


\subsection{The equations} \label{ss+3}
The important result of this subsection is stated in next theorem \ref{2.3}.
The important case is essentially that which concerns the electromagnetic field.

We begin by recalling and clarifying what are the potentials used and the assumptions that concern them:

In the context of the active potential without electromagnetism, the cell is \\ $\C=I\times \mathscr U\times
S^1(\delta)\times W$ and the potential metric satisfies: \\ $g_\mathcal P=g_O-2vX^\flat_1\otimes X^\flat_1$ ~~~ (prop.
\ref{p1.2}). \\ The function $a$ is written $a=\phi\beta$ ~ where ~ $\phi:I\times\mathscr U\times S^1(\delta)\rightarrow\R$ ~ and ~
$\beta\in E_W(\mu)$.

We consider \textbf{the following hypothesis $H_{2,N}$} (compare with hypothesis $H_{1,N}$ section \ref{s2.9}):
\begin{enumerate}
 \item $S_{g_\mathcal P}=S_{g_0}$.
 \item $v$ is a function defined on $\mathscr U$.
 \item $X_1$ is a vector field defined on $I\times W$, ~~ $D_{g_0}X_1=0$, ~ and ~ $X_1$ vanish on $E_W(\mu)$ ~~
(ie. ~ for any $ \beta\in E_W(\mu) ~~X_1(\beta)=0$)
\end{enumerate}
($v$ and $X_1$ can be considered defined on $\C$).
\bigskip

In the context of the specific electromagnetic potential introduced in the previous section \ref{ss2.15+}, the cell is $\C=\Theta\times
S^1(\delta)\times S^3(\rho)\times V$ and the potential metric satisfies: \\ $g_\mathcal P=g_O+sym(\Upsilon^\flat\otimes
X^\flat_2)$ ~~(proposition \ref{p1.3}) \\. The function $a$ is written $a=\phi\beta$ ~ where ~ $\phi:\Theta\times
S^1(\delta)\times S^3(\rho)\rightarrow\R$ ~ and ~ $\beta\in E_V(\nu)$.

We consider \textbf{the following hypothesis $H_{2,E}$} (compare with hypothesis $H_{1,E}$ section \ref{s2.9}):
\begin{enumerate}
 \item $S_{g_\mathcal P}=S_{g_0}$.
 \item $\Upsilon$ is a vector field defined on $\Theta\times S^3(\rho)$.
 \item $X_2$ is a vector field defined on $S^1(\delta)\times V$, ~~ $D_{g_0}X_2=0$, ~ and ~ $X_2$ vanish on
$E_V(\nu)$.
 \item The electric field defined by the potential $A$ is zero. 
\end{enumerate}
The remark that follows the assumption statements $H_{1,N}$ and $H_{1,E}$ is also valid for the assumptions $H_{2,N}$ and
$H_{2,E}$.
\begin{thme} \label{2.3}
 We consider a domain of type "elementary oscillating metric with spin in a potential". Then, for the
three cases of possible potentials, the canonical function $a_c$ satisfies the following equations:
\begin{enumerate}
 \item \textbf{In neutral potential}.
 \vspace{-3mm}
 \begin{eqnarray}\label{F47}
  \Box_\Theta a_c+ M^2a_c=0 
 \end{eqnarray}
 \vspace{-7mm}
 
where ~~ $\Box_\Theta=\frac{\partial^2}{(\partial t)^2}-\sum_{k=1}^3\frac{\partial^2}{(\partial x^k)^2}$ and $M$ is the
mass frequency.

\item \textbf{In a potential without electromagnetism under the assumption $H_{2,N}$}.
\vspace{-3mm}
\begin{eqnarray}\label{F48}
\Box_\Theta a_c+ M^2a_c-2v\frac{\partial^2a_c}{(\partial t)^2}=0
\end{eqnarray}
\vspace{-7mm}

where  $v$ is the potential function (def. \ref{def:11})
\item \textbf{In an electromagnetic potential under the assumption $H_{2,E}$}.
\vspace{-3mm}

\begin{eqnarray}\label{F49}
\hspace{-2cm}\sum_{j=0}^3\varepsilon_j(i\frac{\partial}{\partial x^j}+Q^+\Upsilon^j)^2a_c+M^2a_c+\varrho
Q^+\sum_{k=1}^3B^k\hat S_k(a_c)+{Q^+}^2(\textstyle{\frac{\varrho}{2}})^2\rho^2|B|^2a_c=0
\end{eqnarray}
\vspace{-5mm}

where ~~ $\varepsilon_j=g_{0jj}$ ~~ ie: ~~ $\varepsilon_0=-1$ ~~ and ~~ $\varepsilon_1=\varepsilon_2=\varepsilon_3=+1$, ~~
$|B|^2:=\sum_	{k=1}^3{B^k}^2$, ~ $\rho$ is the radius of the sphere $S^3(\rho)$, ~ $\varrho$ is the 
gyromagnetic constant, ~ the $\hat S_k$ are the canonical isomorphisms of the complexified eigenspaces of $S^3(\rho)$
~ (def. \ref{d2.22}) ~~ and ~~ $\hat S_k(a_k):\Theta\rightarrow E^{\CC}_p$ is defined for any $ x\in \Theta$ ~ by ~ $\hat
S_k(a_c)(x):=\hat S_k(a_c(x))$.
\end{enumerate}
\end{thme}
The proof of this theorem is detailed in Annex \ref{a3.7}.
\bigskip

Equations \ref{F47}, \ref{F48}, are identical to those "without spin" of theorem \ref{2.1}, but here \textbf{the
canonical function $a_c$ has values in $E^{\CC}_p$}. The additional terms of equation \ref{F49}
associated to those given by theorem \ref{2.1} express the "spin effect" into the electromagnetic potential.
\bigskip

When the electric charge is well defined, equations \ref{F47}, \ref{F48} and \ref{F49} are expressed, via the
definition \ref{d2.16}, in terms of state function $\varPsi$. These give, in approximation, the standard Schrödinger (or Pauli) equations. We write here only the result corresponding to the electromagnetic potential.

\begin{coro} \label{c2.1}
Under the hypotheses of theorem \ref{2.3}, when the electric charge is well defined (def. \ref{d2.14}) and in the
case of the specific electromagnetic potential, the state function $\varPsi$ satisfies the equation:
\begin{eqnarray}\label{F50}
 2iM\frac{\partial\varPsi}{\partial t}=\sum_{j=0}^3\varepsilon_j(i\frac{\partial}{\partial
x^j}+Q\Upsilon^j)^2\varPsi+\varrho Q\sum_{k=1}^3B^k\hat S_k(\varPsi)+{Q}^2(\textstyle{\frac{\varrho}{2}})^2\rho^2|B|^2\varPsi
\end{eqnarray}
\end{coro}
This result is obtained quickly by replacing $a_c$ expressed as a function of $\varPsi$ (def. \ref{d2.16}) in
 equation \ref{F49}.

When the $\varepsilon$-approximations (specified in the following paragraph) will be valid, the last term of
 equation \ref{F50} can be "neglected" (as well as the term for $j=0$ of the sum $\sum$) to give exactly
"Pauli's equation" of classical quantum physics.

\begin{rmq} \label{r2.10}
 Considering the definitions of "oscillating metric domains in a potential" for which
the metric $g$ is of the form $g=|a|^{4/n-2}g_\mathcal P$, it seems natural to consider that the function $a$ has
 the characteristics attributed to particles in classical physics (mass, electric charge, spin, etc.) and $g_\mathcal P$ has the characteristics of the "potential alone". However, this interpretation is faulty in the case of oscillating metrics with
spin in an electromagnetic field since the gyromagnetic constant $\varrho$ has been introduced into the specific potential
 $g_\mathcal P$ whereas it is classically rather associated with "particles". In fact, I do not think,
 with the point of view of physics presented here, that it is necessary to separate the two "objects" $a$ and $g_\mathcal P$, the
domain (metric $|a|^{4/n-2}g_\mathcal P$) should be considered as a "whole".

The constant $\varrho$ is here called "gyromagnetic constant" because it corresponds, in the equations
obtained, to the gyromagnetic constant of standard quantum physics ($\varrho\simeq2$ in the case of the electron). An "interpretation" would be to pose $\varrho=2$ in the definition of $g_\mathcal P$ and to consider that domains of type "elementary oscillating metric with spin $1/2$" correspond only to the notion of "electrons" (or some other fermions) in classical physics. The domains describing other particles would have a
more complex form for the function $a$ than $a=\phi\beta$ (which would correspond to composite particles). In this case, the equation \ref{F49} given by theorem \ref{2.3}, would only be an "approximation" but with $\varrho\neq2$. $\varrho$ could then be considered as a caracteristic of the function $a$ and not of the metric $g_\mathcal P$.
\end{rmq}


\subsection{The $\varepsilon$-approximations}
\subsubsection{The $\varepsilon$-approximations for the state function}
We write them identical to those of the definition \ref{d2.17}, but here $\varPsi$ is with values in
$E^{\CC}_{S^3(\rho)}(\gamma)$.
\vspace{-5mm}
\subsubsection{The $\varepsilon$-approximations for potentials}
We take again the conditions given in definition \ref{d2.18} to which we add a condition on the fields
$B^1$, $B^2$, $B^3$ because these are involved in definition of $g_\mathcal P$ as well as in
 equation \ref{F49} of theorem \ref{2.3}. Conditions on ($\rho B^k$) are chosen identical to those written
for $\Upsilon_0$ since these two terms occur in the same way in equation \ref{F49} of theorem
\ref{2.3}.

It is recalled that the sole interest of the $\varepsilon$-approximations (on state functions or on potentials) is to
give precise conditions which make "negligible" some terms of equations \ref{F47}, \ref{F48},
\ref{F49}. These last equations, after the deletion of these negligible terms, become identical to the equations obtained
in classical quantum physics.


\subsection{Probability of presence of a singularity in a domain of "elementary oscillating metric 
with spin in a potential "\label{ss2.13.6}}
We go fast here, taking into account the "spin", which has been exposed in section \ref{s2.12}. 
Modifications are minor, they consist only in taking into account the compact manifold
$S^3(\rho)$ and associated eigenspaces $E_{S^3(\rho)}(\gamma)$. The canonical function $a_c$ associated with the function
$a$ as well as the state function $\varPsi$ are now with values in $E_{S^1(\delta)}(\lambda)\otimes
E_{S^3(\rho)}(\gamma)$ identified with $E^{\CC}_{S^3(\rho)}(\gamma):=E^{\CC}_p$.

$E_{S^3(\rho)}(\gamma)$ is naturally equiped with the scalar product: \\
$$\langle\alpha_1,\alpha_2 \rangle_{L^2}:=\int_{S^3(\rho)}\alpha_1\alpha_2 dv_{S^3}$$ \\
The complexification $E^{\CC}_{S^3(\rho)}(\gamma)$ is then equiped with the Hermitian product:
\begin{eqnarray}\label{F50'}
 \langle\alpha,\alpha'\rangle:=\int_{S^3(\rho)}\alpha\overline{\alpha'} dv_{S^3}
\end{eqnarray}
where $\alpha=\alpha_1+i\alpha_2$ ~ and ~ $\alpha'=\alpha'_1+i\alpha'_2\in E^{\CC}_{S^3(\rho)}(\gamma)$ \\
In particular:
\begin{eqnarray}\label{F51'}
 |\alpha|^2:=\langle\alpha,\alpha\rangle=\int_{S^3(\rho)}(\alpha_1^2+\alpha_2^2)
dv_{S^3}=|\alpha_1|^2_{L^2}+|\alpha_2|^2_{L^2}
\end{eqnarray}

The fact that no law governs the singularities in a domain $\D\subset\M$ has been mathematically expressed in 
section \ref{s2.12} by \ref{F35} and \ref{F36}. Recall that, if $\HH_{2_{(t,u)}}\subset\HH_{1_{(t,u)}}$ are
spacelike submanifolds
of dimension $(n-2)$, then, when $\varsigma$ is an elementary singularity in $\HH_{1_{(t,u)}}$,
the probability that it is in $\HH_{2_{(t,u)}}$ is: \\

$p_{(t,u)}=\int_{\HH_{2_{(t,u)}}}dv_{g_{\HH_{1_{(t,u)}}}}/\int_{\HH_{1_{(t,u)}}}dv_{g_{\HH_{1_{(t,u)}}}}$
\bigskip

We consider a chart of the observation atlas for which the cell $(\C,g)$ is that of an oscillating metric
 with spin in a potential, for which: \\
$\C=\Theta\times S^1(\delta)\times S^3(\rho)\times V$ ~~ ($\Theta$ is of the
form $I\times \mathscr U\subset \R\times\R^3$) ~ and ~ $g=|a|^{4/n-2}g_{\mathcal P}$.

We fix  $(t,u)\in \R\times S^1(\delta)$. \\
Let $\HH_{2_{(t,u)}}\subset\HH_{1_{(t,u)}}\subset\C$ be of the form:

$\HH_{2_{(t,u)}}=\{t\}\times \omega\times \{u\}\times S^3(\rho)\times V$

$\HH_{1_{(t,u)}}=\{t\}\times \Omega\times \{u\}\times S^3(\rho)\times V$

with $\omega\subset\Omega\subset\R^3$

Potentials are assumed to satisfy the $\varepsilon$-approximations presented in the previous paragraph.

We obtain then (cf \ref{F39}):
\begin{eqnarray}\label{F52'}
 p_{(t,u)}\simeq\int_{\HH_{2_{(t,u)}}}(*)/\int_{\HH_{1_{(t,u)}}}(*)
\end{eqnarray}
where $(*)=a^2_{(t,u)}\sqrt{detg_0|_{\HH_{1_{(t,u)}}}}dv_{g_{0|_{\mathscr U\times S^3(\rho)\times V}}}$

Consider the case of an elementary oscillating metric \textbf{of order 2} for which, by definition:

$a=\phi\beta$

where $\phi:\Theta\times S^1(\delta)\times S^3(\rho)\rightarrow\R$ satisfies: $\phi=\phi_1\cos(Q^+u)+\phi_2\sin(Q^+u)$,

$\phi_1$ and $\phi_2$ are two real functions defined on $\Theta\times S^3(\rho)$,

for any $ x\in \Theta$ ~~~ $\phi_{1,x}(.)$ and $\phi_{2,x}(.)\in E_{S^3(\rho)}(\gamma)$,

$\beta\in E_V(\nu)$

(Recall that the canonical function $a_c$ is here defined by $a_c=\phi_1+i\phi_2$, then: for any $ x\in \Theta$
~~~ $a_{c,x}(.)\in E^{\CC}_{S^3(\rho)}(\gamma)$).

We obtain, using \ref{F52'} and "simplifying" by $\int_V\beta^2 dv_{g_{0|_V}}$: \\
$p_{(t,u)}\simeq\int_{\omega\times S^3(\rho)}(**)dx^idv_{g_{0|_{S^3}}}/\int_{\Omega\times
S^3(\rho)}(**)dx^idv_{g_{0|_{S^3}}}$ \\
where
$(**)=\cos^2(Q^+u){\phi_1}^2_{(t,x^i,s)}+\sin^2(Q^+u){\phi_2}^2_{(t,x^i,s)}+2\cos(Q^+u)\sin(Q^+u)\phi_1{\phi_2}_{(t,x^i,
s)}$ \\
and $(x^i):=(x^1,x^2,x^3)\in \mathscr U\subset \R^3$, ~~ $dx^i:=dx^1dx^2dx^3$. \\

Considering the "average" on $S^1(\delta)$ (cf \ref{F41}), we write: \\
$p_{(t)}\simeq\int_{\omega\times S^3(\rho)}(\phi^2_1+\phi^2_2)_{(t,x^i,s)}dx^idv_{g_{0|_{S^3}}}/\int_{\Omega\times
S^3(\rho)}(\phi^2_1+\phi^2_2)_{(t,x^i,s)}dx^idv_{g_{0|_{S^3}}}$. \\

Then: \\
$p_{(t)}\simeq\int_\omega|a_{c(t,x^i)}|^2dx^i/\int_\Omega|a_{c(t,x^i)}|^2dx^i$. \\
Here, $|.|$ is the norm in $E^{\CC}_{S^3(\rho)}(\gamma)$ defined in \ref{F51'}. \\

When the elementary oscillating metric has a well-defined electric charge (def. \ref{d2.14}) and taking into account the
definition of the state function $\psi:\Theta\rightarrow E^{\CC}_{S^3(\rho)}(\gamma)$, we get: \\

$p_{(t)}\simeq\int_\omega|\psi_{(t,x^i)}|^2dx^i/\int_\Omega|\psi_{(t,x^i)}|^2dx^i$. \\

We recover a standard result of classical quantum physics, but here the
denominator\\
$\int_\Omega|\psi_{(t,x^i)}|^2dx^i$ may depend on $t$ (see comments on this in section
\ref{s2.12}).


\subsection{Some examples}


\subsubsection{Example 1 - The spin $1/2$} \label{ss+1}
By definition, an "oscillating elementary metric with spin in a potential" has a spin
$1/2$ when the eigenspace $E_{S^3(\rho)}(\gamma)$ corresponds to the space $E_1$ of the classification given in
\ref{F45}. In this case, $\gamma=3\rho^{-2}$ and $dim E_1=4$.

A natural basis of this eigenspace is obtained by taking on $S^3(\rho)$ the restrictions
$(\alpha_1,\alpha_2,\alpha_3,\alpha_4)$ of the coordinate functions of $\R^4$ (homogeneous harmonic polynomials
of degree ~ 1)
which will be written here: $(x_1,x_2,x_3,x_4)$ (see annex \ref{sa3.7}).

Since the spin index is half-integer, we know that eigenfunctions $(\alpha_1,\alpha_2,\alpha_3, \alpha_4)$
do not come from the sphere $S^2$ by the Hopf fibration.

We write $M_1$, $M_2$, $M_3$ the matrix of $S_1$, $S_2$, $S_3$ (def. \ref{d2.21}) associated to the basis 
$(\alpha_1,\alpha_2,\alpha_3, \alpha_4)$. These are the same as matrix of differential operators $L_1$,
$L_2$, $L_3$
associated to the basis\\
$((x_1), (x_2), (x_3), (x_4))$ since $(L_k(x_l))|_{S^3}=(L_k|_{S^3})(x_l|_{S^3})$.

A very fast calculation, starting from the expressions of $L_1$, $L_2$, $L_3$, gives:

\begin{center}
$ M_1 = \begin{pmatrix}
  0&0&+1&0\\
  0&0&0&+1\\
  -1&0&0&0\\
  0&-1&0&0
\end{pmatrix}$ ~~~
 $ M_2 = \begin{pmatrix}
  0&0&0&+1\\
  0&0&-1&0\\
  0&+1&0&0\\
  -1&0&0&0
\end{pmatrix}$ ~~~
$M_3 = \begin{pmatrix}
 0&-1&0&0\\
 +1&0&0&0\\
 0&0&0&+1\\
 0&0&-1&0
\end{pmatrix}$
\end{center}

The three matrix $\hat M_k$ of the endomorphisms $\hat S_k$ of $E^{\CC}_1$ are therefore $\hat M_k={\textstyle{\frac{i}{2}}}M_k$ according to the
definition (\ref{d2.22}).

(We can verify that these three matrix $\hat M_k$ have the commutative properties of Pauli matrix which are special cases of the commutative properties of the observables of kinetic momentum (with factor $\hbar$ deleted)).
\bigskip

Under corollary assumptions \ref{c2.1}, the state function $\varPsi$ satisfies equation \ref{F50}. Denote
$\varPsi^1$, $\varPsi^2$, $\varPsi^3$, $\varPsi^4$ the four complex functions defined on $\Theta$, components
of $\varPsi$ in the basis $(\alpha_1,\alpha_2,\alpha_3,\alpha_4)$: $\varPsi=\sum_{j=1}^4\varPsi^j\alpha_j$.

Equation \ref{F50} is decomposed into four equations:
 \begin{eqnarray}\label{F51}
  2iM\frac{\partial\varPsi^1}{\partial t}=(\alpha)\varPsi^1-{\textstyle{\frac{\varrho}{2}}} Qi(-B^1\varPsi^3-B^2\varPsi^4+B^3\varPsi^2)
 \end{eqnarray}
\begin{eqnarray}\label{F52}
 2iM\frac{\partial\varPsi^2}{\partial t}=(\alpha)\varPsi^2-{\textstyle{\frac{\varrho}{2}}} Qi(-B^1\varPsi^4+B^2\varPsi^3-B^3\varPsi^1)
\end{eqnarray}
\begin{eqnarray}\label{F53}
 2iM\frac{\partial\varPsi^3}{\partial t}=(\alpha)\varPsi^3-{\textstyle{\frac{\varrho}{2}}} Qi(B^1\varPsi^1-B^2\varPsi^2-B^3\varPsi^4)
\end{eqnarray}
\begin{eqnarray}\label{F54}
 2iM\frac{\partial\varPsi^4}{\partial t}=(\alpha)\varPsi^4-{\textstyle{\frac{\varrho}{2}}} Qi(B^1\varPsi^2+B^2\varPsi^1+B^3\varPsi^3)
\end{eqnarray}
where $(\alpha)$ is the operator: ($\sum_{j=0}^3\varepsilon_j(i\frac{\partial}{\partial
x^j}+Q\Upsilon^j)^2+({\textstyle{\frac{\varrho}{2}}})^2Q^2\rho^2|B|^2$).

It is then interesting to define the four functions:

$\varphi^1:=\frac{1}{\sqrt{2}}(\varPsi^3-i\varPsi^4), ~~~\varphi^2:=\frac{1}{\sqrt{2}}(\varPsi^2+i\varPsi^1)$

$\varphi^3:=-\frac{1}{\sqrt{2}}(\varPsi^1+i\varPsi^2),~~~\varphi^4:	
=\frac{1}{\sqrt{2}}(-\varPsi^4+i\varPsi^3)$

The coefficient $\frac{1}{\sqrt{2}}$ is taken so
that:
$\sum_{j=1}^{4}|\varphi^j|^2=\sum_{j=1}^{4}|\varPsi^j|^2=|\varPsi|^2_{E^{\CC}_1}$.

These four complex functions are the components of the function $\varPsi$ in the basis of
$E^{\CC}_1$:
 $$\beta_1:=\frac{1}{\sqrt{2}}(\alpha_3+i\alpha_4),~~  \beta_2:=\frac{1}{\sqrt{2}}(\alpha_2-i\alpha_1)$$
 \begin{eqnarray}\label{F+6}
 \beta_3:=\frac{1}{\sqrt{2}}(i\alpha_2-\alpha_1)
\end{eqnarray}
In this basis, the three matrix $\hat M_k$ of the endomorphisms $\hat S_k$ of $E^{\CC}_1$ are:
 
\begin{center}
 $\hat M_1 = -1/2\begin{pmatrix}
  0&1&0&0\\
  1&0&0&0\\
  0&0&0&1\\
  0&0&1&0
\end{pmatrix}$ ~~~
 $\hat  M_2 = 1/2\begin{pmatrix}
  0&+i&0&0\\
  -i&0&0&0\\
  0&0&0&+i\\
  0&0&-i&0
\end{pmatrix}$ ~~~
$\hat M_3 =1/2 \begin{pmatrix}
 -1&0&0&0\\
 0&+1&0&0\\
 0&0&-1&0\\
 0&0&0&+1
\end{pmatrix}$
\end{center}

Equations \ref{F51}, \ref{F52}, \ref{F53}, \ref{F54} show that the functions ($\varphi^j$) satisfy:
 \begin{eqnarray}\label{F55}
  2iM\frac{\partial\varphi^1}{\partial t}=(\alpha)\varphi^1-{\textstyle{\frac{\varrho}{2}}} Q(B^1\varphi^2-iB^2\varphi^2+B^3\varphi^1)
 \end{eqnarray}
 \begin{eqnarray}\label{F56}
  2iM\frac{\partial\varphi^2}{\partial t}=(\alpha)\varphi^2-{\textstyle{\frac{\varrho}{2}}} Q(B^1\varphi^1+iB^2\varphi^1-B^3\varphi^2)
 \end{eqnarray}
 \begin{eqnarray}\label{F57}
  2iM\frac{\partial\varphi^3}{\partial t}=(\alpha)\varphi^3-{\textstyle{\frac{\varrho}{2}}} Q(B^1\varphi^4-iB^2\varphi^4+B^3\varphi^3)
 \end{eqnarray}
 \begin{eqnarray}\label{F58}
  2iM\frac{\partial\varphi^4}{\partial t}=(\alpha)\varphi^4-{\textstyle{\frac{\varrho}{2}}} Q(B^1\varphi^3+iB^2\varphi^3-B^3\varphi^4)
 \end{eqnarray}

 The pair of equations (\ref{F55},  \ref{F56}) is \textbf{identical} to the pair (\ref{F57},  \ref{F58}) when
$\varphi^1$ becomes $\varphi^3$ and $\varphi^2$ becomes $\varphi^4$.\\
The function $\varphi$ can be considered as the sum of two functions $\varphi'$ and $\varphi''$
corresponding to the superposition of two oscillating metrics, the first function having
$(\varphi^1,\varphi^2,0,0)$ for components in the basis $(\beta_1,\beta_2,\beta_3,\beta_4)$ and the second
$(0,0,\varphi^3,\varphi^4)$. In classical physics this corresponds to the sending of two particles of the same mass and electric charge, but possibly with  different initial conditions for the spin state.

\textbf{Each of these pairs of equations corresponds exactly to Pauli's equations of classical quantum physics
} when the last term as well as that in $j=0$ (in the sum $\sum$) of the operator ($\alpha$) are
deleted. These are actually "neglected" when the $\varepsilon$-approximations are valid.

We let the reader interpret the results of a "Stern-Gerlach" experiment,
for a spin
$1/2$, in terms of "space-time deformation" specified by the pairs of equations (\ref{F55},
~ \ref{F56}) and (\ref{F57}, ~ \ref{F58}), however this will be developed in section \ref{+2.1} which will allow us
to study the phenomena of "quantum entanglement" in section \ref{+2.2}.


\subsubsection{Example 2 - The spin $1$}
The spin $1$ is associated with the domain for which the eigenspace $E_{S^3(\rho)}(\gamma)$ corresponds to the classified space
$E_2$ (see \ref{F45}). In this case $\gamma=8\rho^{-2}$. We know that spaces $E_q$ whose index $q$ is even,
contain a subspace $E'_q$ that identifies with the eigenspace $E_{S^2(\rho/2)}(\gamma)$ of the sphere
$S^2(\rho/2)$ (see annex \ref{sa3.7}). We could make the general description of "spin $q$", when $q$ is
pair, considering the spaces $E_q$ (and not $E'_q$), but we will present, as an example, only the particular case
where the
state function associated with oscillating metric satisfies $\varPsi_x(.)\in E'^{\CC}_2$. Dimension of spaces
$E'_q$ is $q+1$, it corresponds to the considered dimension in standard quantum physics for "spin
integers ".

We choose as basis of $E'_2$, the restrictions to $S^3(\rho)$ of the three homogeneous harmonic polynomials $\tilde{P}_k$ of degree 2, which are written in the form $\tilde{P}_k=P_k\circ\pi$ where $\pi:\R^4\rightarrow\R^3$ is
the map that defines the Hopf fibration and the polynomials $P_k$ are the "coordinate functions" of $\R^3$:
 ($y_1$), ($y_2$), ($y_3$).
 
We have (see annex \ref{sa3.7}):
\begin{center}
\begin{tabular} {rcl}
$\tilde{P}_1(x_1,x_2,x_3,x_4)$ & $=$ & $x_1x_3+x_2x_4$\\[0.7em]
$\tilde{P}_2(x_1,x_2,x_3,x_4)$ & $=$ & $x_1x_4-x_2x_3$\\[0.7em]
$\tilde{P}_3(x_1,x_2,x_3,x_4)$ & $=$ & $\frac{1}{2}(x^2_3+x^2_4-x^2_1-x^2_2)$\\[0.7em]
\end{tabular}
\end{center}

We write $(\beta_1,\beta_2,\beta_3):=(\tilde{P}_1|_{S^3(\rho)},\tilde{P}_2|_{S^3(\rho)},\tilde{P}_3|_{S^3(\rho)})$ the
basis chosen in $E'_2$.

The matrix $M_1$, $M_2$, $M_3$ of $S_1$, $S_2$, $S_3$ (def. \ref{d2.21}) associated to the basis
$(\beta_1,\beta_2,\beta_3)$ are the same as the matrix of the differential operators $L_1$, $L_2$, $L_3$ associated to
the basis $(\tilde{P}_1$, $\tilde{P}_2$, $\tilde{P}_3)$.

A quick calculation gives:
\begin{center}
$ M_1 =2 \begin{pmatrix}
0&0&+1\\
0&0&0\\
-1&0&0
\end{pmatrix}$ ~~~
$ M_2 = 2\begin{pmatrix}
0&0&0\\
0&0&+1\\
0&-1&0
\end{pmatrix}$ ~~~
$ M_3 = 2\begin{pmatrix}
0&+1&0\\
-1&0&0\\
0&0&0
\end{pmatrix}$
\end{center}

The three matrix $\hat{M}_k$ of the endomorphisms $\hat{S}_k$ of ~ $E'^{\CC}_2$ are therefore ~ $\hat{M}_k={\textstyle{\frac{i}{2}}}M_k$.
\bigskip

Under the corollary assumptions \ref{c2.1}, the state function $\varPsi$ satisfies  equation \ref{F50}.\\
Write $\varPsi_1$, $\varPsi_2$, $\varPsi_3$ the three complex functions defined on $\Theta$, components of $\varPsi$
in the basis $(\beta_1,\beta_2,\beta_3)$.

Equation \ref{F50} is decomposed into three equations:
\begin{center}
\begin{tabular} {rcl}
$2iM\displaystyle\frac{\partial\varPsi^1}{\partial t}$ &$=$& $(\alpha)\varPsi^1+\varrho
Qi(-B^3\varPsi^2+B^1\varPsi^3)$\\[1em]
$2iM\displaystyle\frac{\partial\varPsi^2}{\partial t}$ &$=$& $(\alpha)\varPsi^2+\varrho
Qi(-B^3\varPsi^1+B^2\varPsi^3)$\\[1em]
$2iM\displaystyle\frac{\partial\varPsi^3}{\partial t}$ &$=$& $(\alpha)\varPsi^3+\varrho
Qi(-B^1\varPsi^1-B^2\varPsi^2)$\\[1em]
\end{tabular}
\end{center}
where $(\alpha)$ is the operator: $\sum_{j=0}^3\varepsilon_j(i\frac{\partial}{\partial
x^j}+Q\Upsilon^j)^2+({\textstyle{\frac{\varrho}{2}}})^2Q^2\rho^2|B|^2$

Here again, we let the reader interpret the results of a "Stern-Gerlach" experiment for
a spin $1$, in terms of "space-time deformation" specified by the three previous equations.

\subsection{The particular case of spin 1/2 for which $\varrho=2$}\label{ssAv}

In the specific case of the electromagnetic field the vector field $\Upsilon$ is written here:
\begin{eqnarray}
\Upsilon=\sum_{i=1}^3A^i\partial_i+\sum_{k=1}^3B^k{L_k}|_{S^3} 
\end{eqnarray}
Equation \ref{F49} is written:
\begin{eqnarray}\label{F49Av}
\sum_{j=0}^3\varepsilon_j(i\frac{\partial}{\partial x^j}+Q^+\Upsilon^j)^2a_c+M^2a_c+2
Q^+\sum_{k=1}^3B^k\hat S_k(a_c)=-{Q^+}^2(\textstyle{\frac{\varrho}{2}})^2\rho^2|B|^2a_c
\end{eqnarray}
The right side of this equation has coefficient $\rho^2$ and an evaluation of this radius $\rho$, introduced in section \ref{++2,3},
shows that, for values of $|B|$ corresponding to standard experiments, this term is negligible compared to the left side terms of the equation. It is therefore assumed in the following lines that solutions $a_c$ of the equation \ref{F49Av} are close to solutions of the following equation:
\begin{eqnarray}\label{F50Av}
\sum_{j=0}^3\varepsilon_j(i\frac{\partial}{\partial x^j}+Q^+\Upsilon^j)^2a_c+M^2a_c+2
Q^+\sum_{k=1}^3B^k\hat S_k(a_c)=0
\end{eqnarray}
\textbf{We will verify that this last equation can be written in factorized form and brings up the standard Dirac equation of the electron}. For this, we consider the three endomorphisms $\hat{S}_k$ of ~ $E^{\CC}_1$ whose matrix associated to the basis 
$(\beta_1, \beta_2, \beta_3, \beta_4)$, defined in \ref{F+6}
, are $\hat M_1, \hat M_2, \hat M_3$ specified in the previous subsection.\\
We define the two endomorphisms of $E^{\CC}_1$, denoted $\sigma$ and $\hat S_0$, whose matrix $\hat M_{\sigma}$ and 
$\hat M_0$ associated to the basis $(\beta_1, \beta_2, \beta_3, \beta_4)$ are:

\begin{center}
 $\hat M_{\sigma} = \begin{pmatrix}
  0&0&+1&0\\
  0&0&0&+1\\
  -1&0&0&0\\
  0&-1&0&0
\end{pmatrix}$ ~~~
 $\hat  M_0 = 1/2\begin{pmatrix}
  -1&0&0&0\\
  0&-1&0&0\\
  0&0&+1&0\\
  0&0&0&+1
\end{pmatrix}$ ~~~
\end{center}
We then quickly verify that the four endomorphisms of $E^{\CC}_1$, denoted $S^{\sigma}_0, S^{\sigma}_1, S^{\sigma}_2, S^{\sigma}_3$ defined by 
$S^{\sigma}_j:=2\hat S_j\circ\sigma$ satisfy:\\
$S^{\sigma}_j\circ S^{\sigma}_k=-S^{\sigma}_k\circ S^{\sigma}_j$ if $k\neq j,~~(S^{\sigma}_1)^2=(S^{\sigma}_2)^2=(S^{\sigma}_3)^2=-I$,~~$(S^{\sigma}_0)^2=I$.\\
Moreover: $S^{\sigma}_1\circ S^{\sigma}_2=2i\hat S_3,~~S^{\sigma}_1\circ S^{\sigma}_3=-2i\hat S_2,~~S^{\sigma}_2\circ S^{\sigma}_3=2i\hat S_1$.\\
The four matrices $M^{\sigma}_j$ of endomorphisms $S^{\sigma}_j$, associated to the basis $(\beta_1,\beta_2,\beta_3,\beta_4)$, are exactly the standard Dirac matrix.\\
We then obtain the following operators factorization (according to the definitions, the operators $S^{\sigma}_k$ commute with the $\frac{\partial}{\partial x^j}$):\\
$\sum_{j=0}^3\varepsilon_j(i\frac{\partial}{\partial x^j}+Q^+\Upsilon^j)^2+2
Q^+\sum_{k=1}^3B^k\hat S_k+M^2=$\\

\hspace{4cm} $(M+\sum_{j=0}^3(i\frac{\partial}{\partial x^j}+Q^+\Upsilon^j)\circ S^{\sigma}_j)\circ(M-\sum_{j=0}^3(i\frac{\partial}{\partial x^j}+Q^+\Upsilon^j)\circ S^{\sigma}_j)$\\
Equation \ref{F50Av} is then written:\\
$(M+\sum_{j=0}^3(i\frac{\partial}{\partial x^j}+Q^+\Upsilon^j)\circ S^{\sigma}_j)\circ(M-\sum_{j=0}^3(i\frac{\partial}{\partial x^j}+Q^+\Upsilon^j)\circ S^{\sigma}_j)a_c=0$\\
In classical quantum physics the Dirac equation for the electron is written:\\
$(M-\sum_{j=0}^3(i\frac{\partial}{\partial x^j}+Q^+\Upsilon^j)\circ M^{\sigma}_j)\phi=0$\\
where the four components of the state function $\phi$ correspond to the components of the canonical function $a_c$ of our theory associated to the basis $(\beta_1,\beta_2,\beta_3,\beta_4)$.\\
We deduce that the solutions of Dirac's equation are solutions of equation \ref{F50Av} and are therefore very close to  solutions of equation \ref{F49Av} taking into account the assumed conditions.\\
One of the interests of these factorizations is to make the magnetic field $B$ "disappear" in equations (\ref{F49Av})  (except the term in $|B|^2$ in general negligible) and (\ref{F50Av}). These equations, however, are not invariant by Poincaré transformations on $\Theta$ (see remark \ref{r7}).


\section{Zero mass oscillating metrics \label{s2.14}}
Notion of mass has been defined for elementary oscillating metrics associated with an 
eigenspace $E_{\lambda,\mu}$ (def. \ref{d2.3}). The mass frequency $M$ is (def. \ref{d2.8}) the non negative constant
 $\sqrt{S+\mu-\lambda}$ ~ where ~ $S:=\frac{n-2}{4(n-1)}S_{g_0}$, ~ $S_{g_0}$ is the constant scalar curvature of
the metric $g_0$.\\
\textbf{We therefore consider in this section, the case for which $M=0$,~~ ie: $S=\lambda-\mu$}.\\
Note that, if $\lambda=0$, this equality can be satisfied only if the scalar curvature is non positive. (As written in the remark \ref{++r7} in Chapter 1, the value of the scalar curvature of $(W,g_W)$ can achieve equality $S=\lambda-\mu$). \\

Electric charge frequency has been defined (def. \ref{d2.7}) by $Q^+:=\sqrt{\lambda}$, however, when
 $M=0$, the relative electric charge is not defined (def. \ref{d2.14}) and it is the same,
therefore, as the state function $\psi$.

Constant $Q^+$ (which may be zero) will remain an important feature of "zero mass oscillating metrics", but the term "electric charge frequency "will no longer be appropriate if we want to stay close to the language
of standard physics.

\textbf{Equations satisfied by the canonical function $a_c$ of a zero mass elementary oscillating metric of order 1 or 2  were given by theorems \ref{2.1} and \ref{2.3} when we set $M=0$}.

We rewrite here equations obtained for elementary oscillating metrics of order 1 with zero mass (for those
of order 2 see thm. \ref{2.3}), only the case of the electromagnetic potential is different because we take into account the
notion of "spin", moreover, $a_c$ has values in $E^{\CC}_{S^3(\rho)}(\gamma)$ and not in $\CC$).

\begin{enumerate}
 \item In a neutral potential.
 \begin{eqnarray}\label{F59}
 \Box_\Theta a_c=0  
 \end{eqnarray}
 where ~~ $\Box_\Theta=\frac{\partial^2}{(\partial t)^2}-\sum_{k=1}^3\frac{\partial^2}{(\partial x^k)^2}$
 \item In an active potential without electromagnetism.
\begin{eqnarray}\label{F60}
 \Box_\Theta a_c=2v\frac{\partial^2a_c}{(\partial t)^2}
\end{eqnarray}
\item In an electromagnetic potential.
  \begin{eqnarray}\label{F61}
 \sum_{j=0}^3\varepsilon_j(i\frac{\partial}{\partial x^j}+Q^+\Upsilon^j)^2a_c=0
\end{eqnarray}
where ~~ $\varepsilon_j=g_{0jj}$
\end{enumerate}
In the case of neutral potential, the canonical function $a_c$ is therefore a solution of the classical wave equation.
In the context of our theory, standard electromagnetic waves can be seen as particular zero-mass oscillating metrics, however, it is possible that some electromagnetic wave fall into a very different category of oscillating metrics as is specified in the subsection (\ref{4-ss-3}).\\
Equations \ref{F60} and \ref{F61} describe the influence of a potential on the zero mass oscillating metric.
That such an influence exists is not surprising and already appears in classical general relativity where we know that,
for example, an electromagnetic field distorts space-time via its energy-momentum tensor.
\subsection{Zero-mass oscillating metric close to a potential}\label{4-ss-3}
If one considers a domain with "constant scalar curvature, conform to a potential" for which the tensor $g={\left|a\right|}^{4/n-2}g_0$ has a function $a$ "close" to 1, ie, of the form $a=1+b$ where $b:\Theta \times W\rightarrow \R$ satisfies $b<<1$, equation (\ref{F0'}) is written:
\begin{eqnarray}
\dfrac{4(n-1)}{n-2}\Box_{g_{\mathcal P}}b+S_{g_{\mathcal P}}(1+b)=S_{g_\mathcal P}{(1+b)}^{n+2/n-2}\simeq S_{g_{\mathcal P}}({1+\frac{n+2}{n-2}b)}
\end{eqnarray}

In other words:
\begin{eqnarray}\label{4-F0}
\Box_ {g_{\mathcal P}}b-\frac{1}{n-1}S_{g_{\mathcal P}}b\simeq 0  
\end{eqnarray}
This linear approximation, satisfied only when $a=1+b$ with $b<<1$, is fundamentally different from that given by
\ref{F1}, in particular the sign in front of the scalar curvature is here negative (here the expression "oscillating metric" no longer matches definition \ref{d2.3}).\\
Assuming, by analogy with definition \ref{d2.3}, that for any $ x\in\Theta~~ a_x(.)\in E_{\lambda,\mu}$, that $g_{\mathcal P}=g_0$ (neutral potential) and $\mu-\lambda-\frac{1}{n-1}S_{g_0}=0$, equation \ref{4-F0} is written:\\
$\Box_\Theta a=\Box_\Theta b=0$\\
Then, the caracteristic function $a_c:\Theta\rightarrow E_{\lambda,\mu}$ (defined for $x\in\Theta$ by $a_c(x)=a_x(.)\in E_{\lambda,\mu}$) satisfies standard wave equation. (In this case the constant $\mu-\lambda-\frac{1}{n-1}S_{g_0}$ corresponds to the notion of mass).\\
An electromagnetic wave can be defined by such a method (by specifying $E_{\lambda,\mu}$) but then, the definition is very different from one given above. In one case it is a small disturbance of the metric $g_0$ (or $g_{\mathcal{P}}$) (described by the function $b$) on the considered whole space-time domain, in the other case it is a strong disturbance of the metric $g_0$ (or $g_{\mathcal{P}}$) but only on localised domains of space-time (cf (\ref{4-1}), "density" of these domain in the considered space-time domain can be very low. It is important to note that the concept of photon only becomes important,
 in our theory, in the context of "strong disturbances" because the localizations by the singularities (cf \ref{d2.19}) use the equiprobability principle specified in \ref{s2.12}. When the metric $g_0$ is weakly disturbed by the function $b$, the singularities remain quasi-equiprobably distributed relative to the metric $g_0$ and it is not a count of these which is interesting in term of measurement (cf \ref{s2.16}), but the physical effects (which remain to be clarified) associated to the function $b$ which will be important.\\
 It is conceivable that, with the point of view of our theory, what is classically considered as "electromagnetic waves" separates into these two categories which eventually "superset" (the "low frequencies" wave would be classified in the framework of the weak disturbances by the function $b$, and visible light in the context of strong localized disturbances). The zero-mass oscillating metrics just described may also correspond to other notions than those classified in electromagnetic wave according to the caracteristics of the spaces $E_{\lambda,\mu}$, we can cite, for example, the gravitational waves classified in the framework of weak disturbances of $g_0$ (those classified in the framework of the strong disturbances could introduce, by the singularities, the notion of graviton).
In this section, zero-mass oscillating metrics have been presented from a "quantum" point of view (with
our point of view). However, we can also be interested in these notions in terms of "lightlike fluid",
but then, it's not the same kind of approximations that we use and this corresponds to the first chapter of
this paper in which we have excluded the "quantum phenomena" since we made, from the beginning, an average of the
metric $g$, in particular on the circles $S^1_x(\delta)$ (see Lemma \ref{l1}). The study, as a lightlike fluid, can then be applied to much larger domain of space-time since, in this case, definitions of the important objects are not associated to a reference metric $g_0$, whereas on the other hand, the study of zero-mass oscillating metrics that we have just presented is limited to restricted domains of space-time.
\section{Influence of singularities on a neutral potential metric.\\
The cosmological constant}
A neutral potential metric has been defined in \ref{ss1.1.1} on a cell $\C:=\Theta\times S^1(\delta)\times W$ by $g_0=g_\Theta\times(-g_{S^1(\delta)})\times g_W$ where $g_\Theta$ is the Minkovski metric.\\
Notion of "mass frequency" M has been introduced by definition \ref{d2.8} for elementary oscillating metric associated with $E_{\lambda,\mu}$ (cf \ref{s2.5}) because it is within this framework that comparison with classical theories has been made (see theorems \ref{2.1} and \ref{2.2}). We defined: $M=(S+\mu-\lambda)^{1/2}$ where $S:=\frac{n-2}{4(n-1)}S_{g_0}$,~~ $S_{g_0}$ being the scalar curvature of the metric $g_0$. In fact, notion of mass frequency (as well as the other constants, such as electric charge frequency) can be clearly defined identically for all the functions $a:\C\rightarrow \R$ which satisfy: for any $ x\in\Theta~~~ a_x(.)\in E_{\lambda,\mu}$ and without any other condition. These functions may be called {\it elementary functions associated with $E_{\lambda,\mu}$}.\\
To make the link with classical theories and to return to standard units, one defined the {\it mass} of an elementary oscillating metric associated with $E_{\lambda,\mu}$ (in kg) by (definition \ref{d2.10}) $m:=\hbar c^{-1}M$. It is important to note, when one makes the link with standard theories of particles (by theorems \ref{2.1} and \ref{2.2} for example) that this notion correponds to the mass of \textbf{one} particle whose characteristic constants correspond to those of the oscillating metric.\\
When measuring positions, the notion of "particle" corresponds, within the framework of our theory, to the notion of elementary singularity (definition \ref{d2.20}). More precisely, to a particle at time "t" in a domain $\omega\subset\R^3$ corresponds an elementary singularity found at the double time $(t,u)$ in the spacelike domain $\{t\}\times \omega\times\{u\}\times W$ (of dimension $n-2$). To $k_t$ particles which are at time "t" in $\omega$ we will associate $k_t$ elementary singularities found in $\omega\times W$ where $k_t:=\frac{1}{2\pi\delta}\int_{S^1(\delta)}k(t,u)du$ and $k(t,u)$ denotes the number of elementary singularities found in $\{t\}\times\omega\times\{u\}\times W$.\\
Given the preceding considerations, it is natural to define \textbf{the mass at time "t" of an elementary function $a$ associated with $E_{\lambda,\mu}$ in a domain $\omega\subset\R^3$} by: 
\begin{eqnarray}\label{4-F2}
m_{t,\omega}=\hbar c^{-1}k_t(S+\mu-\lambda)^{1/2}
\end {eqnarray}
where $k_t$ is the number of elementary singularities in $\omega\times W$ at time "t". This definition applies in particular to elementary oscillating metric associated with $E_{\lambda,\mu}$ (via the function $a$) as well as to  neutral potential domain for which the associated function $a$ is the constant $1$. It is the latter case that interests us here. The eigenspace $E_{\lambda,\mu}$ is reduced to the space of constant functions and the eigenvalues $\lambda$ and $\mu$ are zero (zero electric charge frequency, zero spin, etc.).\\
\textbf{The mass at time "t" of a domain $\omega$ for a neutral potential of metric $g_0$} is, therefore \ref{4-F2}:
\begin {eqnarray} \label{4-F3}
m_{t,\omega,g_0}=\hbar c^{-1}k_tS^{1/2}
\end {eqnarray}
In the language of classical physics we can express above considerations by saying that an "empty domain" $\omega\subset \R^3$ (defined in our theory by a metric of neutral potential $g_0$ on $\R\times\omega\times S^1\times W$) contains particles each having a mass $m= \hbar c^{-1}S^{1/2}$. These particles may be called "vaccum particles" and, since $\lambda$ and $\mu$ are zero, these particles can only be detected by means of their masses. An interpretation on this consists in that the "mass of vaccum" (dark energy?) intervenes in the expansion of space-time and is, in standard theory of general relativity, expressed by the cosmological constant in the Einstein equation. From the point of view of our theory, we can consider that these elementary singularities intervene by "averaging" on the value of the energy density function $\mu$ that we defined in the first chapter in \ref{ss1.1}(3.) when we neglected the "quantum effects" by other types of averaging (cf section \ref{s1.2}).\\
\textbf{Some orders of magnitude}\\
Denote $\Omega_t$ the spacelike and dimension 3 domain at time "t" on which the metric is "quasi"-Euclidean (it comes from $g_0$) and $D_t$ the density of "vaccum particles" in $\Omega_t$.\\
We have, according to (\ref{4-F3}):\\
$D_t=k_t(vol\Omega_t)^{-1}=m_{t,\Omega_t,g_0}\hbar^{-1}cS^{-1/2}(vol\Omega_t)^{-1}$\\
If we assume that, in the considered space-time domain, $D_t$ does not depend on time (then denoted by $D$) and that the cosmological constant is estimated (according to the standard theory) at:\\
$\Lambda\simeq10^{-9}J/m^3~~~ie~~~\Lambda\simeq10^{-26}kg/m^3$,\\
if, as an example, we suppose that the scalar curvature $S_{g_0}$ is given only by the sphere $S^3(\rho)$, then $S=\frac{n-2}{4(n-1)}S_{g_0}$ is of order $1/\rho^2$.\\
If $\rho$ is assumed to be of the order of $10^{-17}m$ (cf end of section \ref{++2,3}) then $D=\Lambda\hbar^{-1}cS^{-1/2}$ is of order of a few "particles" per $m^3$.\\
Of course, we should not give much credit to these estimates, it seems much more likely that in the "universe expansion" intervene "expansions or retractions" of the compact Riemannian manifold $(W,g_W)$ as it is specified at the end of annex \ref{a3.10+}. In particular, large variations of the radius $\rho$ of the sphere $S^3(\rho)$ render the previous calculations inappropriate ($\Lambda$ is no longer a constant).\\

In this section, we focused on singularities in a neutral potential domain for which the sole non zero characteristic constant is given by the sole scalar curvature $S_{g_0}$ of the domain, which then determines the notion of mass. In the following subsection, we focus on elementary oscillating metrics for which only the notion of mass intervenes.


\section [The Higgs field?] {A very basic oscillating metric \label{s2.15}}
Consider an "elementary oscillating metric domain in a neutral potential associated with $E_{\lambda,\mu}$" (cf definition
\ref{d2.3}) \textbf{for which $\lambda=0$ and $\mu=0$}. Eigenfunctions involved in the spectral decomposition
 of the function $a$ are all constant. Electric charge frequency $Q^+$ is zero, 
"spin" is null, etc. The canonical function $a_c$ is therefore reduced to a real function defined on $\Theta\subset\R^4$
(the equivalent of a real scalar field in Q.F.T).

In the framework of the linear approximation, fundamental equation \ref{F1} is written:

$\Box_{g_\mathcal P}a_c+Sa_c=0$ ~ where ~ $S=\frac{n-2}{4(n-1)}S_{g_0}$ ~ and ~ $\Box_\Theta=\frac{\partial^2}{(\partial
t)^2}-\sum_{k=1}^3\frac{\partial^2}{(\partial x^k)^2}$.

Assuming that $S_{g_0}$ is positive, the mass frequency is none other than:
\bigskip

$M=(\frac{n-2}{4(n-1)}S_{g_0})^{1/2}$.
\bigskip

This oscillating metric has the \textbf{simplest} form that can be written in a neutral potential
 and its mass depends only on the scalar curvature of the neutral potential. It should be noted that this oscillating metric is rather to be classified within the framework of the strong perturbations of $g_0$, very localized in the considered space-time.

We have already noticed, when defining the mass (definition \ref{d2.10}), that a 
 positive scalar curvature of the neutral potential "gives" mass to elementary oscillating metrics
associated with $E_{\lambda,\mu}$ for which $\mu-\lambda$ is negative (when $S+\mu-\lambda$ is
positive).

It seems natural to link this elementary oscillating  metric to the Higgs field presented in
Q.F.T. We can consider that singularities in such a domain correspond to Higgs bosons, however it is
likely that these singularities associated with this precise domain are "undetectable". Experiments that allow
to actually detect presence of Higgs bosons certainly are no parts of the oscillating metric that we have just presented,
but describe complex interaction phenomena for which the notion of
"mass" does not have the meaning given here and the notion of "life time" becomes important.
These last points are briefly discussed in section \ref{s2.17}.


\section [Magitudes and measures] {Magnitudes and measures} \label{s2.16}
This section is conceptually very important, it specifies the notion of "magnitude measurement" which turns out to be fundamental for the description of quantum phenomena since, for them, "measuring instruments" are necessarily an integral part of the physical system studied (read about it in Annex \ref{a3.10+}). \\
In standard quantum physics, the notion of measurement is introduced directly into the axiomatic system (and the imprecision of this notion ends up raising some conceptual problems). \\
In the theory presented in this paper, the point of view is very different since the notion of "measurement of a magnitude" is introduced by giving only a precise definition of what is an "idealized measuring instrument" for this quantity.  Considered quantities will be only characteristics of the oscillating metrics (and position measurements for singularities) and will not be in any case quantities associated with "macroscopic objects" (Schrödinger's cat for example). \\

It is important to note that, throughout this paper, the notion of\begin{it}
momentum \end{it} did not occur, whereas this notion is fundamental in all standard physical theories. In the chapter on quantum phenomena the sole notion of \begin{it} measure \end{it} that we have
considered up to this point is
that of \begin{it} position \end{it} and this has been defined only for \begin{it} singularities \end{it} (cf
\ref{s2.11}). This is
sufficient to describe qualitatively and quantitatively all the standard experiments of quantum physics
(diffraction, Young's slits, potential deviations, Stern-Gerlach experiment, quantum entanglement, etc.).
We now address more complex phenomena (currently treated by the Q.F.T). There is no reason, with our
point of view on physics, to introduce notions of \begin{it} momentum \end{it},
\begin{it} energy \end{it}, etc. However, to be able to compare our theoretical results
 with those of standard physical theories, we will introduce in this section the \begin{it} magnitudes \end{it} notions such as
\begin{it} momentum \end{it}, \begin{it} energy \end{it}, etc. (which will be associated with metrics conform to a
potential and not
singularities), as well as the concept \begin{it} of measuring instrument of these quantities \end{it}. In particular, we will recover in this study, inequalities similar to those of the \begin{it} uncertainty principles \end{it} of standard quantum physics.

It should be noted that the concept of \begin{it} velocity \end{it} (which allows to introduce the 
\begin{it} momentum \end{it} notion) was only defined
in the particular case of homogeneous oscillating metrics (see \ref{s2.6}), this will serve as a basis for 
what we will introduce.

In all that concerns this section the cell is of the form $\C=\Theta\times S^1(\delta)\times W$ where
$\Theta=I\times \Omega\subset\R\times\R^3$ and the reference metric is $g_0=g_\Theta\times
(-g_{S^1(\delta)})\times g_W$.


\subsection{An example} \label{ssn2.1}
Before giving definitions associated to the measurement notion in a general context, we begin by considering
the particular case of an oscillating metric domain of order 1 in a neutral potential which will serve as a
reference example in this section.

The metric $g$ is written: $g=|a|^{4/n-2}g_0$ and the function $a$ is of the particular form:

$a=\beta\sum_{j=1}^p\varphi_j$ ~~ where ~~, for any $ j$ from 1 to p, ~~ $\varphi_j:\Theta\times S^1(\delta)\rightarrow\R$
satisfies:

$\varphi_j=C_j \cos(M_jt+Qu-\sum_{k=1}^3\lambda_{jk}x^k)+C'_j\sin(M_jt+Qu-\sum_{k=1}^3\lambda_{jk}x^k)$.

$C_j$ and $C'_j$ are here constants, ~~ $\beta\in E_W(\mu)$, ~~ $Q$ is the relative electric charge.

According to definition \ref{d2.12} each function $a_j:=\beta\varphi_j$ corresponds to an \textbf{homogeneous} elementary oscillating metric of
 order 1 and the velocity vector of this oscillating metric is given by:

$\overrightarrow{v_j}=(1/M_j)(\lambda_{j1},\lambda_{j2},\lambda_{j3})$.

The corresponding momentum is then naturally defined by:

$\overrightarrow{\lambda_j}:=M_j\overrightarrow{v_j}=(\lambda_{j1},\lambda_{j2},\lambda_{j3})$.

$M_j$ is the "relativistic" mass frequency, it is associated to the mass "at rest" $M_0$ (which we consider
do not depend
of j in this example) by (see \ref{s2.6}):

$M_0=(1-|\overrightarrow{v_j}|^2)^{1/2}M_j$.

It is easy to verify that $\Box_{g_0}a+Sa=0$ since $M_0^2=S+\mu-Q^2$. The oscillating metric given by
$g=|a|^{4/n-2}g_0$ is therefore  \begin{it} elementary \end{it} of order 1 in a neutral potential.

Note (see \ref{s2.8}) that the associated canonical function satisfies:

$a_c=\sum_{j=1}^p(C_j+iC'_j)e^{-i(M_jt-\sum_{k=1}^3\lambda_{jk}x^k)}$ ~~ when $Q>0$, and is equal to the conjugate when
$Q<0$.

The state function is:

$\varPsi=\sum_{j=1}^p(C_j+iC'_j)e^{-i((M_j-M_0)t-\sum_{k=1}^3\lambda_{jk}x^k)}$.

In the "classical" language, the study of the field just described (when one
considers that the state function $\Psi$ is that of standard quantum physics) is expressed from the
following way: \\
Let a particle of mass $m$ and electric charge $q$ in a domain $\Omega\subset\R^3$ such that the
state function $\varPsi$ is the one just mentioned. So, according to the classical quantum physics principles
 (delicate in this case because $\varPsi$ is not normalizable), one can conclude that, during a momentum measurement of the particle at time $t_0$, the probability of getting
$\overrightarrow{q_j}=(\lambda_{j1},\lambda_{j2},\lambda_{j3})$ is $({C_j}^2+{C'_j}^2)/\sum_{j=1}^p({C_j}^2+{C'_j}^2)$.

For us, the process will be fundamentally different and we are going to use a precise definition of
concept of \begin{it} magnitudes \end{it} then that of \begin{it} measuring instruments of these magnitudes \end{it} which,
used especially on
the previous example, will give similar results to those of standard quantum physics. Although the
definitions start from the same principles, we will present in two different sub-sections the 
\begin{it} magnitudes \end{it} notions associated to apparent space $\Theta$ and those associated to compact manifolds of the cell
 $\C$. The
non-compactness of apparent space will require a more elaborate presentation than that concerning compact manifolds.


\subsection{Magnitudes associated to apparent space $I\times\Omega\subset\R\times\R^3$ and their measurements}
As the actual functions defined on $\R^p$ that we are going to consider will not necessarily belong to
$L^2(\R^p)$, we introduce the following definition where $B_L$ is the cube in $\R^p$ whose coordinates of
vertex are $(\pm L/2,...,\pm L/2)$.

\begin{dfn} \label{d2.25}
 A family $(f_i)_{i\in \mathcal{G}}$ of continuous real functions on $\R^p$ is \textbf{B-orthonormal}
if \begin{enumerate}
\item for any $i\in\mathcal{G}$ ~~ $\frac{2}{L^p}\int_{B_L}f_i^2\rightarrow 1$ when $L$ goes to infinity.
\item for any $i,j\in\mathcal{G}$ such that $i\neq j$ ~~ $\frac{2}{L^p}\int_{B_L}f_if_j\rightarrow 0$ when $L$ goes to
infinity.

\end{enumerate}

\end{dfn}

$\frac{2}{L^p}\int_{B_L}f_if_j$ will sometimes be written $<f_i,f_j>_{B_L}$. The choice of factor "2" will allow
simplify some coefficients in the following calculations, essentially when the functions $f_i$ will be
trigonometric functions.

The \begin{it} magnitude \end{it} notion that we introduce now will be of interest only associated to the concept of
\begin{it} measure instrument \end{it} that will be defined next.
\begin{dfn} \label{d2.26}
 Consider the apparent space $I\times\Omega\subset\R\times\R^3$. \\
\textbf{A magnitude defined on $\Omega\subset\R^3$ (resp. $I\subset\R$)} is a family $(h_i)_{i\in\mathcal{G}}$ of
 $C^\infty$ functions defined on $\R^3$ (resp. $\R$) with values in $\R^m$ such that the "union" of the $m$ components families
 $(h_{1,i})_{i\in\mathcal{G}},...,(h_{m,i})_{i\in\mathcal{G}}$ form a B-orthonormal family
(def. \ref{d2.25}).
\end{dfn}
\begin{rmq}
 In the examples that follow $m$ will be 2 and in this case we could say that the functions $h_i$ are with
values in $\CC$ rather than in $\R^2$. Introduction of $\CC$ makes it possible to simplify calculations
in the particular cases of the domains for which the canonical function (cf \ref{s2.8}) is defined. However,
the following definitions apply in a more general context.
\end{rmq}

The two basic examples are as follows.


\subsubsection{Example 1}
The \textbf{ momentum magnitude} is defined by the family
$(h_{\overrightarrow{q}})_{\overrightarrow{q}\in\mathcal{G}}$

where $\mathcal{G}=\R^{3^*}$ and for any $\overrightarrow{q}=(q_1,q_2,q_3)\in\mathcal{G}$ ~~
$h_{\overrightarrow{q}}:\R^3\rightarrow\R^2$ ~ satisfies: \\ for any $ (x_1,x_2,x_3)\in\R^3$ ~~
$h_{\overrightarrow{q}}(x_1,x_2,x_3)=(\cos(\sum_{k=1}^3q_kx^k),\sin(\sum_{k=1}^3q_kx^k))$.
\subsubsection{Example 2}
The \textbf{energy magnitude} is defined by the family
$(h_e)_{e\in\mathcal{G}}$ ~~ where $\mathcal{G}=\R^*$ \\ and, for any $ e\in\mathcal{G}$ ~~ $h_e:\R\rightarrow\R^2$
satisfies: for any $ t\in\R$ ~~ $h_e(t)=(\cos(et),\sin(et))$.
\begin{rmq}
 Usually, the denomination of the considered \begin{it} magnitude \end{it} is rather given to the elements of the
set $\mathcal{G}$: we speak about \begin{it} the momentum \end{it} $\overrightarrow{q}=(q_1,q_2,q_3)$ or the 
\begin{it} energy \end{it} $e$.

\end{rmq}
Note that the set $\R^3$ which corresponds to the \begin{it} momentum \end{it} magnitude represents the tangent space to
$\Omega$ in
each of its points (which canonically identifies with $\R^3$ by the choice of the cell).
\bigskip

To clarify the connection between physical reality and the notion of \begin{it} measure of a magnitude \end{it} associated with a "metric conform to a potential type domain", we present in the form of
\begin{it} definitions \end{it} the concept of
\begin{it}magnitude measuring instrument\end{it}. These definitions only concerns the magnitudes
that we have
just introduced associated with the apparent space, but they will adapt without difficulty to the magnitudes associated with the
compact manifolds that will be specified later.

We consider a  "metric conform to a potential" in a cell $(\C,g)$ where
$\C=I\times\Omega\times S^1(\delta)\times W$ ~ and ~ $g=|a|^{4/n-2}g_\mathcal{P}$, $g_\mathcal{P}$ being a potential metric  
(this domain can be in particular a domain "with constant scalar curvature" (def. \ref{d2.1}) or of 
type "oscillating metric" (def. \ref{d2.2})).

Let a magnitude $(h_i)_{i\in\mathcal{G}}$ (def. \ref{d2.26}) associated with this domain. The three important characteristics of a measuring instrument of this magnitude are as follows.
\begin{enumerate}
 \item Its \textbf{spectrum} which, by definition, is a finite subfamily of the family $(h_i)_{i\in\mathcal{G}}$,
that is, a family $(h_i)_{i\in Sp}$ where $Sp$ is a finite part of $\mathcal{G}$.
\item Its \textbf{measurement domain} $B_L$ which is a cube of $\R^3$ whose vertex coordinates are\\
$(\pm L/2,\pm L/2,\pm L/2)$ (with Euclidean volume $L^3$). For simplicity, and without
restrict the generality, it will be assumed that $0\in\Omega$ and that for $L$ sufficiently small, $B_L\subset\Omega$,
which can always be obtained by a "translation" of the coordinates $(x^1,x^2,x^3)$ in the cell.
\item Its \textbf{measurement duration T} which will correspond to the time measured by the observer associated with the  cell.
\end{enumerate}
\begin{dfn} \label{d2.27}
 \textbf{An instrument for measuring the magnitude $(h_i)_{i\in\mathcal{G}}$ defined on $\Omega\subset\R^3$
(\ref{d2.26})} whose characteristics are $Sp,B_L$ and $T$, is a physical system that has the following properties:
\begin{enumerate}
 \item From a time $t_0\in I$, it transforms the "subdomain" $(\C_{B_L}, g)$ of $(\C,g)$ into a
"domain"
$(\C',g')$ for which:
\begin{enumerate}
 \item $\C_{B_L}=I\times B_L\times S^1(\delta)\times W$
 \item $\C'=]t_0 , t_0+T[\times(\bigcup_{i\in Sp}B_i)\times S^1(\delta)\times W$
 
 where $]t_0 , t_0+T[\subset I$ and $B_i$ are 2 to 2 disjointed cubes of $\R^3$ (or at least which have zero measure for their 2 to 2 intersections), each isometric to $B_L$.
\item for any $ i\in Sp$ ~~ $g'_{B_i}=|a_i|^{4/n-2}g_0$ ~ where ~ $a_i:]t_0 , t_0+T[\times B_i\times S^1(\delta)\times
W\rightarrow\R$ is defined by:
$$a_i=\sum_{l=1}^m <a,h_{il}>_{B_L}(h_{il}\circ\sigma_i)$$

Here $\sigma_i$ is the isometry between $B_i$ and $B_L$ ($h_{il}\circ\sigma_i$ can be considered defined on
$\C_i:=]t_0 , t_0+T[\times B_i\times S^1(\delta)\times W$). It is recalled that
$<f_1,f_2>_{B_L}:=\frac{2}{L^3}\int_{B_L}f_1f_2$.
\end{enumerate}
\item For each $i\in Sp$, the measuring instrument estimates the average number of elementary singularities found
in $ B_i$ during the time $(t,u)\in ]t_0 , t_0+T[\times S^1(\delta)$. (When $n_i(t,u)$ is the number of
elementary singularities that are, at time $(t,u)$, in $\HH_i:=\{t\}\times B_i\times\{u\}\times W$, the
average number is defined by $\bar{n}_i:=\frac{1}{2\pi\delta T}\int_{S^1(\delta)}\int_{t_0}^{t_0+T}n_i(t,u)dudt)$).
\end{enumerate}

\end{dfn}
The properties required of such a measuring instrument can be succinctly summarized by saying that it "separates
in space" the components that interest us of the metric on which we make the measurements, then
analyzes the created domain (separated in space) by counting the elementary singularities found there. An example
of such measuring instrument is given by a prism which breaks down the light consisting of p separate momentums of
same direction by spreading it over p domains (rainbow). In this case the considered oscillating metrics are of
zero mass (see section \ref{s2.14}).
\bigskip

When the magnitude $(h_i)_{i\in\mathcal{G}}$ is defined on $I\in\R$, the definition of a measuring instrument of
this magnitude is analogous to the previous definition \ref{d2.27}, only the condition $1.(c)$ is different.

\begin{dfn} \label{d2.28}
\textbf{An instrument for measuring the magnitude $(h_i)_{i\in\mathcal{G}}$ defined on $I\subset\R$ (def. \ref{d2.26})} whose
the characteristics are $Sp,B_L$ and $T$, is a physical system that has the following properties:
\begin{enumerate}
 \item From a time $t_0\in I$, it transforms the "subdomain" $(\C_{B_L}, g)$ from $(\C,g)$ into a
"domain"
$(\C',g')$ for which:
\begin{enumerate}
\item Same as Definition \ref{d2.27}
\item Same as Definition \ref{d2.27}
\item for any $ i\in Sp$ ~~ $g'_{B_i}=|a_i|^{4/n-2}g_0$ ~ where ~ $a_i:\C_i:=]t_0 , t_0+T[\times B_i\times S^1(\delta)\times
W\rightarrow\R$ is defined by:
$$a_i=(\sum_{l=1}^m(\frac{2}{T}\int_{t_0}^{t_0+T}a_{\C_L}h_{il}dt)h_{il})\circ\sigma_i$$
Here $\C_L:=]t_0 , t_0+T[\times B_L\times S^1(\delta)\times W$ and it is recalled that, in the context of this definition,
$h_i$ is a function defined on $\R$. $\sigma_i$ is an isometry between $B_i$ and $B_L$ naturally extended by
"Identity" between $\C_i$ and $\C_L$.
\end{enumerate}
\item Same as Definition \ref{d2.27}
\end{enumerate}
\end{dfn}
It is important to note that the measuring instruments just presented measure magnitudes
associated to the \textbf{metric} conform to a potential but have no connection with the singularities found in
the space-time "before" the measure.


\subsubsection{Examples of momentum and energy measurements}\label{4-sss-1}
Once again, we use as an example the elementary oscillating metric of order 1 presented at the beginning of this section:

The metric $g$ is written: $g=|a|^{4/n-2}g_0$ and the function $a$ is of the particular form:

$a=\beta\sum_{j=1}^p\varphi_j$ ~~ where ~~, for any $ j$ from 1 to p, ~~ $\varphi_j:\Theta\times S^1(\delta)\rightarrow\R$
satisfies:

$\varphi_j=C_j \cos(M_jt+Qu-\sum_{k=1}^3\lambda_{jk}x^k)+C'_j\sin(M_jt+Qu-\sum_{k=1}^3\lambda_{jk}x^k)$.

$C_j$ and $C'_j$ are here constants, ~~ $\beta\in E_W(\mu)$, ~~ $Q$ is the relative electric charge.
\begin{enumerate}
 \item We consider a \textbf{momentum measurement} performed on this oscillating metric using the measuring instrument
described in the definition \ref{d2.27}. \\ As we will see, the essential conditions in order that
the measuring instrument "correctly separates into space" each of the $\beta\varphi_j$ are as follows:
 \begin{eqnarray}\label{Fn62}
  \forall j\neq j'\in\{1..p\}~~~~ L|\overrightarrow{\lambda_j}-\overrightarrow{\lambda_{j'}}|>>1 ~~~~ and~~~~
L|\overrightarrow{\lambda_j}|>>1
 \end{eqnarray}
where $|\overrightarrow{\lambda_j}|:=sup(|\lambda_{j_1}|,|\lambda_{j_2}|,|\lambda_{j_3}|)$ and $L$ is the length of
sides of the cube that defines the measurement domain. \\
The fact that for any $ j \in\{1..p\} ~~~~ L|\overrightarrow{\lambda_j}|>>1$ allows us to write the following result: \\
for any $ j\in\{1..p\},~~ \forall \overrightarrow{q}\in Sp$
\begin{eqnarray}\label{Fn63}
 <\varphi_j,h_{1\overrightarrow{q}}>_{B_L}h_{1\overrightarrow{q}}+<\varphi_j,h_{2\overrightarrow{q}}>_{B_L}h_{
2\overrightarrow{q}}\simeq\varphi_j ~~~~ if ~~ L|\overrightarrow{\lambda_j}-\overrightarrow{q}|<<1
\end{eqnarray}
\begin{eqnarray}\label{Fn64}
 <\varphi_j,h_{1\overrightarrow{q}}>_{B_L}h_{1\overrightarrow{q}}+<\varphi_j,h_{2\overrightarrow{q}}>_{B_L}h_{
2\overrightarrow{q}}\simeq 0 ~~~~ if ~~ L|\overrightarrow{\lambda_j}-\overrightarrow{q}|>>1
\end{eqnarray}
where $|\overrightarrow{\lambda_j}-\overrightarrow{q}|:=sup(|\lambda_{j_1}-q_1|,|\lambda_{j_2}-q_2|,|\lambda_{j_3}-q_3|)$.

This is achieved by using equalities of the following form written here for simplicity in "dimension 1" (they are in "dimension 3"): for any $ \lambda\in\R$, ~~ for any $ q\in\R$

$\frac{2}{L}\int_{-L/2}^{L/2}\cos(\lambda x)\cos(qx)dx=\frac{2}{L(\lambda+q)}\sin\frac{L}{2}(\lambda+q)+\frac{2}{
L(\lambda-q)}\sin\frac{L}{2}(\lambda-q)$ ~~ if $\lambda\neq q$

\hspace{3.7cm} $=1+\frac{1}{L\lambda}\sin(L\lambda)$ ~~ if $\lambda=q$

$\frac{2}{L}\int_{-L/2}^{L/2}\cos(\lambda x)\sin(qx)dx=0$ ~~ and we have:

$\frac{2}{L(\lambda-q)}\sin\frac{L}{2}(\lambda-q)\simeq 1$ ~~ if $L|\lambda-q|<<1$

$\frac{2}{L(\lambda-q)}\sin\frac{L}{2}(\lambda-q)\simeq 0$ ~~ if $L|\lambda-q|>>1$
\bigskip

We deduce from \ref{Fn62},~ \ref{Fn63} and \ref{Fn64} that: for any $ j\in\{1..p\} ~~ \forall \overrightarrow{q}\in Sp$

$<a,h_{1\overrightarrow{q}}>_{B_L}{h_{1\overrightarrow{q}}}_{B_L}+<a,h_{2\overrightarrow{q}}>_{B_L}{h_{
2\overrightarrow{q}}}_{B_L}\simeq{\beta\varphi_j}_{\C_L}$ ~~ if $L|\overrightarrow{\lambda_j}-\overrightarrow{q}|<<1$

\hspace{7.3cm} $\simeq 0$ \hspace{0.9cm} if $L|\overrightarrow{\lambda_j}-\overrightarrow{q}|>>1$

Then (def. \ref{d2.27}) the domain $(\C',g')$ created by the measuring instrument is such that $g'=|a'|^{4/n-2}g_0$ where the
function $a'$ satisfies the following properties: for each $j\in\{1..p\}$

$a'_{\C'_{\overrightarrow{q}}}\simeq\beta\varphi_j\circ\sigma_{\overrightarrow{q}}$ ~~ if $\overrightarrow{q}\in Sp$
satisfies $L|\overrightarrow{\lambda_j}-\overrightarrow{q}|<<1$

\hspace{0.85cm} $\simeq 0$ \hspace{1.6cm} if $\overrightarrow{q}\in Sp$
satisfies $L|\overrightarrow{\lambda_j}-\overrightarrow{q}|>>1$

where $\C'_{\overrightarrow{q}}:=]t_0 , t_0+T[\times B_{\overrightarrow{q}}\times S^1(\delta)\times W$.

This means that the function $a'$ is $\simeq 0$ on the domains $\C'_{\overrightarrow{q}}$
when $\overrightarrow{q}$ is "far" from all $\overrightarrow{\lambda_j}$ (ie. when
$L|\overrightarrow{\lambda_j}-\overrightarrow{q}|>>1$) and that $a'$ is "very close" to a function representing an \textbf{homogeneous} elementary oscillating metric on $\C'_{\overrightarrow{q}}$ when $\overrightarrow{q}$ is
"very close" to a $\overrightarrow{\lambda_j}$ (ie. when $L|\overrightarrow{\lambda_j}-\overrightarrow{q}|<<1$).

It is now assumed that the spectrum is sufficiently "rich" so that one can choose for each
$j\in\{1..p\}$, a $\overrightarrow{q_j}\in Sp$ such that $L|\overrightarrow{\lambda_j}-\overrightarrow{q_j}|<<1$.

We can then write: for any $ j\in\{1..p\}$

$a'_{\C'_{\overrightarrow{q_j}}}\simeq\beta\varphi_j\circ\sigma_{\overrightarrow{q_j}}$

where $\varphi_j(t,x,u)=C_j\cos(M_jt+Qu-\sum_{k=1}^3\lambda_{jk}x^k)+C'_j\sin(M_jt+Qu-\sum_{k=1}^3\lambda_{jk}x^k)$

and $(\C'_{\overrightarrow{q_j}},|a'|^{4/n-2}g_0)$ is isometric to $(\C_L,|\beta\varphi_j|^{4/n-2}g_0)$.

It is deduced from the results of section \ref{s2.12} that, if at a given time $t\in ]t_0, t_0+T[$, an elementary singularity "$\varsigma$" is in $\bigcup_{j=1}^p\HH_j(t)$ where $\HH_j(t):=\{t\}\times
B_{\overrightarrow{q_j}}\times S^1(\delta)\times W$, then the probability that it is in a $\HH_j(t)$ is
$({C_j}^2+{C'_j}^2)/\sum_{j=1}^p({C_j}^2+{C'_j}^2)$ since the $B_{\overrightarrow{q_j}}$ are 2 to 2 disjointed (or
at least have their 2 to 2 intersections of zero measure). This probability does not depend on $t$.

This result is comparable to that given by standard quantum physics (specified in \ref{ssn2.1}).
However, the interpretation associated to physics that we present is, as we have just seen, profoundly
different. It is particularly important to note that the elementary singularities observed in
$\HH_j(t)$ and whose probability of presence is known, have nothing to do with the singularities found in the
domain $I\times\Omega\times S^1(\delta)\times W$ for $t<t_0$ in other words, before the measurement.


\subsubsection{Comparison between conditions \ref{Fn62} and "momentum-position" uncertainty relation}
As already mentioned, the momentum $\overrightarrow{\lambda_j}$ is given by
$\overrightarrow{\lambda_j}=M_j\overrightarrow{v_j}$ where $M_j$ is the mass frequency and
$\overrightarrow{v_j}=(v_1,v_2,v_3)$ is the velocity at which the homogeneous elementary oscillating metric moves.
In the SI unit system $M_j=m_jc/\hbar$ where $m_j$ is the mass (in kg) (def. \ref{d2.10}) and
$\overrightarrow{\text{v}_j}=(\text{v}_1,\text{v}_2,\text{v}_3)$ (in $ms^{-1}$) is given by
$\overrightarrow{\text{v}_j}=c\overrightarrow{v_j}$.

The conditions \ref{Fn62} are therefore also written, when we write the "standard" momentum
$\overrightarrow{p_j}=m_j\overrightarrow{\text{v}_j}=\frac{m_0}{(1-|\overrightarrow{\text{v}_j}|^2/c^2)^{
1/2}}\overrightarrow {\text{v}_j}$:
\begin{eqnarray}\label{Fn65}
 \forall j\neq j'\in\{1..p\}~~~~ L|\overrightarrow{p_j}-\overrightarrow{p_{j'}}|>>\hbar ~~~~ and ~~~~
L|\overrightarrow{p_j}|>>\hbar
\end{eqnarray}
Inequality $L|\overrightarrow{p_j}-\overrightarrow{p_{j'}}|>>\hbar$ is to be compared with the uncertainty relation of 
quantum physics that links the measurements of \begin{it} position \end{it} and \begin{it} momentum \end{it} of a
particle.
The interpretation of inequalities \ref{Fn65} is nevertheless quite different. Recall that the notion of \begin{it} measures
of position \end{it} is, for us, defined for \textbf{elementary singularities} while that of \begin{it}
measures of momentum \end{it}, which has just been described, concerns \textbf{metrics conform to a potential}, the uncertainty relation of standard quantum physics loses, for us, all its meaning. Inequalities \ref{Fn65} say
simply that, if one measures momentums using the measuring instrument (def. \ref{d2.27}) and that one wishes
that this one "correctly separates in space" the initial metric in homogeneous elementary oscillating metric
(for which the notion of momentum is naturally defined), then, the measurement domain $B_L$ must be
sufficiently "large" relative to the considered momentum (moreover, inequality $L|\overrightarrow{p_j}|>>\hbar$
specifies that one can not measure, without the criteria which one has just specified, a "too small" momentum).

\item We now consider an \textbf{energy measurement} performed on this same oscillating metric using
the measuring instrument described in the definition \ref{d2.28}.

The essential condition for the measuring instrument to "correctly separate in space" each of the
$\beta\varphi_j$ is the following:
\begin{eqnarray}\label{Fn66}
 \forall j\neq j'\in \{1..p\} ~~~~ T(M_j-M_{j'})>>1 ~~~~et ~~~~ TM_j>>1 
\end{eqnarray}
Analogous calculations (faster here) than those of the previous part show that: for any $ j\in\{1..p\}$ ~~ for any $
e\in Sp$

$(\frac{2}{T}\int_{t_0}^{t_0+T}a_{\C_L}h_{1e}dt)(h_{1e})_{]t_0,
t_0+T[}+(\frac{2}{T}\int_{t_0}^{t_0+T}a_{\C_L}h_{2e}dt)(h_{2e})_{]t_0, t_0+T[}\simeq\beta\varphi_{j_{\C_l}}$

if $T|M_j-e|<<1$

(Recall that $h_e(t)=(h_{1e}(t),h_{2e}(t))=(\cos (et),\sin (et))$)

Then (def. \ref{d2.28}) the domain $(\C',g')$ created by the measuring instrument is such that $g'=|a'|^{4/n-2}g_0$ where the
function $a'$ satisfies the following properties: for each $j\in\{1..p\}$

$a'_{\C'_e}\simeq\beta\varphi_j\circ\sigma_e$ ~~ if $e\in Sp$
satisfies $T|M_j-e|<<1$

\hspace{0.6cm} $\simeq 0$ \hspace{1.3cm} if $e\in Sp$
satisfies $T|M_j-e|>>1$

where $\sigma_e$ is the isometry between $B_e$ and $B_L$ and $\C'_e:=]t_0,t_0+T[\times B_e\times S^1(\delta)\times W$.

This means that the function $a'$ is $\simeq 0$ on the domains $\C'_e$ when $e$ is "far" from all
$M_j$ (ie. when
$T|M_j-e|>>1$) and that $a'$ is "very close" to a function representing an homogeneous
elementary oscillating metric on $\C'_e$ when $e$ is "very close"
to $M_j$ (ie. when $T|M_j-e|<<1$).

It is now assumed that the spectrum is sufficiently "rich" so that one can choose for each
$j\in\{1..p\}$, a $e_j\in Sp$ such that $T|M_j-e_j|<<1$.

We can then write: for any $ j\in\{1..p\}$

$a'_{\C'_{e_j}}\simeq\beta\varphi_j\circ\sigma_{e_j}$

where $\varphi_j(t,x,u)=C_j\cos(M_jt+Qu-\sum_{k=1}^3\lambda_{jk}x^k)+C'_j\sin(M_jt+Qu-\sum_{k=1}^3\lambda_{jk}x^k)$

and $(\C'_{e_j},|a'|^{4/n-2}g_0)$ is isometric to $(\C_L,|\beta\varphi_j|^{4/n-2}g_0)$.

It is deduced from the results of the section \ref{s2.12} that if, at a given time $t\in ]t_0, t_0+T[$, an elementary singularity "$\varsigma$" is in $\bigcup_{j=1}^p\HH_j(t)$ where $\HH_j(t):=\{t\}\times
B_{e_j}\times S^1(\delta)\times W$, then the probability that it is in a $\HH_j(t)$ is
$({C_j}^2+{C'_j}^2)/\sum_{j=1}^p({C_j}^2+{C'_j}^2)$ since $B_{e_j}$ are 2 to 2 disjointed (or
at least have their 2 to 2 intersections of zero measure). This probability does not depend on $t$.

Again this result is comparable to that given by standard quantum physics.


\subsubsection{Comparison between conditions \ref{Fn66} and "time-energy" uncertainty relation}
The conditions \ref{Fn66} written in SI units become:
\begin{eqnarray}\label{Fn67}
 \forall j\neq j'\in\{1..p\}~~~~\mathcal{T}|m_jc^2-m_{j'}c^2|>>\hbar ~~~~ and~~~~\mathcal{T}m_jc^2>>\hbar
\end{eqnarray}
where $\mathcal{T}$ (unit: second) is defined by $c\mathcal{T}=T$ and since $M_j=m_jc/\hbar$.

Inequality $\mathcal{T}|m_jc^2-m_{j'}c^2|>>\hbar$ is to be compared with the uncertainty relation
\begin{it} time-energy \end{it} from
classical quantum physics, but, as for the \begin{it} position-momentum uncertainty relation \end{it}, the interpretation is
for us,
quite different. (It should be noted that the conditions \ref{Fn66} do not involve the "dimension" of $B_L$ which,
here, is of little importance).
\end{enumerate}


\subsection{Magnitudes associated to compact manifolds and their measurements}
The considered cell is of the form $\C=I\times\Omega\times V_1\times V_2$ where $V_1$ and $V_2$ are two compact manifolds. We are interested in the "magnitudes" defined on $V_1$.

The definitions of \begin{it} magnitudes \end{it}, \begin{it} measuring instruments of these magnitudes \end{it} that we are going to specify, are based on the same principles that were used for the magnitudes defined on apparent space, but now the compactness of
$V_1$ (which allows to use the spectral theorem (section \ref{s2.4})) greatly simplifies things.

We consider an eigenspace $E_{V_1}(\mu)$ of the Laplacian operator $\Delta_{g_{V_1}}$. The following scalar product is naturally defined:
\begin{eqnarray}\label{Fn68}
 <\beta_1,\beta_2>:=(vol V_1)^{-1}\int_{V_1}\beta_1\beta_2
\end{eqnarray}
\begin{dfn} \label{d2.29}
\textbf{A magnitude defined on $E_{V_1}(\mu)$} is a family $(h_i)_{i\in\mathcal{G}}$ of functions defined on
$V_1$ with
values in $\R^m$ ($m\geqslant 1$) such that the "union" of the $m$ component families
$(h_{1,i})_{i\in\mathcal{G}},...,(h_{m,i})_{i\in\mathcal{G}}$ forms an orthonormal family (for the scalar product
\ref{Fn68}) of eigenfunctions belonging to $E_{V1}(\mu)$.

According to the spectral theorem this family is finite.
\end{dfn}
Here again, the denomination of the quantity (spin magnitude for example) will often be assigned to the element $i$ of
$\mathcal{G}$.
\begin{dfn} \label{d2.30}
\textbf{An instrument for measuring the magnitude $(h_i)_{i\in\mathcal{G}}$ defined on $V_1$ which three
characteristics
are $Sp$, $B_L$ and $T$} is a physical system that has the following properties:
\begin{enumerate}
 \item From a time $t_0\in I$, it transforms the "subdomain" $(\C_{B_L}, g)$ from $(\C,g)$ into a
"domain"
$(\C',g')$ for which:
\begin{enumerate}
 \item $\C_{B_L}=I\times B_L\times V_1\times V_2$
 \item $\C'=]t_0 , t_0+T[\times(\bigcup_{i\in Sp}B_i)\times V_1\times V_2$
 
 where $]t_0 , t_0+T[\subset I$ and $B_i$ are cubes in $\R^3$ 2 to 2 disjointed (or at least have their 2 to 2 
intersections of zero measurement) each isometric to $B_L$.
\item for any $ i\in Sp$ ~~ $g'_{B_i}=|a_i|^{4/n-2}g_0$ ~ where ~ $a_i:\C_i:=]t_0 , t_0+T[\times B_i\times V_1\times
V_2\rightarrow\R$ is defined by:
$$a_i=\sum_{l=1}^m (<a_{\C_{B_L}},h_{il}>h_{il})\circ\sigma_i$$

Here $\sigma_i$ is the isometry between $B_i$ and $B_L$ naturally extended by "identity product" into an isometry
between $\C_i$ and $\C_{B_L}$.
\end{enumerate}
\item For each $i\in Sp$, the measuring instrument estimates the average number of elementary singularities found
in $ B_i$ during the time $(t,u)\in ]t_0 , t_0+T[\times S^1(\delta)$. (When $n_i(t,u)$ is the number of
elementary singularities in $\HH_i:=\{t\}\times B_i\times\{u\}\times W$ at time $(t,u)$, the
average number is defined by $\bar{n}_i:=\frac{1}{2\pi\delta T}\int_{S^1(\delta)}\int_{t_0}^{t_0+T}n_i(t,u)dudt)$).
\end{enumerate}
\end{dfn}

Of course, the measuring instruments described by definitions \ref{d2.27}, \ref{d2.28} and \ref{d2.30} are
"idealized". In fact, the measurements of \begin{it} magnitudes \end{it} that have been studied in this section and which are
involved
in quantum phenomena, are delicate. It can be considered that a bubble chamber (wire or drift) is an
"association" of measuring instruments as defined (but very imperfect).

The (fictitious) measurement instruments described in this section have, for operating mode, a
"separation in space" of the different "components" of the considered metric conform to a potential.
However, other types of (real) measuring instruments are used and some have, for operating mode, the use of
\begin{it} resonance phenomena \end{it} (magnetic resonance for example). They can not
 be assimilated to those presented here. \\

The following section is dedicated to an important example of magnitude measurement associated to compact manifolds: Measurement
of spin. We detail this study because it is with this particular notion that we will present the "quantum entanglement 
phenomenon" in section \ref{+2.2}.
\newpage

\section{Spin measurement} \label{+2.1}
\subsection{Oscillating metric beam and "spin state"}
In classical physics many experiments use sending "particle beams". This notion will become for us that of "oscillating metric beams".

The considered oscillating metrics are the homogeneous elementary oscillating metrics moving at a velocity
$\overrightarrow{v}$ in a neutral potential, defined in section \ref{s2.6}. These are the ones for which the
notion of velocity (and momentum) is well defined. Here we limit ourselves to 2nd order oscillating metrics with
spin 1/2 (cf \ref{ss+1}) but generalization does not cause any difficulties. Recall (def. \ref{d2.12} and
generalization \ref{ss+2}) that for an homogeneous elementary metric of order 2 and of spin 1/2 defined on the cell
$\C=\Theta\times S^1(\delta)\times S^3(\rho)\times V$, the corresponding metric $g=|a|^{4/n-2}g_0$ is such
that:
\begin{eqnarray}\label{F+1}
a=\beta(\sum_{l=1}^4(C_l\cos(M't-\sum_{k=1}^3\lambda_kx^k+Qu)+C_l'\sin(M't-\sum_{k=1}
^3\lambda_kx^k+Qu))\alpha_l
\end{eqnarray}
where $(\alpha_1,\alpha_2,\alpha_3,\alpha_4)$ is an orthonormal basis of the eigenspace $E_{S^3(\rho)}(\gamma)$ (and can also be considered as a basis of $E^{\CC}_{S^3(\rho)}(\gamma)$), ~~ $\beta\in E_V(\nu)$, ~~
$\overrightarrow{v}=(1/M')(\lambda_1,\lambda_2,\lambda_3)$ ~~ and ~~ $M'=(1-|\overrightarrow{v}|^2)^{-1/2}M$.

The function $a$ satisfies the fundamental equation \ref{F1}: ~~ $\Box_{g_0}a+Sa=0$.

Associated canonical function: $a_c:\Theta\rightarrow E^{\CC}_{S^3(\rho)}(\gamma)$ and state function
$\varPsi:\Theta\rightarrow E^{\CC}_{S^3(\rho)}(\gamma)$ then satisfy (see \ref{ssn2.1} extended to order 2):
\begin{eqnarray}\label{F+2}
a_c=\sum_{l=1}^4(C_l+iC'_l)e^{-i(M't-\sum_{k=1}^3\lambda_{k}x^k)}\alpha_l
\end{eqnarray}
when $Q>0$, and is equal to the
conjugate when $Q<0$.
\begin{eqnarray}\label{F+3}
\varPsi=\sum_{l=1}^4(C_l+iC'_l)e^{-i((M'-M)t-\sum_{k=1}^3\lambda_{k}x^k)}\alpha_l
\end{eqnarray}
for any $Q$.

The term "beam" specifies the fact that the considered oscillating metrics are defined on a "tube" of the
cell $\C=\Theta\times S^1(\delta)\times S^3(\rho)\times V$ ~ where ~ $\Theta=I\times \mathscr U\subset \R\times\R^3$.  Canonical coordinates of $\mathscr U$ will be written more conventionally $(x,y,z)$. For example, we consider (see Figure\ref{f+1}) the "tube" 
$T(r_1,y_1,y_2):=D(r_1)\times]y_1,y_2[\times S^1(\delta)\times S^3(\rho)\times V\subset \mathscr U\times S^1(\delta)\times S^3(\rho)\times V$ where ~ $D(r_1)$ is the Euclidean disk of radius $r_1$ centered in $O$ in the
plane $(\overrightarrow x, \overrightarrow z)$ and $]y_1,y_2[$ is an interval of the axis $\overrightarrow y$. We have
choosen here $\overrightarrow y$ as axis of the tube, the disc $D(r_1)$ being perpendicular to $\overrightarrow y$. Of course, we
extend, if necessary, this definition to a tube of any axis.
\newpage
We give the following definition.
\begin{dfn} \label{d+1}
 \textbf{A domain of type "oscillating metric beam in a neutral potential with  spin 1/2, of radius $r_1$ and
axis
$\overrightarrow y$"} is a domain $(\C,g)$ where $g=|a|^{4/n-2}g_0$ and the function $a:\C\rightarrow \R$ (assumed
regular) satisfies: \begin{enumerate}
                    \item $a/_{I\times T(r_1)}$ is of the form given in (\ref{F+1}) when $\sum_{k=1}^3\lambda_kx^k$
is replaced by ($\lambda y$) to specify that the velocity $\overrightarrow v$ has direction $\overrightarrow y$.
\item $a/_{\C-(I\times T(r_2))}=C^{te}$ ~~ where ~~ $0\leq C^{te}\leq 1$ ~~ and ~~ $r_2=r_1+\varepsilon$.

($\varepsilon$ is generally chosen $<r_1$ and is only used to allow the regularity of the function ($a$) on $\C$).
                   \end{enumerate}

\end{dfn}
(The second condition does not have a great importance in the following, it means that, in the considered domain, the metric
will be approximately that of a neutral potential).
\begin{figure}[]
\begin{center}
    \label{fig:1}
    \begin{tikzpicture}[yscale=0.8, xscale=0.9]
      \draw (0, 0)  ellipse (0.5 and 1);
      \draw (7, 0) ellipse (0.5 and 1);
      \draw (8, 0) ellipse (0.5 and 1);
      \draw     (0, 1) -- (8, 1);
      \draw[->] (-2, 0) -- (12, 0);
      \draw     (0, -1) -- (8, -1);
      \draw[->] (-1, -3) -- (-1, 3);
      \draw[->] (0, 1.5) -- (-2.4, -2);
  
      \draw ( 9, -3) -- ( 9,  2);
      \draw ( 9,  2) -- (11,  3);
      \draw (11,  3) -- (11, -2);
      \draw (11, -2) -- ( 9, -3);
      
      \draw (-1, -0.01) node {\tiny$\bullet$};
      \draw (0, -0.01) node {\tiny$\bullet$};
      \draw (7, -0.02) node {\tiny$\bullet$};
      \draw (8, -0.02) node {\tiny$\bullet$};
      \draw (10, -0.02) node {\tiny$\bullet$};
      
      \draw[thick, ->] (0, 0) node[below right]{$~~~~P$} -- (1.5, 0);
      \draw (10, -0.02) node {\tiny$\bullet$};

      \draw ( 9, -2.7) node[above right]{$E$};
      \draw (-1, -0.01) node[below right]{$0$};
      \draw (0, -0.01) node[below ]{$y_1$};
      \draw (7, -0.01) node[below ]{$y_2$};
      \draw (8, -0.01) node[below ]{$y_3$};
      \draw (10, -0.01) node[below ]{$y_4$};
      \draw (11.5, 0) node[below ]{$\overrightarrow{y}$};
      \draw (-2.8,-1.5) node[right ]{$\overrightarrow x$};
      \draw (-1,2.5) node[right ]{$\overrightarrow z$};
   
      \draw[->] (7.15, -0.7) -- (7.35, 0.7);
      \draw[->] (7.35, -0.7) -- (7.55, 0.7);
      \draw[->] (7.55, -0.7) -- (7.75, 0.7);
      \draw (7.25, 0.7) node{$\overrightarrow{B}$};
    \end{tikzpicture}
  \end{center}
  \caption{The Stern-Gerlach instrument} \label{f+1}
\end{figure}
 
It is easer for the rest to introduce the notion of "spin state" of an homogeneous oscillating metric in
a neutral potential, so we give the following definition.
\begin{dfn} \label{d+2}
\textbf{The spin state} of an homogeneous oscillating metric of order 2 with spin 1/2 is the element $\zeta\in
E^{\CC}_{S^3(\rho)}(\gamma)$ obtained from the state function $\varPsi$ specified in (\ref{F+3} ) by setting:
\begin{eqnarray}\label{F+4}
 \zeta=\frac{\sum_{l=1}^4(C_l+iC'_l)\alpha_l}{(\sum_{l=1}^4|C_l+iC'_l|^2)^{1/2}}
\end{eqnarray}
\end{dfn}
Of course, this definition may be extended without difficulty to "spin" other than 1/2.

\subsection{Idealized Stern-Gerlach instrument. \\ Spin measurement}

The principle of the Stern-Gerlach experiment (Figure \ref{f+1}), described in classical language, is as follows:

We send particles (or atoms) with a well-defined velocity in an area where is present a inhomogeneous magnetic field
with an orthogonal direction to the initial velocity of the particles. We measure the possible deviation
of the particles by this magnetic field by looking at the impacts on a screen orthogonal to the
direction of the particle flow.
The magnetic field $\overrightarrow B$ is zero on the domain for which $y\in]y_1,y_2[\cup]y_3,y_4[$ and
is directed in an orthogonal direction to $\overrightarrow y$ and has a nonzero gradient on the domain for which
$y\in]y_2,y_3[$. \\
For us, in terms of oscillating metrics, Stern-Gerlach experiment is as follows. \\
We consider the cell
$\C=]t_1,t_2[\times\omega\times S^1\times S^3\times V$ ~~ where  $\omega$ is the union of a tube $T_1(r_1,y_1,y_2)$ of 
$\overrightarrow y$-axis and radius $r_1$ (on which $\overrightarrow B=0$), of a tube $T_2(r_2,y_2,y_3)$ (on which
$\overrightarrow B\neq0$) with radius $r_2>r_1$, and a tube $T_3(r_2,y_3,y_4)$ (on which $\overrightarrow B=0$). \\
The metric $g=|a|^{4/n-2}g_\mathcal P$ defined on the cell $\C$ is an elementary oscillating metric
of order 2, which will be assumed of spin 1/2. Its state function is written $\varphi$. \\
The cell $\C_{y_1,y_2}=]t_1,t_2[\times T_1(r_1,y_1,y_2)\times S^1\times S^3\times V$ with the metric
$|a|^{4/n-2}g_0$
defines a domain of type "oscillating metric beam" (def. \ref{d+1}). On this cell, $\varphi$ satisfies:
\begin{eqnarray}\label{F+5}
 \varphi_{|_{\C_{y_1,y_2}}}=\sum_{l=1}^4(C_l+iC'_l)e^{-i((M'-M)t-\lambda y)}\alpha_l
\end{eqnarray}
The associated spin state (def. \ref{d+2}) is $\zeta=\frac{\sum_{l=1}^4(C_l+iC'_l)\alpha_l}{(\sum_{l=1}^4|C_l+iC'_l|^2)^{1/2}}$. Equality (\ref{F+5}) is
therefore a form of "initial condition" for the function $\varphi$.

According to corollary \ref{c2.1} of section \ref{ss+3}, the state function $\varphi$ is solution of the equation
\ref{F50}. We choose the orthonormal basis $(\beta_1,\beta_2,\beta_3,\beta_4)$ of $E^{\CC}_{S^3(\rho)}(\gamma)$ for
which the four components of $\varphi$ satisfy equations (\ref{F55}), (\ref{F56}), (\ref{F57}) and (\ref{F58}).

As we are only interested in the "spin effect", we will neglect, in the differential operator
$(\alpha)$, the terms that have coefficient $Q$, which, in the language of classical physics, means that 
we are only interested in the spin-deviations and not possible deviations from electric charge in an
 electromagnetic potential.

Given the fact that the component $(B^2)$ of the electromagnetic field is here zero, the four equations
 (\ref{F55}) (\ref{F56}) (\ref{F57}) (\ref{F58}) are written :
$$ -2iM\frac{\partial\varphi^1}{\partial t}=\Box_\Theta\varphi^1+{\textstyle{\frac{\varrho}{2}}} Q(B^1\varphi^2+B^3\varphi^1)$$
$$ -2iM\frac{\partial\varphi^2}{\partial t}=\Box_\Theta\varphi^2+{\textstyle{\frac{\varrho}{2}}} Q(B^1\varphi^1-B^3\varphi^2)$$
$$ -2iM\frac{\partial\varphi^3}{\partial t}=\Box_\Theta\varphi^3+{\textstyle{\frac{\varrho}{2}}} Q(B^1\varphi^4+B^3\varphi^3)$$
$$ -2iM\frac{\partial\varphi^4}{\partial t}=\Box_\Theta\varphi^4+{\textstyle{\frac{\varrho}{2}}} Q(B^1\varphi^3-B^3\varphi^4)$$

The first two equations are identical to the last two equations where $\varphi^1$ becomes $\varphi^3$ and $\varphi^2$
becomes $\varphi^4$.

The state function $\varphi$ can be considered as the sum of two state functions $\varphi'$ and $\varphi''$
corresponding to the superposition of two oscillating metrics, the first state function having\\
$(\varphi^1,\varphi^2,0,0)$ for components in the basis $(\beta_1,\beta_2,\beta_3,\beta_4)$ and the second
$(0,0,\varphi^3,\varphi^4)$. In classical physics this corresponds to the sending of two particles of the same mass and electric charge, but possibly with a different spin state. The result of the measurement of the Stern-Gerlach instrument
gives, when an impact takes place on the screen, the probability that it is in one of the two disjoint domains of
the screen (for spin 1/2). This result will be easily obtained for the superposition of the two oscillating metrics if
it is obtained for only one because the calculations are identical for $\varphi'$ and
$\varphi''$. Indeed, this is consequence of the two following facts:
\begin{enumerate}
 \item $E_1^{\CC}:=E^{\CC}_{S^3(\rho)}(\gamma)=E_1'^{\CC}\oplus E_1''^{\CC}$ where $E_1'^{\CC}$ is the linear subspace
 generated by $\beta_1$ and $\beta_2$ and $E_1''^{\CC}$ that generated by $\beta_3$ and $\beta_4$
\item $E_1'^{\CC}$ and $E_1''^{\CC}$ are invariant by the endomorphisms $\hat{S}_1,\hat{S}_2,\hat{S}_3$  in view of
the shape of their matrix in the basis of $\beta_k$ given in \ref{ss+1}.
\end{enumerate}
The following terminology is then introduced.
\begin{dfn} \label{d+3}
The space $E_1'^{\CC}$ is called \textbf{the restricted eigenspace} (It is of complex dimension 2 and corresponds to
the space of spin states 1/2 of standard quantum physics).

Endomorphisms $\hat{S}_1,\hat{S}_2,\hat{S}_3$ induced in $E_1'^{\CC}$ written $\hat{S'}_1,\hat{S'}_2,\hat{S'}_3$
are called \textbf{restricted canonical endomorphisms}.
\end{dfn}
Matrix in the basis $(\beta_1,\beta_2)$ of restricted canonical endomorphisms $\hat{S'}_1$ and $\hat{S'}_3$
(the only ones that will interest us because $B^2=0$) are:

\begin{center}
 $\hat M'_1 = 1/2\begin{pmatrix}
  0&-1&\\
  -1&0&
 \end{pmatrix}$ ~~
$\hat M'_3 = 1/2\begin{pmatrix}
 -1&0&\\
 0&1&
\end{pmatrix}$
\end{center}
We have, of course, the same results if we consider $E_1''^{\CC}$ and $\hat{S''}_1$, $\hat{S''}_3$.

\textbf{In the following of this section we will therefore consider that state functions have values in the restricted eigenspace $E_1'^{\CC}$. Spin states (def. \ref{d+2}) of the homogeneous oscillating metrics of order 2 with
spin 1/2 that will be considered will be elements of $E_1'^{\CC}$}.

On the cell $\C_{y_2,y_3}$ of the idealized Stern-Gerlach instrument, the magnetic field is of the form
$\overrightarrow B=\mathcal{B}\overrightarrow s $ ~~ where $\overrightarrow s$ is an unit vector in the plane
$(\overrightarrow x,\overrightarrow z)$ and $\mathcal B$ is a non-constant function of the variable
written $s=(\cos \theta)z+(\sin \theta)x$ and we will assume that $\frac{\partial\mathcal B}{\partial s}>0$ so that there is 
no ambiguity on the direction of $\overrightarrow s$. The angle $\theta$ is specified in the following definition.
\begin{dfn} \label{d+4}
We call \textbf{measurement angle} of the Stern-Gerlach instrument, the angle $\theta\in]-\pi,\pi]$ given by
$\cos\theta=(\overrightarrow s/\overrightarrow k), ~~~ \sin\theta=(\overrightarrow s/\overrightarrow i)$ ~~ where
$\overrightarrow i$ and $\overrightarrow k$ are the unit vectors of the axis $\overrightarrow x$ and $\overrightarrow
z$.
\end{dfn}
The magnetic field $\overrightarrow B$ is therefore written $\overrightarrow B=\mathcal B(\cos\theta\overrightarrow
k+\sin\theta\overrightarrow i)$ and the two components of $\varphi'$ (now considered as values in
$E_1'^{\CC}$) in the basis
$(\beta_1,\beta_2)$: $(\varphi^1,\varphi^2)$, satisfy, considering the
already specified approximations, the two equations:
\begin{eqnarray}\label{F+8}
  -2iM\frac{\partial\varphi^1}{\partial t}=\Box_\Theta\varphi^1+{\textstyle{\frac{\varrho}{2}}} Q\mathcal
B((\sin\theta)\varphi^2+(\cos\theta)\varphi^1)
\end{eqnarray}
\begin{eqnarray}\label{F+9}
-2iM\frac{\partial\varphi^2}{\partial t}=\Box_\Theta\varphi^2+{\textstyle{\frac{\varrho}{2}}} Q\mathcal
B((\sin\theta)\varphi^1-(\cos\theta)\varphi^2)
\end{eqnarray}
We define the endomorphism $\hat S_\theta$ of $E_1'^{\CC}$ by setting:
\begin{eqnarray}\label{F+10}
 \hat S_\theta:=(\sin\theta) 2\hat S'_1+(\cos\theta)2\hat S'_3
\end{eqnarray}
This one has a matrix associated to the basis $(\beta_1,\beta_2)$:

\begin{center}
 $\hat M_\theta = \begin{pmatrix}
  {-\cos\theta}&\sin\theta&\\
  \sin\theta&\cos\theta&
 \end{pmatrix}$
\end{center}
The two preceding equations can then be written in the following form:
$$-2iM\frac{\partial\varphi'}{\partial t}=\Box_\Theta\varphi'-{\textstyle{\frac{\varrho}{2}}} Q\mathcal B\hat S_\theta(\varphi')$$
We easily verify that the following pair:
\begin{eqnarray}\label{F+11}
 \beta'_1:=(\cos\theta/2)\beta_1+(\sin\theta/2)\beta_2 ~~~,~~~ \beta'_2:=(\sin\theta/2)\beta_1-(\cos\theta/2)\beta_2
\end{eqnarray}
forms a basis of eigenvectors of $\hat S_\theta$, orthonormed in $E_1'^{\CC}$, this for respective eigenvalues
 $(-1)$ and $(+1)$.

The following rotation matrix $R_{\theta/2}$:

\begin{center}
 $R_{\theta/2} = \begin{pmatrix}
  {\cos\theta/2}&\sin\theta/2&\\
  \sin\theta/2&-\cos\theta/2&
 \end{pmatrix}$
\end{center}
then satisfies the equality: $R_{\theta/2}\hat M_\theta=JR_{\theta/2}$ where $J$ is the matrix:
\begin{center}
$J = \begin{pmatrix}
 -1&0&\\
 0&1&
\end{pmatrix}$
\end{center}
It is deduced that the two components $\varphi'^1$ and $\varphi'^2$ of $\varphi'$ in this basis, satisfy
the two equations \textbf{now independent}:\
\begin{eqnarray}\label{F+12}
  -2iM\frac{\partial\varphi'^1}{\partial t}=\Box_\Theta\varphi'^1+{\textstyle{\frac{\varrho}{2}}} Q\mathcal
B\varphi'^1
\end{eqnarray}
\begin{eqnarray}\label{F+13}
  -2iM\frac{\partial\varphi'^2}{\partial t}=\Box_\Theta\varphi'^2-{\textstyle{\frac{\varrho}{2}}} Q\mathcal
B\varphi'^2
\end{eqnarray}
Each of these equations is identical to that which makes it possible to describe, in classical quantum physic, the distribution of the impacts on the screen when the "particles" are under the influence of a potential that here is given by the function ${+\textstyle{\frac{\varrho}{2}}} Q\mathcal B$ for the equation \ref{F+12} and function $-{\textstyle{\frac{\varrho}{2}}} Q\mathcal B$ for the equation
\ref{F+13}. It can be considered that ${\textstyle{\frac{\varrho}{2}}} Q\mathcal B$ (resp. $-{\textstyle{\frac{\varrho}{2}}} Q\mathcal B$) represents an electric potential, the corresponding electric field is then $\overrightarrow E={\textstyle{\frac{\varrho}{2}}} Q\frac{\partial\mathcal B}{\partial
s}\overrightarrow s$ (resp. $\overrightarrow E=-{\textstyle{\frac{\varrho}{2}}} Q\frac{\partial\mathcal B}{\partial s}\overrightarrow s$), 
which, in classical physics, "deflects" the particles in a direction along the direction $\overrightarrow s$ for
the equation (\ref{F+12}) and in the other direction for the equation (\ref{F+13}). \\
We do not detail this study here, but give the qualitative results in the following proposition (without
proof):
\begin{prop} \label{p+1}
 In a domain for which $y_3<y<y_4$ (ie. "after" the magnetic field $\overrightarrow B$), the function
 $a$, characteristic of the metric $g=|a|^{4/n-2}g_0$, is approximately zero outside two disjointed domains
 $\mathscr D_1=]t_1,t_2[\times\omega_1\times S^1\times S^3\times V$ and
$\mathscr D_2=]t_1,t_2[\times\omega_2\times S^1\times S^3\times V$ ~ where  $\omega_1$ and $\omega_2$ are two
disjoint "tubes" included in the part limited by $y_3$ and $y_4$ (Fig. \ref{f+1}). The two "impact spots"
on the screen, classic for the spin 1/2, are limited by the intersection of these two tubes with the screen. The axis  $\overrightarrow y_1$ and $\overrightarrow y_2$ of the
tubes $\omega_1$ and $\omega_2$ intersect the axis
$\overrightarrow y$. The state function of oscillating metric in $\mathscr D_1$ is 
$\varphi'_1\beta'_1|_{]t_1,t_2[\times\omega_1}$, it is similar to that of an oscillating metric beam with axis 
$\overrightarrow
y_1$ (def. \ref{d+1}). It is the same for the state function of the oscillating metric in $\mathscr D_2$.
\end{prop}
The study of this phenomenon is identical to that made by standard quantum physics since, if one
uses $\varepsilon$-approximations, equations (\ref{F+12}) and (\ref{F+13}) are similar to those satisfied by
the state function of the latter. Moreover, the probabilistic analysis made in terms of "singularities" for us,
gives results that also correspond to those of standard quantum physics as we saw in the section
\ref{s2.12}. The precise results are given in proposition \ref{p+2} which will follow in the context of an
idealized Stern-Gerlach instrument. For this purpose, we continue the study as follows.

On the domain for which $y\in]y_1,y_2[$ (Fig. \ref{f+1}), the state function $\varphi'$ is that of the homogeneous oscillating metric moving at a velocity $\overrightarrow v=\lambda\overrightarrow j$ where $\lambda$ is positive and
$\overrightarrow j$ is the unit vector of the axis $\overrightarrow y$. Restricted to this domain, we will denote it
$\varphi'_0$.

According to (\ref{F+2}) and considering restriction to $E_1'^{\CC}$, the function $\varphi'_0$ is written:
\begin{eqnarray}\label{F+14}
 \varphi'_0=e^{-i((M'-M)t-\lambda y)}(z_1\beta_1+z_2\beta_2)
\end{eqnarray}
where $z_1$ and $z_2\in \CC$.

According to definition \ref{d+2}, the spin state of the beam entering the Stern-Gerlach instrument is~~
$(z_1\beta_1+z_2\beta_2)/(|z_1|^2+|z_2|^2)^{1/2}$.

If we consider the function $\epsilon:\Theta\rightarrow\CC$ defined by:
\begin{eqnarray}\label{F+15}
 \epsilon:=e^{-i((M'-M)t-\lambda y)}
\end{eqnarray}
then, on the domain where $y\in]y_1,y_2[$, the two components of $\varphi'_0$ in the basis $(\beta_1,\beta_2)$, that
are $\varphi_0^1$ and $\varphi_0^2$, are written:
$$\varphi_0^1=z_1\epsilon ~~~~\varphi_0^2=z_2\epsilon$$
The components of $\varphi'_0$ in the basis $(\beta'_1,\beta'_2)$ specified in (\ref{F+11}) (and on the domain where
$y\in]y_1,y_2[$) are therefore:
\begin{eqnarray}\label{F+16}
 {\varphi'_0}^1=((\cos\theta/2) z_1+(\sin\theta/2)z_2)\epsilon
\end{eqnarray}
\begin{eqnarray}\label{F+17}
 {\varphi'_0}^2=((\sin\theta/2) z_1-(\cos\theta/2)z_2)\epsilon
\end{eqnarray}
They are, of course, solutions of equations (\ref{F+12}) and (\ref{F+13}) in the domain where $y\in]y_1,y_2[$ since the
magnetic field $\overrightarrow B$ is zero. They can be considered as "initial conditions" which
determine $\varphi'^1$ and $\varphi'^2$ on the "whole" domain, in particular on the domain where $y\in]y_3,y_4[$ which gives
the "measure of spin". \\
To obtain a precise result on this "measure of spin" by the Stern-Gerlach instrument, we specify
now the data on the magnetic field $\overrightarrow B$ in the area where $y\in]y_2,y_3[$.
It has already been assumed that $\overrightarrow B=\mathcal B\overrightarrow s$ where $\mathcal B$ depends on the variable
$s=(\cos\theta)z+(\sin\theta)x$. We will assume now that $\mathcal B$ is an odd function. (We can, 
for example, consider that, by a linear approximation, $\mathcal B=\left(\frac{\partial\mathcal B}{\partial s}(0)\right)s$.
Since, under this condition, $\mathcal B(0)=0$, the field $\overrightarrow B$ is null on the axis $\overrightarrow
y$ and this can be seen as a way to neglect "deviations" (in classical physics language) other than
those associated to "spin". Of course, $\frac{\partial \mathcal B}{\partial s}(0)\neq0$. \\
The choice of an odd function $\mathcal B$ will allow to use the following property:
$$\mathcal B\circ\sigma=-\mathcal B$$
where $\sigma$ is the isometry of $\R^3$ defined by $\sigma(x,y,z)=(-x,y,-z)$.
\begin{prop} \label{p+2}
 It is assumed that the following condition is realized. \\ The domain where $y\in]y_1,y_2[$ (Fig. \ref{f+1}) is a domain of type "oscillating metric beam in a neutral potential with spin 1/2 and axis $\overrightarrow y$"
(def. \ref{d+1}), its spin state is $(z_1\beta_1+z_2\beta_2)/(|z_1|^2+|z_2|^2)^{1/2}$ (def. \ref{d+2}) where $(\beta_1,\beta_2)$ is the basis of
$E_1'^{\CC}$ in which equations \ref{F+8} and \ref{F+9} are satisfied. \\
The measurement angle of the Stern-Gerlach instrument (def. \ref{d+4}) is written $\theta$.
Then: \\ If at time $t\in]t_1,t_2[$, $\varsigma$ is an elementary singularity in $\HH(t):=\{t\}\times
(\omega_1\cup\omega_2)\times S^1\times S^3\times V$, where $\omega_1$ and $\omega_2$ are specified in proposition
(\ref{p+1}), the probability that this one will be in $\HH_1(t):=\{t\}\times\omega_1\times S^1\times S^3\times V$
is:
\begin{eqnarray}\label{F+18}
p_1(t)=\frac{|(\cos\theta/2)z_1+(\sin\theta/2)z_2|^2}{|z_1|^2+|z_2|^2}
\end{eqnarray}
(Of course, the probability that it is in $\HH_2(t):=\{t\}\times\omega_2\times S^1\times S^3\times V$
is: $p_2(t)=\frac{|(\sin\theta/2)z_1-(\cos\theta/2)z_2|^2}{|z_1|^2+|z_2|^2}=1-p_1(t)$). \\
We have chosen here for $\omega_1$ (resp. $\omega_2$) the domain associated with the component $\varphi'_1$ (resp. $\varphi'_2$)
of the state function and therefore to the eigenvalue ($-1$) (resp. ($+1$)) of the endomorphism $\hat S_\theta$ defined by
equality (\ref{F+10}). \\
The probability given by (\ref{F+18}) is also written in the following way:
\begin{eqnarray}\label{F+19}
 p_1(t)=|\langle\zeta,\zeta'_1\rangle|^2
\end{eqnarray}
where $\langle,\rangle$ denotes the hermitian product of the restricted eigenspace $E_1'^{\CC}$ (def.\ref{d+3}), ~ $\zeta$
is the spin state of the oscillating metric beam defined on the domain for which $y\in]y_1,y_2[$, ~
$\zeta'_1:=(\cos\theta/2)\beta_1+(\sin\theta/2)\beta_2$ is the eigenvector of $\hat S_\theta$ defined by
 equality (\ref{F+11}).
\end{prop}
Results given by this proposition are conform to results of standard quantum physics (and
with experimental results).
\bigskip

\textbf{Proof}: according to the imposed initial condition, the restriction of the state function $\varphi'$ to
domain for which $y\in]y_1,y_2[$ is written: $\varphi'_0=\varphi_0'^1\beta'_1+\varphi_0'^2\beta'_2$ ~~~ where $\varphi_0'^1$
and $\varphi_0'^2$ satisfy (\ref{F+16}) and (\ref{F+17}) and $\beta'_1$, $\beta'_2$ are defined by equalities
(\ref{F+11}). \\
Denote $\varPsi:\Theta\supset\omega\rightarrow{\CC}$ the solution of equation (\ref{F+12}) that satisfies
$\varPsi=\epsilon$ on the domain for which $y\in]y_1,y_2[$, where $\epsilon$ is defined by (\ref{F+15}). So,
since
$\epsilon\circ\sigma=\epsilon$ and $\mathcal B\circ\sigma=-\mathcal B$, \\ $\varPsi\circ\sigma$ is solution of equation
(\ref{F+13}) and $\varPsi\circ\sigma=\epsilon$ on the domain for which $y\in]y_1,y_2[$. Moreover,
$\omega_1=\sigma(\omega_2)$ where $\omega_1$ and $\omega_2$ are the two disjoint tubes involved in proposition
\ref{p+1}. \\ From the linearity of the equations we deduce that: \\
\hspace*{1cm} $\varphi'_1=((\cos\theta/2)z_1+(\sin\theta/2)z_2)\varPsi$ ~ and ~ $\varphi'_2=((\sin
\theta/2)z_1-(\cos\theta/2)z_2)(\varPsi\circ\sigma)$\\
According to section \ref{s2.12} and its extension in \ref{ss2.13.6}, if, at time $t$, an elementary singularity
$\varsigma$ is in $\HH(t)$, so the probability that it is in $\HH_1(t)$ is:
$$p_1(t)=\frac{\int_{\omega_1}|\varphi'_1|^2(t,x^i)dx^i}{\int_{\omega_1}|\varphi'_1|^2(t,x^i)dx^i+\int_{\omega_2}
|\varphi'_2|^2(t , x^i)dx^i}$$
So,
$$p_1(t)=\frac{|(\cos\theta/2)z_1+(\sin\theta/2)z_2|^2}{
|(\cos\theta/2)z_1+(\sin\theta/2)z_2|^2+|(\sin\theta/2)z_1-(\cos\theta/2)z_2|^2}$$
because: $\int_{\omega_1}|\varPsi(t,x^i)|^2dx^i=\int_{\omega_2}|\varPsi(t,\sigma(x^i)|^2dx^i$

Which completes the proof of equality (\ref{F+18}) since:
$$|(\cos\theta/2)z_1+(\sin\theta/2)z_2|^2+|(\sin\theta/2)z_1-(\cos\theta/2)z_2|^2=|z_1|^2+|z_2|^2$$
Equality (\ref{F+19}) follows immediately since:
$$(|z_1|^2+|z_2|^2)^{1/2}\langle\zeta,\zeta'_1\rangle=\langle
z_1\beta_1+z_2\beta_2,(\cos\theta/2)\beta_1+(\sin\theta/2)\beta_2\rangle=(\cos\theta/2)z_1+(\sin\theta/2)z_2$$


\subsubsection{Particular case where there exists $\theta'\in]-\pi, \pi]$ such that the spin state is an eigenvector of the endomorphism $\hat S_{\theta'}$.}
\begin{dfn} \label{d+5}
 The angle $\theta'\in]-\pi,\pi]$ such that the spin state $\zeta$ is an eigenvector of $\hat S_{\theta'}$ will be called
\textbf{the spin state angle}. \\
The spin state is then written $\zeta_{\theta'}$.
\end{dfn}
Denote $E_{\theta'}(-1)$ and $E_{\theta'}(+1)$ the eigenspaces of the endomorphism $\hat S_{\theta'}$ associated to
eigenvalues $(-1)$ and $(+1)$. Eigenvectors of norm 1 for $\langle,\rangle$ in $E_{\theta'}(-1)$ are
 of the form: $\zeta_{\theta',k}=e^{ik}((\cos\theta'/2)\beta_1+(\sin\theta'/2)\beta_2)$ ~~ where ~~ $k\in\R$. \\
So:
\begin{eqnarray}\label{F+20}
 \zeta_{\theta',k}^\bot=\zeta_{\theta'-\pi,k}=e^{ik}((\sin\theta'/2)\beta_1-(\cos\theta'/2)\beta_2)
\end{eqnarray}
is an eigenvector of $E_{\theta'}(+1)$.
\begin{dfn} \label{d+6}
 $ \zeta_{\theta',k}^\bot$ given by equality (\ref{F+20}) will be called \textbf{the canonical orthogonal of
$\zeta_{\theta',k}$}.
\end{dfn}
(Any other orthogonal $\zeta_{\theta',k}$ with norm 1 is of the form $e^{ik'} \zeta_{\theta',k}^\bot$ ~~ where
~~ $k'\in\R$). \\
In the case where the spin state of the oscillating metric beam is of the form $\zeta_{\theta',k}$, proposition
\ref{p+2} has a conclusion whose writing is very simple, specified in the following corollary.
\begin{coro} \label{c+2}
 Under the assumptions of proposition \ref{p+2} and when the spin state $\zeta$ has the form $\zeta_{\theta',k}\in
E_{\theta'}(-1)$ defined above, the probability $p_1(t)$ given by (\ref{F+18}) is written:
\begin{eqnarray}\label{F+21}
 p_1(t)=\cos^2(\frac{\theta-\theta'}{2})
\end{eqnarray}
\end{coro}
Indeed, according to (\ref{F+19}):
$$p_1(t)=|\langle
e^{ik}((\cos\theta'/2)\beta_1+(\sin\theta'/2)\beta_2,(\cos\theta/2)\beta_1+(\sin\theta/2)\beta_2\rangle|^2$$
$$=(\cos\theta'/2\cos\theta/2+\sin\theta'/2\sin\theta/2)^2=\cos^2(\frac{\theta-\theta'}{2})$$


\subsection{The idealized Stern-Gerlach instrument as a measuring instrument corresponding to the definition
(\ref{d2.30})}

\noindent Here we specify the correspondence between a measuring instrument based on the Stern-Gerlach experiment and the notion
of "measuring instrument" given by the definition (\ref{d2.30}) of section \ref{s2.16}. \\
We limit ourselves to "spin 1/2", although generalization does not cause conceptual difficulties. \\
The measurement angle of the idealized Stern-Gerlach instrument (def. \ref{d+4}) is written $\theta$. \\
The first notions are: \\
The cube $B_L$ introduced into 1 (a) of the definition (\ref{d2.30}) is a cube contained in the domain $]y_1,y_2[\times
D(r_1)$ of the Stern-Gerlach instrument (Fig. \ref{f+1}). \\
We denote $\beta_{\theta,-1}$ a normalized eigenvector of the endomorphism $\hat S_\theta$ defined by equality
(\ref{F+10}) associated to the eigenvalue (-1) and $\beta_{\theta,+1}$ its canonical orthogonal (def. \ref{d+6}) (then
associated to the eigenvalue (+1)). \\
The "spin magnitude" corresponding to definition (\ref{d2.29}) is the family $(h_i)_{i\in\mathcal G}$ of functions
defined on $S^1(\delta)\times S^3(\rho)$ with values in $\R^2$ ~ where ~ $\mathcal G$ (here equal to the spectrum $Sp$) is
the two-element set that will be written $\{-1,+1\}$. Functions $h_{-1}$ and $h_{+1}:S^1(\delta)\times
S^3(\rho)\rightarrow\R^2$ are defined by:
$$h_{-1}:=(\eta_1\beta_{\theta,-1},\eta_2\beta_{\theta,-1}),
~~~~h_{+1}:=(\eta_1\beta_{\theta,+1},\eta_2\beta_{\theta,+1})$$
where $(\eta_1,\eta_2)$ is the orthonormal basis of $E_{S^1(\delta)}(\lambda)$ given by 
$\eta_1(u)=\sqrt2\cos(Qu)$ and $\eta_2(u)=\sqrt2\sin(Qu)$. \\
The two cubes $B_{-1}$ and $B_{+1}$ introduced in 1 (b) of definition (\ref{d2.30}) are two cubes such that
$B_{-1}\subset\omega_1$ and $B_{+1}\subset\omega_2$ where $\omega_1$ and $\omega_2$ are defined in proposition
\ref{p+1}. They will be considered isometric to $B_L$ by $\sigma_{-1}$ and $\sigma_{+1}$. \\
The oscillating metric beam "entering" into the Stern-Gerlach instrument (def. \ref{d+1}), on which we
perform the spin measurement, have for metric $g=|a|^{4/n-2}g_0$ where the function $a$ is of the form $a=\beta\Phi$.  According to (\ref{F+1}), $\Phi$ satisfies:
$$\Phi=\sum_{k=1}^2(C_k\cos(M't-\lambda y+Qu)+C_k'\sin(M't-\lambda y+Qu))\alpha_k$$
Functions $\alpha_1$ and $\alpha_2$ which form an orthonormal basis of $E_1'^{\CC}$ are chosen such that \\
$\alpha_1=\beta_{\theta,-1}$ and $\alpha_2=\beta_{\theta,+1}$. \\
Functions $a_{-1}$ and $a_{+1}$ defined in 1 (c) of definition (\ref{d2.30}) are then:

$a_{-1}=\frac{\beta}{\sqrt2}((C_1\cos(M't-\lambda y)+C'_1\sin(M't-\lambda y))\eta_1$

\hspace{6cm} $+(C'_1\cos(M't-\lambda y)-C_1\sin(M't-\lambda y))\eta_2)\beta_{\theta,-1}\circ\sigma_{-1}$
\vspace{5mm}

$a_{+1}=\frac{\beta}{\sqrt2}((C_2\cos(M't-\lambda y)+C'_2\sin(M't-\lambda y))\eta_1$

\hspace{6cm} $+(C'_2\cos(M't-\lambda y)-C_2\sin(M't-\lambda y))\eta_2)\beta_{\theta,+1}\circ\sigma_{+1}$ \\
It follows that the associated canonical functions are:
\begin{eqnarray}\label{F+22}
 a_{-1,c}=z_{\theta,-1}e^{-i(M't-\lambda y)}
\end{eqnarray}
\begin{eqnarray}\label{F+23}
  a_{+1,c}=z_{\theta,+1}e^{-i(M't-\lambda y)}
\end{eqnarray}
where $z_{\theta,-1}:=\frac{1}{\sqrt2}(C_1+iC'_1)$ ~~ and ~~ $z_{\theta,+1}:=\frac{1}{\sqrt2}(C_2+iC'_2)$.

As stated in proposition \ref{p+1}, the canonical function $a_{-1,c}$ (resp. $a_{+1,c}$) corresponding to the
oscillating metric in the domain where $\omega_1\subset B_{-1}$ (resp. $\omega_2\subset B_{+1}$) is that given by
(\ref{F+22})
(Resp. (\ref{F+23})). We therefore deduce that: if, at a given time $t$, a singularity $\varsigma$ "is seen" in
$\omega_1\cup\omega_2$, the probability of seeing it in $\omega_1$ is (cf \ref{F42}):
$$p(t)=\frac{|z_{\theta,-1}|^2}{|z_{\theta,-1}|^2+|z_{\theta,+1}|^2}$$
\textbf{This result corresponds to that obtained in proposition \ref{p+2}}. Indeed, according to (\ref{F+11}): \\
\hspace*{0,4cm} $\beta_{\theta,-1}=\beta'_1=(\cos\theta/2)\beta_1+(\sin\theta/2)\beta_2$ ~~~ and ~~~
$\beta_{\theta,+1}=\beta'_2=(\sin\theta/2)\beta_1-(\cos\theta/2)\beta_2$ \\
Then: \\
\hspace*{1,2cm} $z_{\theta,-1}=(\cos\theta/2)z_1+(\sin\theta/2)z_2$ ~~~ and
~~~ $z_{\theta,+1}=(\sin\theta/2)z_1-(\cos\theta/2)z_2$ \\
(We have here, for simplicity, considered a singularity $\varsigma$ at a given time $t$, and not an average number of
singularities during time $(t,u)\in]t_0,t_0+T[\times S^1(\delta)$ as presented in 2, of definition
\ref{d2.30}).
\newpage


\section{Quantum entanglement} \label{+2.2}

\begin{figure}
  \begin{center}
    \label{fig:2}
    \begin{tikzpicture}[yscale=0.8, xscale=0.9]
      \draw (0, 0)  ellipse (0.5 and 1);
      \draw (10, 0) ellipse (0.5 and 1);
  \draw (0.8, 0)  ellipse (0.5 and 1);
  \draw (9.3, 0)  ellipse (0.5 and 1);
  \draw (4.7, 0)  ellipse (0.5 and 1);
  \draw (5.3, 0)  ellipse (0.5 and 1);
      \draw     (0, 1) -- (10, 1);
      \draw[->] (-3, 0) -- (13, 0);
      \draw     (0, -1) -- (10, -1);
        \draw[->] (5, -3) -- (5, 3);
      \draw[->] (8, 0.7) -- (2.4, -0.6);
      
      \draw (-2, -2) -- (-2,  1.5);
      \draw (-2,  1.5) -- ( 0,  2);
      \draw ( 0,  2) -- ( 0, -1.5);
      \draw ( 0, -1.5) -- (-2, -2);
  
      \draw ( 10, -2) -- ( 10,  1.5);
      \draw ( 10,  1.5) -- (12,  2);
      \draw (12,  2) -- (12, -1.5);
      \draw (12, -1.5) -- ( 10, -2);
      
   \draw (9.3, -0.02) node[below] {$y_{2R}$} node {\tiny$\bullet$};   
  \draw (0.8, -0.02) node[below] {$y_{2L}$} node {\tiny$\bullet$};
  \draw (11, -0.02) node[below] {$y_{4R}$}node {\tiny$\bullet$};
  \draw (-1, -0.02) node[below] {$y_{4L}$}node {\tiny$\bullet$};
      \draw (-0.2, -0.02) node[below] {$y_{3L}$}node {\tiny$\bullet$};
      \draw[thick, <-] (2.5, 0) node[above right]{$~~~~~~P_L$} -- (5, 0);
      \draw (5, -0.01) node[below] {O} node {\tiny$\bullet$};
      \draw[thick, ->] (5, 0) node[above right]{$~~~~~~~~~P_R$} -- (7.5, 0);
      \draw (10.2, -0.02) node[below] {$y_{3R}$}node {\tiny$\bullet$};
    \draw (4.5, -0.01) node[below] {$y_{1L}$} node {\tiny$\bullet$};
      \draw (5.5, -0.01) node[below] {$y_{1R}$} node {\tiny$\bullet$};
      \draw (-1, -2.7) node[above right]{$E_L$};
      \draw ( 9, -2.7) node[above right]{$E_R$};
  
      \draw[->] (0.25, -0.35) -- (0.05, 0.35);
      \draw[->] (0.45, -0.35) -- (0.25, 0.35);
      \draw[->] (0.65, -0.35) -- (0.45, 0.35);
      \draw (0.75, 0.7) node{$\overrightarrow{B_L}$};
    \draw (12.5, 0) node[below ]{$\overrightarrow{y}$};
    \draw (5, 2.5) node[right]{$\overrightarrow{z}$};
    \draw (3, -0.4) node[below ]{$\overrightarrow{x}$};
    
      \draw[->] (9.35, -0.35) -- (9.55, 0.35);
      \draw[->] (9.55, -0.35) -- (9.75, 0.35);
      \draw[->] (9.75, -0.35) -- (9.95, 0.35);
      \draw (9.25, 0.7) node{$\overrightarrow{B_R}$};
    \end{tikzpicture}
  \end{center}
  \caption{Stern-Gerlach dual instrument} \label{f+2}
\end{figure}

This section is devoted to the study of a physical phenomenon that is correctly described with the axiomatic system of the
classical quantum physics but which cannot be described in the context of a more classical particle theory. The actual observation of this phenomenon has been one of the arguments that makes
 quantum theory more studied than other classical possible theories. The reader will be able to
obtain more information on this subject by consulting the specialized works of quantum mechanics
(for example \cite{basd}). However, we will show that this phenomenon can be correctly described in the context of
the
theory presented in this paper. \\
For the sake of simplicity, we will focus in this section on the following experiment described in
classical language: In a tube presented in the figure \ref{f+2} for which $y\in]y_{1L},y_{1B}[$,
are "created" two particles $P_L$ and $P_R$, each with spin 1/2. $P_L$ is thrown into a 
Stern-Gerlach instrument (see section \ref{+2.1}) on the left, $P_R$ to a Stern-Gerlach instrument on the right. These
instruments measure the spin of particles and each contains a magnetic field $\overrightarrow B$
($\overrightarrow B_L$ for the left one, $\overrightarrow B_R$ for the right one) inhomogeneous but
approximately made up of parallel vectors whose direction can be chosen in the plane
($\overrightarrow x,\overrightarrow z$), then a screen $E$ ($E_L$ and $E_R$). \\
Each of the particles, after passing through the magnetic field $\overrightarrow B$ of the corresponding instrument,
mark an impact on the screen $E$. The spin of the particles being "1/2", the impacts take place only inside two zones
for each relevant screen, which gives the (binary) spin measure of the particle. \\
Denote $(+1_L)$, $(-1_L)$ ~~ (resp. $(+1_R)$, $(-1_R)$) the two possible measurements for the particle $P_L$
~~ (resp. $P_R$). \\
We are interested in the probability of obtaining pairs of simultaneous measurements:\\
$(+1_L,+1_R)$, $(+1_L,-1_R)$, $(-1_L,+1_R)$, $(-1_L,-1_R)$, and this according to the orientation of $\overrightarrow B_L$ relative to $\overrightarrow B_R$.
Experimentally, we perform a statistical study, one carries out several times an identical experiment with
the one just described. \\
The remarkable fact, deduced from experiment, is that the statistical result obtained on the $4$ possible 
 "pairs of measures" that we have just written, can only be explained if \textbf{the measurement of the spin
of one particle has an influence on the measurement of the spin of the other.} This excludes any interpretation in terms of "two
independent particles" in classical sense. \\ We summarize in the following lines the arguments that lead to
 this conclusion. For this, we introduce the \textbf{correlation function $E(\overrightarrow B_L,
\overrightarrow B_R)$.} This function is equal, for the orientations of $\overrightarrow B_L$ and $\overrightarrow
B_R$, to the average value of the product of the measurement results ($(+1)$ or $(-1)$) of the Stern-Gerlach instruments
located on the left and right. We have necessarily $|E(\overrightarrow B_L,\overrightarrow B_R)|\leq1$. \\
We perform the pair of measurements for two different directions of the magnetic field on the left ($\overrightarrow
B_L,\overrightarrow B'_L$) and two directions of the magnetic field on the right ($\overrightarrow B_R,\overrightarrow
B'_R$). \\
We define the quantity:
\begin{eqnarray}\label{F+24}
S= E(\overrightarrow B_L,\overrightarrow B_R)+E(\overrightarrow B_L,\overrightarrow B'_R)+E(\overrightarrow
B'_L,\overrightarrow B_R)-E(\overrightarrow B'_L,\overrightarrow B'_R)
\end{eqnarray}
(We sum three of the values and subtract the fourth).

Then, for any "hidden variable theory" whose definition is specified for example
in \cite{basd}, the following inequality, called "Bell inequality", is satisfied (there are several Bell inequalities 
of that type):
\begin{eqnarray}\label{F+25}
 |S|\leq2
\end{eqnarray}
The axiomatic system of classical quantum physics (which does not correspond to an hidden variable theory)
gives, for the experiment just described and with a good choice of the four directions of measurements
$\overrightarrow B_L,\overrightarrow B'_L,\overrightarrow B_R,\overrightarrow B'_R$, the theoretical result $|S|=2\sqrt2$
which contradicts Bell's inequality (\ref{F+25}).

The experiment, carried out at Orsay \cite{aspe} by A. Aspect, P.Grangier and G.Roger (which will be written later A.G.R),
gets a result for which $|S|$ is very close to $2\sqrt2$ (the
small
difference with this value is justified by the characteristics of the instruments used). \\
This experiment uses pairs of photons emitted by an atomic cascade of calcium atoms caused by
lasers, it is not therefore "spin 1/2" that we measured but "polarities", however the
transposition of the theory from one notion to another does not cause difficulties and gives identical results. \\
We begin by quickly recalling the principles of standard quantum mechanics that make possible to obtain the
results cited above. The analysis of the phenomenon with the point of view of our theory will be presented in subsection
\ref{ss+4}, this will be the essential point of section \ref{+2.2}.
\subsection{The type A.G.R experiment \cite{aspe} as seen by standard quantum physics}
The axiomatic data of the standard quantum physics that we will use are the following:

\begin{enumerate}
 \item The space of spin states 1/2 of a particle is a complex dimension 2 Hilbert space $(\mathcal{E},\langle,\rangle)$  that can be likened to $\CC^2$ with the standard Hermitian product.
\item The three observables of spin $\hat S_1$, $\hat S_2$, $\hat S_3$ associated to directions $\overrightarrow x$,
$\overrightarrow y$, $\overrightarrow z$, are hermitian endomorphisms whose matrix are:
\begin{center}
 $\hat M_1= \begin{pmatrix}
  0&1&\\
  1&0&
 \end{pmatrix}$ ~~~
$\hat M_2= \begin{pmatrix}
  0&-i&\\
  i&0&
 \end{pmatrix}$ ~~~
 $\hat M_1= \begin{pmatrix}
  1&0&\\
  0&-1&
 \end{pmatrix}$

\end{center}

(Coefficient $\hslash/2$ that usually appears, where $\hslash$ is the Planck constant, has been removed here
without any consequences on the studied phenomena). \\
More generally, the spin observable $\hat S_{\overrightarrow u}$ associated to the direction $\overrightarrow u=
\sin\theta\cos\phi \overrightarrow i+ \sin\theta\sin\phi \overrightarrow j+ cos\theta \overrightarrow k$ has for matrix:
\begin{center}
 $\hat M_{\overrightarrow u}= \begin{pmatrix}
  \cos \theta&\sin\theta e^{-i\phi}&\\
  \sin \theta e^{i\phi}&-\cos \theta&
 \end{pmatrix}$
 \end{center}
Note that the eigenvalues of $\hat S_{\overrightarrow u}$ are $(-1)$ and $(+1)$.
\item The probability of obtaining $(-1)$ ~~ (resp. $(+1)$) for the spin measurement associated to the direction $\overrightarrow u$
of a particle, whose normalized spin state is $\xi\in\mathcal{E}$, is given by:
$$p=|\langle\xi_{\overrightarrow u},\xi\rangle|^2$$
where $\xi_{\overrightarrow u}$ is a normalized vector of the endomorphism $\hat S_{\overrightarrow u}$ for the
eigenvalue $(-1)$ ~~ (resp. $(+1)$).
\item For a pair of particles $(P_1,P_2)$, the spin state space is: $\mathcal{E}\otimes\mathcal{E}$. As
the dimension of $\mathcal{E}$ is $2$ for the spin 1/2, every element of $\mathcal E\otimes\mathcal{E}$ goes under the
form $\Psi=\xi^1\otimes\xi^2+\xi^3\otimes\xi^4$. When $\Psi$ is of maximum rank 2, the concept of "entangled states" appears.
\item For a pair of particles $(P_1,P_2)$ whose normalized state function $\Psi\in\mathcal E\otimes\mathcal E$,
when performing a spin measurement on the particle $P_1$ associated to the direction $\overrightarrow u_1$ and a measurement
of spin on the particle $P_2$ associated to the direction $\overrightarrow u_2$, the probability of obtaining $(-1,-1)$
~~ (resp. $(-1,+1), (+1,-1),(+1,+1)$) is:
$$p=|\langle\xi_{\overrightarrow u_1}\otimes\xi_{\overrightarrow u_2}, \Psi\rangle_{\mathcal
E\otimes\mathcal{E}}|^2$$
where $\xi_{\overrightarrow u_1}$ is a normalized eigenvector of $\hat S_{\overrightarrow u_1}$ for the eigenvalue
$(-1)$ (resp. $(+1)$) and $\xi_{\overrightarrow u_2}$ a normalized eigenvector of $\hat S_{\overrightarrow u_2}$
for the eigenvalue $(-1)$ (resp. $(+1)$).
\item If, for a pair $(P_1,P_2)$ of particles whose normalized state function is $\Psi\in\mathcal
E\otimes\mathcal E$, we only measure the spin on (for example) $P_1$ associated to the direction $\overrightarrow
u_1$ (without being concerned about $P_2$), the probability of obtaining the measurement $(-1)$ ~~ (resp. $(+1)$) is:
$$p=\sum_{k=1}^2|\langle\xi_{\overrightarrow u_1}\otimes\xi_k,\Psi\rangle|^2$$
where $\xi_{\overrightarrow u_1}$ is a normalized vector of $\hat S_{\overrightarrow u_1}$ for the eigenvalue
$(-1)$ ~~ (resp. $(+1)$) and $(\xi_1,\xi_2)$ is an orthonormal basis of $\mathcal{E}$ (the result does not depend on the choice of the basis).
\end{enumerate}
\begin{rmq}
The reader will draw the parallel between the space of spin state given in 1. (resp. the spin observables given in 2.) with
the restricted eigenspace (resp. the restricted endomorphisms) given by the definition \ref{d+3}. It is important to notice
that the first are abstract axiomatic data of classical quantum physics whereas those defined
in the point of view of our theory are precise data defined from the geometry of $S^3$.
\end{rmq}

In the context of the A.G.R \cite{aspe} experiment, the normalized state function of the pair of particles
$(P_1,P_2)$ is assumed to be of the form: $\Psi=\frac{1}{2}(\xi\otimes\xi+\xi^\bot\otimes\xi^\bot)$ ~~ and is of
rank 2. \\
Here $\xi$ (norm 1) in $\mathcal{E}$ and $\xi^\bot$, well chosen, satisfy
$\langle\xi^\bot,\xi\rangle=0$.

When axiom 5 is applied, the probability of obtaining $(-1,-1)$\\
(resp. $(-1,+1), (+1,-1), (+1,+1)$) as a spin measure of
$(P_1,P_2)$ associated to the direction $(\overrightarrow u_1,\overrightarrow u_2)$ is:
$$p=|\langle\xi_{\overrightarrow u_1}\otimes\xi\rangle\langle\xi_{\overrightarrow
u_2}\otimes\xi\rangle+\langle\xi_{\overrightarrow u_1}\otimes\xi^\bot\rangle\langle\xi_{\overrightarrow
u_2}\otimes\xi^\bot\rangle|^2$$
where $\xi_{\overrightarrow u_1}$ ~~ (resp. $\xi_{\overrightarrow u_2}$) is a normalized eigenvector of $\hat
S_{\overrightarrow u_1}$ ~~ (resp. $\hat S_{\overrightarrow u_2}$) for the eigenvalue $(-1)$ ~~ (resp. $(+1)$) . \\
It is then possible to calculate the correlation function $E(.,.)$ for two
directions ${\overrightarrow u_1}$ and ${\overrightarrow u'_1}$ associated to $P_1$ and two directions ${\overrightarrow
u_2}$ and ${\overrightarrow u'_2}$ associated to $P_2$ and then calculating the magnitude $S$ defined by (\ref{F+24}). We
then verifies, by a good choice of measurement directions ${\overrightarrow u_1}$, ${\overrightarrow u'_1}$,
${\overrightarrow u_2}$, ${\overrightarrow u'_2}$, that $|S|>2$, which contradicts the Bell inequality (\ref{F+25}).
\begin{rmq}
 In classical quantum physics, the fact that, in a physical system, there is a precise number of particles
is a well-defined notion. The state space is then axiomatically given by the tensor product of the spaces
states of each type of particle. In the context of our theory, things are profoundly different, only
the oscillating metrics are perfectly defined. The singularities (see section \ref{s2.11}) which, during a measurement
(see section \ref{s2.16}), give the equivalent of the presence of particles in classical quantum physics, are managed only
by the equiprobability principle relative to the metric $g$ (see section \ref{s2.12}) and only a "probabilistic estimation"
makes sense to us. It is then difficult to draw a mathematical parallel between the axiomatic system of the
classical quantum physics (which deals with several particles using the notion of tensor product) and
computational procedures associated to our theory. This will be clarified in the next subsection.
\end{rmq}


\subsection{The quantum entanglement described by our theory} \label{ss+4}
For the sake of clarity, we start by defining some important notions.
\begin{dfn} \label{d+7}
 \textbf{A double oscillating metric} is an oscillating metric of order 2 in an electromagnetic potential with
metric $g=|a|^{\frac{4}{n-2}}g_\mathcal P$, it is defined on the domain determined by the "double instrument of
Stern-Gerlach "described by the figure \ref{f+2}. On the domain for which $y\in]y_{1R}, y_{3R}[$ ~~ (resp. $y\in]y_{1L},
y_{3L}[$) the metric is written $g_R=|a_R|^{\frac{4}{n-2}}g_{\mathcal P_R}$
~~ (resp. $g_L=|a_L|^{\frac{4}{n-2}}g_{\mathcal P_L}$). The domain for which $y\in]y_{1R}, y_{2R}[$
~~ (resp. $y\in]y_{1L}, y_{2L}[$) is that of an oscillating metric beam in a neutral potential with spin 1/2,
axis $\overrightarrow y$, velocity $\overrightarrow v$ in the same direction as $\overrightarrow y$ ~~ (resp. $-\overrightarrow
y$).
\end{dfn}
 The concepts introduced in the following definitions are essential in the description of the phenomenon
of "quantum entanglement". We will only use a particular case of definition (\ref{d+8}) which follows,
this will be useful for a more general study of entanglement phenomena. The particular case that we are going to
use here will be presented in definition (\ref{d+9}). \\
We consider an elementary oscillating metric of order 2 in an electromagnetic potential
$g=|a|^{\frac{4}{n-2}}g_\mathcal P$, as well as an isomorphism $\sigma$ of space
$E_{\lambda,\gamma}:=E_{S^1}(\delta)\otimes
E_{S^3}(\gamma)$. This space will be considered as real vector space of dimension 8 for which we will take as a basis
$(\nu_{ij})\\
i\in\{1,2\} ~~j\in\{1,2,3,4\}$ ~~~ where ~~~ $\nu_{1j}(u,s):=(\cos u)\alpha_j(s)$ ~ and ~$\nu_{2j}(u,s):=(\sin
u)\alpha_j(s)$. \\ (Here $(\alpha_1,\alpha_2,\alpha_3,\alpha_4)$ is the basis of $E_{S^3}(\gamma)$ defined in \ref{ss+1}
example 1). \\
The function $a$ is of the form $\phi\beta$ (def. \ref{d2.9}) ~~ where ~~ $\phi:\Theta\times S^1(\delta)\times
S^3(\gamma)\rightarrow\R$ is written: \\ for any $ x\in\Theta, ~~~\phi_x(.)=\sum_{i,j}\phi^{ij}(x)\nu_{ij}(.)$. \\
Denote $\phi_\sigma:\Theta\times S^1(\delta)\times S^3(\gamma)\rightarrow\R$ the function defined by: \\ for any $
x\in\Theta, ~~~{\phi_\sigma}_x(.)=\sum_{i,j}\phi^{ij}(x)\sigma(\nu_{ij})(.)$ ~~ and ~~ $a_\sigma:=\phi_\sigma\beta$
\begin{dfn} \label{d+8}
 \textbf{The $\sigma$-transformed} of the metric $g=|a|^{\frac{4}{n-2}}g_\mathcal P$ is the metric
$g_\sigma=|a_\sigma|^{\frac{4}{n-2}}g_\mathcal P$.
\end{dfn}
Of course, $g_\sigma$ is still an elementary metric of order 2 with spin 1/2. \\

We limit ourselves now to the particular case of an homogeneous oscillating metric of order 2 with spin 1/2 which will be,
in the following, the restriction to the domain for which $y\in]y_{1R},y_{2R}[$  (resp. $y\in]y_{1L},y_{2L}[$) of a 
double oscillating metric (def. \ref{d+7}). We only consider the case where the spin state $\zeta_R\in E_1'^{\CC}$ ~~ (
resp. $\zeta_L$) is an eigenvector of an endomorphism $\overrightarrow S_{\theta_R}$
~~ (resp. $\overrightarrow S_{\theta_L}$) defined by equality (\ref{F+10}).

One consider $\theta_T\in]-\pi, +\pi]$ ~~ and ~~ $g=|a|^{\frac{4}{n-2}}g_0$ an homogeneous oscillating metric of order 2 with spin
1/2 whose spin state is of the form $\zeta_{\theta,k}=e^{ik}((\cos\theta/2)\beta_1+(\sin\theta/2)\beta_2)$.
\begin{dfn} \label{d+9}
 \textbf{The $\theta_T$-transformed oscillating metric $g$, written $g_{\theta_T}$}, is the $\sigma$-transformed of
$g$ (def. \ref{d+8}) for which the spin state is:
$\zeta_{\theta-\theta_T,k}=e^{ik}((\cos\frac{\theta-\theta_T}{2})\beta_1-(\sin\frac{\theta-\theta_T}{2})\beta_2)$.
\end{dfn}
\textbf{The $\theta_T$-transformations are therefore those that cause an angle rotation $\theta_T/2$ of the spin state} ~~ (The reader will be able to explain the isomorphisms $\sigma$ (in fact, rotations) of
$E_{S^1}(\lambda)\otimes E_{S^3}(\gamma)$ that correspond to these $\theta_T$-transformations).

The last notion introduced is that of "successive oscillating metrics", it express in terms of
oscillating metric what in classical physics language we would call "successive particle waves".  
Particles can be of different nature depending on the considered time interval. \\
We therefore consider a time interval $]t_1,t_m[$ decomposed into successive intervals $I_1,I_2 \dots I_m$ \\ where
~ $I_k:=]t_k,t_{k+1}[$. \\
A domain of the form $\C=]t_1,t_m[\times\Omega\times S^1(\delta)\times S^3(\rho)\times V$ is decomposed in a 
disjunct union of subdomains $:\C_k:=I_k\times\Omega\times S^1(\delta)\times S^3(\rho)\times V$. \\
We give the following definition:
\begin{dfn} \label{d+10}
 \textbf{A domain of type "successive oscillating metrics"} is a domain $(\C,g)$ ~ where ~ $\C=\cup_{k=1}^m\C_m$
is described in the previous lines. For even indices ($k=2l$) (for example), the metric $g|_{\C_k}$ is
that of an oscillating metric that is written $|a_k|^{\frac{4}{n-2}}g_\mathcal P$. For odd numbers ($k=2l+1$) the
metric $g|_{\C_k}$ is of the form $c_kg_\mathcal P$ ~ where ~ $c_k$ is a non negative constant.
\end{dfn}
This notion of "successive oscillating metrics" is, of course, idealized. In the following, constants $c_k$
can be considered $<<1$ (when $c_k=0$ the metric is degenerate).


\subsubsection{Analysis of the A.G.R experiment \cite{aspe} in terms of oscillating metrics}
We consider the Stern-Gerlach double instrument described in the figure \ref{f+2}. We are referring to
the experiment \cite{aspe}, where the spin is measured and not the polarity of
photons. Indeed we have not developed, in this paper, this last notion associated for us to  
zero mass oscillating metrics (see \ref{s2.14}). \\
In standard quantum physics, the  polarization state of a pair of photons emitted by an "atomic cascade"
of calcium atoms in a domain where $y\in]y_{1L},y_{1R}[$ is assumed to be the next "entangled state"
(cf figure \ref{f+2}):
\begin{eqnarray}\label{F+26}
 \frac{1}{2}(\xi\otimes\xi+\xi^\bot\otimes\xi^\bot)
\end{eqnarray}
We will study the more general case for which the polarization (or spin) state is of the form:
 \begin{eqnarray}\label{F+27}
  \frac{1}{2}(\xi_1\otimes\xi_2+\xi_1^\bot\otimes\xi_2^\bot)
 \end{eqnarray}

\textbf{In terms of oscillating metrics of order 2 with spin 1/2, this is replaced by the data presented in the
following paragraph}.
\subsubsection{Successive double oscillating metrics (def. \ref{d+7} and def. \ref{d+10}) created in the domain where
$y\in]y_{1L},y_{1R}[$}
These oscillating metrics will be defined by the succession of double oscillating metrics which can be of two
different types, these types being equiprobably distributed over the concerned time intervals. We start by giving the following details. \\
The total time interval of the experiment $I=]t_1,t_m[$ is decomposed as $I=\cup_{k=1}^{m-1}I_m$. \\
The subdomain of $\C=I\times\omega\times S^1\times S^3\times V$ for which $y\in]y_i,y_j[$ is written $\C_{y_i,y_j}$.
Double oscillating metrics are defined from the following four functions $a_1^L,a_2^L,a_1^R,a_2^R$
($L$ and $R$ refer to the left and right part of the instrument described in Figure \ref{f+2}):
\begin{enumerate}
 \item $a_1^L:\C_{y_{4L},y_{1L}}\rightarrow\R$ is such that $g_1^L:=|a_1^L|^{\frac{4}{n-2}}g_\mathcal P$ is an
elementary oscillating metric of order 2 with spin 1/2, its restriction to $\C_{y_{2L},y_{1L}}$ is a beam of
oscillating metric (def. \ref{d+1}) whose spin state has angle $\theta_1^L$ and is written $\gamma_{\theta_1^L}$
(def. \ref{d+5}), its velocity $\overrightarrow v_1^L$ has a direction opposite to $\overrightarrow y$.
\item $a_1^R:\C_{y_{1R},y_{4R}}\rightarrow\R$ has the same properties as $a_1^L$, the spin state of the restriction of
$g_1^R$ to $\C_{y_{1R},y_{2R}}$ is written $\gamma_{\theta_1^R}$, its velocity is $\overrightarrow v_1^R=- \overrightarrow
v_1^L$ (we assume that $|\overrightarrow v_1^L|= |\overrightarrow v_1^L|$ to simplify
conditions on the dimensions of the instrument for simultaneous measurements in $L$ and $R$).
\item $a_2^L:\C_{y_{4L},y_{1L}}\rightarrow\R$ has the same properties as $a_1^L$ but the spin state
$\gamma_{\theta_2^L}$ of the restriction of $g_2^L$ to $\C_{y_{2L},y_{1L}}$ is the canonical orthogonal
$\gamma_{\theta_1^L}^\bot$ from $\gamma_{\theta_1^L}$ (def. \ref{d+6}).
\item $a_2^R:\C_{y_{1R},y_{4R}}\rightarrow\R$ has the same properties as $a_1^R$ but the spin state
$\gamma_{\theta_2^R}$ of the restriction of $g_2^R$ to $\C_{y_{1R},y_{2R}}$ is the canonical orthogonal
$\gamma_{\theta_2^R}^\bot$ from $\gamma_{\theta_2^d}$.
\end{enumerate}
Since the concepts are fundamentally different, the reader can draw a parallel between the conditions that
we have just given on the spin states of the metrics associated with the four preceding functions and the state of
polarization (\ref{F+27}) given in classical quantum physics. \\
The effective oscillating metrics that will allow to obtain results identical to those of 
classical quantum physics are specified in the following lines. \\
We write $\C_{y_i,y_j}^k$ the subdomain of $\C_{y_i,y_j}$ for which $t\in I_k$. \\
For an index $k=2l+1$, the metric on $\C$ is of the form $g=c_kg_\mathcal P$ ~ where ~ $c_k$ 
 is a non negative constant. \\
For an index $k=2l$ ~~ the metric is a double oscillating metric (def. \ref{d+7}), it is
\textbf{equiprobably} one of the following two metrics:
\begin{enumerate}
 \item $g_{\theta_{T_1}}$: The $\theta_{T_1}$-transformed (def. \ref{d+9}) of the metric
$g=|a_1|^{\frac{4}{n-2}}g_\mathcal P$ ~ where the function $a_1$ defined on $\C_{y_{4L},y_{1L}}^{2l}$ is $a_1^L$ and
defined on $\C_{y_{1R},y_{4R}}^{2l}$ is $a_1^R$.
\item $g_{\theta_{T_2}}$: The $\theta_{T_2}$-transformed metric
$g=|a_2|^{\frac{4}{n-2}}g_\mathcal P$ ~ where the function $a_2$ defined on $\C_{y_{4L},y_{1L}}^{2l}$ is $a_2^L$ and
defined on $\C_{y_{1R},y_{4R}}^{2l}$ is $a_2^R$.
\end{enumerate}
Here are some examples for the choices of $\theta_{T_1}$ and $\theta_{T_2}$ which give the results (identical to those
of standard quantum mechanics) presented in proposition \ref{p+3} that will follow. A comment on these choices is
given after the presentation of the examples. \\
We write $\theta^L$ the measurement angle of the Stern-Gerlach instrument on the left and $\theta^R$ that of the instrument
located on the right (Figure \ref{f+2}).
\vspace{3mm}

\textbf{Example 1}: $\theta_{T_1}=\theta^L-\theta_1^L$ ~~~ and ~~~ $\theta_{T_2}=\theta^L-\theta_2^L$. \\
(Here, we choose the left part of the instrument but we can symmetrically choose the right part by setting

$\theta_{T_1}=\theta^R-\theta_1^R$ ~~~ and ~~~ $\theta_{T_2}=\theta^R-\theta_2^R$). \\
The left-right dissymmetry of Example 1 may seem unjustified but note that this happens in the A.G.R experiment 
(\cite{aspe}), since one selects by their frequencies the photons going towards the
left
and those going to the right. \\
The following example has a left-right symmetry.
\vspace{3mm}

\textbf{Example 2}: There is equiprobability between the four following possibilities: \\
$\theta_{T_1}=\theta^L-\theta_1^L$, ~~ $\theta_{T_1}=\theta^R-\theta_1^R$, ~~ $\theta_{T_2}=\theta^L-\theta_2^L$,
~~ $\theta_{T_2}=\theta^R-\theta_2^R$. \\

But we can present other examples like the one for which: \\
$|\theta_{T_1}|=min(|\theta^L-\theta_1^L|,|\theta^R-\theta_1^R|)$ ~~ and
~~ $|\theta_{T_2}|=min(|\theta^L-\theta_2^L|,|\theta^R-\theta_2^R|)$. \\

In fact, a choice of $\theta_{T_1}$ and $\theta_{T_2}$ can only be specified by a study of the phenomenon of
"creation" of the double oscillating  metrics (study of the atomic cascade of calcium atoms for example), but this
paper does not yet address this complex subject. Otherwise, it is possible that the assumptions that have been made to
characterize successive double oscillating metrics are significantly modified in the following. The goal of
this section is mainly to show that the phenomena of "entanglement" studied by quantum mechanics
can be described qualitatively and quantitatively with our theory. \\

(We can briefly summarize the process by saying that, systematically, a "rotation" in $S^1(\delta)\times
S^3(\rho)$ of the double oscillating metric $g=|a_1|^{\frac{4}{n-2}}g_\mathcal P$ makes coincide: either the angle of
the spin state linked to $a_1^L$ with the measurement angle $\theta_L$, or the angle of the spin state linked to $a_1^R$ with
the measurement angle $\theta_R$). \\

It may seem surprising at first sight that the choices of $\theta_{T_1}$ and $\theta_{T_2}$, given in the
previous examples and which make it possible to obtain the results of proposition \ref{p+3}, depend on the 
measurement angles $\theta^L$ and $\theta^R$ defined from the magnetic fields $\overrightarrow B_L$ and $\overrightarrow B_R$.
These fields are, in fact, located at a $g_0|_\Theta$ -distance of the "creation" domain of oscillating metrics
which can be 
important \\ \textbf{However, we must not forget that all notions introduced on the spin (or
 the polarity) only concern the compact manifold $S^1(\delta)\times S^3(\rho)$ and that it appears
in the Cartesian product with apparent space $\Theta$ (objects associated to the notion of spin are 
not "included" in apparent space). This profoundly changes habits from in 
classical physics where a similar procedure is not usable. Otherwise, we can notice that,
during the time intervals $I_{2l+1}$ (immediately before creating a double oscillating metric) the
metric $g$ is of the form $cg_\mathcal P$ and, if the constant $c$ is $<<1$, the effective distance, which is the
$g|_\Theta$-distance, can be as small as desired compared to the ${g_0}|_\Theta$-distance, this
 last is the one measured by an observer outside the field of experiment}. \\
"The ghost link", which Einstein was ironically talking about two entangled particles, is supported by, with the point of view
 of our theory, the compact manifold $S^1(\delta)\times S^3(\rho)$. \\

The important result of this section is:
\begin{prop} \label{p+3}
We consider the instrument described in Figure \ref{f+2} \\
We denote $\theta^L$ ~ (resp. $\theta^R$) the measurement angle
(def. \ref{d+4}) of the Stern-Gerlach instrument on the left (resp. on the right). \\
Double oscillating metrics are those described in the previous chapter. We consider the
examples presented in this section. \\
We denote $\omega_{1L}$ and $\omega_{2L}$ ~ (resp. $\omega_{1R}$ and $\omega_{2R}$) the two domains of $\R^3$ presented
in propositions \ref{p+1} and \ref{p+2} for the Stern-Gerlach instrument located on the left (resp. on the right). \\
We suppose that, at time $t\in]t_1,t_2[$, a couple of elementary singularities $(\varsigma_L,\varsigma_R)$ is such
as $\varsigma_L\in\omega_{1L}\cup\omega_{2L}$ ~ and ~ $\varsigma_R\in\omega_{1R}\cup\omega_{2R}$. \\
We denote $(+,+)$ the situation where $\varsigma_L\in\omega_{1L}$ and $\varsigma_R\in\omega_{1R}$, $(+,-)$ the situation
where $\varsigma_L\in\omega_{1L}$ and $\varsigma_R\in\omega_{2R}$, likewise for $(-,+)$ and $(-,-)$. \\

Then:
\begin{enumerate}
       \item The probabilities of obtaining $(+,+)$, $(+,-)$, $(-,+)$, $(-,-)$ are:
  \begin{eqnarray}\label{F+28}
   p(+,+)=p(-,-)=\frac{1}{2}\cos^2\frac{(\theta^L-\theta^R)-(\theta_1^L-\theta_1^R)}{2}
  \end{eqnarray}
    \begin{eqnarray}\label{F+29}
  p(+,-)=p(-,+)=\frac{1}{2}\sin^2\frac{(\theta^L-\theta^R)-(\theta_1^L-\theta_1^R)}{2}   
    \end{eqnarray}
    
  where $\theta_1^L$ (resp. $\theta_1^R$) is the spin state angle of the metric $g_1^L$ ~ (resp. $g_1^R$) defined
in the previous chapter. \\
These probabilities can also be written as follows:
 \begin{eqnarray}\label{F+30}
p(+,+)=p(-,-)=\frac{1}{2}|\langle\gamma_{\theta^L},\gamma_{\theta_1^L}\rangle\langle\gamma_{\theta^R},\gamma_{\theta_1^R
}\rangle+\langle\gamma_{\theta^L},\gamma_{\theta_1^L}^\bot\rangle\langle\gamma_{\theta^R},\gamma_{\theta_1^R}
^\bot\rangle|^2
 \end{eqnarray}
 \begin{eqnarray}\label{F+31}
 p(+,-)=p(-,+)=\frac{1}{2}|\langle\gamma_{\theta^L},\gamma_{\theta_1^L}\rangle\langle\gamma_{\theta^R}^\bot,\gamma_{
\theta_1^R
}\rangle+\langle\gamma_{\theta^L},\gamma_{\theta_1^L}^\bot\rangle\langle\gamma_{\theta^R}^\bot,\gamma_{\theta_1^R}
^\bot\rangle|^2 
 \end{eqnarray}
where $\gamma_{\theta^L}$ ~ (resp. $\gamma_{\theta^R}$) is the eigenvector of $\hat S_{\theta^L}$ ~ (resp.
$\hat S_{\theta^R}$) defined by (\ref{F+11}), ~~ $\gamma_{\theta_1^L}$ ~ (resp. $\gamma_{\theta_1^R}$) is
the spin state of the metric $g_1^L$ ~ (resp. $g_1^R$), and $\gamma_{\theta_1^L}^\bot=\gamma_{\theta_2^L}$ ~ (resp.
$\gamma_{\theta_1^R}^\bot=\gamma_{\theta_2^R}$) is the canonical orthogonal of $\gamma_{\theta_1^L}$ ~ (resp.
$\gamma_{\theta_1^R}$) ~ (def. \ref{d+6}).
\item If one denotes $(+)_L$ ~ (resp. $(+)_R$) the situation for which $\varsigma_L\in\omega_{1L}$ ~ (resp.
$\varsigma_R\in\omega_{1R}$) and $(-)_L$ ~ (resp. $(-)_R$) that for which $\varsigma_L\in\omega_{2L}$ ~ (resp.
$\varsigma_R\in\omega_{2R}$), then: \\ the probabilities of obtaining situations $(+)_L$ and $(-)_L$ (independently of 
situations $(+)_R$ and $(-)_R$) are equal to $\frac{1}{2}$. Similarly for those to obtain situations $(+)_R$ and
$(-)_R$.
\end{enumerate}
\end{prop}
Here we have obtained results identical to those obtained by classical quantum mechanics from the state
"entangled" given by (\ref{F+27}). \\

In the particular case that corresponds to the A.G.R experiment (\cite{aspe}), the spin state of the couple of
particles
is given, in classical quantum mechanics, by (\ref{F+26}). For us, in 
proposition \ref{p+3} this hypothesis is written $\theta_1^L=\theta_1^R$. \\ It follows that: \\
$p(+,+)=p(-,-)=\frac{1}{2}\cos^2\frac{(\theta^L-\theta^R)}{2}$ ~~ and ~~
$p(+,-)=p(-,+)=\frac{1}{2}\sin^2\frac{(\theta^L-\theta^R)}{2}$. \\
The correlation function ~ $E_{\theta^L,\theta^R}:=p(+,+)+p(-,-)-p(+,-)-p(-,+)$ therefore satisfies: \\
 $E_{\theta^L,\theta^R}=\cos^2\frac{(\theta^L-\theta^R)}{2}-\sin^2\frac{(\theta^L-\theta^R)}{2}
=\cos(\theta^L-\theta^R)$. \\

If we choose as angles measurement: \\
($\theta^L=0$ ~ or ~ $\theta'^L=\pi/2$), ~~~~ ($\theta^R=\pi/4$ ~ or ~ $\theta'^R=-\pi/4$),\\ the four measurements give: \\
$E_{0,\pi/4}=\frac{\sqrt2}{2}$, ~~~ $E_{0,-\pi/4}=\frac{\sqrt2}{2}$, ~~~ $E_{\pi/2,\pi/4}=\frac{\sqrt2}{2}$,
~~~ $E_{\pi/2,-\pi/4}=-\frac{\sqrt2}{2}$. \\

Quantity $S:=E_{0,\pi/4}+E_{0,-\pi/4}+E_{\pi/2,\pi/4}-E_{\pi/2,-\pi/4}$ is then $2\sqrt2$ which contradicts
Bell's inequality and corresponds to the experimental result \cite{aspe} (the theory presented here
is not an hidden variable theory). \\
\textbf{Proof of proposition \ref{p+3}} \\
We consider Example 1, the other examples are treated in the same way.
\begin{enumerate}
 \item It is assumed that, at time $t\in I_{2l}$ of the measurement, the double oscillating metric is $g_{\theta_{T_1}}$. \\
Proposition \ref{p+2} (corollary \ref{c+2}) says that the probability for $\varsigma_L$ to be seen in $\omega_{1L}$ is: \\
$\cos^2\frac{(\theta^L-\theta_1^L)-(\theta^L-\theta_1^L)}{2}=1$ and probability that it will be seen in $\omega_{2L}$
is zero. \\
Probability for $\varsigma_L$ to be seen in $\omega_{1R}$ is:
$\cos^2\frac{(\theta^L-\theta_1^L)-(\theta^R-\theta_1^R)}{2}$ and  probability that it will be seen in $\omega_{2R}$
is: $\sin^2\frac{(\theta^L-\theta_1^L)-(\theta^R-\theta_1^R)}{2}$.
Probability of obtaining the situation $(+,+)$ is therefore: $\cos^2\frac{(\theta^R-\theta_1^R)-(\theta^L-\theta_1^L)}{2}$,
that of obtaining $(+,-)$ is:\\
$\sin^2\frac{(\theta^R-\theta_1^R)-(\theta^L-\theta_1^L)}{2}$, ~~ probabilities of obtaining
$(-,+)$ and $(-,-)$ are null.
\item It is assumed that, at time $t\in I_{2l}$ of the measurement, the double oscillating metric is $g_{\theta_{T_2}}$. \\
Proposition \ref{p+2} shows here that the probability for $\varsigma_L$ to be seen in $\omega_{1L}$ is: \\
$\cos^2\frac{(\theta^L-\theta_1^L)-(\theta^L-\theta_2^L)}{2}=0$ since $\theta_2^L$ is the angle of the canonical orthogonal
 $\gamma_{\theta_1^L}^\bot$ and therefore satisfies, according to (\ref{F+20}): $\theta_2^L=\theta_1^L-\pi$. Probability
so that $\varsigma_L$ is seen in $\omega_{2L}$ is then equal to $1$. \\
Probability for $\varsigma_R$ to be seen in $\omega_{1R}$ is: \\
$\cos^2\frac{(\theta^R-\theta_1^R)-(\theta^L-\theta_2^L)}{2}=\sin^2\frac{(\theta^R-\theta_1^R)-(\theta^L-\theta_1^L)}{2
}$ and probability to be seen in $\omega_{2R}$
is: $\cos^2\frac{(\theta^R-\theta_1^R)-(\theta^L-\theta_1^L)}{2}$. \\
Probabilities of obtaining the situations $(+,+)$ and $(+,-)$ are therefore zero, probability of obtaining $(-,+)$ is:
$\sin^2\frac{(\theta^R-\theta_1^R)-(\theta^L-\theta_1^L)}{2}$, that of obtaining $(-,-)$ is:
$\cos^2\frac{(\theta^R-\theta_1^R)-(\theta^L-\theta_1^L)}{2}$.
\end{enumerate}
Since, by hypothesis, the presented situations are the only possible ones and are equiprobable, we deduce  
equalities (\ref{F+28}) and (\ref{F+29}). \\
Equalities (\ref{F+30}) and (\ref{F+31}) are quickly proven since:
$$\gamma_{\theta^L}=((\cos\frac{\theta}{2})\beta_1+(\sin\frac{\theta}{2})\beta_2))e^{ik_L}
~~~~\gamma_{\theta^L}^\bot=((\sin\frac{\theta}{2})\beta_1-(\cos\frac{\theta}{2})\beta_2))e^{ik_L}$$
$$\gamma_{\theta_1^L}=((\cos\frac{\theta_1}{2})\beta_1+(\sin\frac{\theta_1}{2})\beta_2))e^{ik_{L_1}}
~~~~\gamma_{\theta_1^L}^\bot=((\sin\frac{\theta_1}{2})\beta_1-(\cos\frac{\theta_1}{2})\beta_2))e^{ik_{L_1}}$$
where $k_L$ and $k_{L_1}$ are real numbers. \\
Likewise for $\gamma_{\theta^R}$, $\gamma_{\theta^R}^\bot$, $\gamma_{\theta_1^R}$, $\gamma_{\theta_1^R}^\bot$. \\
Part 2. of proposition \ref{p+3} is a quick consequence of the symmetry of the data.


\section [Generalization of potential type domains \\ The fine structure constant \\ The anomalous magnetic moment of the electron] {Generalization of potential type domain. \\ The fine structure constant. \\ The anomalous magnetic moment of the electron.} \label{++2,3}
  
As announced in section \ref{s1.4} we will now abandon some simplifying assumptions that
 was assumed during the presentation of potential type domains. Recall that the fundamental hypothesis
given in \ref{ss1.3} specifies that the endomorphism field $\leftidx{^e}{h}$ is nilpotent. This hypothesis allowed
greatly to simplify calculations associated to curvature since the "inverse" $g^{-1}$ of the metric tensor
was obtained quickly (see expression (\ref{F1.2+})), as its determinant $|g|$.
These two objects are those which, technically, greatly complicate all calculations. Results
obtained using this hypothesis (for example theorems \ref{2.1} and \ref{2.3}) show that this hypothesis constitutes a good approximation
for many experimental results, they are very close to those obtained by classical quantum physics.
However, some experimental results of a high precision like, for example, the measurement of the magnetic moment of
the electron, suggest to modify the definitions of potential type domains by refining them and by no longer assuming
necessarily the hypothesis of nilpotence. The purpose of these changes is to correctly describe the 
extremely accurate results which are currently only addressed by the Q.F.T. \\
The new hypotheses that we will propose for a finer notion of "electromagnetic potential domain" are, of course, inspired by that chosen in \ref{ss1.3}. The choice of the vector field $X_2$
(see proposition \ref{p1.3}) will be different. This vector field will not necessarily satisfy the properties of proposition \ref{p1.3}
since the field $\leftidx{^e}{h}$ will no longer be assumed nilpotent, it will still be defined on $S^1(\delta)\times W$, ie. it will not depend on the variables of $\Theta$.


\subsubsection{The new specific form of electromagnetic potential}
As before, the cell $\C$ is of the form $\C=\Theta\times S^1(\delta)\times S^3(\rho)\times V$
 (see section \ref{ss2.15+}) but here we assume that $V=S^1(r)\times Z$ and $g_0|_V=g|_{S^1(r)}\times g|_Z$ where $g|_{S^1(r)}$ is the standard Riemannian circle $S^1(r)$ of radius $r$  and $g|_Z$ a Riemannian metric on $Z$.\\
The metric $g_\mathcal{P}$ of the electromagnetic potential has the same form as the one given by proposition \ref{p1.3} (here $X$ replaces $X_2$):
\begin{eqnarray}\label{F0.++}
 g_\mathcal P=g_0+sym(\Upsilon^\flat\otimes X^\flat)
\end{eqnarray}
where (see section \ref{ss2.15+}):
\begin{eqnarray}\label{F1.++}
\Upsilon=A+\sum^3_{k=1}B^kL_k|_{S^3}
\end{eqnarray}
since we assume here that $\varrho=2$ which corresponds, in standard theories, to the case of fermions with spin $1/2$. \\
We now assume that the field $X$ is \textbf{timelike} (and not lightlike as was $X_2$) and that it represents a \textbf{fluid} defined on $S^1(\delta)\times S^1(r)$. We can draw parallels with the notion of fluid defined on $\R\times\R$ in standard relativity theory, but care must be taken here that "time u" is defined on the \textbf{compact} manifold $S^1(\delta)$. This requires that the velocity $v$ of a point of the fluid (represented by $X$) is, in the standard coordinates of $S^1(\delta)\times S^1(r)$, an \textbf{integer multiple} of $r/\delta$: $v=pr/\delta$ where $p\in\varmathbb{Z}$. Indeed, coordinates $(u,r(u))$ of a fluid point must satisfy $r(u+2\pi\delta)=r(u)$. \\
Field $X$ is written:
\[X=c_{1}\delta_u+c_{2}\delta_r\]
The components $c_{1}$ and $c_{2}$ are assumed constant and $c_{1}^2-c_{2}^2=1$ since the field $X$ is timelike. \\
The velocity $v=c_{2}/c_{1}$ can be positive or negative (remember that $S^1(\delta)$ and $S^1(r)$ have a canonical orientation (see section \ref{ss1.1.1})). \\
In particular, we deduce:
\[c_{1}=\pm(1-v^2)^{-1/2}~~ and~~~~ c_{2}=\pm v(1-v^2)^{-1/2}\]
We choose the simplest possible field $X$: we assume that $v$ is the minimum non-zero speed, ie. $v=r/\delta$ (a choise of an integer multiple of this last would not fundamentaly change anything).

Then:
\begin{eqnarray}\label{F2.++}
 X=(1-(r/\delta)^2)^{-1/2}\partial_u+(r/\delta)(1-(r/\delta)^2)^{-1/2}\partial_r
\end{eqnarray}
Of course, we can consider that the vector field $X$ is defined on the cell $\C$ but that it depends only on the variables of $S^1(\delta)\times S^1(r)$ and that, for any $x\in\C$, $X(x)$ is a vector of the subspace  $T_x(S^1(\delta)\times S^1(r))$ of $T_x(\C)$. \\
The field $X$ will be much more important than the field $X_2$ given by proposition \ref{p1.3}.
Since the metric $g_\mathcal{P}$ given by \ref{F0.++} will not characterize only the electromagnetic field given by $\Upsilon$, we will say later that $g_\mathcal{P}$ is the metric of the \textbf{X-electromagnetic potential}. \\
This choice of X (and therefore of g) allows, as we will see, to precisely describe the experimental results on the magnetic moment measurement of the electron and other fermions (proposition \ref{PROP}). This choice may seem artificial because it introduces a circle $S^1(r)$ in the structure of the compact manifold $V$. \textbf{However, such a choice should be considered as a first attempt for a description of $V$ and $X$. It will then be necessary to modify these choices and, perhaps, to modify the expression of the geometric type chosen to obtain a result even closer to the experimental result than that which we will obtain here}. In the following definition, we present the geometric type that corresponds, in particular, in standard quantum theory, to "electrons in a magnetic field".
\begin{dfn} \label{d2.++1}
 \textbf{A domain of type "fermions in an electromagnetic field"} is a domain of type "elementary oscillating metric of order 2 in a potential" (definition \ref{d2.9}) for which the pseudo-Riemannian metric $g_\mathcal{P}$ of the X-electromagnetic potential is given by (\ref{F0.++}) where $\varUpsilon$ satisfies (\ref{F1.++}). The metric $g$ is of the form: $|a|^{4/{n-2}}g_\mathcal{P}$ with $a=\phi\beta$ ~~ where ~~ $\beta\in E_V(\mu)$ ~~ and the function $\phi:\Theta\times S^1(\delta)\times S^3(\rho)\rightarrow\R$ is such that, for any $x\in\Theta$, $\phi_x(.)\in E_{S^1(\delta)}(\lambda)\otimes E_1$ . In addition, the potential $g_\mathcal{P}$ is neutral on $Z$ (Def. \ref{d2.4}).
\end{dfn}
\begin{rmq} \label{r1+++}
 A possible alternative is to use the sphere $S^3(\rho)$ instead of the circle $S^1(r)$. In particular the field $X$ is then defined on $S^1(\delta)\times S^3(\rho)$ which slightly modifies the given definitions since $S^3(\rho)$ is a parallelizable manifold (as is $S^1(r)$) and just replace "$r$" with "$\rho$" in ( \ref{F2.++}). However, some additional terms in the resulting equations are difficult to estimate. However, we can prove that, under reasonable assumptions, these terms remain negligibles when using $S^1(r)$. Choosing $S^3(\rho)$ in place of $S^1(r)$ would have the advantage of not introducing a new compact manifold to compose $W$. Conversely, the choice of $S^1(r)$ gives more possibilities for the construction of interesting oscillating metrics (remember that "$\rho$" is involved in the definition of the mass frequency, "$r$" is then a constant independent of "$\rho$"). In fact, a more elaborate (but "natural") choice of $X$ defined on other manifold than $S^1(\delta)\times S^1(r)$ or $S^1(\delta)\times S^3(\rho)$ may be necessary to obtain a result as precise as that given by Q.F.T. It would be very interesting to make a link between the calculation procedures used in Q.F.T (used to obtain the result on the magnetic moment of the electron) and the calculations necessary to obtain this same result with the geometric theory presented here with a good choise of $X$. However, the great technical complexity of the Q.F.T procedures (which require, among other things, the use of renormalization procedures which are difficult to justify mathematically) makes this project difficult to approach. The choice of an $X$ that satisfies \ref{F2.++} must not use more data than those requested by the quantum electrodynamics (Q.E.D), the data that we use in this first test is determined by the only choice of the radius $r$ of $S^1(r)$ since $\delta=Q^{-1}$ is determined by the electric charge of the electron. This choice may be compared to the choice of the interaction term in the Lagrangian of the Q.E.D.
\end{rmq}

\subsubsection{Physical Interpretation - The parallel with Q.F.T}
\begin{enumerate}
 \item In definition \ref{d2.++1}, the function $a$ characterizes the fact that one represents fermions with spin  1/2. 
\item The metric $g_\mathcal{P}$ of the X-electromagnetic potential characterizes the electromagnetic field by $\Upsilon$ and involves the "new data" $X$ defined on $S^1(\delta)\times S^1(r)$. The field $X$ can be compared to a field of the Q.F.T (here defined on $S^1(\delta)\times S^1(r)$ and not on $\Theta$) whose quantas correspond to "virtual" particles introduced by the Feynmann diagrams.
\item The choice of $X$, for which the speed of a point of the fluid is $v = r/\delta$ (the smallest possible, non-zero), can be compared to the choice of the interaction term in the Lagrangian of Q.E.D since, as we will see in the following lines, $r/\delta$ will be chosen so that $\pi r^2/\delta^2=\alpha$ where $\alpha$ is the fine structure constant ($\alpha\simeq1/137 $). The coupling constant in the Lagrangian of Q.E.D is then $2\pi r/\delta$.
\end{enumerate}

\textbf{The fine stucture constant takes for us a very precise physical sense, it is a very simple characteristic of the geometry of the compact manifold $V$ (for this first test $\alpha=\pi r^2/\delta^2$}).\\
Given the conceptual importance of these notions, we give the following definition (see also section \ref{ss1.1.1} for details on the notion of "physical constants" in the context of our theory).

\begin{dfn}
 \textbf{The structural constant}, written $\Xi$, is the quotient of the radius of the circle $S^1(r)$ by the radius of the circle $S^1(\delta)$ that characterises the cell $(\C,g_0)$ specified at the beginning of the section:
 \[\Xi=r/\delta\]
 (If one chose $S^3(\rho)$ in place of $S^1(r)$ (see remark \ref{r1+++}) one would define $\Xi=\rho/\delta$.)
\end{dfn}
 If the unique motivation to modify the definitions of electromagnetic potential type domains was to theoretically obtain a value very close to the experimental value of the magnetic moment of the electron, it would be much simpler to keep the specific form of the electromagnetic potential presented in section \ref{ss2.15+} and to consider that, in this one, the constant $\varrho$ is not exactly equal to 2 but to a value deduced from the experimental result. The factorization in equation (\ref{F49}) presented in subsection \ref{ssAv}, would not be longer possible. This last point may be compared with the fact that, in standard quantum theory, the Dirac theory imposes to take $\varrho=2$. In fact, the spectacular result of the Q.F.T on this subject is the precise description (particularly difficult to obtain) of a link between the fine structure constant $\alpha$ and the gyromagnetic constant obtained experimentally. It is the precision of this result obtained by the Q.F.T which suggests that, in the new theory that we present here, it is not the constant $\varrho$ which is different from 2 but the potential type domain which has a more elaborate expression than that given in section \ref{ss2.15+}. The first attempt suggested here is clarified by the choice of $g_\mathcal{P}$ and the result is stated in proposition \ref{PROP} which follows.

\subsubsection{The fundamental equation of the domain of type "electrons (or other fermions) in an electromagnetic field"}
In the framework of the simplifying assumptions of section \ref{s2.13}, where the field $X_2$ replaces $X$, the fundamental equation is given by (\ref{F49}) (where $\varrho=2$ and the considered spin is 1/2):
\begin{eqnarray}\label{eq2++}
\sum_{j=0}^3\varepsilon_j(i\frac{\partial}{\partial x^j}+Q^+\Upsilon^j)^2a_c+M^2a_c+2
Q^+\sum_{k=1}^3B^k\hat S_k(a_c)+{Q^+}^2(\rho/2)^2|B|^2a_c=0
\end{eqnarray}
Difference between different kind of fermions is made by the choice of the mass frequency $M$. \\
\textbf{We will assume in the following that the electric charge frequency $Q^+$ is the non-zero minimum frequency  $Q^+=1/\delta$ and that it corresponds to the domains of type "electrons"}. \\
The new fundamental equation is given by the following proposition, it is a simple consequence of the chosen geometrical type where $X$ satisfies \ref{F2.++}, it is expressed as a small perturbation of equation (\ref{eq2++}).

\begin{prop} \label{PROP}
 We consider a domain of type "fermions in an electromagnetic field" (definition \ref{d2.++1}). We no longer assume here the part 3 in hypothesis $H_{2E}$ of theorem \ref{2.3} that we replace by \ref{F2.++}.\\
 Then the canonical function $a_c$ (definition \ref{d2.15}) satisfies the following equation:
 \begin{multline}\label{F4.++}
 \sum_{j=0}^3\varepsilon_j(i\frac{\partial}{\partial x^j}+Q^+\Upsilon^j)^2a_c+M^2a_c+2(1-\Xi^2)^{-1/2}Q^+\sum_{k=1}^3B^k\hat S_k(a_c)+\\
{Q^+}^2(\rho/2)^2|B|^2a_c+(*)=0
 \end{multline}
where the disruptive term $(*)$ is specified in the lines that follow. \\
\end{prop}
\noindent \textbf{Proof}. We do not detail here the proof since it is similar to the proof of the theorem \ref{2.3} part 3, presented in section \ref{a3.7}, but we evaluate the additional terms which appear in equation (\ref{F4.++}), associated to the choice of the field $X$ instead of $X_2$. These are very difficult to give in detail because they come from new expressions, now complex, of the scalar curvature, the determinant and the "inverse" of the tensor $g_\mathcal{P}$. The "inverse" of this tensor can in fact be computed without great difficulties (see the following lemma) and makes it possible to express very precisely the "Landé factor" $2(1-\Xi^2)^{-1/2}$ which is the essential element of this section.
\begin{rmq}
 The inverse $g^{-1}$ of the metric can be obtened as the limit of a series: \\
 $g^{ij}=g_0^{ij}+h^{ij}+\dots+(-1)^ph_{k_1}^ih_{k_2}^{k_1}\dots h^{k_pj}+\dots$
\end{rmq}
\begin{lem}
 We recall that we defined $g_\mathcal{P}$ by:
 $g_{\mathcal{P}ij}=g_{0ij}+h_{ij}$ ~~ where ~~ $h_{ij}=\Upsilon_iX_j+X_i\Upsilon_j$ \\
 So:
 \begin{eqnarray}
  g^{ij}_\mathcal{P}=g_0^{ij}+a_X(\Upsilon^iX^j+X^i\Upsilon^j)+b_XX^iX^j+c_X\Upsilon^i\Upsilon^j
 \end{eqnarray}
where: \\
$a_X=-(1+X^i\Upsilon_i)/(1+X^i\Upsilon_i)^2-(X^iX_i)(\Upsilon^i\Upsilon_i)$ \\
$b_X=\Upsilon^i\Upsilon_i/(1+X^i\Upsilon_i)^2-(X^iX_i)(\Upsilon^i\Upsilon_i)$ \\
$c_X=X^iX_i/(1+X^i\Upsilon_i)^2-(X^iX_i)(\Upsilon^i\Upsilon_i)$
\end{lem}
This result is easily obtained, we can verify for example that $g_\mathcal{P}^{ij}g_{\mathcal{P}jk}=\delta_k^i$. \\
The important term $2(1-\Xi^2)^{-1/2}Q^+\sum_{k=1}^3B^k\hat S_k(a_c)$ of equation (\ref{F4.++}) comes from the calculation of the terms $|g_\mathcal{P}|^{-1/2}\partial_4(|g_\mathcal{P}|^{1/2}X^4\Upsilon^i\partial_i\phi)$ and $|g_\mathcal{P}|^{-1/2}\partial_i(|g_\mathcal{P}|^{1/2}X^4\Upsilon^i\partial_4\phi)$ (see annex \ref{a3.7}). However, since $g_\mathcal{P}$ does not depend on the variable "$u$", we have $\partial_4|g_\mathcal{P}|=0$. Then, the
term $2(1-\Xi^2)^{-1/2}Q^+\sum_{k=1}^3B^k\hat S_k(a_c)$ comes from $2X^4\Upsilon^i\partial_i\partial_4\phi$ where
 $X^4=(1-\Xi^2)^{-1/2}$ and it follows that the coefficient 2 which appear in the "simplified" equation \ref{F4.++} is now $2(1-\Xi^2)^{-1/2}$ which corresponds to the Landé factor of classical quantum physics denoted $l_e$ in the case of the electron. \\
Term $2(1-\Xi^2)^{-1/2}Q^+\sum_{k=1}^3B^k\hat S_k(a_c)$ is the preponderant term in which the magnetic field $B$ appear, the other terms that are ${Q^+}^2(\rho/2)^2|B|^2a_c$ and those which appear in $(*)$ will have a negligible influence relative to the order of precision obtained for the value of the Landé factor. This will only be satisfied under some conditions specified in the following lines. \\

\noindent \textbf{Magnitude orders of some fundamental constants} \\

The radius $\delta$ of the circle $S^1(\delta)$ impose the possible electric charge frequencies (see definition \ref{d2.6}). The smallest possible non-zero frequency is then $1/\delta$. \textbf{We suppose that this corresponds to the electric charge frequency of the electron which is written $Q_e$} ~~ ($Q_e=1/\delta$). \\
Then (definition \ref{d2.7}) $\delta=1/Q_e=\hbar/|q|$ \\
where $q$ is the electric charge (in Coulomb) of the electron. \\
Which gives: $\delta\simeq6,24149~~~10^{-16}m$ ~~ and ~~ $Q_e\simeq1,51927~~~10^{15}m^{-1}$ \\
The mass frequency of the electron is $M_e=m_ec/\hbar$ (definition \ref{d2.8}) where $m_e$ is the mass of the electron (in kg). \\
Which gives: $M_e\simeq2,58961~~~10^{12}m^{-1}$ \\
It follows: $Q_e/M_e\simeq5,86679~~~10^2$ \\
An estimate of the radius $r$ of $S^1(r)$ is associated to the choice of the vector field $X$, which gives  $r^2=\alpha\delta^2/\pi$ where $\alpha$ is the fine structure constant: $\alpha\simeq1/137$ (see end of this section). It is then reasonable to consider that $r$ is of the order of $10^{-17}m$. \\
An estimate of the radius $\rho$ of $S^3(\rho)$ is currently more difficult to obtain but it is not unreasonable to consider that it is of the same order of $r$ (see also remark \ref{r1+++}): $10^{-17}m$. \\

\noindent \textbf{Magnitude orders of disruptive terms} \\

Given the orders of magnitude given in the preceding paragraph, the coefficient 
$Q_e^{+2}(\rho/2)|B|^2$ is of the order of $10^{-29}M^2$ if $B$ is of the order of 1 Tesla. It is then absorbed by $M^2$ which comes from the term $M^2a_c$ of equation (\ref{F4.++}). \\
Calculation of the scalar curvature $S_{g_\mathcal{P}}$ (of the same type as for proposition \ref{p1.5}) introduces the term $\Delta_{g_\mathcal{P}}(\Upsilon^i\Upsilon_i)a_c$. As $\Upsilon^i\Upsilon_i=|A|^2+(\rho/2)^2|B|^2$ it will be assumed that $\Delta_{g_\mathcal{P}}(\Upsilon^i\Upsilon_i)<<M^2$ so that it has only negligible influence. It should be noted that these perturbations can be interpreted as a small change in the mass frequency $M$. \\
Other disruptive terms, derived from the difference between the determinant of $g_\mathcal{P}$ and that of $g_0$, also occur but may be considered negligible under reasonable assumptions about measurement experiments. \\

\noindent \textbf{The Landé factor obtained by our theory for the choice of the field $X$ given by \ref{F2.++}}\\

Since $X=(1-(r/\delta)^2)^{-1/2}\partial_u+(r/\delta)(1-(r/\delta)^2)^{-1/2}\partial_r$, the component on $S^1(\delta)$ is\\ $X^4=(1-\Xi^2)^{-1/2}$ and, as we have seen, the Landé factor $l_e$ is then: 
$$2(1-\Xi^2)^{-1/2}=2(1+\frac{1}{2}\Xi^2+\frac{3}{8}\Xi^4+\dots)$$
The spectacular result of Q.F.T is the link obtained between the Landé factor $l_e$ and the fine structure constant $\alpha$ which is written in the form of a sequence:
$$l_e=2(1+\frac{1}{2}\alpha/\pi+A_2(\alpha/\pi)+\dots+A_n(\alpha/\pi)^n+\dots)$$
The comparison between the sequence obtained by our theory and that of the T.Q.C suggests choosing $\Xi^2:=r^2/\delta^2=\alpha/\pi $ which determines the value of the radius $r$ of $S^1(r)$ since $\delta=Q^{-1}=\hbar/|q|$.\\
Taking into account the experimental value of $\alpha$, we obtain with our theory $l_e=2,002326$ while the experimental result at the same order is $l_e(exp)=2,002319$.
\begin{rmq}
 If we choose the field $X$ (defined on $S^1(\delta)\times S^1(r)$) of the form $X=\sum_{k=1}^m X_k$ where the vector fields $X_k$ are timelike  (ie. satisfy $g_0|_{S^1(\delta)\times S^1(r)}(X_k,X_k)=-1$), we can show that several choices of finite families $\{X_k\}$ can give a result \textbf{as close as desired} to the experimental result (this, with the constraints on the fluids  velocities $v_k$ associated with the fields $X_k$: $v_k=p_k r/\delta$ where $p_k\in\varmathbb{Z}$). Fields $X_k$ can be compared to T.Q.C fields whose quantas correspond to "virtual" particles introduced by the Feynmann diagrams.
 However, a choice of \textbf{several} fields $X_k$ would be very artificial if the unique objective is to give a result very close to the experimental result and I do not see any "natural procedure" which would justify a good choice of these fields, unlike the choice of the field X specified in \ref{F2.++} whose data is to be compared, as we have already said, with the data of the Lagrangian interaction term of the Q.E.D with its coupling constant.
\end{rmq}
The aim of this section was, first of all, to show that the physical theory presented in this paper makes it possible to deal precisely, in geometric terms, with the subtle phenomena that only the Q.F.T is currently addressing.
\section{Study of phenomena which, in classical physics,\\
correspond to interactions between particles}\label{4-S-2}
\hspace{10mm} \textbf{- We propose here a comparative study between our theory and the Q.F.T -}
\subsection{Our theory}
The geometric type chosen to represent particles in a vaccum (with possible interactions between them) is given by a domain with constant scalar curvature conform to $g_0$ (definition \ref{d2.1}). The metric $g$ defined on $\C=\Theta\times V$ is thus of the form $g=|a|^{4/{n-2}}g_0$ and $S_g=S_{g_0}$.\\
The function $a$ satisfies (see equation \ref{F0'}):
\begin{eqnarray}\label{4-F-7}
 \Box_{g_0}a+Sa=S|a|^{4/{n-2}}a
\end{eqnarray}
where $S=\frac{4(n-1)}{n-2}S_{g_0}$.\\
To be able to carry out a sufficiently detailled study of this equation we will limit the space-time domain $\Theta$ to the subdomain $]t_1~~t_2[\times B_L\subset\R\times\R^3$ where $B_L$ is the box $]${-L/2}$~~$L/2$[^3\subset\R^3$ \\
(L sufficiently large). This will allow to decompose the function $a_t(.)$ in Fourier series.\\
We start by looking at the simplest case for which, in the considered space-time, the particles have a mass $m_0$ but have neither electric charge nor spin, which correspond to real scalar fields of the Q.F.T. (In particular, there is no creation of other types of particles.)\\
The function $a:\Theta\times V\rightarrow\R$ is here considered to depend only on the variables $(t,\vec{x})$ of $\Theta$ and not on those of $V$.\\
For each $t\in]t_1~~t_2[$,~~~$a_t(\vec{x})=\sum_{j=0}^{\infty}\beta_j(t)e^{i\vec{k_j}.\vec{x}}+\overline{\beta_j(t)}e^{-i\vec{k_j}.\vec{x}}$~~~~ since $a=\bar{a}$\\
where $j\in J_L:=(\frac{2\pi}{L} \mathbb{Z})^3$,~~  $\vec{k_j}=(k_{1_j},k_{2_j},k_{3_j})$,~~ $\vec{x}=(x^1,x^2,x^3)$,~~ $\vec{k}.\vec{x}=\sum_{l=1}^3k_lx^l$.\\
We can then write, more interestingly, $a_t(\vec{x})$ in the form:
\begin{eqnarray}\label{4-F-3}
 a_t(\vec{x})=\sum_{j=0}^{\infty}\lambda_j(t)e^{i(M_jt+\vec{k_j}.\vec{x})}+\overline{\lambda_j(t)}e^{-i(M_jt+\vec{k_j}.\vec{x})}
\end{eqnarray}
by setting $\lambda_j(t):=e^{-iM_jt}\beta_j(t)$~~where~$M_j:=(M_0^2+|\vec k_j|^2)^{1/2}$\\
($M_0$ is the mass frequency of the type of particles).\\
Each term $ a_j(t,\vec{x}):=\lambda_j(t)e^{i(M_jt+\vec{k_j}.\vec{x})}+\overline{\lambda_j(t)}e^{-i(M_jt+\vec{k_j}.\vec{x})}$
of the series can then be seen as the function of an homogeneous oscillating metric $g=|a_j|^{4/{n-2}}g_0$ in the vaccum, with momentum $\vec{k_j}$ (of mass $m_0$, without electric charge and without spin) (cf \ref{4-ss-1}).\\
It is recalled that a momentum measuring instrument, as defined in section \ref{s2.16}, can be considered (to simplify) as "selecting", for each $t\in]t_1~~t_2[$ each of the terms $a_j(t,.)$ in the spectral decomposition of the function $a$, this on a domain $\Theta'(t)\times V$ of the same $g_0$-volume as $\Theta(t)\times V$ where $\Theta(t)\subset B_L$, on which $g=|a_j|^{4/{n-2}}g_0$.\\

We now calculate the probability of getting a number $m$ of singularities at time $t$ for a momentum $\vec{k_j}$, that is, when the measuring instrument has created a domain $\Theta'(t)\times V$ with similary $g_0$-volume that $\Theta(t)\times V$ on which $g=|a_j|^{4/{n-2}}g_0$ where  $ a_j(t,\vec{x}):=\lambda_j(t)e^{i(M_jt+\vec{k_j}.\vec{x})}+\overline{\lambda_j(t)}e^{-i(M_jt+\vec{k_j}.\vec{x})}$.\\
The $g_j$-volume of $\Theta'(t)\times V$ is $v_j(t)=\int_{\Theta(t)\times V}|a_j(t,x)|^2dv_{g_0}\simeq2|\lambda_j(t)|^2v_{0,t}$ where $v_{o,t}$ is the $g_0$-volume of $\Theta(t)\times V$ ~~~(see \ref{4-sss-1}).\\
Then, according to (\ref{F44'}) proven in subsection \ref{4-1}, the probability of obtaining $m$ singularities for a momentum $\vec{k_j}$ at time $t$ is:
$$ p_t(m,j)=\frac{(2Dv_{0,t}|\lambda_j|^2)^m}{m!}e^{-2Dv_{0,t}|\lambda_j|^2}$$
where $D$ is the "density of elementary singularities" (see \ref{4-1}).\\
Note that the $g$-volume of $\Theta(t)\times V$ is ~~$c_t=\int_{\Theta(t)\times V}|a|^2dv_{g_0}\simeq2\sum_{j=0}^\infty|\lambda_j(t)|^2v_{0,t}$.\\
We define $K_t=Dc_t~~(=2Dv_{0,t}\sum|\lambda_j|^2)$ and $\mu_j=|\lambda_j|/{(\sum_{j=0}^\infty |\lambda_j|^2)^{1/2}}$.\\
Then:\\
$p_t(m,j)=\frac{(K_t\mu_j^2)^m}{m!}e^{-K_t\mu_j^2}$ ~~~~(we have $\sum_{j=0}^\infty \mu_j^2=1$).\\
$K_t:=Dc_t$ would represent, in classical theory, the number of particles found in $\Theta(t)$ at time $t$. We suppose that $K_t$ does not depend on time $t$ (no creation and annihilation, only "shocks" of particles) and $K_t$ will be written $K$.\\
With our theory, for any momentum, the probability that no elementary singularity is "observed" by the measuring instrument is (for $m=0$):\\
$p_0=\prod_{j=0}^\infty e^{-K\mu_j^2}=e^{-K}$\\
More generaly, the probability that $m_1$ singularities are "observed" with a momentum $\vec{k_{j_1}}$, $m_2$ singularities with a momentum $\vec{k_{j_2}}$, .. , $m_l$ singularities with a momentum $\vec{k_{j_l}}$, and no other singularity observed by the measuring instrument is: (by ordering $j_1<j_2<..<j_l$)\\
$p_{j_1,..,j_l,m_1,..m_l}=\prod_{j\neq j_1,..,j_l}e^{-K\mu_j^2}\prod_{k=1}^l\frac{(K\mu_{j_k}^2)^{m_k}}{m_k!}e^{-K\mu_{j_k}^2}=\prod_{j=0}^\infty e^{-K\mu_j^2}\prod_{k=1}^l\frac{(K\mu_{j_k}^2)^{m_k}}{m_k!}$\\
Then:
\begin{eqnarray}\label{4-F-4}
 p_{j_1,..,j_l,m_1,..m_l}=e^{-K}\prod_{k=1}^l\frac{(K\mu_{j_k}^2)^{m_k}}{m_k!}
\end{eqnarray}
\subsection{Study of this same phenomenon by the Q.F.T}
In Q.F.T, the evolution over time of this physical system is given by the state function $\psi:\R\rightarrow\mathcal{F}$ where $\mathcal{F}$ is the Fock space of the system. Here $\mathcal{F}=\oplus_{k=0}^\infty\mathcal{H}_k$ ~~where $\mathcal{H}_k=\otimes_s^k\mathcal{H}_1$, ~~ $\otimes_s$ is the symetric tensor product, ~~ $\mathcal{H}_0=\mathbb C$, ~~ $\mathcal{H}_1$ is the space of one particle states (here a boson) whose basis is given by the family of functions $(e^{i\vec{k_j}.\vec{x}})_{j\in J_L}$~~ where 
$J_L=(\frac{2\pi}{L} \mathbb{Z})^3$.\\
For $t\in\R$, ~ $\psi(t)$ is written in the form:\\
$\psi(t)=\sum_J\gamma_J(t)\otimes_s^{m_1}e^{i\vec{k_{j_1}}.\vec{x}}\otimes_s^{m_2}e^{i\vec{k_{j_2}}.\vec{x}}..\otimes_s^{m_l}e^{i\vec{k_{j_l}}.\vec{x}}$\\
where $J$ is the multi-index $j_1,..,j_1,....,j_l,..,j_l$.\\
According to the axiomatic system of the Q.F.T, during a measurement at time $t$, the probability of obtaining $m_1$ particles with momentum $\vec{k_1}$, $m_2$ particles with momentum $\vec{k_2}$, etc. is:
\begin{eqnarray}\label{4-F-5}
 \|(\psi(t)/e_J)\|^2=|\gamma_J(t)|^2
\end{eqnarray}
where $e_J=\otimes_s^{m_1}e^{i\vec{k_{j_1}}.\vec{x}}\otimes_s^{m_2}e^{i\vec{k_{j_2}}.\vec{x}}..\otimes_s^{m_l}e^{i\vec{k_{j_l}}.\vec{x}}$.\\
The state function $\psi$ is a solution of the equation:
\begin{eqnarray}\label{4-F-6}
 i\hbar\frac{\partial\psi}{\partial t}=H(\psi)
\end{eqnarray}
where $H$ is the Hamiltonian of the system (with interaction) and is obtained in Q.F.T by the process of "standard quantization" from the Lagrangian density $\mathcal L$ for which:\\
$$\mathcal L(\Phi)=1/2(g_o^{ij}\partial_i\Phi\partial_j\Phi-m^2\Phi^2)+\lambda\Phi^4$$\\
where $g_0$ is the Minkovski metric, $m$ is the particle mass, $\lambda$ is a coupling constant (the term $\lambda \Phi^4$ is the term that specifies the interaction). Functions $\Phi$ are, for the physical system considered here, real scalar fields.
\subsection{Comparative study}
If we want that our results, given in terms of presence probabilities during a measurement (specified by our theorie by \ref{4-F-4}) coincides with those of the Q.F.T (given by \ref{4-F-5}) on this physical system, the link between the function $a$ (which caracterizes the physical system with our theory) and the state function $\Psi$ (of the Q.F.T) must obligatorily be the following one:\\
$a=(\sum_{j=0}^{\infty}\lambda_j(.)e^{i(M_j(.)+\vec{k_j}.\vec{(.)})}+\overline{\lambda_j(.)}e^{-i(M_j(.)+\vec{k_j}.\vec{(.)})})~~~\longmapsto~~~\Psi_a=\sum_J\gamma_Je_J$\\
where $|\gamma_J|^2=e^{-K_a}\prod_{k=1}^l\frac{(K_a\mu_{j_k}^2)^{m_k}}{m_k!}$\\
$J$ is the multi-index defined previously.\\
Here $K_a=K_{a,t}$ is supposed to do not depend on $t$ and is defined, as we have seen, by $K_a=Dc_t$ ~ where $c_t=vol_g(\Theta(t)\times V)$.\\
Since $c_t$ does not depend on $t$ (the density $D$ does not depend on it), the value of $K_a$ is determined by the "initial condition" specified for a $t<<0$ (if, for a $t<<0$, the function $a$ is multiplied by $k$ then $K_a$ is multiplied by $k^2$, ie., in classical vision, the number of particles "sent" is multiplied by $k^2$).\\
Let us note that, to a function $a$ of our theory corresponds (at least) a function $\Psi_a$ of the Q.F.T, but that the data of a function $\Psi:\R\rightarrow\mathcal{F}$ does not always come from a function $a$. The condition is very restrictive. Functions $\Psi$ that come from a function $a$ will be said to be "admissible by our theory".\\
The link between the function $a$ and $\Psi$ is obviously \textbf{not linear}.\\
Equation \ref{4-F-7} for $a$, given by our theory, is \textbf{not linear}.\\
Equation \ref{4-F-6} for $\Psi$ of the Q.F.T is \textbf{linear}, but the operator $H$ is given by formal series (or formal integrals if $B_L$ becomes $\R^3$) and is not correctly mathematically defined. Calculations lead to "infinites" and it is not known at this stage how to suppress them by principles of "renormalization". We will only be able to do this if we restrict the problem to "asymptotic" results calculations (we assume that the function $\Psi$ is known for $t<<0$ and we only look for information on $\Psi$ for $t>>0$) and, by using the perturbation theory, in this case, the principles of renormalization can work. A particulary interesting result would be to prove that, under standard initial conditions, the solutions $\Psi$ of the equation \ref{4-F-6}, admissible by our theory (hense of the form $\Psi_a$) correspond to functions $a$ which are solutions of equation \ref{4-F-7}. But, the complexity of the Q.F.T procedure (mainly associated to the principles of renormalization) makes this study particularly delicate. Note that the exponent chosen in Q.F.T in the Lagrangian density for the "interaction" term is the integer $4$. The Euler equation associated with the "action" of this Lagrangian density gives an equation close to the equation of curvature \ref{4-F-7} of our theory and is of the form $\Box_{g_0}\Phi+M^2\Phi=\lambda\Phi^3$~~(where $M^2$ corresponds well here to the scalar curvature $S$ taking into account the assumtions made). However, the exponent of the nonlinear term is here 3 and does not correspond to that of equation \ref{4-F-7} because $n$ is $>4$. It is recalled that the axiomatic system of the Q.F.T imposes the fact that the exponent of the interaction term is an even integer $>2$ and that the theory can not be renormalizable if $n>4$, which imposes the choice of the value $4$ (which gives the exponent $3$ in the Euler equation). Moreover, the term $|a|^{4/{n-2}}a$ is only imposed by the "geometric type"  chosen in our theory, which guarantees that the scalar curvature remains constant. Since, for $n>6$,~ $2n/{(n-2)}$ is not an integer, this exponent can not replace $4$ in Q.F.T (It is difficult to "quantify" $\Phi^k$ if $k$ is not integer).\\
\indent Equation \ref{4-F-7} of our theory is expressed, when one writes the function $a$ in the form \ref{4-F-3}, by a family of differential order $2$ equations, non linear, concerning the functions $\lambda_j$. This system is, of course, unmanageable, it will be useful to restrict the study, as in Q.F.T, by seeking information on the $\lambda_j(t)$ only when $t>>0$, this for initial conditions given for $t<<0$.\\
\indent An interesting point of view for our theory is to consider that the "good" exponent of the term corresponding to the interaction-Lagrangian in the Q.F.T, is not $4$ but $2n/{(n-2)}$ (non integer and $<4$) and that the quantization process appropriate for this exponent was not found. The renormalization processes may be only tricks that can rectify the error on this exponent. Of course, this remains speculative as long as actual results, confirmed by experiment, have not been found using equation \ref{4-F-7} with the term $|a|^{4/{n-2}}a$ in the context of our theory.


\section [Some remarks for the rest] {Some remarks on the future description of 
more complex phenomena} \label{s2.17}
In the second chapter of this paper we basicaly limited ourselves (except section \ref{4-S-2}) to the study of physical phenomena that are
involved in the generic experiments of the standard quantum mechanics (diffraction, Young's slits,
potential deviations, Stern-Gerlach experiment, quantum entanglement, etc.). With our point of view on physics we
considered only the domains of type "oscillating metric in a potential" (remaining within the framework of
the linear approximation (cf section \ref{s2.3})), and elementary singularities (cf section \ref{s2.11}) played 
 a secondary role by intervening only during "measurements" (position measurements or other (cf section \ref{s2.16})).
The more complex phenomena we will study are those that, in the language of
standard physics, involve "interactions" of particles with each other (in a more general case than that presented in section \ref{4-S-2}). \\
The following two facts will be called into question:
\begin{enumerate}
 \item The fact of having used the linear approximation and therefore equation \ref{F1} (while equation that considers a form of interaction between particles is the nonlinear equation \ref{F0'}).
\item The fact that singularities have only a negligible influence on oscillating metrics.

\end{enumerate}
We will present very briefly in the following some elements of study in which the previous two facts
 are not necessarily true, this by describing very specific examples only to
show the encountered difficulties.


\subsection{Where we abandon the linear approximation but we continue to neglect the possible influence of
singularities on the equations that govern metrics conform to a potential} \label{ssn2.2}
We describe an experiment which, in the language of standard physics, corresponds to the sending of two particle flows
thrown at opposite velocities (with modulus $\lambda$) and which meet at time $t=0$ on the surface
$x^1=0$.

For this we consider two domains in $\R^4$ written $\Theta_1$ and $\Theta_2$ defined by:

$\Theta_1=\{(t,x)\in \R^4 ~~~ /t<t_0<0,~~~ x^1<\lambda t<0,~~~ (x^2,x^3)\in\omega\subset\R^2\}$

$\Theta_2=\{(t,x)\in \R^4 ~~~ /t<t_0<0,~~~ x^1>-\lambda t>0,~~~ (x^2,x^3)\in\omega\subset\R^2\}$

It is assumed that the domain $\C_1:=\Theta_1\times S^1(\delta)\times S^3(\rho)\times V$ is a domain with
constant scalar curvature conform to $g_0$ (def. \ref{d2.1}) whose metric $g_1=|a_1|^{4/n-2}g_0$ is such that the
spectral decomposition of the function $a_1$ is of the form (cf subsection \ref{ssn1}): \\
$$a_{1}(x,.)=\sum_{i=1}^{\infty}\varphi_i(x)\alpha_i(.)$$ where $(\alpha_i)$ is an 
orthonormal hilbertian basis of eigenfunctions of $(S^1(\delta)\times S^3(\rho)\times V,g_0)$.

It is assumed that this spectral decomposition has a main term which represents an \textbf{homogeneous} elementary oscillating metric
 moving at a velocity $\vec{v}$ along the $x^1$-axis and $|\vec v|=\lambda$, the
other terms of the spectral decomposition (the "harmonics") being negligible relatively to the main term
(It should be noted that a function $a_1$ corresponding to an homogeneous elementary oscillating metric cannot satisfy
the non-linear equation \ref{F0'} since a solution of such an equation necessarily contains an infinity of
non-zero terms in the decomposition).

It is assumed that the domain $\C_2:=\Theta_2\times S^1(\delta)\times S^3(\rho)\times V$ satisfies the same properties with
$g_2=|a_2|^{4/n-2}g_0$, the velocity corresponding to the main term being equal to $-\vec v$ (it is not assumed that the main term necessarily represents an oscillating metric of the same type as that defined on $\C_1$). The
velocity $\vec v$ is chosen so that the "encounter" of the two oscillating metrics takes place for $t=0$. The
metric is assumed to be close to $g_0$ on the domain $\C_{t_0}-(\C_1 \cup \C_2)$ where $\C_{t_0}:=\{t<t_0\}\times\R^3\times
S^1(\delta)\times S^3(\rho)\times V$.

We then seek to determine the function $a$ (or at least to have information on it) defined on
$\C=\R^4\times S^1(\delta)\times S^3(\rho)\times V$, corresponding to a metric with constant scalar curvature such
that $g=|a|^{4/n-2}g_0$ (ie. which satisfies  equation \ref{F0'} with $g_\mathcal{P}=g_0$) this one being equal to
the function that has been specified on $\C_{t_0}$. This last point is the "boundary condition" imposed on the
solution of equation \ref{F0'}. (In fact, it may be necessary to use a more complex expression for the metric $g$ than that of the form $g=|a|^{4/{n-2}}g_0$).

This technical study is very difficult and we will probably limit ourselves the evaluation of the main terms of
the spectral decomposition of the function $a$. These terms may represent oscillating metrics of different types from those involved in the boundary condition for $t<t_0$. It is conceivable that terms appear
representing oscillating "life-limited" metrics that will be presented later in this section. We can also limit ourselves to the evaluation of the function $a$ only for large values of time $t$ (study to compare with that of the scattering matrix in Q.F.T).
An early comparative study with Q.F.T was proposed in the section \ref{4-S-2}.

\subsection{Singularities and their influence}
Until now we have not given precise expressions of metrics in a neighborhood of singularities (cf
\ref{s2.11}) because this was not necessary in the context of the phenomena we have studied. The behavior of
metrics in a neighborhood of singularities will become an important ingredient in the precise description of notions that,
in language of the standard physics, concern the "composite particles", "the nuclei of atoms", "the
atoms", etc. Of course, for us, these notions will remain within the framework of the theory already presented in terms of
"n-dimensional space-time deformation" and more specifically, in the study of particular domains with 
constant scalar curvature conform to a potential, but precise results can only be obtained
with a sufficiently detailed description of the asymptotic behavior of the metric in a neighborhood of a
singularity.

Examples of asymptotic behaviors have been presented in the
manuscript \cite{vaugon-2}, this for singularities of different "dimensions" (relative to the corresponding submanifolds). We will restrict ourselves here to some "naive" examples.

In the following the considered cell is of the form $\C=I\times\Omega\times W$ where $W=S^1(\delta)\times S^3(\rho)\times V$. We
write $(t,x^1,x^2,x^3)$ the coordinates of $I\times\Omega$, we assume that $0\in\Omega$ and we set
$r=(\sum_{k=1}^3(x^k)^2)^{1/2}$.\\
\textbf{Example 1}\\
Consider a metric with constant scalar curvature conform to a neutral potential $g=|a|^{4/n-2}g_0$ where the function
$a$ is of the following form:
\begin{eqnarray}\label{Fn69}
 a=\frac{c}{r}a_1+a_2
\end{eqnarray}
where $c$ is a constant.

It is assumed that the metrics $|a_1|^{4/n-2}g_0$ and $|a_2|^{4/n-2}g_0$ are those of an elementary oscillating metric
 (def. \ref{d2.3}) and therefore satisfy $\Box_{g_0}a_i+Sa_i=0$. When the function $a_1$ does not depend on
variables $(x^1,x^2,x^3)$ we deduce, since $\Box_{g_0}(\frac{1}{r})=\Delta(\frac{1}{r})=0$, that $\Box_{g_0}a+Sa=0$
on $\C-\mathcal{S}$ where $\mathcal{S}=I\times\{0\}\times W$.

For $(t,u)\in I\times S^1(\delta)$ the singularity
$\mathcal{S}(t,u)$ is therefore of the form $\mathcal{S}(t,u)=\{t\}\times\{0\}\times\{u\}\times S^3(\rho)\times V$.

The function $a$ is here that of an oscillating metric with a singularity (stationary in $0\in\Omega$).

It is assumed that functions $a_1$ and $a_2$ are of the same order of magnitude. As an example, and to make the link with the
standard physics, it is assumed that the oscillating metric $g:=|a|^{4/n-2}g_0$ is that associated with a proton
(stationary). Outside a domain very close to the singularity we have $a\simeq a_2$. The function $a_2$ is therefore
that which characterizes the oscillating metric "proton". When we "neglect" the singularity, it is the function
$a_2$ which carries, as we have seen, the important characteristics: mass, electric charge, spin, etc. and that one
used in the description of standard experiments (diffraction, etc.). The metric $g_1:=|\frac{c}{r}a_1|^{4/n-2}g_0$
becomes dominant only when $\frac{c}{r}>1$, so that the function $a_2$ becomes negligible compared with
$\frac{c}{r}a_1$. the function $a_1$ can therefore be considered as that which describes the \textbf{internal constitution of the
proton}. The spectral decomposition of the function $a_1$ can contain in its main terms the notion of quark and significant caracteristic quantities probabily comme from 
compact manifold $V$ in decomposition of $W$ in the form $S^1(\delta)\times S^3(\rho)\times V$ (it is possible
that the group $SU(3)$ is naturally associated with $V$ as $SU(1)$ was associated with $S^1(\delta)$ (diffeomorph) and
$SU(2)$ to $S^3(\rho)$). It should be noted that the metric $g_1$ is written $g_1:=|ca_1|^{4/n-2}g'_0$ where
$g'_0:=\frac{1}{r^{4/n-2}}g_0$ and that $g'_0$ is an "inverted" metric with respect to $g_0$ (when $n=6$, the
metric $g'_0=\frac{1}{r}g_0$ is, on $\R^{3*}$, the one that corresponds to the inversion
$\varphi:\R^{3*}\rightarrow\R^{3*}$ defined by $\varphi(x)=\frac{x}{r^2}$ where $x=(x^1,x^2,x^3)$ and one has
$\varphi^*g_0=g'_0$).

We have, in this example, used the linear approximation since the considered functions are solutions of
 equation \ref{F1}. When functions $a_1$ and $a_2$ are $<<1$ and $\frac{c}{r}$ is $<1$ this can be considered
 a good approximation of the general case that uses equation \ref{F0'}, however, this is no longer the case when
$\frac{c}{r}>>1$, ie. in the description of the internal constitution of the proton.

\begin{rmq}
 It can be considered that the term $\frac{c}{r}a_1$ in the expression of the function $a$ which describes the
singularity comes from a "collapse" of an oscillating metric characterized by the function $a_2$, but this
is pure speculation.
\end{rmq}
\noindent\textbf{Example 2 (Naive descrition of an atom)}\\
On the cell $\C$ we consider the metric with constant scalar curvature conform to a potential $g:=|a|^{4/n-2}g_0$
where the function $a$ is of the following form:
\begin{eqnarray}\label{Fn70}
 a=\frac{c}{r}(\varepsilon \cos(\lambda_1r)a_1+\sin(\lambda_2r)a_2)+a_3
\end{eqnarray}
where $c$, $\lambda_1$, $\lambda_2$ are constants. \\
The following conditions are imposed:

Functions $a_1$, $a_2$ and $a_3$ do not depend on the variables $x^1,x^2,x^3$ of $\Omega$, which makes the calculations
very simple (and gives in particular a spherical symmetry to the function $a$). Functions $a_1$ and $a_2$ satisfy
$\Box_{g_0}a_1=k_1a_1$ and $\Box_{g_0}a_2=k_2a_2$ (and are therefore not directly associated with elementary oscillating metrics).

The metric $|a_3|^{4/n-2}g_0$ is that of an elementary oscillating metric (and thus has the important characteristics
: mass, electric charge, spin, etc., as in the previous example). It satisfies equation \ref{F1}.

Constants are chosen so that $c\lambda_2<\varepsilon<<1$, functions $a_1$, $a_2$ and $a_3$ are
 "of the same order of magnitude" and are $<<1$.

We deduce: \\ - when $r<<c\varepsilon$, the terms in $a_2$ and $a_3$ are negligible compared with the term
$\frac{c\varepsilon}{r}\cos(\lambda_1r)a_1$. \\
-When $c\varepsilon<r<<c$, the terms in $a_1$ and $a_3$ are negligible compared with the term
$\frac{c}{r}\sin(\lambda_2)a_2$. \\
-When $r>>c$, the terms in $a_1$ and $a_2$ are negligible compared with the term
$a_3$. \\

Since $\Delta(\frac{1}{r}\cos(\lambda_1r))=\lambda_1^2(\frac{1}{r}\cos(\lambda_1r))$ and
$\Delta(\frac{1}{r}\sin(\lambda_2r))=\lambda_1^2(\frac{1}{r}\sin(\lambda_2r))$ ($\Delta$ "geometric"), when
$k_1+\lambda_1^2+S=k_2+\lambda_2^2+S=0$ the function $a$ satisfies $\Box_{g_0}a+Sa=0$ and so is that of an oscillating metric
 with singularity (stationary), just as in the first example, but here the singularity is more
complex. The term $\frac{c\varepsilon}{r}\cos(\lambda_1r)a_1$ characterizes the nucleus of an atom (predominant for
$r<<c\varepsilon$). The term $\frac{c\varepsilon}{r}\sin(\lambda_2r)a_2$ characterizes the electronic domain of the atom
(predominant for $c\varepsilon<r<<c$). The term $a_3$ characterizes the "oscillating metric" of the atom
when we neglect the singularity.

Function $\sin(\lambda_2r)$ specifies the "layered distribution" of the electronic domain of the atom, function
$\cos(\lambda_1r)$ that of the structure of the nucleus of the atom.

Of course, example of the atom that we have just presented is very "naive", in particular, spherical symmetry
 on $\Omega$ for the function $a$ that we assumed, is very restrictive, moreover, just as in
the previous example, we are in the context of the linear approximation (but this allowed us to give
"explicit formulas"). The use of nonlinear equation \ref{F0'} will certainly deeply modify
things, at least in a neighborhood of singularities, in other words, in the study of the structure of the nucleus of the atom.


\subsubsection{Particle Shocks}
Oscillating metrics with singularities that have been described in the two previous examples, for which 
functions $a$ satisfy equation \ref{F0'} or \ref{F1}, have a "stationary" singularity on $\Omega$
relative to the cell $\C=\Theta\times W$ where $W=S^1(\delta)\times S^3(\rho)\times V$.

We consider an "extended" Poincaré transformation 
where $\C'=\Theta'\times W$ which comes from a standard Poincaré transformation
$\Lambda:\Theta'\rightarrow\Theta$ (see remark \ref{r7}). Using the function $a\circ\sigma$ on $\C'$ the
oscillating metric type domain with a stationary singularity relative to $\C$ becomes a domain
of oscillating metric type with a singularity moving at a velocity $\vec v$ relative to $\C'$ if one has
chosen the Poincaré transformation $\Lambda$ which reflects the fact that an observer associated with $\C'$ moves at
velocity $-\vec v$ relative to an observer associated with $\C$.

We now summarize the experiment described at the beginning of this section (see \ref{ssn2.2}). We consider the two domains
$\C_1$ and $\C_2$ on which the metrics are respectively $g_1=|a_1|^{4/n-2}g_0$ and $g_2=|a_2|^{4/n-2}g_0$, but
we now take for $a_1$ and $a_2$ functions associated with oscillating metrics with a singularity
moving at a velocity $\vec v$ (respectively $-\vec v$) along the $x^1$-axis and $|\vec v|=\lambda$. This
choice of "position" of the singularities in $\R^3$ (not necessarily on the same
axis) is important. This constitutes boundary conditions of the considered experiment, which, in
 classical physics, corresponds to a "shock" of particles. For example, we can consider that the function
$a_1$ characterizes a proton (see example 1) moving at a velocity $\vec v$ and that the function $a_2$ characterizes an
anti-proton (for which, in particular, $Q$ becomes $-Q$) moving at a velocity $-\vec v$.

We then study the function $a$ defined on $\C$ corresponding to a metric with constant scalar curvature
 (which thus satisfies non-linear equation \ref{F0'} with eventually $g_\mathcal{P}=g_0$) such that the "limit" conditions
 are those which have just been specified. The difficulty of this study is still greater than that
presented in \ref{ssn2.2} since metrics with singularities are involved. It is conceivable that, in the spectral decomposition of the function $a$, appear "principal terms"
representing oscillating metrics that do not exist in boundary conditions (particle creation) and
some terms have a limited "life". These last notions are presented in the next subsection.


\subsection{Lifetime, creation, annihilation}
We start with the following remark: On a cell $\C=\Theta\times S^1(\delta)\times V$, we consider
the metric $g=|a|^{4/n-2}g_0$ where, for $t>0$ and when $K$ is a positive constant:
\begin{eqnarray}\label{Fn71}
 a=e^{-Kt}a_0 
\end{eqnarray}
It is assumed that the function $a_0$ defined on $\C$ does not depend on the variable $t$ and satisfies \\
$\Box_{g_0}a_0=\nu a_0$ (we can take, for example, $a_0=C \beta\cos(Qu-\sum_{k=1}^3x^k)$ where $\beta\in E_V(\mu)$
in which case $\nu=-Q^2+\sum_{k=1}^3x^k+\mu$). Then $\Box_{g_0}a+Sa=(K^2+\nu+S)a$.

When $K^2+\nu+S=0$, the function $a$ defined by \ref{Fn71} satisfies the same equation as that corresponding to an
elementary oscillating metric.

The important properties of this oscillating metric are the constants $Q$ and $\mu$ ($\mu$
which can be decomposed and introduce the notion of spin for example). The constant $K$ appears as a characteristic
analogous to the mass frequency (although here the term "frequency" is no longer appropriate) but with an
opposite sign.\\

Recall that the mass frequency $M>0$ has been defined (def. \ref{d2.8}) by $M^2=\mu+S-Q^2$ for elementary 
oscillating metrics because it corresponds exactly to the usual notion of mass as shown by the equations
 given by theorems \ref{2.1} and \ref{2.2}. But these theorems only make sense for oscillating metrics
for which mass frequency is well defined and do not concern metrics of the form
$g=|a|^{4/n-2}g_0$ where the function $a$ satisfies \ref{Fn71}. These metrics (with "pseudo-mass" $K$) can
possibly appear via important terms in the spectral decomposition of functions
solutions of equation \ref{F0'} that describe results of experiments of the type that we have just
presented. \nopagebreak\\

Note that if a measuring instrument creates a domain $\C=]t_0, t_0+T[\times B\times S^1(\delta)\times W$ (see section
\ref{s2.16}) such that the metric is of the form $g=|a|^{4/n-2}g_0$ where $a=Ce^{-Kt}\beta\cos(Qu-\sum_{k=1}^3x^k)$, the
probability that a singularity is at time $t$ in $\HH_t=\{t\}\times B\times W$ (after an "average"
on $u\in S^1(\delta)$) is of the form $C^2e^{-2Kt}$. This describes a phenomenon \textbf{of annihilation} of the "pseudo
oscillating metric "(whose characteristic constants are $Q$ and $\mu$) for which the average life time
is $1/2K$. (Note that there is no reason to assume that $K^2+\nu+S=0$ because $e^{-Kt}a_0$ is just a particular term
 of a spectral decomposition).

By analogy, one can consider metrics of the form $g=|a|^{4/n-2}g_0$ where, for $t>0$, $a=(1-e^{-Kt})a_0$. This
 describes a "pseudo oscillating metric" \textbf{creation} phenomenon.

Of course, examples of "pseudo oscillating metrics" that we have just presented are very particular and we
can define, in a more general context, "pseudo oscillating  metrics" which have analogous properties:
characteristic constants, lifetime (associated to the pseudo-mass), etc.

\bigskip

To be continued ...

\chapter{Annexes}


\section{Proof of Theorem \protect \thmunun \label{a3.1}}
\begin{enumerate}
 \item \textbf{We prove that an eigenspace $E$ of the endomorphism $\leftidx{^e}{G}|_{H_x}$ which is of dimension 1 and
 timelike, is necessarily unique in $H_x$}.

Suppose that $X_0$ and $X'_0$ are two timelike eigenvectors of $\leftidx{^e}{G}|_{H_x}$, with respective eigenvalues
 $\lambda$ and $\lambda'$. We have:

$G_x(X_0,X'_0)=\lambda g(X_0,X'_0)$

$G_x(X'_0,X_0)=\lambda' g(X'_0,X_0)$

hence ~ $\lambda=\lambda'$ because, according to the signature of $g|_{H_x}$, necessarily $g(X_0,X'_0)\neq0$. \\ $X'_0$ is therefore
an eigenvector for the eigenvalue $\lambda$. As the eigenspace is of dimension 1, $X'_0$ is proportional to
$X_0$.
\item Let us denote $(\sigma_s)_{s\in\R}$ the 1-parameter group of diffeomorphisms of the field
$Y$. \textbf{We prove that the field $X_0$, the functions $\mu$ and $\rho$ are invariant by $\sigma_s$}.

By definition, for any $s\in\R$, ~~~ $\sigma_{s_*}(Y)=Y$.

On the other hand, for any $x\in\D$, ~ for any $ Z\in T_x(W_x)$, \\
$\sigma_{s_{*x}}(Z)\in T_{\sigma_s(x)}(S^1_{\sigma_s(x)})\oplus T_{\sigma_s(x)}(W_{\sigma_s(x)})$.

This point can be proved, for example, by considering the Jacobian matrix of $\sigma_{s_*}$ associated to a chart aroud $x$    of the observation atlas.

It follows that, for any $s\in\R$, \\
$\sigma_{s_*}(T_x(S^1_x)\oplus T_x(W_x))= T_{\sigma_s(x)}(S^1_{\sigma_s(x)})\oplus
T_{\sigma_s(x)}(W_{\sigma_s(x)})$.

(The fact that $\sigma_s$ are isometries was not used to obtain this equality).

Since $\sigma_x$ are isometries, we deduce that:

for any $s\in\R$, ~~~ $\sigma_{s_*}(H_x)=H_{\sigma_s(x)}$.

Given the uniqueness of the eigenspace proven in 1., we have:

for any $s\in\R$, ~ for any $x\in\D$, ~~~ $\sigma_{s_x}(X_{0_x})=\pm X_{0_{\sigma_s(x)}}$.

We now consider a chart around $x$ of the observation atlas $(\mathscr V,\varphi)$.

For any $x'\in\mathscr V$, ~~~ $g(X_{0_x'}, \varphi^*(\frac{\partial}{\partial t}))<0$ ~ according to the selected orientation,
So: for any $s\in\R$, \\
(*) ~~~ $g(X_{0_{\sigma_s(x)}},\varphi^*(\frac{\partial}{\partial t}))=\pm
g(\sigma_{s_x}(X_{0_x},\varphi^*(\frac{\partial}{\partial t}))<0$

But for $s=0$ ~~ $\sigma_s=I_d$ ~~ and ~~ $g((X_{0_x},\varphi^*(\frac{\partial}{\partial t}))<0$.

By continuity in "$s$" of $g(\sigma_{s_x}(X_{0_x},\varphi^*(\frac{\partial}{\partial t}))$ which does not cancel, we
deduce that the sign in (*) is necessarily "$+$". Then, for any $s\in\R$, for any $x \in\D$,
~~~ $\sigma_{s_x}(X_{0_x})=X_{0_{\sigma_{s_x}}}$.

Since $G$, $Y$, $X_0$ are invariant by the $g$-isometries $\sigma_s$, it is now easy to deduce that
functions $\mu$ and $\rho$ are too, in other words: $Y(\mu)=Y(\rho)=0$.
\item \textbf{Results that we have just obtained allow us to write the following equalities}:

\begin{enumerate}
 \item $D_YX_0=D_{X_0}Y$, ~~ this since the field $X_0$ is invariant by the isometries $\sigma_s$ and therefore: $\mathscr
L_Y X_0=D_Y X_0-D_{X_0}Y=0$. We also deduce, since the torsion is zero, that $[X_0~~ Y]=0$.
\item (i) ~ $D_YY=0$

(ii) ~ $\nabla\cdotp  Y=0$

(iii) ~ $\leftidx{^e}{F}(X_0)=-2D_YX_0$

Indeed: $Y$ is a Killing field, so $\nabla_iY_j=-\nabla_jY_i$, especially
$\nabla_iY^i=\nabla\cdotp  Y=0$.

As $Y^iY_i=-1$, ~~ $0=\nabla_j(Y^iY_i)=2Y^i\nabla_jY_i=-2Y^i\nabla_iY_j$.

Moreover:
\begin{eqnarray}\label{F62}
 F_{ij}=\nabla_iY_j-\nabla_jY_i=2\nabla_iY_j=-2\nabla_jY_i 
\end{eqnarray}

In other words ~ $\leftidx{^e}{F}(X_0)=-2D_{X_0}Y=-2D_YX_0$ according to (a).
\end{enumerate}
\item \textbf{Proof of 1. in Theorem \protect \thmunun}.

$G^{ij}=\mu X^iX^j+\sigma Y^iY^j+P^{ij}=\mu X_0^iX_0^j+\rho(X_0^iY^j+Y^iX_0^j)+(\sigma+\rho^2/\mu)Y^iY^j+P^{ij}$.

So, according to Bianchi's second identity:
\begin{multline}\label{F63}
\nabla_iG^{ij}=X_0^j\nabla_i(\mu X_0^i)+\mu X_0^i\nabla_iX_0^j+\nabla_i(\rho(X_0^iY^j+Y^iX_0^j))+\nabla_iP^{ij}\\
=X^j\nabla_i(\mu X^i)+\mu X^i\nabla_iX^j+\nabla_iP^{ij}=0
\end{multline}
because $\nabla_i((\sigma+\rho^2/\mu)Y^iY^j)=\nabla_i(\sigma Y^iY^j)=0$, indeed:

$\sigma$ and $\sigma+\rho^2/\mu$ are invariant functions by the $\sigma_s$ of the 1-parameter group of
diffeomorphisms associated to $Y$, so $Y(\sigma+\rho^2/\mu)=0$, then:

$\nabla_i((\sigma+\rho^2/\mu)Y^iY^j)=Y(\sigma+\rho^2/\mu)Y^j+(\sigma+\rho^2/\mu)(Y^i\nabla_iY^j+Y^j\nabla_iY^i)=0$

according to (i) and (ii).

According to \ref{F63} and since $X_{0j}\nabla_iX_0^j=0$ because $X_{0j}X_0^j=-1$, we can write:

${X_0}_j \nabla_iG^{ij}=-\nabla_i(\mu X_0^i)+{X_0}_j\nabla_i(\rho(X_0^iY^j+Y^iX_0^j))+{X_0}_j\nabla_iP^{ij}=0$.

But ${X_0}_j\nabla_i(\rho(X_0^iY^j+Y^iX_0^j))=0$, this last point is easily obtained by developing and using
2., 3. (a), (i), (ii), the fact that $X_0^i{X_0}_i=Y^iY_i=-1$ and $X_o^iY_i=0$.

So:
\begin{eqnarray}\label{F64}
 \nabla_i(\mu X_0^i)=X_{0_j}\nabla_iP^{ij}
\end{eqnarray}

which is the first result sought since in addition:

$\nabla_i(\mu X^i)=\nabla_i(\mu X_0^i)+\nabla_i(\rho Y^i)=\nabla_i(\mu X_0^i)$.

On the other hand, reusing \ref{F63}:

$X_j\nabla_iG^{ij}=X_jX^j\nabla_i(\mu X^i)+\mu X^iX_j\nabla_iX^j+X_j\nabla_iP^{ij}=0$

but ~ $X^jX_j=-(1+\rho^2/\mu^2)$, ~ then ~ $X_j\nabla_iX^j=\frac{1}{2} \nabla_i(X^jX_j)=-\rho/\mu\nabla_i(\rho/\mu)$.

We can write:

$\rho X^i\nabla_i(\rho/\mu)=-(1+\rho^2/\mu^2)\nabla_i(\mu X^i)+X_j\nabla_iP^{ij}$

and according to \ref{F64}:

$\rho X^i\nabla_i(\rho/\mu)=-(1+\rho^2/\mu^2)X_{0_j}\nabla_iP^{ij}+X_{0_j}\nabla_iP^{ij}+(\rho/\mu )Y_j\nabla_iP^{ij}$

\hspace{1.8cm} $=-(\rho^2/\mu^2)X_{0_j}\nabla_iP^{ij}+(\rho/\mu) Y_j\nabla_iP^{ij}$

which gives the second result sought.
\item \textbf{Proof of 2. in Theorem \protect \thmunun}.

According to \ref{F63}: $Y_j\nabla_i(\rho(X_0^iY^j+Y^iX_0^j))+Y_j\nabla_iP^{ij}=0$ ~ since $Y_jX_0^j=0$.

But: $Y_j\nabla_i(\rho(X_0^iY^j+Y^iX_0^j))=-\nabla_i(\rho X_0^i)+Y_jY^i\nabla_i(\rho X_0^j)$ ~ according to (i) and (ii).

Moreover:

$Y_jY^i\nabla_i(\rho X_0^j)=\rho Y_jY^i\nabla_iX_0^j=\rho Y_jX_0^i\nabla_iY^j=0$ ~~ from 3. (a) and since
$Y_j\nabla_iY^j=0$

We thus obtain:

$\nabla_i(\rho X_0^i)=Y_j\nabla_iP^{ij}$.

Which is the first result sought knowing that:

$\nabla_i(\rho X^i)=\nabla_i(\rho X_0^i+(\rho^2/\mu) Y^i)=\nabla_i(\rho X_0^i)$ ~~ according to 2. and (ii).

On the other hand:

$\nabla_i(P^{ij}Y_j)=Y_j\nabla_iP^{ij}+P^{ij}\nabla_iY_j$.

But ~~ $P^{ij}\nabla_iY_j=\frac{1}{2} P^{ij}F_{ij}=0$ ~~ because $P$ is symmetric and $F$ antisymmetric.

Then ~~ $\nabla_i(P^{ij}Y_j)=Y_j\nabla_iP^{ij}$ ~~ in other words:
$\nabla\cdotp  (\leftidx{^e}{P}(Y))=g(Y,\nabla\cdotp  P)$.

\item \textbf{Proof of 3 in Theorem \protect \thmunun}.

According to \ref{F63} ~ $\mu X^i\nabla_iX^j=-\nabla_iP^{ij}-X^j\nabla_i(\mu X^i)$

In other words:

$\mu D_XX=-\nabla\cdotp  P-g(X_0,\nabla\cdotp  P)X$ ~~ since according to 1.,
$\nabla_i(\mu X^i)=g(X_0,\nabla\cdotp  P)$.

Which proves 3. (a).

On the other hand:

$D_XX=D_{X_0+(\rho/\mu)Y}(X_0+(\rho/\mu)Y)=D_{X_0}X_0+D_{X_0}((\rho/\mu)Y)+(\rho/\mu)D_YX_0$

$=D_{X_0}X_0+2(\rho/\mu)D_YX_0+X_0(\rho/\mu)Y$ ~~~ since $D_{X_0}Y=D_YX_0$, ~ $Y(\rho/\mu)=0$, ~ $D_YY=0$.

Since $X(\rho/\mu)=X_0(\rho/\mu)=\frac{1}{\mu}g(Y,\nabla\cdotp  P)-(\rho/\mu^2)g(X_0,\nabla\cdotp  P)$

we can deduce:

$\mu D_{X_0}X_0=\mu D_XX-2\rho D_YX_0-\mu X_0(\rho/\mu)Y$, ~~ then:

$\mu D_{X_0}X_0=-\nabla\cdotp  P-g(X_0,\nabla\cdotp  P)X-2\rho
D_YX_0-g(Y,\nabla\cdotp  P)Y+(\rho/\mu)g(X_0,\nabla\cdotp  P)Y$

\hspace{1.3cm} $=-\nabla\cdotp  P+\rho\leftidx{^e}{F}(X_0)-g(X_0,\nabla\cdotp  P)X_0-g(Y,\nabla\cdotp  P)Y$

\hspace{1.3cm} $=\rho\leftidx{^e}{F}(X_0)-pr_{\mathcal{T}^\bot}(\nabla\cdotp  P)$

Which proves 3. (b).
\item \textbf{Proof of 4 in Theorem \protect \thmunun}.

According to \ref{F62} $\nabla_iF^{ij}=-2\nabla_i\nabla^jY^i$.

But:
\begin{eqnarray}\label{F65}
 \nabla_i\nabla^jY^i=\nabla^j\nabla_iY^i+R^j_kY^k=R^j_kY^k	
\end{eqnarray}

since, according to (ii): $\nabla_iY^i=0$.

We can deduce:

$\nabla_iF^{ij}=-2R^j_kY^k$.

On the other hand, by definition: $2R^j_k=G^j_k+S_g\delta^j_k$

(where $S_g$ is the scalar curvature and $\delta^j_k$ the Krönecker symbols)

Using the fact that ~ $G^j_k=\mu X_0^jX_{0_k}+\rho(X_0^jY_k+Y^jX_{0_k})+(\sigma+\rho^2/\mu)Y^jY_k+P^j_k$ ~ we writte,
since ~ $X_{0_k}Y^k=0$:

$G^j_kY^k=-\rho X_0^j-(\sigma+\rho^2/\mu)Y^j+P^j_kY^k$.

We thus obtain:

$\nabla_iF^{ij}=\rho X_0^j+(\sigma+\rho^2/\mu-S_g)Y^j-P^j_kY^k$.

So according to \ref{F62}:

$2Y_j\nabla_i\nabla^jY^i=-Y_j\nabla_iF^{ij}=\sigma+\rho^2/\mu-S_g$ ~~ since ~ $P(Y,Y)=P^j_kY^kY_j=0$, ~ $Y^jY_j=-1$,
~ $X_0^jY_j=0$.

As $Y^j\nabla_jY^i=0$, ~ one can write:

$0=\nabla_i(Y^j\nabla_jY^i)=(\nabla_iY^j)(\nabla_jY^i)+Y^j\nabla_i\nabla_jY^i=-\frac{1}{4}F_{ij}F^{ij}
+Y^j\nabla_i\nabla_jY^i$.

So: $\sigma+\rho^2/\mu-S_g=\frac{1}{2}F_{ij}F^{ij}$ ~ and finally:

$\nabla_iF^{ij}=\rho X_0^j+\frac{1}{2}(F_{ij}F^{ij})Y^j-P^j_kY^k$.

Which prove 4.
\end{enumerate}

\newpage
\section [Proof of properties 1, 2 and 3] {Proof of properties 1, 2 and 3 on the metrics representing
active potential
(paragraph \ref{ss1.3} B.) and lemma \ref{l2} \label{a3.2}}
\begin{enumerate}
\item Proof of Properties 1. 2. 3.
\begin{enumerate}
 \item Symmetry of $\leftidx{^e}{h}_x$.
 
 $g_0(\leftidx{^e}{h}_x(X),Y)=g_{0_{ij}}h^i_kX^kY^j=h_{jk}X^kY^j=h_{kj}X^kY^j=g_{0_{ki}}h^i_jY^jX^k=g_0(X,\leftidx{^e}{h
}_x(Y))$.
\item When an endomorphism is nilpotent, there is always a basis in which its matrix is 
strictly superior triangular, therefore with zero trace. This is the case, regardless $q\in\N^*$, for $(\leftidx{^e}{h}_x)^q$.
\item For any $x\in\C$:

$(g_0)^{-1}_B(g)_B=(g_0)^{-1}_B(g_0)_B+(\leftidx{^e}{h})_B=I+(\leftidx{^e}{h})_B$

where $(~)_B$ means "the matrix of" in the basis $B$.

In a basis $B$ where the matrix of $\leftidx{^e}{h}_x$ is strictly superior triangular we have:

$det(g_0)^{-1}_B det(g)_B=det(I+(\leftidx{^e}{h}_x)_B)=1$

therefore ~ $det(g)_B=det(g_0)_B$.

If we denote by $P$ the passing matrix from basis $B'$ to the basis $B$, we obtain:

$(detP)^2det(g)_{B'}=(detP)^2det(g_0)_{B'}$

hence the result.
\end{enumerate}
\item Proof of lemma \ref{l2}

For $x\in\C$, let us denote ~ $A$ ~ the linear subspace of $T_x(\C)$ generated by
 $(Y$, $X_0)$, and $B$ ~ the linear subspace generated by
$(Y,\leftidx{^e}{h}(Y),\dots,\leftidx{^e}{h}^{p-1}(Y),X_0,\dots,\leftidx{^e}{h}^{p-1}(X_0))$. The metric $g_0$ is non-
degenerate by hypothesis since of signature $(-,+,+,+,-,+,\dots,+)$.

As $Y$ and $X_0$ are timelike and $g_0$-orthogonal, $g_0|_A$ is of signature $(-,-)$. ~~ $g_0|_
{A^\perp}$ is therefore of signature $(+,+,\dots,+)$ and is a scalar product.

Since ~ $B^\perp\subset A^\perp$, ~ $g_0|_{B^\perp}$ ~ is still a scalar product. As ~ $(\leftidx{^e}{h})^p=0$,
~~ $B$ is invariant by $\leftidx{^e}{h}$ and, since $\leftidx{^e}{h}$ is $g_0$-symmetrical, $B^\perp$ is also invariant
by $\leftidx{^e}{h}$.

$\leftidx{^e}{h}|_{B^\perp}$ is an endomorphism of $B^\perp$ ~ $g_0|_{B^\perp}$-symmetric and, as
~ $g_0|_{B^\perp}$ ~ is a scalar product, ~ $\leftidx{^e}{h}|_{B^\perp}$ is diagonalizable therefore identically zero
since nilpotent.

It has been proven that ~ $\leftidx{^e}{h}|_{B^\perp}=0$ ~. As ~ $T_x(\C)=B\oplus B^\perp$ ~ because ~ $B^\perp$ and (therefore) $B$
are $g_0$-regular, lemma \ref{l2} is proved.
\end{enumerate}

\newpage
\section{Proof of proposition \ref{p1.4} \label{a3.3}}
Here we use the notation used in the determination of geodesics in section \ref{ss1.4}.

In a standard coordinate system of the cell $\C$, denote:

$\Gamma^k_{ij}$ ~ (resp. $\tilde \Gamma^k_{ij})$ the Christoffel symbols of $g$ ~ (resp. $g_0$).

 Denote $T^k_{ij}$ the coordinates of the \textbf{tensor} ($\Gamma^j_{ij}-\tilde\Gamma^k_{ij})$.

We have, when $g=g_0+h$:
\begin{eqnarray}
 T^k_{ij}=\frac{1}{2}g^{kl}(\nabla_ih_{jl}+\nabla_jh_{il}-\nabla_lh_{ij}) \label{F65'}
\end{eqnarray}

\textbf{where the $\nabla_i$ are associated to $g_0$}.

Moreover, when we denote $R_{ij}$ (resp. $\tilde R_{ij}$) the Ricci curvature of $g$ (resp. $g_0$), we prove
easily that:
\begin{eqnarray}
 R_{ij}=\tilde R_{ij}+\nabla_kT^k_{ij}-T^k_{li}T^l_{kj} \label{F65*}
\end{eqnarray}

\textbf{Here $h=-2vX_1\otimes X_1$ and, in the following calculations of this proof, $X_1$ will be denoted $X$
to simplify writing}.

The expression of $T^k_{ij}$ was given by \ref{F1.5} but it is simpler here thanks to the hypothesis $DX=0$ (stronger
than the fact that $X$ is a Killing field), we obtain:
\begin{eqnarray}
 T^k_{ij}=-(X_iX^k\nabla_jv+X_jX^k\nabla_iv-X_iX_j\nabla^kv) \label{F66}
\end{eqnarray}

where properties $X^kX_k=0$ ~ and ~ $X(v)=0$ were also used.

\begin{enumerate}
 \item Calculation of $R_{ij}$.
 
With \ref{F66} we get:

$T^k_{li}T^l_{kj}=0$ ~~ and
~~ $\nabla_kT^k_{ij}=-(X_iX^k\nabla_k\nabla_jv+X_jX^k\nabla_k\nabla_iv-X_iX_j\nabla_k\nabla^kv)$

Hence, since $X(v)=0$:

$\nabla_kT^k_{ij}=-(\Delta_{g_0}v)X_iX_j$ ~~~ (where ~ $\Delta_{g_0}:=-\nabla^k\nabla_k$).

Which finally gives for the  Ricci curvature according to \ref{F65*}:

$R_{ij}=\tilde R_{ij}-(\Delta_{g_0}v)X_iX_j$.
\item Calculation of $S_g$.

$S_g=g^{ij}R_{ij}=(g_0^{ij}+2vX^iX^j)(\tilde R_{ij}-(\Delta_{g_0}v)X_iX_j)$.

By developing and since $X^kX_k=0$:

$S_g=S_{g_0}+2vR_{icc_{g_0}}(X_1,X_1)$.
\item $g(X_1,X_1)=g^{ij}X_iX_j=(g_0^{ij}+2vX^iX^j)X_iX_j=0$.

On the other hand:

$\nabla_{g,i}X^j=\partial_iX^j+\Gamma^j_{il}X^l=\nabla_iX^j+T^j_{il}X^l=T^j_{il}X^l$.

And using \ref{F66}:

$\nabla_{g,i}X^j=0$.
\item $\nabla_{g,i}Y^j=\partial_iY^j+\Gamma^j_{il}Y^l=\nabla_iY^j+T^j_{il}Y^l$.

$\nabla_iY^j=\partial_iY^j+\tilde\Gamma^j_{il}Y^l=\tilde\Gamma^j_{il}Y^l=\frac{1}{2}g_0^{mj}(\partial_ig_{0km}
+\partial_kg_{0im}-\partial_mg_{0ik})Y^k$.

But $Y^4=1$ and $Y^k=0$ for $k\neq4$, therefore, given the definition of $g_0$: \\

\vspace{-5mm}

$\nabla_iY^j=0$.

On the other hand according to \ref{F66}, since $Y(v)=0$ ~ and ~ $Y^jX_j=0$:

$T^j_{il}Y^l=0$.

Finally:

$\nabla_{g,i}Y^j=0$.
\end{enumerate}


\section{Proof of proposition \protect \propuncinq \label{a3.4}}
Expressions \ref{F65'} and \ref{F65*} are given in the previous section.

\textbf{Here $h=\Upsilon^\flat\otimes X^\flat_2+X^\flat_2\otimes\Upsilon^\flat$ and, in the following calculations of this
proof, $X_2$ will be simply denoted $X$ to simplify writing}.

Calculation of $T^k_{ij}$.
\begin{eqnarray}
 T^k_{ij}=\frac{1}{2}(g_0^{kl}-h^{kl}+h^{km}h^l_m)(*) \label{F67}
\end{eqnarray}

where we denote:

$(*)=\nabla_j(\Upsilon_iX_l+\Upsilon_lX_i)+\nabla_i(\Upsilon_jX_l+\Upsilon_lX_j)-\nabla_l(\Upsilon_iX_j+\Upsilon_jX_i)$
\bigskip

We have, since $DX=0$:

$(*)=X_l(\nabla_j\Upsilon_i+\nabla_i\Upsilon_j)+X_i(\nabla_j\Upsilon_l-\nabla_l\Upsilon_j)+X_j(\nabla_i\Upsilon_l-\nabla
_l\Upsilon_i)$.

That is, since $F=d\Upsilon^\flat$:

$(*)=X_l(\nabla_j\Upsilon_i+\nabla_i\Upsilon_j)+X_iF_{jl}+X_jF_{il}$

So:

$g_0^{kl}(*)=X^k(\nabla_j\Upsilon_i+\nabla_i\Upsilon_j)+X_iF_j^{~~k}+X_jF_i^{~~k}$.

$\nabla_j\Upsilon_i=\partial_j\Upsilon_i-\tilde\Gamma^l_{ij}\Upsilon_l=\partial_j\Upsilon_i-\frac{1}{2}g_0^{lm}
(\partial_ig_{0jm}+\partial_jg_{0im}-\partial_mg_{0ij})\Upsilon_l$.

It is deduced that $\nabla_j\Upsilon_i$ does not depend on the variables of $\Theta$ under the assumption $H_E$ and the variables of
$\Theta\times S^3(\rho)$ under the assumption $H'_E$. It is the same for $F_{ij}$ and one deduces therefrom
$X^lF_{il}=0$.

Since more $X^lX_l=0$, we have:

$h^{kl}(*)=(\Upsilon^kX^l+\Upsilon^lX^k)(*)=\Upsilon^lX^k(X_iF_{jl}+X_jF_{il})$.

And:

$h^{km}h^l_m(*)=\Upsilon^m\Upsilon_mX^kX^l(*)=0$.

So, according to \ref{F67}:
\begin{eqnarray}
 2T^k_{ij}=X^k(\nabla_j\Upsilon_i+\nabla_i\Upsilon_j)+X_iF_j^{~~k}+X_jF_i^{~~k}-\Upsilon^lX^k(X_iF_{jl}+X_jF_{il})
\label{F68}
\end{eqnarray}

\begin{enumerate}
 \item 1- Calculation of $R_{ij}$.
 
 According to the dependence on the variables deduced from the preceding lines, we have:
 
 $X^k\nabla_k(\nabla_j\Upsilon_i+\nabla_i\Upsilon_j)=0$ ~~ and ~~ $X^k\nabla_k(\Upsilon^l(X_iF_{jl}+X_jF_{il}))=0$.
 
 So:
 
 $\nabla_kT^k_{ij}=\frac{1}{2}(X_i\nabla_kF_j^{~~k}+X_j\nabla_kF_i^{~~k})$.
 
 On the other hand, according to \ref{F68}, developing and using the properties already mentioned:
 
 $4T^k_{li}T^l_{kj}=F_l^{~~k}F_k^{~~l}X_iX_j$.
 
 Hence:
 \begin{eqnarray}
  R_{ij}=\tilde R_{ij}+\frac{1}{2}(X_i\nabla_kF_j^{~~k}+X_j\nabla_kF_i^{~~k})-\frac{1}{4}F_l^{~~k}F_k^{~~l}X_iX_j
\label{F69}
 \end{eqnarray}
 
Which gives the result 1. of proposition \ref{p1.5} since $F_j^{~~k}=-F^k_{~~j}$.
\item 2- Calculation of $S_g$.

$S_g=g^{ij}R_{ij}=(g_0^{ij}-h^{ij}+h^i_kh^{kj})R_{ij}$.

According to \ref{F69} and properties already used:

$g_0^{ij}R_{ij}=S_{g_0}$.

On the other hand, since $g_0$ is a product metric and according to the assumptions $H_E$ and $H'_E$:

$h^{ij}R_{ij}=(\Upsilon^iX^j+\Upsilon^jX^i)R_{ij}=2\tilde R_{ij}\Upsilon^iX^j=0$.

$h^i_kh^{kj}R_{ij}=\Upsilon^k\Upsilon_kR_{ij}X^iX^j=\Upsilon^k\Upsilon_k\tilde R_{ij}X^iX^j$.

Finally:

$S_g=S_{g_0}+(\Upsilon^k\Upsilon_k)R_{icc_{g_0}}(X_2,X_2)$.
\item 3- $g(X_2,X_2)=g^{ij}X_iX_j=(g_0^{ij}-\Upsilon^iX^j-\Upsilon^jX^i+\Upsilon^k\Upsilon_kX^iX^j)X_iX_j=0$.

On the other hand:

$\nabla_{g,i}X^j=\partial_iX^j+\Gamma^j_{il}X^l=\nabla_iX^j+T^j_{il}X^l=T^j_{il}X^l$.

And using \ref{F68} and the properties already mentioned:
$\nabla_{g,i}X^j=0$.

\item 4-can be quickly deduced from the following lemma which gives a more general result.
\begin{lem} \label{al1}
 We consider a cell $\C=\Theta\times S^1(\delta)\times W$ (coordinate system
$(x^0,x^1,\dots,x^{n-1})$). Let $g$ the metric tensor transported on $\C$ from the tensor $g_\M$ and $Y$ the vector field  defined as before ($Y^i=0$ if $i\neq4$, $Y^4=1$). So:

$Y$ is a Killing field if and only if, in the coordinate system,
 $\partial_4g_{ij}=0$ for any $i$ and $j$  (that is, the matrix terms $(g_{ij})$ do not depend
on $x^4$).
\end{lem}

\textbf{Proof of the lemma}.

$Y$ is a Killing field if and only if $\nabla_iY_j+\nabla_jY_i=0$. We have:

$\nabla_iY_j=\partial_iY_j-\Gamma^k_{ij}Y_k$.

So:

 $\nabla_iY_j+\nabla_jY_i=\partial_iY_j+\partial_jY_i-2\Gamma^k_{ij}Y_k$
 
 In the chosen coordinate system, $Y_j=g_{ij}Y^i=g_{4j}$ ~ therefore:
 
 $\nabla_iY_j+\nabla_jY_i=\partial_i(g_{4j})+\partial_j(g_{4i})-2\Gamma^k_{ij}g_{4k}$, ~~ but:
 
 $\Gamma^k_{ij}g_{4j}=\frac{1}{2}\delta^l_4(\partial_ig_{jl}+\partial_jg_{il}-\partial_lg_{ij})=\frac{1}{2}
(\partial_ig_{4j}+\partial_jg_{4i}-\partial_4g_{ij})$.

So:

$\nabla_iY_j+\nabla_jY_i=\partial_4g_{ij}$

hence the result.
\end{enumerate}

\section [A very simple example of approximation] {A very simple example of approximation of solutions of a non linear  equation by those of a linear equation
 \label{ea3.4}}
The only interest of the example that we will present is to help understand the process of approximation of solutions
of the fundamental equation \ref{F0'} by solutions of linear equation \ref{F1} associated with domains of type
"Oscillating metric in a potential".

On the interval $[0 ~1]\subset\R$, we consider the differential equation:
$$(*) \text{~~~}y'-y=y^2$$
(Equation $(*)$ is here of order 1 so that calculations are very simple, equation \ref{F0'} is of course
of order 2).

 The solution of this differential equation which satisfies boundary condition  ~~ $y(0)=\varepsilon$ ~~ (which one
will then assume $\ll 1$) is, as can be quickly proved:
$$y(t)=\frac{\varepsilon e^t}{1+\varepsilon(1-e^t)}$$

We consider the associated linear equation:
$$(**) \text{~~~}y'-y=0$$
The solution $y_1$ of this equation which satisfies the same boundary condition $y_1(0)=\varepsilon$ is:
$$y_1(t)=\varepsilon e^t$$

As we can see, when $\varepsilon \ll 1$, this solution is "very close" to that of $(*)$ since
$y_1(t)/y(t)=1+\varepsilon(1-e^t)$ ~~ and ~ $-2<(1-e^t)<0$ ~ because we have limited $t$ to
the interval $[0 ~1]$.


\section{Proof of proposition \protect \propdeuxun \label{a3.5}}
\begin{enumerate}
 \item Let $\varphi:\Theta\times S^1(\delta)\times V_1\rightarrow\R$ ~~ and ~~ $\beta\in E_{V_2}(\mu_2)$.
 
 $\Box_{g_\mathcal P}(\varphi\beta)=-\nabla^i\nabla_i(\varphi\beta)=\beta~\Box_{g_\mathcal
P}(\varphi)+\varphi~\Box_{g_\mathcal P}(\beta)-2\nabla^i\varphi\nabla_i\beta$

where $\nabla_i$ are relative to $g_\mathcal P$.

$\nabla^i\varphi\nabla_i\beta=g_\mathcal
P^{ij}\nabla_j\varphi\nabla_i\beta=g_0^{ij}\nabla_j\varphi\nabla_i\beta+2v X_1^iX_1^j\nabla_j\varphi\nabla_i\beta=0$

since ~~ $g_0^{ij}\nabla_j\varphi\nabla_i\beta=0$ ~~ and ~~ $X_1^i\nabla_i\beta=X_1(\beta)=0$ ~~ taking into account
hypotheses.

On the other hand:
\begin{eqnarray}
 \Box_{g_\mathcal P}\beta=|g_0|^{-\frac{1}{2}}\partial_i(|g_0|^{\frac{1}{2}}g_\mathcal P^{ij}\partial_j\beta)
\label{F70}
\end{eqnarray}

since $|g_\mathcal P|=|g_0|$.

But ~ $g_\mathcal P^{ij}\partial_j\beta=g_0^{ij}\partial_j\beta+2vX_1^iX_1^j\partial_j\beta=g_0^{ij}\partial_j\beta$.

So:

$\Box_{g_\mathcal P}\beta=\Box_{g_0}\beta=\mu_2\beta$.

Finally:

$\Box_{g_\mathcal P}(\varphi\beta)=(\Box_{g_\mathcal P}(\varphi)+\mu_2\varphi)\beta$.

\item Here ~~
$\nabla^i\varphi\nabla_i\beta=g_0^{ij}
\nabla_j\varphi\nabla_i\beta-(\Upsilon^iX_2^j+\Upsilon^jX_2^i)\nabla_j\varphi\nabla_i\beta+\Upsilon^k\Upsilon_kX_2^iX_2^
j\nabla_j\varphi\nabla_i\beta$

but, given the assumptions:
\begin{eqnarray}
 \Upsilon^i\nabla_i\beta=\Upsilon(\beta)=0 ~~~\text{and} ~~~X_2^i\nabla_i\beta=X_2(\beta)=0 \label{F71}
\end{eqnarray}

hence:

$(\Upsilon^iX_2^j+\Upsilon^jX_2^i)\nabla_j\varphi\nabla_i\beta=0$.

So, as in 1.:

$\nabla^i\varphi\nabla_i\beta=0$.

on the other hand:

$g_\mathcal
P^{ij}\partial_j\beta=g_0^{ij}
\partial_j\beta-(\Upsilon^iX_2^j+\Upsilon^jX_2^i)\partial_j\beta+\Upsilon^k\Upsilon_kX_2^iX_2^
j\partial_j\beta=g_0^{ij}\partial_j\beta$.

Hence, according to \ref{F70} and \ref{F71}:

$\Box_{g_\mathcal P}\beta=\Box_{g_0}\beta=\mu_2\beta$.

Finally:

$\Box_{g_\mathcal P}(\varphi\beta)=(\Box_{g_\mathcal P}(\varphi)+\mu_2\varphi)\beta$.
\end{enumerate}


\section{Proof of Theorem \ref{2.1} \label{a3.6}}
\begin{enumerate}
 \item \textbf{In a neutral potential}.
 
 The function $a$ satisfies: $a=\varphi\beta$ ~~ where ~~ $\varphi:\Theta\times S^1(\delta)\rightarrow\R$ ~~ and ~~ $\beta\in
E_W(\mu)$.

It is known that: $\Box_{g_0}a+Sa=0$ ~~ where ~~ $S=\frac{n-2}{4(n-1)}S_{g_0}$.

So:

$\beta(\Box_{\Theta\times S^1}\varphi+\mu\varphi+S\varphi)=0$ ~~ where ~~ $\Box_{\Theta\times
S^1}=\frac{\partial^2}{\partial t^2}+\frac{\partial^2}{\partial u^2}-\sum_{j=1}^3\frac{\partial^2}{\partial (x^j)^2}$

Then:

$\Box_\Theta\varphi+(\mu+S-\lambda)\varphi=0$ ~~ where ~~ $\Box_\Theta=\frac{\partial^2}{\partial
t^2}-\sum_{j=1}^3\frac{\partial^2}{\partial (x^j)^2}$,

in other words:

$\Box_\Theta\varphi+M^2\varphi=0$.

 So, for any $x\in \Theta$, ~~ $\CC_\lambda((\Box_\Theta\varphi+M^2\varphi)_x(.))=0$
 
 where ~ $\CC_\lambda$ ~ is the isomorphism defined in section \ref{s2.8}.
 
 Since $\CC_\lambda((\Box_\Theta\varphi)_x(.))=\Box_\Theta(\CC_\lambda(\varphi_x(.)))=\Box_\Theta a_c$
 we can deduce:
 
 $\Box_\Theta a_c+M^2a_c=0$.
 \item \textbf{In a potential without electromagnetism}.
 
 The cell is $\C=\Theta\times S^1(\delta)\times W$ ~~ where ~~ $\Theta=I\times \mathcal U\subset\R\times\R^3$,
~~ $a=\varphi\beta$ ~~ with ~~ $\varphi:\Theta\times S^1(\delta)\rightarrow\R$ ~~ and ~~ $\beta\in E_W(\mu)$, ~ but
now:

$\Box_{g_\mathcal P}a+Sa=0$ ~~ where ~~ $g_\mathcal P=g_0+h$ ~~ and ~~ $h=-2vX_1\otimes X_1$.

(In the following calculations $g_\mathcal P$ will be denoted $g$ to simplify the writing).

In a standard coordinate system of the cell $\C=\Theta\times S^1(\delta)\times W$, we have:

$\Box_ga=-|g|^{-\frac{1}{2}}\partial_i(g^{ij}|g|^{\frac{1}{2}}\partial_ja)$ ~~ where ~~ $|g|:=detg=detg_0$ ~ (cf
\ref{ss1.3}) \\ and ~~ $g^{ij}=g_0^{ij}+2vX^iX^j$.

So:

$\Box_ga=\Box_{g_0}a-2|g|^{-\frac{1}{2}}\partial_i(v|g|^{\frac{1}{2}}X^iX^j\partial_ja)$.

According to the hypotheses of theorem: $X^1=X^2=X^3=0$ ~~ and $v$ only depends on the variables $x^1,x^2,x^3$.

Therefore:

$|g|^{-\frac{1}{2}}\partial_i(v|g|^{\frac{1}{2}}X^iX^j\partial_ja)=v|g|^{-\frac{1}{2}}\partial_i(|g|^{\frac{1}{2}}
X^iX(a))$.

But:

$|g|^{-\frac{1}{2}}\partial_i(|g|^{\frac{1}{2}}X^iX(a))=|g|^{-\frac{1}{2}}(\partial_i(|g|^{\frac{1}{2}}
X^i))X(a)+X^i\partial_i(X(a))$

\hspace{3.2cm} $=(\nabla\cdotp  X)X(a)+X(X(a))$.

Moreover, since $X(\beta)=0$, ~ $X$ ~ is tangent to $\R\times W$ ~ and ~ $X^0=-1$ it follows that:

$X(a)=\beta X(\varphi)+\varphi X(\beta)=-\beta\frac{\partial\varphi}{\partial t}$ ~~ and
~~ $X(X(a))=\beta\frac{\partial^2\varphi}{\partial t^2}$.

So, since $DX=0$:

$|g|^{-\frac{1}{2}}\partial_i(v|g|^{\frac{1}{2}}X^iX^j\partial_ja)=\beta(\frac{\partial^2\varphi}{\partial
t^2}-(\nabla\cdotp  X)\frac{\partial\varphi}{\partial t})=\beta\frac{\partial^2\varphi}{\partial t^2}$.

Hence, as ~ $\Box_{g_0}a+Sa=\beta(\Box_\Theta\varphi+M^2\varphi)$ ~~ (as in 1.):

$0=\Box_{g}a+Sa=\beta(\Box_\Theta\varphi+M^2\varphi-2v\frac{\partial^2\varphi}{\partial t^2})$.

Taking the image by the isomorphism $\CC_\lambda$, we obtain equation \ref{F16}.

\item \textbf{In an electromagnetic potential}.

Function $a$ is written: $a=\varphi\beta$ ~~ where ~~ $\varphi:\Theta\times S^1(\delta)\rightarrow\R$ ~~ and
~~ $\beta\in E_W(\mu)$.

We have:
\begin{eqnarray}
 \Box_{g_\mathcal P}a+Sa=0 \text{~~where ~~}g_\mathcal P=g_0+h\text{~~and
~~}h=\Upsilon^\flat\otimes X_2^\flat+X_2^\flat\otimes\Upsilon^\flat \label{F72}
\end{eqnarray}

(In the following calculations $g_\mathcal P$ will be denoted $g$ to simplify writing, $\Upsilon$ will be denoted $A$ and
corresponds to the vector field "electromagnetic potential" defined on $\Theta$ and $X_2$ will be denoted simply
$X$).

The vector field $A$ is tangent to $\Theta$ and depends only on $\Theta$ variables.

In a standard coordinate system of the cell $\C=\Theta\times S^1(\delta)\times W$, we have:

$\Box_ga=-|g|^{-\frac{1}{2}}\partial_i(g^{ij}|g|^{\frac{1}{2}}\partial_ja)$ ~~ where ~~ $|g|:=detg=detg_0$

and $g^{ij}=g_0^{ij}-h^{ij}+h^i_kh^{kj}$ ~~ with ~~ $h^{ij}=A^iX^j+A^jX^i$.

Given assumptions:

~~ $A^i=0$ ~ if ~ $i>3$, ~ $X^j=0$ ~ if ~ $j<4$ ~~ and ~ $X^4=-1$ ~ since ~ $g_0(X_2,Y)=1$

So:
\begin{eqnarray}
 \Box_ga+Sa=\Box_{g_0}a+Sa+(*)+(**) \label{F73}
\end{eqnarray}

where ~ $(*):=|g_0|^{-\frac{1}{2}}\partial_i(|g_0|^{\frac{1}{2}}h^{ij}\partial_ja)$ ~ and
~ $(**):=-|g_0|^{-\frac{1}{2}}\partial_i(|g_0|^{\frac{1}{2}}h^i_kh^{kj}\partial_ja)$.

Since $g_0$ is a product metric on $\Theta\times S^1(\delta)\times W$ we have (see 1. and 2.):
\begin{eqnarray}
 \Box_{g_0}a+Sa=\beta(\Box_\Theta\varphi+M^2\varphi) \label{F74}
\end{eqnarray}

\begin{enumerate}
 \item Study of (*).
 
 $(*)=(*_1)+(*_2)$ ~~ where ~~ $(*_1):=|g_0|^{-\frac{1}{2}}\partial_i(|g_0|^{\frac{1}{2}}h^{ij}\varphi\partial_j\beta)$
 
 and ~~ $(*_2):=|g_0|^{-\frac{1}{2}}\partial_i(|g_0|^{\frac{1}{2}}h^{ij}\beta\partial_j\varphi)$.
 
 $(*_1)=|g_0|^{-\frac{1}{2}}(\partial_i(|g_0|^{\frac{1}{2}}A^iX^j\varphi\partial_j\beta)+
\partial_j(|g_0|^{\frac{1}{2}}A^iX^j\varphi\partial_i\beta))=0$

because $X^j\partial_j\beta=0$ ~~ and ~~ $A^i\partial_i\beta=0$ ~ according to (i) since $\partial_i\beta=0$ ~ if ~ $i\leq4$.

On the other hand:

$(*_2):=(*_{2_1})+(*_{2_2})$ ~~ where
~~ $(*_{2_1})=|g_0|^{-\frac{1}{2}}\partial_i(|g_0|^{\frac{1}{2}}A^iX^j\beta\partial_j\varphi)$

and ~~ $(*_{2_2})=|g_0|^{-\frac{1}{2}}\partial_j(|g_0|^{\frac{1}{2}}A^iX^j\beta\partial_i\varphi)$.

$(*_{2_1})=-|g_0|^{-\frac{1}{2}}\partial_i(|g_0|^{\frac{1}{2}}A^i\beta\partial_4\varphi)$ ~~ since
$\partial_j\varphi=0$ if $j>4$ ~~ and ~~ $X^4=-1$ ~~ ($\partial_4\varphi:=\frac{\partial\varphi}{\partial u}$).

So, since $A^i=0$ ~ for ~ $i\geqslant4$:

$(*_{2_1})=-\beta(|g_0|^{-\frac{1}{2}}(\partial_i(|g_0|^{\frac{1}{2}}
A^i))\partial_4\varphi+A^i\partial_i\partial_4\varphi)$

\hspace{8mm} $=	-\beta((\nabla_{g_0}\cdotp  A)\partial_4\varphi+A(\partial_4\varphi))$.

On the other hand, since $X^4=-1$:

$(*_{2_2})=-|g_0|^{-\frac{1}{2}}\partial_4(|g_0|^{\frac{1}{2}}\beta
A(\varphi))+\sum_{j>4}|g_0|^{-\frac{1}{2}}\partial_j(|g_0|^{\frac{1}{2}}X^j\beta A(\varphi))$.

But, since $|g_0|$, $\beta$, and $A^i$ do not depend on $x^4=u$:

$(*_{2_2})=-\beta A(\partial_4\varphi)$.

Finally:
\begin{eqnarray}
 (*)=(*_1)+(*_{2_1})+(*_{2_2})=-\beta((\nabla_{g_0}\cdotp  A)\partial_4\varphi+2A(\partial_4\varphi)) \label{F75}
\end{eqnarray}

\item Study of $(**)$.

$(**)=-|g_0|^{-\frac{1}{2}}\partial_i(|g_0|^{\frac{1}{2}}h^i_kh^{kj}\partial_ja)$

\hspace{7mm} $=-|g_0|^{-\frac{1}{2}}\partial_i(|g_0|^{\frac{1}{2}}(A^iX_k+A_kX^i)(A^kX^j+A^jX^k)\partial_ja)$

So, since $A^kX_k=0$ ~~ and ~~ $X^kX_k=0$:

$(**)=-|g_0|^{-\frac{1}{2}}\partial_i(|g_0|^{\frac{1}{2}}X^iX^jA_kA^k\partial_j(\varphi\beta))$

But, since $A_kA^k$ does not depend on the variables of $S^1\times W$:

$(**)=-A_kA^k|g_0|^{-\frac{1}{2}}(\partial_i(|g_0|^{\frac{1}{2}}X^iX^j\beta\partial_j\varphi)+\partial_i(|g_0|^{\frac{1
}{2}}X^iX^j\varphi\partial_j\beta))$.

Since $X^j\partial_j\varphi=-\partial_4\varphi$ since $\partial_j\varphi=0$ for $j>4$ ~ and $X(\beta)=0$, we obtain:

$(**)=A_kA^k|g_0|^{-\frac{1}{2}}\partial_i(|g_0|^{\frac{1}{2}}X^i\beta\partial_4\varphi)$.

We consider that $i=4$ then $i>4$, we have:

$(**)=-A_kA^k(\beta\partial^2_4\varphi-|g_0|^{-\frac{1}{2}}(\sum_{i>4}\partial_i(|g_0|^{\frac{1}{2}}
X^i\beta))\partial_4\varphi)$.

What we can write:

$(**)=-A_kA^k(\beta\partial^2_4\varphi-(\nabla_{g_0}\cdotp  (\beta X))\partial_4\varphi)$.

And, since $\nabla_{g_0}\cdotp  (\beta X))=X(\beta)+\beta\nabla_{g_0}\cdotp  X=0$:
\begin{eqnarray}
 (**)=-A_kA^k\beta\partial^2_4\varphi \label{F76}
\end{eqnarray}

\item End of the proof of 3.

According to \ref{F72}, \ref{F73}, \ref{F74}, \ref{F75} and \ref{F76} we have:

$0=\Box_ga+Sa=(\Box_\Theta\varphi+M^2\varphi-(\nabla_{g_0}\cdotp  
A)\partial_4\varphi-2A(\partial_4\varphi)-A_kA^k\partial^2_4\varphi)\beta$.

This equation is also in the form:
\begin{eqnarray}
 -\sum_{j=0}^3\varepsilon_j(\frac{\partial}{\partial x^j}+A^j\frac{\partial}{\partial u})^2\varphi+M^2\varphi=0
\label{F77}
\end{eqnarray}

where ~~ $\varepsilon_0=-1$, ~~ $\varepsilon_1=\varepsilon_2=\varepsilon_3=1$, ~~ $\frac{\partial}{\partial u}=\partial_4$.

To obtain equation given in the theorem just composes the two hand sides of the equation \ref{F77}
by the isomorphism $\CC_\lambda$. We easily verify that:

$\CC_\lambda\circ((\Box_\Theta\varphi)_x(.))=(\Box_\Theta
a_c)(x)$ ~~~ $\CC_\lambda\circ((\frac{\partial \varphi}{\partial u})_x(.))=-iQ^+a_c(x)$

$\CC_\lambda\circ((\frac{\partial^2 \varphi}{\partial u^2})_x(.))=-(Q^+)^2a_c(x)$.

The equation obtained is then:

$\sum_{j=0}^3\varepsilon_j(i\frac{\partial}{\partial x^j}+Q^+A^j)^2a_c+M^2a_c=0$.
\end{enumerate}
\end{enumerate}


\section [The sphere $S^3$] {The sphere $S^3$, the Hopf fibration, the eigenspaces of the Riemannian Laplacian operator
 \label{sa3.7}}

\subsection{Hopf fibration}
We consider the map ~~ $\tilde{\Pi}:\R^4=\CC^2\rightarrow\R^3$ ~~ defined by:

\begin{eqnarray}\label{F3.1}
 \tilde{\Pi}(x_1,x_2,x_3,x_4):=(x_1x_3+x_2x_4,~~ x_1x_4-x_2x_3, ~~\frac{1}{2}(x_3^2+x_4^2-x_1^2-x_2^2))
\end{eqnarray}

 This can be written, when we set $z_1=x_1+ix_2$ ~and ~$z_2=x_3+ix_4$:
\bigskip

$\tilde{\Pi}(z_1,z_2)=(Re(\bar{z_1}z_2),~~Im(\bar{z_1}z_2),~~\frac{1}{2}(|z_2|^2-|z_1|^2)$
\bigskip

We define:

$S^3(1)=\{(x_1,x_2,x_3,x_4)\in\R^4/~~\sum_{i=1}^4x_i^2=1\}$
\bigskip

$S^2(\frac{1}{2})=\{(y_1,y_2,y_3)\in\R^3/~~\sum_{i=1}^3y_i^2=\frac{1}{4}\}$
\bigskip

It is easy to verify that: $\tilde{\Pi}(S^3(1))=S^2(\frac{1}{2})$

We then consider the following commutative diagram:
$$\begin{CD}
   S^3(1)@>i_1>>\R^4\\
   @V{\Pi}VV @VV{\tilde{\Pi}}V\\
   S^2(\frac{1}{2})@>i_2>>\R^3
  \end{CD}$$
where ~ $i_1$ ~ and ~ $i_2$ are the canonical injections and $\Pi$ the restriction of $\tilde\Pi$ to
$(S^3(1),S^2(\frac{1}{2}))$.

$S^3(1)$ is endowed with its canonical Riemannian metric by the Euclidean metric of $\R^4$ and likewise
$S^2(\frac{1}{2})$ by the Euclidean metric of $\R^3$.

Then, for any point $P$ of $S^2(\frac{1}{2})$, ~~ $\Pi^{-1}\{P\}$ is a great circle of
$S^3(1)$.

$\Pi:S^3(1)\rightarrow S^2(\frac{1}{2})$ is called a "Hopf fibration".

\begin{rmq}
 We have imposed here, by the definition of $\tilde\Pi$, the fact that the image of the sphere $S^3$ of radius $1$ is the sphere
$S^2$ of radius $\frac{1}{2}$. Multiplying $\tilde\Pi$ by a positive real number one can change the radius of the
spheres as desired. The choice given by \ref{F3.1} is justified by the fact that the map: \\
$\Pi:(S^3(1),g_{S^3(1)})\rightarrow(S^2(\frac{1}{2}),g_{S^2(\frac{1}{2})})$ \\ is then a \textbf{Riemannian submersion} , which allows to prove proposition \ref{p3.1} which will follow.
\end{rmq}

The Hopf fibration can also be introduced in the following way: the group $SO(2)$, which is identified with the group of complex numbers with modulus 1~~ $(r_\alpha\leftrightarrow e^{i\alpha})$ ~~ naturally acts on $\R^4=\CC^2$
as \textbf{isometries group} by setting: \\ for any $r_\alpha\in SO(2)
~~~r_\alpha\bullet(z_1,z_2):=e^{i\alpha}(z_1,z_2)$.

It is easy to verify that, under this action, the orbits under $SO(2)$ are exactly great circles
$\Pi^{-1}\{P\}$ of the Hopf fibration defined previously. The quotient manifold $S^3(1)/SO(2)$ is then
diffeomorphic to $S^2(\frac{1}{2})$ by the following commutative diagram:
\vspace{-0.3cm}

\[
  \xymatrix {S^3(1)\ar[r]^\Pi \ar[d]_{\Pi'} & S^2(\frac{1}{2}) \\
  S^3(1)/SO(2) \ar[ur]_f}
\]

where $\Pi'$ is the canonical Riemannian submersion associated with the quotient $S^3(1)/SO(2)$ when $g'$ is the metric on
$S^3(1)/SO(2)$ "quotiented" from $g_{S^3(1)}$. The map $f$ is obviously a diffeomorphism and one can
prove that it is actually an isometry of $(S^3(1)/SO(2),g')$ to $(S^2(\frac{1}{2}),g_{S^2(\frac{1}{2})})$. \\
The Hopf fibration therefore also corresponds to: \\ $\Pi':S^3(1)\rightarrow S^3(1)/SO(2)\sim S^2(\frac{1}{2})$.
\begin{prop} \label{p3.1}
 For any map $\varphi:S^2(\frac{1}{2})\rightarrow\R$ of class $C^2$:
 \bigskip
 
 $\Delta_{S^3(1)}(\varphi\circ\Pi)=(\Delta_{S^2(\frac{1}{2})}\varphi)\circ\Pi$
\end{prop}

where $\Delta_{S^3(1)}$ and $\Delta_{S^2(\frac{1}{2})}$ refer to the standard Riemannian Laplacian operators defined on the
spheres $S^3(1)$ and $S^2(\frac{1}{2})$.
 \bigskip
 
\textbf{Proof}.

This result is basically associated to the fact that, for any point $P$ of $S^2(\frac{1}{2})$,~~ $\Pi^{-1}\{P\}$ is
a great circle of $S^3(1)$, in particular the function $v:S^2(\frac{1}{2})\rightarrow\R$ defined
by $v:=\text{vol}_{g_{\Pi^{-1}\{P\}}}(\Pi^{-1}\{P\})$ is a constant function equal to $2\pi$, and that, moreover, $\Pi$ is
a Riemannian submersion. We do it as follows:

Since $\Pi$ is a Riemannian submersion, we know that:

for any $\varphi$ and $\phi\in C^1(S^2(\frac{1}{2}))$
\begin{eqnarray}\label{F3.2}
 \int_{S^3(1)}\varphi\circ\Pi=\int_{S^2(\frac{1}{2})}\varphi v=2\pi\int_{S^2(\frac{1}{2})}\varphi ~~~\text{and}
\end{eqnarray}
\begin{eqnarray}\label{F3.3}
 (\nabla_{S^2}\varphi\nabla_{S^2}\phi)\circ\Pi= \nabla_{S^3}(\varphi\circ\Pi)\nabla_{S^3}(\phi\circ\Pi)
\end{eqnarray}

For any $\varphi\in C^2(S^2(\frac{1}{2}))$, ~~ $\varphi\circ\Pi$ is invariant by the action of $SO(2)$ on
$S^3(1)$ defined previously, and the same is true for the function $\Delta_{s^3}(\varphi\circ\Pi)$. This one "quotiented" is therefore a function $h:S^2(\frac{1}{2})\rightarrow\R$ such that
$\Delta_{s^3}(\varphi\circ\Pi)=h\circ\Pi$.

It is then necessary to prove that $h=\Delta_{S^2}\varphi$:

According to \ref{F3.2} and \ref{F3.3}, ~ for any $\phi\in C^1(S^2(\frac{1}{2}))$:
  \begin{displaymath}
\begin{tabular} {rcl}
${\displaystyle 2\pi\int_{S^2(\frac{1}{2})}(\Delta_{S^2}\varphi)\phi}$
& = & $\displaystyle 2\pi\int_{S^2(\frac{1}{2})}\nabla_{S^2}\varphi\nabla_{S^2}
\phi=\int_{S^3(1)}(\nabla_{S^2}\varphi\nabla_{S^2}\phi)\circ\Pi$ \\
& = & $\displaystyle \int_{S^3(1)}\nabla_{S^3}(\varphi\circ\Pi)\nabla_{S^3}(\phi\circ\Pi)$
\end{tabular}
\end{displaymath}
Therefore:

 $2\pi\int_{S^2(\frac{1}{2})}(\Delta_{S^2}\varphi)\phi=\int_{S^3(1)}(\phi\circ\Pi)\Delta_{s^3}
(\varphi\circ\Pi)=\int_{S^3(1)}(h\phi)\circ\Pi=2\pi\int_{S^2(\frac{1}{2})}h\phi$

And one deduces: $\Delta_{S^2}\varphi=h$
\bigskip

This proposition shows an important link between the eigenspaces associated to the Laplacian operator
$\Delta_{S^2(\frac{1}{2})}$ and those associated to the Laplacian operator $\Delta_{S^3(1)}$. This is introduced in the following paragraph.
\subsection{Eigenspaces of $(S^3(1),g_{S^3(1)})$ and $(S^2(\frac{1}{2}),g_{S^2(\frac{1}{2})})$}
We begin by recalling the known results on the eigenspaces of $(S^n(1),g_{S^n(1)})$ where $g_{S^n(1)}$ refer to the
standard Riemannian metric on $S^n(1)$.
\begin{prop} \label{p3.2}
~ \\ \vspace*{-1em}
 \begin{enumerate}
  \item Eigenvalues of the (geometric) Laplacian operator of $(S^n(1),g_{S^n(1)})$ are given \textbf{by the sequence
$(\gamma_k)_{k\in\N}$ where $\gamma_k=k(k+n-1)$}.
\item Corresponding eigenspaces $E_k$ are the sets of \textbf{harmonics homogeneous polynomials on $\R^{n+1}$  of degree $k$ restricted to $S^n(1)$}.
 \end{enumerate}
\end{prop}

We quickly deduce that the eigenvalues of $\Delta_{S^n(\rho)}$ are given by the sequence
$(\gamma_k(\rho))_{k\in\N}$ where $\gamma_k(\rho))=\rho^{-2}k(k+n-1)$.

Eigenvalues of $\Delta_{S^3(1)}$ are therefore: $\gamma_k=k(k+2))$ for $k\in\N$ and those of
$\Delta_{S^2(\frac{1}{2})}$
: $\gamma'_l=4l(l+1)$ for $l\in\N$, that is $\gamma'_l=2l(2l+2)$.

On the other hand, proposition \ref{p3.1} shows that if $\varphi$ is an eigenfunction on $S^2(\frac{1}{2})$ for the
eigenvalue $\gamma'_l$ ~ then ~ $\varphi\circ\Pi$ is an eigenfunction on $S^3(1)$ for the same eigenvalue
$\gamma'_l=2l(2l+2)$.

If we denote $F_{2l}$ the eigenspace of $\Delta_{S^2(\frac{1}{2})}$ corresponding to the eigenvalue $\gamma'_l$, then
$E'_{2l}:=\{\varphi\circ\Pi~~/~~\varphi\in F_{2l}\}$ is a linear-subspace of the eigenspace $E_{\gamma'_l}$ and
$F_{2l}$ is naturally isomorphic to $E'_{2l}$.

\textbf{Each eigenspace $E'_{2l}$ of $S^2(\frac{1}{2})$ is therefore identified, by Hopf's fibration, to a 
linear subspace of the even index eigenspace $E_{2l}$ of $S^3(1)$ which corresponds to the eigenvalue $2l(2l+2)$}.


\subsection{Proof that, for any $k$ from 1 to 3, $\nabla_{S^3}.L_{k_{S^3}}=0$ when the
$L_{k_{S^3}}$ are the three vector fields that parallelize $S^3$ \label{ss3.1}}

It is immediate to verify that, for the Euclidean metric $\xi$ of $\R^4$, $\nabla_\xi.L_k=0$.

We consider, at one point $x$ of $S^3$, the four mutualy orthogonal vectors: $L_{1_x}, L_{2_x}, L_{3_x}, N_x$ ~~ where
$N_x$ is the normal vector to $S^3(\rho)$: $N_x=x^1\partial_1+\dots+x^4\partial_4$.

We have, for each $k$ (omitting to write the point $x$ in index):

$0=\nabla_\xi.L_k=\xi(L_1,D_{L_1}L_k)+\xi(L_2,D_{L_2}L_k)+\xi(L_3,D_{L_3}L_k)+\xi(N,D_{N}L_k)$

where ~ $D$ denotes the Euclidean connection of $\R^4$.

But $\xi(N,D_{N}L_k)=\xi(D_N N,L_k)=0$ ~ because $D_N N=N$.

In addition $\xi(L_i,D_{L_i}L_k)=\xi(L_i,\tilde{D}_{L_i}L_k)$ ~ where ~ $\tilde D$ is the Euclidean connection of $\R^4$
induced on $S^3$, ~ this because $\tilde D_{L_i}L_k$ is the orthogonal projection on $T_x (S^3(\rho))$ of $D_{L_i}L_k$.

We can deduce:

$0=\xi(L_1,\tilde D_{L_1}L_k)+\xi(L_2,\tilde D_{L_2}L_k)+\xi(L_3,\tilde D_{L_3}L_k)=\nabla_{S^3}.L_{k_{S^3}}$.


\subsection{Proof of proposition \ref{p2.4} of section \ref{s2.13} \label{ss3.2}}

\subsubsection{Stability of eigenspaces $E_p$ by the action of the three vector fields $L_1, L_2, L_3$}

According to proposition \ref{p3.2} the eigenfunctions of the Laplacian operator which constitute the eigenspace $E_p$ are 
restrictions to $S^3(\rho)$ of the homogeneous harmonic  polynomials of degree $p$ defined on $\R^4$. \\
It is therefore necessary to prove that if $P$ is such a polynomial then, for any $k$ from 1 to 3, ~~ $L_k(P)$ is still an
homogeneous harmonic polynomial of the same degree $p$ since $L_k{|_{S^3(\rho)}}(P|_{S^3(\rho)})=(L_k(P))|_{S^3(\rho)}$.

Given the expression of $L_k$, it is clear that $L_k(P)$ are homogeneous of degree $p$.

The difficulty is to prove that $L_k(P)$ is harmonic knowing that $P$ is harmonic.

For this we develop $\Delta(L_k(P)):=\sum_{i=1}^4{\partial_i}^2(L_k(P))$ ~~ where
~~ $\partial_i:=\frac{\partial}{\partial {x^i}}$.
\bigskip

If $k=1$, ~~~ $L_1(P)= -x^3\partial_1P-x^4\partial_2P+x^1\partial_3P+x^2\partial_4P$.

As $\partial_i^2(x^l\partial_m P)=2\delta_i^l\partial_i\partial_m P+x^l\partial_i^2\partial_m P$, we have:

$\sum_{i=1}^4\partial_i^2(x^l\partial_m)=x^l\sum_{i=1}^4\partial_i^2\partial_m P+2\partial_l\partial_m
P=2\partial_l\partial_m P$

since $\sum_{i=1}^4\partial_i^2\partial_m P=\partial_m\Delta P=0$.

Then:

$\Delta (L_1(P))=2(-\partial_3\partial_1 P-\partial_4\partial_2 P+\partial_1\partial_3 P+\partial_2\partial_4 P)=0$.

We also verify: $\Delta (L_k(P))=0$ for $k=2$ and $3$.

This proves the stability of $E_p$ by the action of the three vector fields $L_1, L_2, L_3$.


\subsubsection{Stability of spaces $E'_q$ by the action of the three vector fields $L_1, L_2, L_3$}
We consider the map $\Pi:\R^4\rightarrow R^3$ which defines the Hopf fibration given by:

\hspace{-5mm} $\Pi(x^1,x^2,x^3,x^4)=(x^1x^3+x^2x^4,x^1x^4-x^2x^3,1/2((x^3)^2+(x^4)^2-(x^1)^2-(x^2)^2))$

\hspace{2cm} $:=(y_1,y_2,y_3)$.
\bigskip
  
Eigenfunctions of $E'_q$ ($q$ even) are restrictions to $S^3(\rho)$ of the homogeneous harmonic polynomials of
degree ${q}$ defined on $R^4$ of the form $\tilde P=P\circ \Pi$ ~ where ~ $P$ is an homogeneous harmonic polynomial of degree
$q/2$ defined on $\R^3$.

It is therefore necessary to show that when $\tilde P$ is of the preceding form, for $k$ from 1 to 3, ~~ $L_k(\tilde
P)$ is of the form $Q\circ \Pi$  where  $Q$ is an homogeneous harmonic polynomial of degree $q/2$ defined on $\R^3$. (We know
already, according to the preceding paragraph, that $L_k(\tilde P)$ is homogeneous harmonic of degree $q$).

For this, we develop $L_k(\tilde P)$:
\bigskip

$\partial_1(\tilde P)=(\frac{\partial P}{\partial y^1}\circ\Pi)\frac{\partial y^1}{\partial x^1}+(\frac{\partial
P}{\partial y^2}\circ\Pi)\frac{\partial y^2}{\partial x^1}+(\frac{\partial P}{\partial y^3}\circ\Pi)\frac{\partial
y^1}{\partial x^1}$

\hspace{1cm} $=(\frac{\partial P}{\partial y^1}\circ\Pi)x^3+(\frac{\partial P}{\partial y^2}\circ\Pi)x^4+(\frac{\partial
P}{\partial y^3}\circ\Pi)(-x^1)$
\bigskip
  
Likewise:
\bigskip

$\partial_2(\tilde P)=(\frac{\partial P}{\partial y^1}\circ\Pi)x^4+(\frac{\partial P}{\partial
y^2}\circ\Pi)(-x^3)+(\frac{\partial P}{\partial y^3}\circ\Pi)(-x^2)$

$\partial_3(\tilde P)=(\frac{\partial P}{\partial y^1}\circ\Pi)x^1+(\frac{\partial P}{\partial
y^2}\circ\Pi)(-x^2)+(\frac{\partial P}{\partial y^3}\circ\Pi)(x^3)$

$\partial_4(\tilde P)=(\frac{\partial P}{\partial y^1}\circ\Pi)x^2+(\frac{\partial P}{\partial
y^2}\circ\Pi)x^1+(\frac{\partial P}{\partial y^3}\circ\Pi)(x^4)$
\bigskip
  
By grouping terms and then simplifying we get, for $k=1$:
\bigskip

$L_1(\tilde P)=(-(x^3)^2-(x^4)^2+(x^1)^2+(x^2)^2)(\frac{\partial P}{\partial y^1}\circ\Pi)+2(x^1x^3+x^2x^4)
(\frac{\partial P}{\partial y^3}\circ\Pi)$

\hspace{1cm} $=2(y^1\frac{\partial P}{\partial y^3}-y^3\frac{\partial P}{\partial y^1})\circ\Pi$

it only remains to verify that $ P_1:=y^1\frac{\partial P}{\partial y^3}-y^3\frac{\partial P}{\partial
y^1}$ is harmonic on $\R^3$.

We have:
\bigskip

$\frac{\partial^2 P_1}{(\partial y^1)^2}=-y^3\frac{\partial^3 P}{(\partial y^1)^3}+y^1\frac{\partial^3
P}{(\partial y^1)^2\partial y^3}+2\frac{\partial^2 P}{\partial y^1\partial y^3}$

$\frac{\partial^2 P_1}{(\partial y^2)^2}=y^1\frac{\partial^3 P}{\partial y^3(\partial y^2)^2}-y^3\frac{\partial^3
P}{\partial y^1(\partial y^2)^2}$

$\frac{\partial^2 P_1}{(\partial y^3)^2}=-y^3\frac{\partial^3 P}{\partial y^1(\partial y^3)^2}-2\frac{\partial^2
P}{\partial y^1\partial y^3}+y^1\frac{\partial^3 P}{(\partial y^3)^3}$

hence:
\bigskip

$\Delta_{R^3} P_1=y^1\frac{\partial\Delta P}{\partial y^3}-y^3\frac{\partial\Delta P}{\partial y^1}=0$
\bigskip

The stability is also satisfied for $L_2$ and $L_3$.

This ends the proof of  proposition \ref{p2.4}.


\section{Proof of Theorem \ref{2.3} \label{a3.7}}

Only the detailed proof of part 3 of the theorem is presented here. The proof of parts 1. and 2.
is very close to that of theorem \ref{2.1} given in annex \ref{a3.6}. Part 1. is obviously only
corollary of part 3. in which it is sufficient to cancel the electromagnetic potential $\Upsilon$. In fact, the
diagram of the proof of part 3 that we will present is the same as that of part 3 of theorem
\ref{2.1}, only a few terms appear in addition to "spin effect".

\textbf{Proof of 3. Theorem \ref{2.3}: "In an electromagnetic potential"}.

The function $a$ satisfies: $a=\phi\beta$ ~~ where ~~ $\phi:\Theta\times S^1(\delta)\times S^3(\rho)\rightarrow\R$ ~~ and
~~ $\beta\in E_V(\nu)$.

We have:
\begin{eqnarray}
 \Box_{g_\mathcal P}a+Sa=0 \text{~~where ~~}g_\mathcal P=g_0+h\text{~~and
~~}h=\Upsilon^\flat\otimes X_2^\flat+X_2^\flat\otimes\Upsilon^\flat \label{F78}
\end{eqnarray}

(In the following calculations $g_\mathcal P$ will be denoted $g$ to simplify the writing and
$X_2$ will be simply denoted $X$).

In a standard coordinate system of the cell $\C=\Theta\times S^1(\delta)\times S^3(\rho)\times V$, we have:

$\Box_ga=-|g|^{-\frac{1}{2}}\partial_i(g^{ij}|g|^{\frac{1}{2}}\partial_ja)$ ~~ where ~~ $|g|:=detg=detg_0$

and $g^{ij}=g_0^{ij}-h^{ij}+h^i_kh^{kj}$ ~~ with ~~ $h^{ij}=\Upsilon^iX^j+\Upsilon^jX^i$.

Given assumptions:

~~ $\Upsilon^i=0$ ~ if ~ ($i=4$ ~ and ~ $i>7$), ~~ $X^j=0$ ~ if ~ ($j<7$ ( ~ and ~ $j\neq4$) ~~ and ~ $X^4=-1$ ~ since
~ $g_0(X_2,Y)=1$

So:
\begin{eqnarray}
 \Box_ga+Sa=\Box_{g_0}a+Sa+(*)+(**) \label{F79}
\end{eqnarray}

where ~ $(*):=|g_0|^{-\frac{1}{2}}\partial_i(|g_0|^{\frac{1}{2}}h^{ij}\partial_ja)$ ~ and
~ $(**):=-|g_0|^{-\frac{1}{2}}\partial_i(|g_0|^{\frac{1}{2}}h^i_kh^{kj}\partial_ja)$.
\vspace{4mm}

Since $g_0$ is a product metric on $\Theta\times S^1(\delta)\times S^3(\rho)\times V$ we have:
\begin{multline}
 \Box_{g_0}a+Sa=(\Box_\Theta+\Box_{S^1\times
S^3}+\Delta_V)(\phi\beta)+S\phi\beta\\=\beta(\Box_\Theta\phi+(\gamma-\lambda+\nu+S)\phi)=\beta(\Box_\Theta\phi+M^2\phi)
\label{F80}
\end{multline}
\begin{enumerate}
 \item \textbf{Study of (*)}.
 
 $(*)=(*_1)+(*_2)$ ~~ where ~~ $(*_1):=|g_0|^{-\frac{1}{2}}\partial_i(|g_0|^{\frac{1}{2}}h^{ij}\phi\partial_j\beta)$
 
 and ~~ $(*_2):=|g_0|^{-\frac{1}{2}}\partial_i(|g_0|^{\frac{1}{2}}h^{ij}\beta\partial_j\phi)$.

$(*_1)=|g_0|^{-\frac{1}{2}}(\partial_i(|g_0|^{\frac{1}{2}}\Upsilon^iX^j\phi\partial_j\beta)+
\partial_j(|g_0|^{\frac{1}{2}}\Upsilon^iX^j\phi\partial_i\beta))=0$

because $X^j\partial_j\beta=0$ according to the assumption $H_{2,E}$ ~~ and ~~ $\Upsilon^i\partial_i\beta=0$ ~ according to (i) since $\partial_i\beta=O$ ~ if ~ $i\leq7$.

On the other hand:

$(*_2):=(*_{2_1})+(*_{2_2})$ ~~ where
~~ $(*_{2_1})=|g_0|^{-\frac{1}{2}}\partial_i(|g_0|^{\frac{1}{2}}\Upsilon^iX^j\beta\partial_j\phi)$

and ~~ $(*_{2_2})=|g_0|^{-\frac{1}{2}}\partial_j(|g_0|^{\frac{1}{2}}\Upsilon^iX^j\beta\partial_i\phi)$.

$(*_{2_1})=-|g_0|^{-\frac{1}{2}}\partial_i(|g_0|^{\frac{1}{2}}\Upsilon^i\beta\partial_4\phi)$ ~~ since
$\partial_j\phi=0$ if $j>7$ ~~ and ~~ $X^4=-1$ ~~ ($\partial_4\phi:=\frac{\partial\phi}{\partial u}$).

So, since $\Upsilon^i=0$ ~ for ~ $i>7$:

$(*_{2_1})=-\beta(|g_0|^{-\frac{1}{2}}(\partial_i(|g_0|^{\frac{1}{2}}
\Upsilon^i))\partial_4\phi+\Upsilon^i\partial_i\partial_4\phi)$

\hspace{8mm} $=	-\beta((\nabla_{g_0}\cdotp  \Upsilon)\partial_4\phi+\Upsilon(\partial_4\phi))$.

On the other hand, since $X^4=-1$:

$(*_{2_2})=-|g_0|^{-\frac{1}{2}}\partial_4(|g_0|^{\frac{1}{2}}\beta
\Upsilon(\phi))+\sum_{j>7}|g_0|^{-\frac{1}{2}}\partial_j(|g_0|^{\frac{1}{2}}X^j\beta \Upsilon(\phi))$.

But, since $\partial_j(\Upsilon(\phi))=0$ ~ if ~ $j>7$:

$\sum_{j>7}|g_0|^{-\frac{1}{2}}\partial_j(|g_0|^{\frac{1}{2}}X^j\beta
\Upsilon(\phi))=\Upsilon(\phi)\nabla_{g_0}\cdotp  (\beta X)=0$.

Then, since $|g_0|$, $\beta$ and $\Upsilon^i$ do not depend on $u$:

$(*_{2_2})=-\beta\Upsilon(\partial_4\phi)$.

Finally:
\begin{eqnarray}
 (*)=(*_1)+(*_{2_1})+(*_{2_2})=-\beta((\nabla_{g_0}\cdotp  \Upsilon)\partial_4\phi+2\Upsilon(\partial_4\phi))
\label{F81}
\end{eqnarray}

\item \textbf{Study of (**)}.

$(**)=-|g_0|^{-\frac{1}{2}}\partial_i(|g_0|^{\frac{1}{2}}h^i_kh^{kj}\partial_ja)$

\hspace{7mm}
$=-|g_0|^{-\frac{1}{2}}\partial_i(|g_0|^{\frac{1}{2}}
(\Upsilon^iX_k+\Upsilon_kX^i)(\Upsilon^kX^j+\Upsilon^jX^k)\partial_ja)$

So, since $\Upsilon^kX_k=0$ ~~ and ~~ $X^kX_k=0$:

$(**)=-|g_0|^{-\frac{1}{2}}\partial_i(|g_0|^{\frac{1}{2}}X^iX^j\Upsilon_k\Upsilon^k\partial_j(\phi\beta))$

But, since $\Upsilon_k\Upsilon^k$ does not depend on the variables of $S^1\times V$:

$(**)=-\Upsilon_k\Upsilon^k|g_0|^{-\frac{1}{2}}(\partial_i(|g_0|^{\frac{1}{2}}
X^iX^j\beta\partial_j\phi)+\partial_i(|g_0|^{\frac{1}{2}}X^iX^j\phi\partial_j\beta))$.

As $X^j\partial_j\phi=-\partial_4\phi$ since $\partial_j\phi=0$ for $j>7$ ~ and $X(\beta)=0$, we obtain:

$(**)=\Upsilon_k\Upsilon^k|g_0|^{-\frac{1}{2}}\partial_i(|g_0|^{\frac{1}{2}}X^i\beta\partial_4\phi)$.

Considering $i=4$ then $i>7$, we have:

$(**)=-\Upsilon_k\Upsilon^k(\beta\partial^2_4\phi-|g_0|^{-\frac{1}{2}}(\sum_{i>7}\partial_i(|g_0|^{\frac{1}{2}}
X^i\beta))\partial_4\phi)$.

What we can write:

$(**)=-\Upsilon_k\Upsilon^k(\beta\partial^2_4\phi-(\nabla_{g_0}\cdotp  (\beta X))\partial_4\phi)$.

And, since $\nabla_{g_0}\cdotp  (\beta X))=X(\beta)+\beta\nabla_{g_0}\cdotp  X=0$:
\begin{eqnarray}
 (**)=-\Upsilon_k\Upsilon^k\beta\partial^2_4\phi \label{F82}
\end{eqnarray}

\item \textbf{End of the proof of 3}.

According to \ref{F78}, \ref{F79}, \ref{F80}, \ref{F81} and \ref{F82} we have:

$0=\Box_ga+Sa=(\Box_\Theta\phi+M^2\phi-(\nabla_{g_0}\cdotp  
\Upsilon)\partial_4\phi-2\Upsilon(\partial_4\phi)-\Upsilon_k\Upsilon^k\partial^2_4\phi)\beta$.

Hence:
\begin{eqnarray}
 0=\Box_\Theta\phi+M^2\phi-(\nabla_{g_0}\cdotp  
\Upsilon)\partial_4\phi-2\Upsilon(\partial_4\phi)-\Upsilon_k\Upsilon^k\partial^2_4\phi \label{F83}
\end{eqnarray}

\begin{enumerate}
 \item Study of the term: $(\nabla_{g_0}\cdotp  \Upsilon)\partial_4\phi+2\Upsilon(\partial_4\phi)$.
 
 $\Upsilon$ was chosen in the form;
 
 $\Upsilon=A +C$ ~~ where ~~ $A =\sum_{i=0}^3A^i\frac{\partial}{\partial x^i}$ ~~ and ~~ $C={\textstyle{\frac{\varrho}{2}}}\sum_{k=1}^3B_kL_{k_{S^3}}$
 
 $\varrho$ is the gyromagnetic constant.
 
 $A$ is a vector field defined on $\Theta$.
 
 $C$ is a vector field tangent to $S^3(\rho)$. ($A$ and $C$ are considered defined on $\Theta\times S^3(\rho)$).
 
 We have:
 
 $(\nabla_{g_0}\cdotp  \Upsilon)\partial_4\phi+2\Upsilon(\partial_4\phi)=\sum_{i\leq7}(|g_0|^{-\frac{1}{2}}
\partial_i(|g_0|^{\frac{1}{2}} \Upsilon^i)\partial_4\phi+2\Upsilon^i\partial_i\partial_4\phi)$.

What we can write, since $\Upsilon^4=0$, in the form:

$\sum_{i=0}^3(|g_0|^{-\frac{1}{2}}\partial_i(|g_0|^{\frac{1}{2}}
\Upsilon^i)\partial_4\phi+2\Upsilon^i\partial_i\partial_4\phi)+\sum_{i=5}^7(|g_0|^{-\frac{1}{2}}
\partial_i(|g_0|^{\frac{1}{2}} \Upsilon^i)\partial_4\phi+2\Upsilon^i\partial_i\partial_4\phi)$.

Hence:

$(\nabla_{g_0}\cdotp  \Upsilon)\partial_4\phi+2\Upsilon(\partial_4\phi)=(\nabla_\Theta\cdotp  
A)\partial_4\phi+2A\partial_4\phi+(\nabla_{S^3}\cdotp  C)\partial_4\phi+2C(\partial_4\phi)$.

But $\nabla_{S^3}\cdotp  C={\textstyle{\frac{\varrho}{2}}}\sum_{k=1}^3B_k(\nabla_{S^3}L_{k_{S^3}})=0$
since ~ $\nabla_{S^3}L_{k_{S^3}}=0$ ~~ (cf \ref{ss3.1}).

Therefore:
\begin{eqnarray}
 (\nabla_{g_0}\cdotp  \Upsilon)\partial_4\phi+2\Upsilon(\partial_4\phi)=(\nabla_\Theta\cdotp  
A)\partial_4\phi+2A\partial_4\phi+2C(\partial_4\phi) \label{F84}
\end{eqnarray}

\item Study of the term $\Upsilon_k\Upsilon^k\partial_4^2\phi$.

$A$ and $C$ being $g_0$-orthogonal we have:

$\Upsilon_k\Upsilon^k=A_kA^k+C_kC^k$.

But, since the three vector fields $L_{j_{S^3}}$ are mutualy
$g_0$-orthogonal\\
and that
$g_{0_{S^3}}(L_{j_{S^3}},L_{j_{S^3}})=\rho^2$ ~ where  $\rho$ is the radius of the sphere $S^3(\rho)$, we have:

$C_kC^k=g_0(C,C)=({\textstyle{\frac{\varrho}{2}}})^2\rho^2\sum_{j=1}^3B^2_j$.

Hence:
\begin{eqnarray}
\Upsilon_k\Upsilon^k\partial_4^2\phi=(A_kA^k+({\textstyle{\frac{\varrho}{2}}})^2\rho^2\sum_{j=1}^3B^2_j)\partial_4^2\phi \label{F85}
\end{eqnarray}

Finally the equation \ref{F83} is written:
\begin{eqnarray}
 0=\Box_\Theta\phi+M^2\phi-(\nabla_{g_0}\cdotp  
A)\partial_4\phi-2A(\partial_4\phi)-A_kA^k\partial^2_4\phi-(\alpha) \label{F86}
\end{eqnarray}

Where $(\alpha):= C(\partial_4\phi)+({\textstyle{\frac{\varrho}{2}}})^2\rho^2(\sum_{k=1}^3B^2_k)\partial^2_4\phi$

In other words:

$(\alpha):={\textstyle{\frac{\varrho}{2}}}\sum_{k=1}^3B_kS_k(\partial_4\phi)+({\textstyle{\frac{\varrho}{2}}})^2\rho^2(\sum_{k=1}^3B^2_k)\partial^2_4\phi$

where the $S_k$ were specified in the definition \ref{d2.21}.

The equation \ref{F86} is also in the form:
\begin{eqnarray}
 0=-\sum_{j=0}^3\varepsilon_j(\frac{\partial}{\partial x^j}+A^j\frac{\partial}{\partial u})^2\phi+M^2\phi-(\alpha)
\label{F87}
\end{eqnarray}

where $\varepsilon_0=-1$ ~~ $\varepsilon_1=\varepsilon_2=\varepsilon_3=1$

To obtain equation for $a_c$ given by the theorem, we consider the isomorphism\\
(cf \ref{s2.8}):
$\CC_{\lambda,\nu}:E_{S^1(\delta)}\otimes E_p\rightarrow E_p^{\CC}$ ~~ where ~~ $E_p=E_{S^3(\rho)}(\nu)$

We easily verify that:
\begin{center}
  \begin{tabular} {rcl}
    $\CC_{\lambda,\nu}\circ((\Box_\Theta\phi)_x(.))$ & $=$ & $(\Box_\Theta a_c)(x)$ \\ [0.7em]
    $\CC_{\lambda,\nu}\circ((\frac{\partial\phi}{\partial u})_x(.))$ & $=$ & $-iQ^+a_c(x)$ \\ [0.7em]
    $\CC_{\lambda,\nu}\circ((\frac{\partial^2\phi}{\partial u^2})_x(.))$ & $=$ & $-{Q^+}^2a_c(x)$ \\ [0.7em]
  \end{tabular}
\end{center}

By composing, for each $x\in\Theta$, each side of equation \ref{F87} by $\CC_{\lambda,\nu}$,\\
we obtain the
 equation \ref{F49} of
theorem \ref{2.3}:
\[\sum_{j=0}^3\varepsilon_j(i\frac{\partial}{\partial x^j}+Q^+\Upsilon^j)^2a_c+M^2a_c+\varrho
Q^+\sum_{k=1}^3B^k\hat S_k(a_c)+{Q^+}^2({\textstyle{\frac{\varrho}{2}}})^2\rho^2|B|^2a_c=0\]
\end{enumerate}
\end{enumerate}

\newpage
\section{The choice of a pseudo-Riemannian manifold}
\label{a3.8}
\begin{enumerate}
\item The structure of manifold.
  \begin{enumerate}
  \item The real numbers $\R$ set. \\
    If we want to study a map defined on a finite set
     and with values in another finite set, these having a
    very large number of elements and being provided, for example, with an
     order relation, we have a priori very few means.
    One method is to "intelligently plug the holes" into
    these two sets by adding virtual elements and then to
    extend the map on these now
    "continuous" two new sets. The first step is treated mathematically 
     in the construction of $\Q$ and then $\R$ from
    $\N$. This non-trivial construction (especially to pass from $\Q$
    to $\R$) introduces notions of "limit",
    "derivative", "differential equation", etc., relating to the functions
    defined from $\R$ to $\R$ and this gives a very large
    power to the process of study of the functions, these being able to
    be, if necessary, restricted to the finished sets of the
    original problem. In fact, the construction of $\R$, purely
    axiomatic, introduces all the tools of what is called
    "the real analysis".

  \item The vector space $\R^n$. \\
    The construction of $\R^n$ as a set of $n$-uplets of real numbers
    is a triviality, but it introduces the important notion of
    \textbf{dimension} by its natural structure of vector space
    on $\R$. The standard topology of $\R^n$ (global) is
    "poor". The differential analysis developed on $\R$ extends
    naturally on $\R^n$.

  \item Topological manifolds. \\
    A topological manifold is, by definition, a topological space
    where each point admits a neighborhood homeomorphic to $\R^n$.The interest
    of this structure is that it keeps the local topological properties of
    $\R^n$ but allows a large spectrum of global topologies.
    There is no canonical algebraic structure on a
    topological manifold unlike on $\R^n$. This is an advantage
    for the description of space-time, the algebraic structure
    of $\R^n$ having proved too rigid to make physic sense
     (which led to the construction of general relativity).
    Of course, all the analytical tools associated with the
    topology can be used (continuity, etc.), however,
    the differential analysis on $\R^n$ (which comes from that of $\R$),
    does not extend to topological manifolds. (One can, locally,
    come back to $\R^n$ by an homeomorphism, but the differential calculation
    defined this way on an open subset of the topological manifold
    completely depends on the choice of the homeomorphism and,
    as a priori none is privileged, this has no interest).
    To correctly introduce differential analysis on a
    topological manifold, it is necessary to give it an 
    additional structure via a "differentiable atlas".

  \item Differential manifolds. \\
    A differential manifold is, by definition, a couple
     (topological manifold,  complete atlas) (the reader is referred to specialized works for
    a precise definition). An atlas is a set of "charts"
    (homeomorphisms from an open subset of the topological manifold to an open subset of
    $\R^n$) and this notion has a very important physical meaning
    . The choice of a chart of this atlas (also called a
    coordinate system) can be seen as the choice of an
    "observer" which express what happens on an open set of the
    manifold by bringing it back to $\R^n$. What is conceptually
    very important is that the data of a complete atlas
    (differentiable, class $C^k$, etc.) on the topological manifold
    allows to reconstruct much of the differential analysis
     on this manifold and this \textbf{regardless of
      any particular choice of a chart in the atlas}. This
    reconstruction begins with that of a "tangent space" in
    each point of the manifold, which is a real vector space 
    whose dimension corresponds to the topological dimension of the
    manifold (this tangent space is independant of a 
    choice of a chart in the atlas). From there, all
    "objects" of differential analysis are redefined without
    difficulties: tensor field, differential map, etc.). All of these
    notions are such that they do not depend of a chart in
    the atlas, in other words, of the "look" of an
    observer on the manifold. In summary, a differential manifold
    is a topological space on which is defined a (complete) set of
     observers so that most of the tools of
    the known analysis on $\R^n$ can be used. In addition, all
     notions of this analysis have their definitions independent
    of an observer choice. There is, a priori, no algebraic structure
     on (the underlying set to) the manifold, so no
     artificial physical notion (as it was the
    case, for example, in $\R^n$ for objects that are:
    straight lines, affine subspaces, origin, etc.), this gives
    great freedom to introduce notions that will have
    a physical sense.
  \end{enumerate}

\item The pseudo-Riemannian manifold structure. \\
  A pseudo-Riemannian manifold is a couple: (differential manifold
   $\M$ , field $g$ of quadratic forms defined
  on this manifold). The field of quadratic forms is the datum, in
  any point $x$ of $\M$, of a symmetric bilinear form $g_x$ on
  the tangent space at this point. We do not impose a 
  restriction on the signature of $g_x$ or even on the fact that it is
  degenerate or not. \textbf{In the theory that we present
    in this paper, all physics is described by the single data
    of a pseudo-Riemannian manifold}. Any usual "objects"
  are defined from $g$. The choice of a field of
  quadratic forms is therefore an important fact, it is
  essentially justified by how easy it becomes to define the concepts
   of "distance" and "time". It is important that a physical theory introduces
  standard concepts of "distance" and "time" associated to an
  observer. Of course, we should not hesitate to
  "try" other objects than a quadratic form field if
  necessary, but, for the moment, this choice seems to be
   perfectly fine.

\item In summary. \\
  Although it is conceivable that "space-time" can be described
  by a finite set provided with a structure ($X$, Struct), it is
  probably much more interesting, for the sole purpose of describing
  this structure, to embed ($X$, Struct) into a pseudo-Riemannian 
  manifold ($\M$, $g$) in order to take advantage from the mathematical analysis developed in the
  latter. I'm not sure, however, that there is much
  interest in specifying a possible couple ($X$, Struct) and its
  embedding in ($\M$, $g$).

  So we chose in this paper to represent the universe by a
  pseudo-Riemannian manifold ($\M$, $g$) of dimension $n>5$.
  As explained in the introduction, this manifold will be
  considered "totally anarchic". The reader may refer to
   chapter 14 of the manuscript \cite{vaugon-1} where it is specified what a
  "totally anarchic manifold" is in the case where
  this one is Lorentzian of dimension 4 (which extends
  naturally to our manifold ($\M$, $g$)). This manifold is
  "dotted" with singularities of different types (big-bangs,
  big-crunches, black-holes, etc.). The singularity "big-bang" that
  it seems to be observed (almost homogeneous and isotropic "with great
  scale") is only a detail of this "anarchic" universe.This
  point of view is quite different from the one commonly accepted in
  cosmology, although it does not change anything fundamental in this
  field. Our point of view excludes any attempt to
  precision in the global description of the universe.
\end{enumerate}

\section{Determinism and approximations}
\label{a3.9}
We consider a domain with a geometric type $(\D, g, \A)$ (cf Mathematical preliminary, it is recalled that a "type" is
 a geometric constraint on $(\D, g)$). It is
natural to say that such a domain is \textbf{totally
  deterministic} if the knowledge of the tensor $g$ on a subdomain
$\D'$ of $\D$ completely determines $g$ on $\D$, that is, if
a couple $(\D, g')$ satisfies the same geometric constraint as the one
which gives the type of $(\D', g)$ and if $g' = g$ on $\D'$ then,
necessarily, $g = g'$ on $\D$. Proof that a domain with a geometric type is
totally deterministic is usually a very complex problem often named "Cauchy problem", but a rigorous presentation of this problem is beyong the scope of this paper.

Although these problems are mathematically interesting, the fact that the manifold $\M$ is of great dimension and is locally diffeomorphic to a product $\Theta \times K$ (where $\Theta$ is an open set of $\R^4$ and $K$ a compact manifold) makes them difficult to solve and we do not adress them in this paper. However we give some approximation methods that make domains that we define "sufficiently deterministic" to be humanely interesting. 
These "approximation methods" are none other than
those commonly used in standard physics and here adapted to
considered spaces.
Take as an example the case of a "fluid" type domain. The
equations given by the theorem \ref{t1.1} are not
"sufficiently deterministic". Differential operators are
associated to the metric tensor $g$ itself. An important particular case
 for which these equations can become usable is
the one where, in a specific coordinate domain $\Theta \times K$,
the metric tensor $g$ is written in the form $g_{ij} = {g_0}_{ij} +
h_{ij}$ where $g_0$ is a neutral potential metric such that $g_0
|_\Theta$ is the Minkovski metric and the functions $h_{ij}$ are
$\ll 1$ (the $\partial_kh_{ij}$ being also controlled). Equations of the theorem \ref{t1.1} can then be
rewritten to approximation by replacing the differential operators
associated to $g$ by those associated to $g_0$. However, we must verify
that the terms deleted by this process (which come from the Christoffel symbols) are "negligible" compared to the terms
remaining. This is not enough to make
usable equations in a classical sense, the terms of
"pressure" give too many "unknowns" compared to the number of
equations. The cases studied are then those where it is assumed, for example, that these
terms are null (or negligible) which corresponds to the case of a
really perfect fluid (definition \ref{def:5}), or, more
generally, those where additional equations are imposed on
the pressure terms (state equations). If now,
it is assumed that the remaining unknown functions do not depend,
in the coordinate system $\Theta \times K$, of the "variables"
of $K$, we recover equations on the fluids (electrically charged or not) of classical physics (this is how the
Newtonian physics appears as an approximation of the theory of standard 
general relativity, itself approximation of physics
presented here).

The approximations just described, associated to the fact that the
metric tensor is "close" to a neutral potential metric
$g_0$, have the advantage of showing that equations of
classical physics (which are sufficiently deterministic)
are easily deduced from those obtained in the theorem \ref{t1.1}. \textbf{The
"approximations" which consist in giving \textbf{explicitly} a
metric tensor $g$ on a domain $\D$ are actually
much more interesting}, they allow precise calculations in a
framework (which does not necessarily assume that $g$ is
"close" to a $g_0$). These are the "approximate types" that
are used in standard general relativity under the name of
"exact solution of the Einstein equation" and these adapt without
difficulties to the domains $\D$ of the n-dimensional manifold $\M$
that we consider here. However, it is especially the "potential" types
that we define in \ref{s1.4} that allow the precise description of
many physical phenomena through their
geodesics as explained in \ref{s1.3} and which give back
simply, in approximation, Newtonian physics and
standard electromagnetism.

It is, in my opinion, continuing in this direction, by the explicit data
of "metrics" $g$ on domains $\D$, which we will study
effectively more complex areas of type "fluid with
pressure or other". Note that 
chapter 2 on quantum phenomena considers only  domains of type
 defined from the tensor $g$.

It is important to note that, technically, it is by imposing
\textbf{group actions invariances} that we define domains sufficiently simple to be interesting. Here is
some examples:
\begin{itemize}
\item The condition \ref{ei12} given in the definition \ref{def:4}
  of a fluid domain, which reflects the fact that we neglect
   quantum effects of electromagnetism, can be written in
  saying that the metric tensor $g$ is invariant by the group of
  diffeomorphisms generated by the vector field $Y$.
\item In various places of this paper, conditions imposed on
  functions defined on a cell $\C = \Theta \times K_1
  \times \cdots \times K_l$ (which often serve to define the
  metric tensor) are stated by saying that these functions 
  do not depend on the variables of $K_i$. This means
  that these functions are invariant by  action of  group that acts
  transitively on $K_i$ (this can also be written in
  making the group acts on the manifold ~ $\M$).
\item Many important examples of domains with geometric type in
  standard general relativity (which is easily adapted to
  $\M$ as we have already seen) are defined by assuming invariances
  by classical groups. We can mention: the domains of
  Schwarzschild, Reissner-Nordström, Friedmann, invariants by
   action group $SO(3)$ as well as the example introduced in
  \ref{s1.7}, Kerr domains, invariants by $SO(2)$, etc.
\end{itemize}

It would be interesting to detail all these approximation methods 
much more than we did in the few lines of this annex but we choose to leave this (sometimes difficult) exercice to the reader.

\section{Summary of proposed theory. \\
Some clarifications on the link between "reality" and "mathematical model"}
\label{a3.10+}
The mathematical model chosen here to represent "reality" is a pseudo-Riemannian manifold $(\M,g)$ of dimension $n>5$. In this model studied domains are described by triplets $(\D,g,\A)$ where $\D$ is a domain of $\M$ and $\A$ a complete g-
observation atlas (cf Mathematical preliminary). Some elements that justify this choice are given in annex \ref{a3.8}. The important fact of the theory presented here is that no other "mathematical object" than the tensor $g$ and the observation atlas $\A$ is introduced into the manifold $\M$, no law, no principle associated to these objects is introduced. All objects of the mathematical model that turn out to be important are thus defined only from the local properties of $\M$ and the tensor $g$, in other words they are only characteristics of the geometry of the pseudo-Riemannian manifold  $(\M,g)$. We will give a little further clarification on the connection between some objects of this mathematical model and the usual "objects" of reality. It is, of course, very important not to confuse what we will continue to call "objects of the mathematical model" and "objects of reality". \\
Following the domains $(\D,g)$ of the manifold $(\M,g)$ equipped with a complete observation atlas $\A$, the geometric characteristics are different and mathematical objects that prove to be important are correctly defined from $(\D,g,\A)$ only when this triplet has some special properties. These particular properties are often obtained after approximations which  consist in neglecting some features of the tensor $g$. Domains $(\D,g,\A)$ that are interesting are those that are sufficiently deterministic (this notion has been specified in the introduction), we have classified them into two categories whose study corresponds to the two chapters of this paper. Nevertheless, some areas of interest can not be entirely situated in one or the other of these categories, they have not been introduced in this paper but we will say some words at the end of this annex. \\
The first category covers classical non-quantum physics, the 2nd category covers quantum physics. To clarify this, it is, of course, necessary to explain the links between mathematical objects described in the chosen mathematical structure and those of reality. \\
We now explain the processes used for this purpose (which differ according to the two categories). \\
\textbf{1st category} -To clarify the link between the mathematical model and reality, a first method consists in bringing back notions of the mathematical model introduced in the 1st chapter, defined on a manifold of dimension $>5$ (domains of fluid and potential types, for example), to notions defined on a manifold of dimension $4$, in order to recover the point of view of standard general relativity and thus to make the link with the "reality". This strategy is supported by the fact that, according to proposition \ref{p0,1++}, $(\D,g,\A)$ can be considered as a fiber manifold $(\M,B, \pi)$ whose base $B$ is a manifold of dimension 4 and fiber type a compact manifold $K$. On the other hand, a way of neglecting quantum effects, and which enters into the definition of this 1st category of domains, is to consider that the tensor $g$ defined on $\M$ is invariant by a group which acts transitively on $K$. This makes it possible to define a quotient metric $\bar{g}$ on the basis $B$, of signature $(-,+++)$ by identifying in particular for each $x\in\M$, the apparent space $H_x$ (cf Mathematical preliminary) to the tangent space $T_{\pi(x)}(B)$. It is important to note that then theorems \ref{t1.1}, \ref{t1.2}, and proposition \ref{p1.1} give results on the fluids very close to those of standard general relativity. \\
In fact, it is simpler, to make the link with reality, to consider a chart of the observation atlas $\A$ (associated to the choice of an observer). "Space-time" is then reduced by this chart to a cell $(\Theta\times K)$ where $\Theta\subset \R^4$ (on which the tensor $g$ transported by the observation diffeomorphism is still denoted $g$). The fact that, to neglect  quantum effects, one assumes that $g$ is invariant by a group acting transitively on $K$ is expressed  here by the fact that  objects of the mathematical model defined on $(\Theta\times K)$ ( which are all constructed from $g$) do not depend on the variables of $K$ and can be considered defined on the classical "space-time" $\Theta\subset\R^4$, and thus we recover  standard notions of: geodesic , energy density function, electric charge density function, etc., and it should be noted that     vector (or tensor) fields are considered, for each $x$, with values in the apparent space $H_x$, identified with $T_{\pi(x)}(B)$, when they are projected. It is from this point of view that results intoduced in sections \ref{s1.5}, \ref{s1.7} and \ref{s1.8} were given. The space-time "seen" by the observer being $\Theta\subset \R^4$ the link with reality is then obtained as in classical physics. \\
Process described above can be summarized by saying that characteristics of the geometry of $(\D,g,\A)$ are transported on $\Theta\subset R^4$ by becoming the standard objects of classical physics, and the simple mathematical properties of $(\D,g,\A)$ give the \textbf{laws} of classical physics. \\
\textbf{2nd category} -The fundamental property that makes it possible to define this category is the fact that, transported on $\Theta\times K$ by a chart of the observation atlas, the tensor $g$ now depends strongly on variables of $K$. As $K$ is a compact manifold, the spectral decompositions on $K$, associated with the variations of the tensor $g$, naturally lead to discrete magnitude definitions, characteristic of the relevant geometric domain $(\D,g,\A)$. It is these magnitudes which make it possible to recover standard quantum physics notions (mass, electric charge, spin, etc.). However, the link with reality is now obtained in a different way compared to domains of the first category, and it is obviously no longer an option to use a "quotient" of the tensor $g$ so as to "bring it back" to $\Theta\subset\R^4$  since it is precisely the variations of $g$ on $K$ that are important. It is necessary here, to link back with reality, to use the notion of "measuring instruments" which is precisely defined in section \ref{s2.16} and which thus becomes a fundamental concept for this category of domains. Within the framework of domains classified in the 1st category we did not specify  notions of "measurement" or "of measuring instruments" because these do not require any originality (just as in classical non-quantum physics) and we can generally  considere that, in the physical systems studied in this first category, measurements do not "disturb" the system. \\
In the framework of the 2nd category domains, measurement instruments are "integral parts" of domains studied when we choose to make a precise link with reality. \\
An important notion that intervenes in the study of domains of this 2nd category is that of \textbf{singularity} that was introduced in section \ref{s2.11} (whereas this notion did not intervene for domains of the 1st category). Singularities describe, in fact, only one characteristic of the tensor $g$ but it is simpler at first to consider them as connected parts (the elementary singularities (cf definition \ref{d2.20})) whose union is a null measure subset in $\D$ (cf definitions \ref{d2.19} and \ref{d2.20}). It is the particular asymptotic behavior of $g$ in the neighborhood of these singularities which makes it possible to "locate" them in the $\Theta$ component of $\Theta\times K$ in the chart of the observation atlas linked to an observer (note that, for the moment, we do not try to analyze, within the framework of our theory, the phenomena which make it possible to "locate" these singularities in $\Theta$). \\
The measurement instruments that have been theorized in section \ref{s2.16} (which are necessary to make the link between the mathematical model and reality) fundamentally use the notion of singularity. The theoretical action of measuring instruments can be summarized very roughly by saying that they make a spectral analysis of the tensor $g$ in the domain $(\D,g)$ associated to the magnitude to be measured, by creating disjoint "space-time" domains $(\D_k,g_k)$ where each $g_k$ describes an elementary component of the spectral decomposition of $g$ (all this is detailed in section \ref{s2.16}). Count of singularities in each $\D_k$ then makes it possible to obtain the "intensity coefficient" of each elementary component of the spectral decomposition of $g$ in $\D$ (and each component $g_k$ is associated with a type of standard particle). It is in this spectral decomposition of $g$ by the measuring instruments that the "position-momentum" or "time-energy" uncertainty relations appear. \\
Link with reality is thus finally obtained (within the framework of domains of this 2nd category) by the identification of singularities in domains created by the measuring instruments. In practice (in reality), measuring instruments are "screens", "bubble chambers", "wire", "drift", etc. (Clarifying the relationship between these actual measuring instruments and those described in section \ref{s2.16}, which are part of the mathematical model, requires some work). \\
As already mentioned in the introduction to Chapter 2, an elementary singularity must not be considered as a "particle" in the usual sense. The only cases for which we can make a precise link between singularities and "standard particles" are those for which considered domains with singularities $(\D,g,\Ss)$ (definition \ref{d2.19}) are such that $g$ describes an elementary oscillating metric of which the classification has been introduced in sections \ref{s2.4} and \ref{s2.5'} and makes it possible to precisely recover notions of electrons, photons, etc. These particular domains with singularities are obtained, for example (from what has just been introduced), as a result of measurements by instruments. \\
It is important to note that, in our mathematical model, we can not associate any standard notion of "particle" with a singularity in a domain for which $g$ has a "complex" spectral decomposition. This is fundamentally different from  classical theories (quantum or not) and it is recalled that, in the theory that we present, physics is entirely described by the sole geometry of domains $(\D,g,\A)$, which has nothing to do with any "entity" built from "elementary bricks" associated with the notion of particles whose laws are sought to determine their behaviors. \\
In the mathematical model that we present, no axiom manages singularities, only a law of equiprobability associated to the tensor $g$ is supposed to apply (cf section \ref{s2.12}) within the framework of generic experiments of quantum physics (diffraction, Young's slits, potential deviations, Stern-Gerlach experiment, etc.) as well as in the design of measuring instruments. (It is conceivable, however, that this equiprobability law may no longer apply for some domains of this 2nd category describing complex situations and discussed in section \ref{s2.17}.) \\
The notion of singularity introduced into our theory (which gives randomness to results of some experiments as the preceding lines specify) certainly does not have the same conceptual importance as the notion of particle in classical theories. The 2nd chapter shows that it is the oscillating metrics (cf definition \ref{d2.2}) that, for us, manage the study of quantum phenomena. This study is "deterministic" in the sense that equations obtained admit an unique solution when  boundary conditions are assumed to be known. Random part on singularities is then specified, when calculations are finished, by the simple law of equiprobability. \\

Let us note now that some interesting domains may not fall into one or the other of the 2 categories. For example, domains which, in the mathematical model that we propose, make it possible to describe phenomena  consistent with all  "\textbf{universe expansion}" measurements. If the considered model $(\D,g,\A)$ is such that $\D=\Theta\times K_1\times\dots\times K_p$ where $\Theta$ of dimension 4 (and can be supposed compact) is associated with the space-time of standard general relativity, it is interesting to study "geometric types"on $(\D,g,\A)$ that give an "evolution" on each $K_i$ (expansion, retraction or stability) and an expansion on $\Theta$ corresponding to the current measurements. Singularities of type "big-bang", "big-crunch" or others, described in standard general relativity can be, in our model of dimension $>5$, of different nature. The influence on $\Theta$ expansion by the evolution of manifolds $K_i$ in the neighborhood of a singularity of "big-bang" type (for example) can be considered as the equivalent of the influence of quantum phenomena in an area close to the "big bang" in standard general relativity (which eventually leads to an "inflationary" expansion). Of course, what is introduced in these last lines is very imprecise and was written only to initiate research on the subject.

\tableofcontents

\label{sec:refs}
\nocite{*}
\renewcommand\bibname{Références}
\bibliographystyle{plain}
\bibliography{references}

\end{document}